\documentclass[a4paper,10pt]{article}
\usepackage{amssymb,amsmath,psfrag,graphicx}
\usepackage{amsthm}
\usepackage{equationarray}
\usepackage{verbatim}
\usepackage{color}
\usepackage{morefloats}
\usepackage{colortbl}
\usepackage{fullpage}
\usepackage{enumitem}
\usepackage{tabularx}
\usepackage{multirow}
\usepackage{booktabs} 
\usepackage{marginnote}
\usepackage{siunitx}
\DeclareSIUnit\mmHg{mmHg}
\usepackage{graphicx} 
\usepackage[x11names]{xcolor}
\usepackage{algorithm,algpseudocode}
\usepackage[hidelinks]{hyperref}
\usepackage{tikz}
\usepackage[edges]{forest}
\usepackage{subcaption}
\usetikzlibrary{arrows.meta,
                positioning,
                shapes.symbols
                }
\tikzset{%
  >={Latex[width=2mm,length=2mm]}, 
            base/.style = {rectangle, rounded corners,
                           minimum width=4cm, minimum height=2.3cm,
                           text centered},
}

\usepackage{mathrsfs}
\usepackage{upgreek}

\usepackage{todonotes}
\usepackage{ulem}
\usepackage[bottom]{footmisc}

\setlist[description]{font=\normalfont\itshape}

\newenvironment{quoteit}
{\begin{quote}\itshape}
{\end{quote}}

\newtheorem{theorem}{Theorem}[section]

\theoremstyle{remark}
\newtheorem{remark}{Remark}[section]

\theoremstyle{definition}
\newtheorem{definition}{Definition}[section]
\graphicspath{{Images/}}

\usepackage{etoolbox} 
\newtheorem{example}{Example}[section]
\AtEndEnvironment{example}{\qed}

\newcommand{\papertitle}{Combining physics--based and data--driven models: advancing the frontiers of research with Scientific Machine Learning}
\newcommand{\paperkeywords}{Scientific Computing, Approximation of PDEs, Machine Learning, Artificial Neural Networks, Scientific Machine Learning}
\definecolor{lightgreen}{RGB}{204,255,204}
\definecolor{darkgreen}{RGB}{0,200,0}
\definecolor{lightorange}{RGB}{255,255,204}
\definecolor{lightred}{RGB}{255,204,204}
\definecolor{tolook}{RGB}{200,50,50}

\newcommand{\argmin}[1]{\underset{#1}{\operatorname{arg}\!\operatorname{min}}\;}

\newcommand{\lifex}{\texorpdfstring{\texttt{life\textsuperscript{{x}}}}{lifex}}

\newcommand{\NNgeneric}{\mathcal{N\!N}}

\newcommand{\actState}{\mathbf{z}_{\textrm{act}}}

\newcommand{\actRHS}{{\boldsymbol\Phi}_{\textrm{act}}}
\newcommand{\actOBS}{{\Psi}_{\textrm{act}}}
\newcommand{\actTension}{T_{\textrm{a}}}
\newcommand{\Cai}{[\textrm{Ca}^{2+}]_i}

\newcommand{\Wratio}{\lambda_{\text{cycle}}}
\newcommand{\Weq}{\lambda_{\text{eq}}}

\newcommand{\prs}{p}
\newcommand{\vol}{V}
\newcommand{\Rtwo}{R\textsuperscript{2}}
\newcommand{\VLV}{V_{\mathrm{LV}}}
\newcommand{\PLV}{p_{\mathrm{LV}}}
\newcommand{\VminLV}{\vol_{\mathrm{LV}}^{\min}}
\newcommand{\VmaxLV}{\vol_{\mathrm{LV}}^{\max}}
\newcommand{\pminLV}{\prs_{\mathrm{LV}}^{\min}}
\newcommand{\pmaxLV}{\prs_{\mathrm{LV}}^{\max}}

\newcommand{\paramSpace} {\mathscr{P}}
\newcommand{\expected}{\mathbb{E}}
\newcommand{\variance}{\mathbb{V}\mathrm{ar}}
\newcommand{\param}  {\mathbf{p}}
\newcommand{\error}{\boldsymbol{\epsilon}}
\newcommand{\qoi}  {\mathbf{q}}
\newcommand{\qoiObs}  {\qoi_{\text{obs}}}
\newcommand{\forward}  {\mathcal{F}}
\newcommand{\forwardRed}  {\widetilde{\mathcal{F}}}
\newcommand{\piPrior}{\pi_{\text{prior}}}
\newcommand{\piPost}{\pi_{\text{post}}}
\newcommand{\NoiseCov}{\boldsymbol{\Sigma}}

\newcommand{\errorROM}{\error_{\text{ROM}}}
\newcommand{\NoiseCovROM}{\NoiseCov_{\text{ROM}}}
\newcommand{\errorEXP}{\error_{\text{exp}}}
\newcommand{\NoiseCovEXP}{\NoiseCov_{\text{exp}}}
\newcommand{\NoiseMagnEXP}{\sigma_{\text{exp}}}

\newcommand{\SobolFirst}[2]{S_{{#1}{#2}}}
\newcommand{\SobolTotal}[2]{S_{{#1}{#2}}^T}

\newcommand{\dynOut}{\mathbf{y}}
\newcommand{\dynState}{\mathbf{z}}
\newcommand{\dynPDEOut}{\mathbf{y}}
\newcommand{\dynPDEState}{\mathbf{z}}
\newcommand{\dynPDEOutVEC}{\mathbf{Y}}
\newcommand{\dynLat}{\mathbf{s}}
\newcommand{\dynInp}{\mathbf{u}}

\newcommand{\dynRHS}{f}
\newcommand{\dynOBS}{g}
\newcommand{\dynRHSop}{\mathcal{F}}
\newcommand{\dynOBSop}{\mathcal{G}}

\newcommand{\dynNumOut}{N_y}
\newcommand{\dynNumState}{N_z}
\newcommand{\dynNumLat}{N_s}
\newcommand{\dynNumInp}{N_u}

\newcommand{\dynNumTimes}{N_T}
\newcommand{\dynNumSamples}{N_S}
\newcommand{\dynNumObservation}{N_{\text{obs}}}

\newcommand{\dynNNdyn}{\NNgeneric_{\text{dyn}}}
\newcommand{\dynWdyn}{\mathbf{w}_{\text{dyn}}}

\newcommand{\dynNNenc}{\NNgeneric_{\text{enc}}}
\newcommand{\dynWenc}{\mathbf{w}_{\text{enc}}}
\newcommand{\dynNNdec}{\NNgeneric_{\text{dec}}}
\newcommand{\dynWdec}{\mathbf{w}_{\text{dec}}}

\newcommand{\dynPDEOutSpace}{\mathcal{Y}}
\newcommand{\dynPDEOutDiscretize}{\mathcal{A}}

\makeatletter
\newenvironment{mydescription}%
  {\list{}{\labelwidth\z@ \itemindent-\leftmargin
           }}%
  {\endlist}

\makeatother

\title{\papertitle}
\author{
    Alfio Quarteroni$^1$, 
    Paola Gervasio$^2$, 
    Francesco Regazzoni$^3$
}
\date{\footnotesize
	$^1$ Politecnico di Milano (Professor Emeritus), piazza Leonardo da Vinci, 32, Milan, 20133, Italy and\\
    Ecole Polytechnique Fédérale de Lausanne (Professor Emeritus), Station 8, Lausanne, 1015, Switzerland \\
    \texttt{alfio.quarteroni@polimi.it}\\[2ex]
	$^2$ DICATAM, University of Brescia, via Branze 43, Brescia, 25123, Italy\\
    \texttt{paola.gervasio@unibs.it}\\[2ex]
	$^3$ MOX, Department of Mathematics, Politecnico di Milano, piazza Leonardo da Vinci, 32, Milan, 20133, Italy\\
    \texttt{francesco.regazzoni@polimi.it}\\[2ex]
    }

\begin{document}
\maketitle

\begin{abstract}
Scientific Machine Learning (SciML) is a recently emerged research field which combines physics--based and data--driven models for the numerical approximation of differential problems.
Physics--based models rely on the physical understanding of the problem at hand, subsequent mathematical formulation, and numerical approximation. Data--driven models instead aim to extract relations between input and output data without arguing any causality principle underlining the available data distribution. In recent years, data--driven models have been rapidly developed and popularized. Such a diffusion has been triggered by a huge availability of data (the so--called big data), an increasingly cheap computing power, and the development of powerful machine learning algorithms.
SciML leverages the physical awareness of physics--based models and, at the same time, the efficiency of data--driven algorithms. With SciML, we can inject physics and mathematical knowledge into machine learning algorithms. Yet, we can rely on data--driven algorithms’ capability to discover complex and non--linear patterns from data and improve the descriptive capacity of physics--based models.
After recalling the mathematical foundations of digital modelling and machine learning algorithms, and presenting the most popular machine learning architectures, we discuss the great potential of a broad variety of SciML strategies in solving complex problems governed by partial differential equations.
Finally, we illustrate the successful application of SciML to the simulation of the human cardiac function, a field of significant socio--economic importance that poses numerous challenges on both the mathematical and computational fronts. The corresponding mathematical model is a complex system of non--linear ordinary and partial differential equations describing the electromechanics, valve dynamics, blood circulation, perfusion in the coronary tree, and torso potential. Despite the robustness and accuracy of physics--based models, certain aspects, such as unveiling constitutive laws for cardiac cells and myocardial material properties, as well as devising efficient reduced order models to dominate the extraordinary computational complexity, have been successfully tackled by leveraging data--driven models.
\renewcommand{\thefootnote}{}
\footnotetext{Published in \emph{Mathematical Models and Methods in Applied Sciences (M3AS)} \url{https://doi.org/10.1142/S0218202525500125}}
\end{abstract}


\pagebreak
\tableofcontents
\pagebreak

\section{Introduction}\label{sec:intro}

The twentieth century has been a pivotal era in scientific progress, driven by a phenomenological approach and focus on a physical understanding of phenomena. This included experimentation, formulating theories, obtaining numerical results through computer simulations, and validating mathematical and numerical models.

\begin{figure}[h!]
\begin{center}
\includegraphics[trim=0 2cm 0 1.cm,  width=0.9\textwidth]{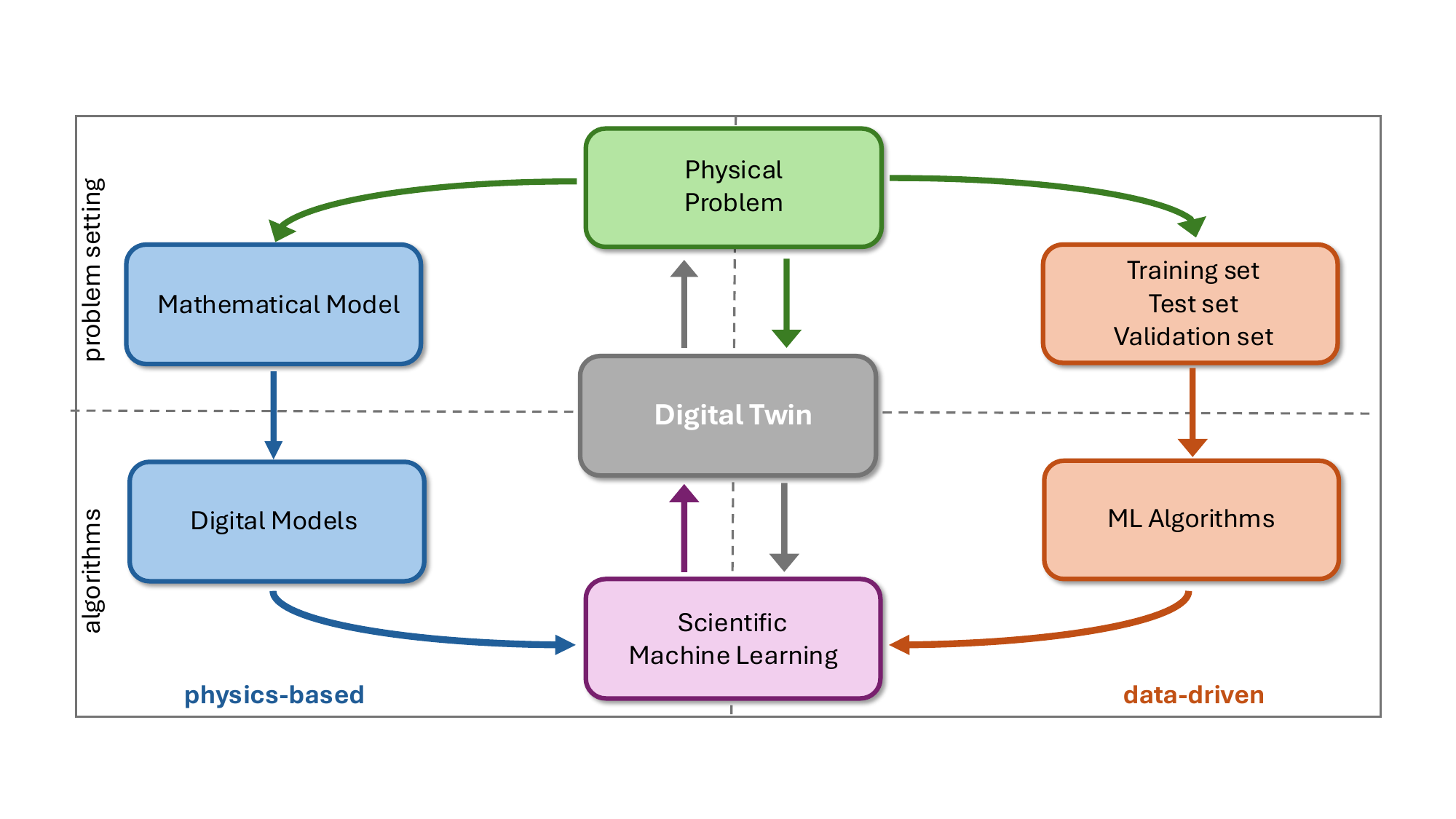}
\end{center}
\caption{The abstract framework.}
\label{fig:the_framework}
\end{figure}

When considering problems from basic sciences (physics, chemistry, biology, etc.) and applied sciences (engineering, medicine, society, etc.), solving them through \emph{physics--based models} consists first of all of formulating the fundamental principles underlying the problem at hand. Subsequently, these principles are translated into \emph{mathematical models} that then need to be analysed. Whenever the mathematical model is well posed, its solution exists and, under certain assumptions, is also unique and continuously dependent on the problem data. However, mathematical models of real--life problems are generally very complex and their solutions can seldom be achieved by paper--and--pen. Resorting to appropriate approximate models is therefore essential. This involves approximating the mathematical models using numerical methods, which in turn are translated into algorithms and then into computer programs. Ultimately, the final computational solution undergoes both verification and validation processes. The physics--based approach, sometimes also called \emph{digital} or \emph{computational modelling} (see the left column of Fig. \ref{fig:the_framework}), is the most widely used strategy suitable to solve complex problems.

On the other hand, \emph{Artificial Intelligence} has gained momentum over the last two decades. In particular, we can refer to algorithms empowered by data--trained \emph{artificial neural networks}. Three triggering factors have contributed to such a rapid development: \emph{(i)} relatively cheap computing power, offered especially by cloud services, and GPUs which can carry out very fast computations; \emph{(ii)} the availability of vast datasets from diverse sources, often referred to as \textit{big data}; \emph{(iii)} the development of powerful \emph{machine learning} algorithms for automatic learning. Among the latter, a prominent role is played by \emph{deep artificial neural networks} inspired by networks of biological neurons. 

As a consequence of this extraordinary development, machine learning algorithms can now be adapted to solve relevant and complex problems arising from both basic and applied sciences. Indeed, a new strategy is emerging as an alternative to physics--based modelling. Physical principles are being replaced by data generated from measurements, clinical images, or a series of experiments. In addition, computational models are being replaced by machine learning algorithms. 

More precisely, first, available data are partitioned into the training and testing sets. Then, artificial neural networks, which represent the core algorithms empowering machine learning (see the right column of Fig. \ref{fig:the_framework}), are built upon these training sets. This strategy, which we have called for the sake of brevity \emph{data--driven}, is meant to build the \emph{input--output function} that maps data to solutions. Data--driven algorithms aim to extract relations between input and output data without inferring any causality principle whatsoever.

These approaches (physics--based and data--driven) could run independently and alternatively, although with a different degree of accuracy, interpretability and trustworthiness.
However, more interestingly, they can be combined jointly, giving rise to a more effective strategy which has recently been called  \emph{Scientific Machine Learning} (SciML). \emph{
Scientific Machine Learning is an interdisciplinary field empowered by the
synergy of physics-based computational models with machine-learning algorithms
for scientific and engineering applications.}
 SciML is a flexible framework that adapts to various goals and requirements, taking different forms as needed. Digital models can be exploited to make machine learning algorithms more accurate, conversely, machine learning algorithms can enhance the computational efficiency of digital models. For instance, we can use digital models to regularize the loss function in artificial neural networks (to avoid overfitting), to augment data if missing or scarce in the training process, or even to learn hidden mathematical operators (sometimes even in symbolic form). On the other hand, machine learning can be used to improve digital models, for instance, to find constitutive laws in the case of new materials, accelerate scientific computing algorithms, estimate and calibrate parameters, solve optimal control problems depending on lots of parameters, or perform sensitivity analysis. 

Generally speaking, SciML allows us to inject physics and mathematical knowledge into machine learning algorithms. At the same time, we can exploit data--driven algorithms' ability to uncover complex patterns from data, to enhance the descriptive capability of physics--based models. In this way, we can open up new paths to leverage the causality principle of physics--based models and the efficiency of data--science algorithms. SciML is also prodromic to the realization of \emph{digital twins} when a real--time dialogue between the physical asset and its digital model is essential.

In this paper, we discuss the great potential of SciML in solving complex problems governed by Partial Differential Equations (PDEs) in general. Before entering this discussion, we provide some basic material to illustrate the mathematical foundation of machine learning algorithms from a mathematical perspective.

Additionally, in the final part of the paper, the previous discussion is adapted to an application of considerable relevance which is the numerical simulation of the \emph{human cardiac function}. The heart is a very sophisticated and complex machine that is the engine of life. Its functioning is the harmonious result of the interaction of different processes. These include the generation and propagation of electrical signals in cardiac tissue, the contraction and relaxation of the myocardium, the motion of blood in the ventricles and atria, the opening and closing of heart valves which are perfectly synchronized with heart mechanics and blood dynamics, and blood supply to the myocardium through perfusion of the coronary arteries. All these processes obey fundamental and general physical principles such as Newton's laws, thermodynamics, chemical kinetics, Navier-Stokes' equations, and Maxwell's equations. The mathematical model that formalizes this complex system is a set of non--linear ordinary and partial differential equations that, once solved, can simulate the human heart function both in physiological and pathological conditions. It is worth noting that solving this highly complex and huge system of equations cannot be achieved simply by pen--and--paper, but requires supercomputers to step into action and play an irreplaceable role.

The \emph{Integrated Heart Model} (IHM) \cite{Fedele2023-comprehensive} allows us to return helpful qualitative and quantitative information to cardiologists and cardiac surgeons to understand, diagnose and treat heart disease. The great potential of mathematical models consists in being universally valid and independent of the particular context to which they are applied. More specifically, the IHM can be used to simulate the cardiac function of any individual, regardless of age, weight, or clinical history (to give some examples), relegating all these peculiarities to the data. 

An essential feature is that data are needed to feed models. Those required by the IHM can be the geometric shape of the patient's heart, which can be obtained using medical images such as computerized tomography or magnetic resonance imaging. Furthermore, the initial and boundary conditions for the differential equations are necessary data, which, however, are seldom available as they would require invasive clinical examinations on patients.  Other data include the parameters characterizing material properties, which are very difficult to recover accurately and even when they are available are often affected by noise. It is worth noting that the greater the complexity of the model, the more parameters are needed to accurately solve the model itself. To compensate for data scarcity and inaccuracy, we can resort to parameter identification, variational techniques for data assimilation, and uncertainty quantification. In particular, the latter is a very powerful tool to understand how uncertainty propagates along the computational process from input to output.

However, despite the robustness and accuracy of physics--based models, certain aspects, such as constitutive laws for cardiac cells and myocardial material properties, remain poorly understood. Data--driven algorithms can bridge these gaps, enabling more accurate and efficient modelling of the heart's intricate processes. 

We believe that the IHM is an ideal and particularly challenging test bed for the development and evaluation of several approaches falling in the general field of SciML.

The outline of the paper is as follows. In Section~\ref{sec:DM}, the basic concepts of digital models are presented. Section~\ref{sec:DDM} introduces the mathematics behind machine learning algorithms, specifically those grounded on artificial neural networks. Section~\ref{sec:SML} presents the basics of SciML and the most popular and emerging approaches in the field. In Section~\ref{sec:SML-IHM}, we briefly illustrate the Integrated Heart Model and how SciML has been called into play. Finally, in Section~\ref{sec:conclusions}, we wrap up the paper with some concluding thoughts and reflections on both Machine Learning and SciML.
\section{Digital models}\label{sec:DM}

In the field of applied sciences, Digital Models (DM) enable solving complex problems that can be formulated using the language of mathematics and which are rarely solvable analytically. Indeed, in general, the mathematical solution of complex problems can seldom be expressed explicitly: this is e.g. the case of non--linear equations or systems (algebraic, differential or integral) for which no solution formulas are known. In other cases, the explicit form, although known, cannot be practically used to determine quantitative values of the solution itself, either because it would require a prohibitive amount of operations, or because it is in turn expressed through other mathematical entities (such as integrals or series developments) for which a quantitative evaluation would be difficult to achieve.

The resolution of problems of applied interest involves several phases:
\begin{enumerate}
\renewcommand\labelenumi{(\alph{enumi})} 
\item the qualitative and quantitative understanding of the problem; 
\item the modelling through mathematical equations, which can be algebraic, functional, differential, or integral, this is the so--called mathematical model; 
\item the identification of numerical analysis methods suitable for approximating this mathematical model; 
\item the implementation of these numerical approximation methods on a computer.
\end{enumerate}

DMs correspond to phases (c) and (d) of the process, and the choice of these approaches requires a thorough understanding of the qualities of the mathematical model's solution and a clear vision of the original problem's phenomenology. Therefore, the use of the term ``digital models'' is justified as an alternative to ``digital methods''. Digital modelling represents an interdisciplinary field of research that combines mathematics, computer science, and applied sciences. Its progress has been remarkable in recent decades thanks to the advanced development of computers and algorithms, along with an increased awareness of the crucial role of scientific computing in simulating complex problems. These problems range from scientific and theoretical domains to industrial, environmental, and social spheres.

We denote by $u_{ph}$ the solution to the original problem, here called the physical problem for brevity, even though problems from all fundamental and applied sciences (not just physics) can be considered, such as chemistry, biological and medical sciences, engineering, economics, etc. 

\subsection{Mathematical models}
Mathematical models translate the first principles of ``physics'' into mathematical terms, which often describe its fundamental characteristics, such as conservation laws (of mass, momentum, energy). These models, expressed in the form of equations, encode an extraordinary and unique body of knowledge, generated by luminaries of the past such as Kepler, Galileo, Newton, Pascal, Bernoulli, Laplace, Einstein, Maxwell, and Schrödinger, to name just a few \cite{Stewart-2012}.

The mathematical model that describes the physical problem can be expressed in the compact form 
\begin{equation}\label{eq:abstract-problem}
F(u,d)=0,
\end{equation}
 where it is assumed that $d$ describes the set of data,  $u$ the solution, and 
$F$ is the functional relation (algebraic, integral, differential, to name a few examples) that links data and solution. Various and diversified mathematical problems can be formulated in this abstract way. Consider, for instance, the definite integration of a function, the problem of finding the roots of an algebraic, trigonometric, or exponential equation, Cauchy problems for ordinary differential equations, boundary value problems for partial differential equations, and linear and non--linear algebraic systems. See, for example, \cite{Atkinson-1989,Eriksson-1996,Golub-1989,Lambert-1991,Quarteroni-NMDP-2017}.

\begin{figure}
\begin{center}
\includegraphics[trim=0 2cm 0 0, width=0.9\textwidth]{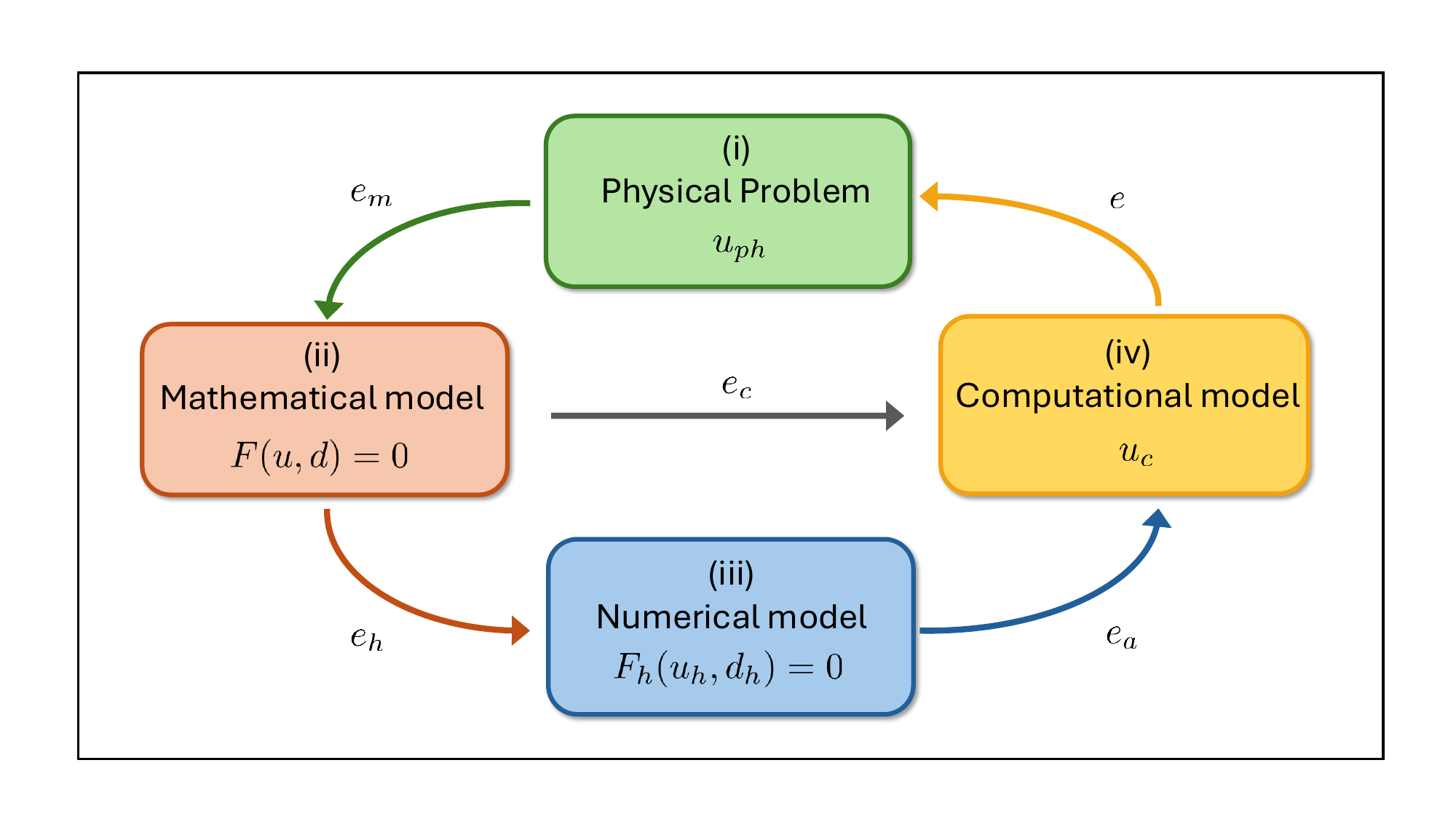}
\end{center}
\caption{Problems, solutions, and errors in digital models}
\label{fig:problemi-soluzioni-errori}
\end{figure}

To fix ideas, an instance of the abstract equation (\ref{eq:abstract-problem}) could be the following Initial Boundary Value Problem (IBVP): let 
$\Omega\subset {\mathbb R}^d$ (with $d=1,\,2,\,3$ denoting the space dimension) an open and bounded set and $(t_0,T)\subset{\mathbb R}$ 
a bounded interval, we look for the function $u=u
(\textbf{x},t):\Omega\times(t_0,T)\to \mathbb{R}^m$ (with $m\geq 1$ denoting the 
number of components of the solution $u$) satisfying 
\begin{eqnarray}\label{eq:IBVP}
\left\{\begin{array}{ll}
\displaystyle\frac{\partial u}{\partial t}+{\cal L}u+{\cal N}(u)=f & \mbox{ in }\Omega\times(t_0,T)\\[2mm]
\mbox{boundary conditions} & \mbox{ on }\partial\Omega\times(t_0,T)\\[2mm]
 u=u_0 & \mbox{ in }\Omega\times\{t_0\},
 \end{array}\right.
\end{eqnarray}
where ${\cal L}$ is a linear differential operator and ${\cal N}$ a non--linear one, both only involving derivatives with respect to space variables, while $f$, $u_0$, and other parameters (or functions) hidden in the definition of the operators ${\cal L}$ and ${\cal N}$ constitute given data.

The IBVP (\ref{eq:IBVP}) can represent the abstract form of several real--life problems. For instance, if we want to simulate the electrophysiology of the heart, i.e., the
result of chemical and electrical processes taking place
from subcellular to the whole organ scale, then $u$ will be a multicomponent function including the transmembrane potential of the myocytes (the cells of the heart) through the myocardium as well as the concentration of some ionic species, such as sodium, potassium, and calcium inside the myocytes.
If instead, we want to simulate the deformations of the cardiac muscle during a heartbeat, $u$ will be the displacement of each cell of the myocardium induced by a force field, in its turn generated by the transmembrane potential.
We could also be interested in simulating both the electrophysiology and the heart mechanics, in this case (\ref{eq:IBVP}) should take into account both the physics and we would speak about a \emph{multiphysics coupled problem}.
For all these examples, the data set $d$ includes suitable external forces, boundary conditions associated with the differential operators and initial conditions for the time--dependent unknowns.

Many other mathematical models are of interest in computational modelling and can be reformulated in the abstract framework (\ref{eq:abstract-problem}). Consider, for example, the problem of representing curves and surfaces, the numerical treatment of images and signals, the minimization of free or constrained functionals, the numerical regularization and resolution of inverse problems, integro--differential problems, as well as the control of differential, integral, or stochastic systems. The relevance of these models in the analysis and simulation of mechanical processes (e.g., in fluid dynamics, structural analysis, bioengineering), physical processes (in neutronics, microelectronics, wave propagation), economic processes (in deterministic or statistical systems for describing macroeconomic systems, or in investment processes with options), etc., is significant and destined to grow over time.

The mathematical model aims to link quantities of physical interest through mathematical relations, often simplified compared to the complexity of the original problem. 
The solution $u$ of the mathematical model will consequently differ from the solution $u_{ph}$ of the physical problem under consideration, and the distance between these to solutions $e_m = \|u_{ph} - u\|$, measured in a suitable norm, will account for both the significance (and reasonableness) of the simplifying assumptions made and the accuracy with which the real physical quantities have been represented by the set $d$ of the model's data. The error $e_m$ is therefore intrinsic to the process and it is named \emph{model error} (see Fig. \ref{fig:problemi-soluzioni-errori}).

\subsection{Numerical models} \label{sec:numerical_models}
The first phase of digital modelling consists of approximating the mathematical model (\ref{eq:abstract-problem}) with a numerical model 
\begin{equation}\label{eq:abstract-numerical-problem}
F_h(u_h,d_h)=0,
\end{equation}
which depends on a parameter $h$ whose meaning varies according to the family of problems considered, but which can always be thought of as referring to the richness of the numerical solution's representation: the smaller $h$, the more representative the numerical solution. 
In this context, $ d_h $ is an approximation 
of $d$, $F_h$ a discretisation of the functional relation $F$, and $u_h$ the numerical solution, approximation of $u$. 

Why do we need to discretise a mathematical problem? 
If we refer to the Initial Boundary Value Problem (IBVP) (\ref{eq:IBVP}), the fundamental reason is that its exact solution $u$ belongs to an infinite dimensional space $V$ and digital computers can approximate such a solution only by functions belonging to suitable finite--dimensional spaces. Indeed, because a digital computer is not able to compute exact derivatives and integrals, and everything that requires passing to a limit, we find ourselves having to approximate these operations.
 Typically, space and time derivatives are approximated separately, although simultaneous space--time discretisation methods could be considered as well.

 \medskip
 \emph{Space approximation.} Space approximation methods aim to approximate space derivatives and reduce the IBVP into a set of Ordinary Differential Equations (ODE) (or Initial Value Problem -- IVP). Many different paradigms can be adopted to achieve this goal, one of them is the Galerkin method \cite{Quarteroni-Valli-1994} which consists of projecting the weak form of the IBVP at any time $t$ on a finite--dimensional subspace $V_h$ of the functional space $V$ which the weak exact solution belongs to. 
 
 If, for instance, we choose the Galerkin Finite Element (FE) method \cite{Quarteroni-Valli-1994,Ern-Guermond,Brenner-Scott}, the first step consists of partitioning the space domain $\Omega$ into elements or cells that are triangles or quadrilaterals when $d=2$, and tetrahedra or hexahedra when $d=3$, and whose characteristic size is $h$ (the so--called mesh--size of the partition). The finite--dimensional space $V_h$ is made of piecewise polynomial functions defined on such a partition and we denote by $N_h$ its finite dimension.
 Denoting by $\{\varphi_j, j=1,\ldots, N_h\}$ a suitable basis of the space $V_h$, for any $t\in (t_0,T)$ we look for a function $u_h(\textbf{x},t)\in V_h$ (approximating $u(\textbf{x},t)$) that takes the form
 \begin{equation}\label{eq:uh-expansion}
     u_h(\textbf{x},t)=\sum_{j=1}^{N_h}u_{h,j}(t)\varphi_j(\textbf{x}).
 \end{equation}
 For any $t\in(t_0,T)$, let $\textbf{u}_h(t)=[u_{h,1}(t),\ldots, u_{h,N_h}(t)]^t$ be an array in ${\mathbb R}^{N_h}$ so that, after taking into account the expansion (\ref{eq:uh-expansion}), we can write the \emph{semidiscrete problem}
 \begin{eqnarray}\label{eq:IBVP-semidiscrete}
\left\{\begin{array}{ll}
\displaystyle M \frac{d\textbf{u}_h}{d t}(t)+L\textbf{u}_h(t)+\textbf{N}(\textbf{u}_h(t))=\textbf{f}(t) & \mbox{ in } (t_0,T)\\[2mm]
 \textbf{u}_h(t_0)=\textbf{u}_{0,h},
 \end{array}\right.
\end{eqnarray}
where $L\in \mathbb{R}^{N_h\times N_h}$ and $\textbf{N}: \mathbb{R}^{N_h}\to  \mathbb{R}^{N_h}$ are suitable linear and non--linear operators related to the operators ${\cal L}$ and ${\cal N}$ introduced in (\ref{eq:IBVP}), respectively, and arising from the Galerkin projection, while 
$M\in \mathbb{R}^{N_h\times N_h}$ such that $M_{ij}=\int_\Omega \varphi_j\, \varphi_i$ is the so--called mass matrix, and $\textbf{u}_{0,h}$ is the array containing the values of $\textbf{u}_{0}$ at the mesh nodes.

Notice that, when the unknown $u$ is a vector field and/or (\ref{eq:IBVP}) collects the solutions of more than one differential equation, the vector function $\textbf{u}_h$ and the formulation (\ref{eq:IBVP-semidiscrete}) will have to take into account all the components of the unknown solution and all the equations of the model.

 \medskip
\emph{Time discretisation.} To discretise the time derivative of the problem (\ref{eq:IBVP-semidiscrete}) we need to introduce a partition of the time interval $[t_0,T]$ into $N_t$ sub--intervals of size $\Delta t$ (for sake of clearness we consider uniform partitions) that induces a set of nodes $t_n=t_0+n\Delta t$ for $n=0,\ldots, N_t$ with the property that $t_{N_t}=T$. 
discretising (\ref{eq:IBVP-semidiscrete}) by classical finite differences schemes means approximating the time derivative by finite ratios, so that, at each $t_n$, 
we look for the real values $\textsf{u}_{h,j}^n$ approximating $u_{h,j}(t_n)$ (for any $j=1,\ldots, N_h$). 
Then we set $\textbf{\textsf u}^n=[\textsf{u}_{h,1}^n,\textsf{u}_{h,2}^n,\ldots,\textsf{u}_{h,N_h}^n]^t$ for any $n=1,\ldots, N_t$. Other time--discretisation methods can be obtained by approximating the integrals (instead of derivatives) arising from the integral form of an ordinary differential equation \cite{Lambert-1991}.

For the sake of simplicity, let us consider the Backward Euler scheme \cite{Quarteroni-CS-2014} and apply it to the problem (\ref{eq:IBVP-semidiscrete}). For $n=0,\ldots, N_t-1$, we look for $\textbf{\textsf u}^{n+1}$ satisfying the \emph{fully discrete system}
\begin{eqnarray}\label{eq:IBVP-discrete}
\left\{\begin{array}{ll}
\displaystyle\frac{1}{\Delta t} M \textbf{\textsf u}^{n+1} +L\textbf{\textsf u}^{n+1}+\textbf{N}(\textbf{\textsf u}^{n+1})=\displaystyle\frac{1}{\Delta t} M \textbf{\textsf u}^{n}+ \textbf{f}(t^{n+1}), & \quad n=0,\ldots,N_t-1\\[2mm]
 \textbf{\textsf u}^0=\textbf{u}_{0,h}.
 \end{array}\right.    
\end{eqnarray}
 The choice of the mesh--size $h$ and time--step--size $\Delta t$ strongly affects the accuracy of the approximation as well as the computational effort: typically the smaller $h$ and $\Delta t$, the more accurate the approximation, but the larger the dimension $N_h$ of the space $V_h$ and the number $N_t$ of time--steps, i.e., the larger the computational cost in computing $\{\textbf{\textsf u}^1,\ldots,\textbf{\textsf u}^{N_t}\} $.

Problems (\ref{eq:IBVP-discrete}) are in general non--linear algebraic systems that can be solved, e.g., by Newton or other iterative methods \cite{Nocedal,Quarteroni-CS-2014}. At each iteration of the Newton method, a linear system of size $N_h$ (whose matrix is given by the sum of the Jacobian of the non--linear operator, or a suitable approximation of it, the linear operator $L$, and the mass matrix divided by $\Delta t$) has to be solved and, typically, the larger $N_h$, the larger the condition number of the system. 
It is extremely common to find ourselves having to solve ill--conditioned systems that may generate inaccurate solutions and, if solved by iterative methods, can require many iterations to reduce the residual even of few orders. 

To overcome these shortbreaking, suitable \emph{preconditioners} and/or scaling matrices must be applied to the linear system. We say that a preconditioner is optimal for a given linear system if the iterative condition number of the preconditioned matrix (that is the ratio between the maximum and minimum moduli of the eigenvalues of the matrix) is quite smaller than the original matrix one and it is independent of both the discretisation parameters. 
A unique optimal preconditioner for all classes of linear systems is not available, although a great effort has been reserved in recent years for the design of preconditioners for matrices arising from the approximation of PDEs \cite{Wathen}.
Another desirable property for a good preconditioner is \emph{scalability}, i.e. the preconditioner should offer excellent performance even on parallel computers and allow solving linear systems in fewer computation times when more cores are used; ideally, the computation time should decrease proportionally with the number of cores. 
We cite, among others, Algebraic and Geometric Multigrid methods \cite{Adams2003,Trottenberg2001,Sundar2015,Xu2017,kronbichler_multigrid,ASGDQ2023} and Additive Schwarz preconditioners \cite{Toselli-Widlund,Xu-Zou-1998,Barker-2010}.

After solving (\ref{eq:IBVP-discrete}), the numerical solution $u_h^n\in V_h$  at time $t^n$ will be reconstructed at each point $\textbf{x}$ of the computational domain $\Omega$ by the expansion
\begin{equation}\label{eq:uN}
    u_h^n(\textbf{x})=\sum_{j=1}^{N_h} \textsf{u}_{h,j}^n \varphi_j(\textbf{x}).
\end{equation}
When (\ref{eq:IBVP}) models multiphysics systems, the size of (\ref{eq:IBVP-semidiscrete}), that is, its number of equations (and unknowns), could be huge, making its resolution by computers in acceptable times practically impossible.
In such cases, the \emph{monolithic formulation}, i.e. the one expressed by a unique system like (\ref{eq:IBVP}) that aims at solving all the equations of the multiphysics model simultaneously \cite{Gerbi-2019-monolithic}, can be replaced by a system of ``simpler'' IBVPs that interact one each other, obtained by applying operator splitting methods and/or domain decomposition techniques \cite{Dede-2020-segregated}.

\emph{Operator splitting} methods allow us to advance in time by taking sub--steps of the complete time--step. At each sub--step only a part of the complex system is solved implicitly, while the others are treated explicitly so that we speak about staggered algorithms. Typically, the staggered parts correspond to the different physics making the multiphysics problem. For instance, if the multiphysics model includes both electrophysiology and mechanics that strongly interact, at each time step, we can first advance the mechanics, and then the electrophysiology.

Operator splitting methods allow us to reduce the computational complexity of the monolithic system, their implementation is often easy, and they respond to the non--intrusive requirement, that is the possibility to use yet available solvers specific for the single sub--problems only with a little effort. 
Operator splitting methods are also very helpful in solving multiphysics problems characterized by different time and spatial scales, by allowing us to choose different finite element spaces and ad--hoc discretisation in both time and space for the specific scale.
However, these methods may introduce errors due to the staggered advance of the sub--problems and be less accurate than the monolithic approach.

\emph{Domain decomposition} techniques allow us to solve multiphysics coupled problems in which different physics govern the phenomena in distinct parts of the computational domain \cite{Quarteroni-Valli-ddm, Toselli-Widlund}. An example is given by the Fluid--Structure--Interaction (FSI) problem modelling the interplay between fluid dynamics (blood flow in a heart chamber) and structural mechanics (displacement of the myocardium). Domain decomposition techniques let us work with non--conforming discretisation at the interfaces between adjacent subdomains \cite{Deparis-2016-fsi,Deparis-2016-internodes, Bernardi-1993-mortar,Belgacem-1999,Wohlmuth-2000-Mortar} and they feature the same computational advantages as operator--splitting methods when the monolithic approach is prohibitive.

\emph{Convergence.}
The fundamental property 
that we ask a numerical model to satisfy is \emph{convergence}, meaning that (in an 
appropriate metric) $\lim_{h \to 0} u_h = u$. A necessary 
condition for this to occur is that $F_h$ (that, in the example below, takes into account both space and time discretisation) correctly translates the 
functional law $F$, that is, $F_h(u, d) \to 0$ as $h \to 0$, with $d$ being an admissible datum for the mathematical model and $u$ the corresponding solution.
This property, known as \emph{consistency}, is not sufficient to guarantee the convergence of the numerical solution to the exact one: for this purpose, the numerical model must also be stable. To define \emph{stability}, we recall that a mathematical model is said well--posed if it admits a unique solution and if small perturbations in the data result in controllable variations in the solution. The same property, rewritten for the numerical problem, will ensure its stability, provided that the ratio between the distance of the solutions and that of the data, in suitable norms, can be uniformly bounded, i.e., with a constant that does not depend on the discretisation parameters (the mesh--size $h$ and the time--step $\Delta t$).

The error introduced by approximating the mathematical solution $u$ by the numerical one $u_h$ is the \emph{numerical error} $e_h=\|u-u_h\|$ and it takes into account all the discretisation processes used in designing our numerical model (see Figure~\ref{fig:problemi-soluzioni-errori}).
The ultimate goal of analysing the error $e_h$ of the numerical model is to demonstrate that it tends to zero as the discretisation parameters $h$ and $\Delta t$ tend to zero. For instance, in the case that the differential problem is time--independent (so that no discretisation in time is required), one seeks estimates of the type:
\begin{equation}\label{eq:numerical-error}
\|u-u_h\|\leq C(u,d)h^{p}
\end{equation}
for a suitable norm depending on the specific problem. The exponent $p$ (i.e., the infinitesimal order of the error with respect to $h$) is named \emph{convergence order of the numerical model} with respect to $h$. Similar arguments apply to time discretisation methods.

The demonstration technique that allows for obtaining error estimates like (\ref{eq:numerical-error}) is based on the analysis of the well--posedness of the mathematical model and the stability of the numerical model. This is called \emph{a--priori analysis}, as it can be performed before actually implementing the numerical model. The constant $C(u,d)$ depends on the problem's data, as well as on $u$ and its derivatives (assuming they exist), up to an appropriate order, dependent on $p $. Quantifying this constant numerically is generally very difficult. Consequently, the above estimate is incomplete: it allows for predicting the error decay trend for vanishing $h$ but does not provide an accurate numerical quantification for a fixed $h$.

It is the task of \emph{a--posteriori analysis} to fill this gap. Based on the solution $u_h$ computed in correspondence with a certain value of the discretisation parameter $h$ and the value of the residual $\|r_h\|=\|F(u_h,d)\|$ (that is an indirect measure of the deviation of $u_h$ from $u$), a--posteriori analysis allows us to establish whether the error is or not below a pre--established tolerance. If this is not true, an adaptive algorithm indicates how to improve the discretisation (providing a criterion for modifying $u_h$) to guarantee a reduction in the error. See, for example, \cite{Quarteroni-NumerMath-2007}.

\medskip
Almost always, mathematical models require the knowledge of physical parameters that strongly affect the solution, but which are difficult or even impossible to observe. For instance, material parameters for the Integrated Heart Model (such as electrical conductivity, or elastic properties of the myocardium) are not easy to obtain for speciﬁc patients due to the difficulty of producing ad hoc measurements.

\emph{Parameters identification} is a critical process aimed at determining the values of such parameters to improve the agreement between the mathematical model output and observed data. Under the hypothesis that the observed data are available and ``exact'', this process requires solving \emph{inverse problems}, namely problems of which we know the output but not some inputs (parameters), by leveraging optimal control theory and numerical optimization techniques \cite{Lions-1971-OCS,Manzoni-Quarteroni-Salsa-book}. 
But very often, observed data are not ``exact'', instead, they are affected by experimental errors (e.g., when a quantity is not correctly measured) or it is not possible to measure some values at all. In this case, we speak about ``Epistemic uncertainty''. Other times, data may depend on the variability between individuals (think, e.g., to patient--specific measurements) or the stochasticity of the considered samples. In this latter case, we speak of ``Aleatory uncertainty''. 

\emph{Uncertainty Quantification} (UQ) mathematical techniques aim to characterize and reduce Epistemic uncertainties in model predictions through two processes: \emph{forward UQ}, which deals with the propagation of the parameters uncertainty on the outputs of the model, and \emph{backward UQ}, which studies how the measurement errors on the outputs aﬀect the estimation of the parameters. 
In particular, in the broad field of UQ, \emph{Sensitivity Analysis} investigates how variations in model inputs affect the outputs. It helps to identify which inputs are most important in determining the output behaviour and is commonly used in model validation, calibration, and optimization to understand and prioritize which inputs to focus on.

\medskip
The parameter identification process is very computationally expensive because it requires the solution of many instances of the same mathematical model which differ one each other for the setting of some parameters. Such an issue is typically faced by \emph{Reduced Order Models} (ROM). Let us denote by $\mu\in P\subset {\mathbb R}^p$ a set of parameters on which the solution of the IVBP (\ref{eq:IBVP}) depends and, after space discretisation, rewrite the semidiscrete form (\ref{eq:IBVP-semidiscrete}) by making the dependence on $\mu$ explicit: given $\mu\in P$, for any $t\in (t_0,T)$ look for the solution $\textbf{u}_h(t;\mu)$ of
\begin{eqnarray}\label{eq:IBVP-ROM-semidiscrete}
\left\{\begin{array}{ll}
\displaystyle M(\mu) \frac{d\textbf{u}_h}{d t}(t;\mu)+L(\mu)\textbf{u}_h(t;\mu)+\textbf{N}(\textbf{u}_h(t,\mu);\mu)=
\textbf{f}(t;\mu) & \mbox{ in } (t_0,T)\\[2mm]
 \textbf{u}_h(t_0;\mu)=\textbf{u}_{0,h}(\mu).
 \end{array}\right.
\end{eqnarray}
Because all the solutions of (\ref{eq:IBVP-ROM-semidiscrete}) corresponding to any possible choice of the parameters $\mu\in P$ belong to a manifold 
$M_h=\{\textbf{u}_h(t;\mu):\ t\in(t_0,T), \mu\in P\}$, the idea of ROM consists in: $(i)$ evaluating a database of discrete solutions in $M_h$, named snapshots, corresponding to different values of the parameters $\mu$ in $P$; $(ii)$ selecting a subset of $N$ linearly independent snapshots that span a linear subspace $V_N$ of $M_h$ (the aim is to take $N\ll N_h$ where $N_h$ is the dimension of the finite element space $V_h$); $(iii)$ given a new parameter $\overline\mu$, computing the corresponding solution $\textbf{u}(t;\overline\mu)$ by the Galerkin projection of the equation (\ref{eq:IBVP-ROM-semidiscrete}) onto the subspace $V_N$. 

If we denote by $W_N$ the matrix whose columns contain the degrees of freedom of the finite element basis functions computed at step $(ii)$ and set
$M_N(\overline\mu)=W_N^T M(\overline\mu) W_N$, $L_N(\overline\mu)=W_N^T L(\overline\mu) W_N$, and 
$\textbf{f}_N(\overline\mu)=W_N^T \textbf{f}(t,\overline\mu)$, the step (iii) may be written as follows: for any $t\in (t_0,T)$ find the solution $\textbf{u}_N(t;\overline{\mu})\in V_N$ of
\begin{eqnarray}\label{eq:IBVP-ROM1-semidiscrete}
\left\{\begin{array}{ll}
\displaystyle M_N(\overline\mu) \frac{d\textbf{u}_N}{d t}(t;\overline\mu)&+L_N(\overline\mu)\textbf{u}_N(t;\overline\mu)\\[2mm]
&+ W_N^T\textbf{N}(W_N\textbf{u}_N(t,\overline\mu);\overline\mu)=\textbf{f}_N(t;\overline\mu) \qquad \mbox{ in } (t_0,T)\\[2mm]
 \multicolumn{2}{l}{\textbf{u}_N(t_0;\overline\mu)=\textbf{u}_0(\overline\mu).}
 \end{array}\right.
\end{eqnarray}
Different strategies can be adopted for step (ii), the most widely used are greedy approaches and Proper Orthogonal Decomposition (POD) methods  \cite{Cohen-2015, Hesthaven-2016, Quarteroni-Manzoni-Negri}.
More precisely, when the problem is stationary, the basis functions are the first $N$ snapshots selected by the greedy algorithm or the first $N$ singular vectors of the snapshots matrix whether POD is used. When the problem is time--dependent, the parameter space can still be sampled by one of the
two techniques mentioned, whereas POD is usually exploited to reduce trajectories of the system over the time interval. 
Typically, step $(i)$ requires the solution of many problems like (\ref{eq:IBVP-semidiscrete}) in the original finite element space and it is very expensive, however, if the differential problem is linear and $M(\mu)$, $L(\mu)$, and $\textbf{F}(t;\mu)$ depend on $\mu$ in affine way, the computation is carried out once for all (offline), jointly with the step $(ii)$. On the contrary, step $(iii)$ is generally very cheap (its dimension $N$ is very small compared to the dimension $N_h$ of the finite element space $V_h$) and is carried out (online) to compute the solution corresponding to each new parameter $\overline\mu$ that comes into play, e.g., in the uncertainty quantification process. Whenever the problem is non--linear or the dependence on the parameters is not affine, the whole process is more involved and hyper--reduction strategies must be adopted \cite{Farhat-2019}.

\null
\emph{Computer implementation.} The numerical model must be implemented on computers using appropriate algorithms. Consequently, instead of the numerical solution $u_h$, we end up with the computational solution $u_c$, and
new errors $e_a$ are generated attributable to the finite arithmetic used by computers to represent real numbers and perform elementary algebraic operations.
The sum of the numerical error $e_h$ with the finite arithmetic error $e_a$ provides the computational error $e_c$, which represents the distance between the computed solution $u_c$ and $u$. The overall error $e$, which measures the difference between the physical solution $u_{ph}$ and the computed solution $u_c$, will result from the combination of the two errors, $e_m$ and $e_c $.\\
Digital modelling should guarantee that the error $e_c$ is small and controllable (\emph{reliability}) with the least possible 
computational cost (\emph{efficiency}). Computational costs refer to 
the amount of resources (computation time, memory usage) needed to 
compute the solution $u_c$. Reliability is a crucial requirement of a 
computational model: the analysis aims to find estimates in a suitable norm of the error 
$e_c$ as a function of the problem's data and the discretisation 
parameters, ensuring that it remains below a certain predetermined 
threshold. To this end, it is reasonable to employ adaptive algorithms 
that use a feedback procedure, based on the results already obtained, 
to adjust the discretisation parameters and improve the quality of the 
solution. Numerical methods that implement an adaptive strategy for 
error control are known as \emph{adaptive computational models}.\\
Under the optimal assumption that the algorithmic error $e_a$ is controllable within the desired tolerance, the computational error $e_c$ will behave like $e_h$ and will tend to zero with respect to $h$ (see Figure~\ref{fig:problemi-soluzioni-errori}).

\section{Data--driven models}\label{sec:DDM}

Data--driven models exploit data from measurements, imaging, or whatever experiments, to obtain a solution to the problem through machine learning algorithms, typically artificial neural networks.

Data--driven discoveries are not a prerogative of the 21\textsuperscript{st} century, just think of Kepler and Newton who derived their milestone laws to predict planets' orbits starting from astronomical data and combining them with the analytical approach \cite{Brunton-Kutz-2022}. The novelty of modern data--driven models consists of employing Machine Learning (ML), instead of analytical approaches, to draw conclusions.

In introducing Digital Models (DM), we highlighted that they are based on computational sciences, and more precisely: $(i)$ on the knowledge of the theory (we think to equations and more in general mathematical models), $(ii)$ the availability of experiments needed to validate the models, and $(iii)$ the numerical simulations. Today, \emph{Artificial Intelligence} (AI) could be reasonably regarded as the fourth pillar of the setting. 

\subsection{Artificial Intelligence}
Many different definitions of AI have been given so far, starting with the very first by John McCarthy in 1955, who defined AI as “the science and engineering of making intelligent machines”. Successively, many definitions have followed the evolution of AI itself, up to the definition most shared today according to which AI is the set of techniques enabling computers to mimic human intelligence. By intelligence, we mean the ability to acquire knowledge and apply it in several contexts: reasoning, learning, problem--solving, perception, linguistic intelligence, social intelligence, decision--making, just to mention the perhaps most relevant.

In the broad realm of AI we find Large Language Models, Natural Language Understanding, Text and Image generation, Dialog Systems, Robotics, and more, which can be applied across a wide range of fields including healthcare, environment, energy, manufacturing, agriculture, finance, retail, transportation, entertainment, security, and many others. 

Over the last thirty years, we have witnessed an exponentially growing interest in artificial intelligence even from scientists from outside the computer science field.
Three main factors stand behind this explosion:
\begin{enumerate}[label=(\roman*)]
    \item the availability of increasingly powerful computers equipped with numerous central processing units (CPUs), graphics processing units (GPUs), and Tensor Processing Units (TPU) capable of executing heavy computations with large amounts of data very quickly and relatively cheaply;
    \item the availability of huge amounts of data from many different sources, stored on either the cloud or our own storage devices at low costs, and transferred at a very high speed;
    \item the availability of powerful Machine Learning (ML) algorithms, especially those based on deep artificial neural networks (NN). 
\end{enumerate}

We can identify three main software tasks that have driven the AI revolution:
\begin{itemize}
    \item \emph{discriminative}: the goal is data categorization, also known as data classification, which is the process of organizing data into specific categories or clusters based on shared characteristics or predefined criteria; 
    \item \emph{generative}: the scope is generating new information (texts, images, sounds, etc). This task is the one that has contributed most to the recent development of AI thanks to Large Language Models (LLM) like GPT--3 and GPT--4 of OpenAI, or T5 and Gemini developed by Google. The list, however, is incomplete and subjected to change very rapidly, due to the tremendous momentum in the LLM field;
    \item \emph{rewards--maximization}: it refers to the process of training an AI agent to make decisions that maximize the cumulative reward over time. This concept stands at the core of various AI applications, including game--playing, robotics, and autonomous systems.
\end{itemize}

For discriminative tasks, whose goal is to construct a mathematical classifier that assigns a label to a given data point, starting from a dataset where the class labels are known, two different approaches are typically followed.

The first one leverages \emph{Statistical Learning} \cite{James-statistical-learning} and is characterized by the total availability of datasets, generated by a prescribed probability distribution, which is used to make the machine learn during a \emph{training phase}. Then, in a successive \emph{testing phase} the machine is put to the test. Two mathematical issues arise in this context: the first one consists of characterizing the \emph{learnability} of the algorithm, i.e., the capability of the algorithm to learn from data and specifically how well it can generalize from a given dataset to unseen data.
The second one consists of characterizing the \emph{computational complexity} of finite learning problems, i.e., how much an algorithm is feasible and efficient in practice. These concepts were introduced in the seminal paper of Valiant in 1984  \cite{valiant1984}.

The second approach, named \emph{Online Learning}, is characterized by a continuous stream of training data generated, e.g., from markets, sensors, social media, etc., so that learning models are updated as soon as new data become available and consequently the train--test paradigm of statistical learning becomes unsuitable.
The mathematical theories modelling this approach are connected to game theory, more precisely to the theory of repeated games started by James Hannan and David Blackwell in the fifties of the last century \cite{hannan1957,blackwell1956}. The learning process behind this approach is based on adversarially chosen sequences of training datasets for which, given a learner and an adversary,  
$(i)$ the adversary selects the next couple of input--output training data $(\mathbf{x}_t,\mathbf{y}_t)$ and reveals the input $\mathbf{x}_t$ to the learner; 
    $(ii)$ the learner outputs a randomized prediction;
    $(iii)$ the adversary reveals the output $\mathbf{y}_t$ to the learner.

\subsection{Machine Learning}

Among the several facets of AI, our interest is above all in \emph{Machine Learning} (ML), for which many definitions have been proposed so far.
Perhaps the simplest one was provided in an apocryphal quote by Arthur Samuel, often attributed to his article on ML for the game checkers\footnote{A. Samuel (1901--1990) created the first checkers program with IBM’s first commercial computer, the IBM 701. IBM’s stock rose 15 points overnight on Wall Street.} \cite{Samuel1959studies}:
\emph{
``Machine Learning is the field of study that gives computers the ability to learn without being explicitly programmed.''}

In 2015, M.I. Jordan and T.M. Mitchell provided the following definition of ML:
\emph{``'Learning' is the process of transforming information into expertise or knowledge; 'machine learning' is automated learning''.}
Here, the ML process takes the form of an algorithm that takes information as input and produces knowledge as output.

A laconic alternative definition was adopted in a more recent report of The Royal Society (of London) of 2017 \cite{RoyalSociety}:
\emph{
``Machine Learning is a technology that allows computers to learn directly from examples and experience in the form of data.''
}

A large subset of ML is represented by \emph{artificial neural networks} (ANNs, often abbreviated simply as NNs), which are computing systems made of layers of artificial neurons and are loosely inspired by biological neural networks that constitute the human brain. Then, in the set of NNs, we point out the subset of \emph{deep neural networks} (DNNs), which are NNs with many hidden layers of neurons.

ML algorithms have been proven very successful in several contexts such as \emph{Autonomous Driving}, \emph{Natural Language Processing}, \emph{Playing Games}, \emph{Image analysis and generation}.
Two are the common elements behind the success of ML algorithms in these fields: first, the availability of large volumes of data, and second the fact that
developers know the truth, that is, they know the context well enough to validate the results, often even during model training. Both these elements are essential to train the algorithm.

However, there are many other sectors in which these two favourable conditions do not ever occur, so ML algorithms do not work so well. For instance, ML is not particularly effective when:
\begin{itemize}[noitemsep]
    \item datasets are sparse and partial. This is usually the case when data acquisition is expensive, like in medicine, when data are generated by clinical images or exams that could be invasive or extremely costly;
    \item training datasets are not available, 
    \item high--consequence decisions require human--interpretable models or outputs.
\end{itemize} 
In these cases, the training of the algorithms might be limited, and what they have learned is not enough to allow them to apply successfully to more general situations.

To better investigate ML, let us start with the definition Tom Mitchell gave in his book in 1997 \cite{Mitchell1997}:
\begin{quoteit}
    ``A computer program is said to learn from experience $E$ with respect to some class of tasks $T$ and performance measure $P$ if its performance at tasks in $T$, as measured by $P$, improves with experience $E$.''
\end{quoteit}

We will clarify this statement through two simple examples.

\begin{example}
\label{ex:classification} {Consider the following task $T$: ``given a picture (identified by the variable $\mathbf x$), write a program returning the result $y=1$ if it portrays a dog, $y=0$ otherwise'', i.e., we want to put a label on the picture. This is an instance of \emph{classification tasks}.

Following a traditional programming approach, to solve this task, we would write a computer program that implements the decision rules imposed by the human programmer.
Instead, a machine learning approach would proceed as follows:
\begin{enumerate}
\item \emph{experience} a collection of training data $(\hat{\mathbf{x}}_i,\hat{y}_i)$ for $i=1,\ldots,N$, i.e., a set of pictures $\hat{\mathbf{x}}_i$ with the corresponding labels $\hat{y}_i\in\{0,1\}$;
\item select a set of candidate \emph{models} 
$y=f(\mathbf{x};\mathbf{w})$ written as functions that depend on an input datum $\mathbf{x}$ (the picture) and some suitable unknown parameters that are the entries of the array $\mathbf{w}\in\mathbb R^M$;
\item define the \emph{loss function} (or \emph{objective function}, or again \emph{cost function}) 
\begin{eqnarray}\label{eq:loss-function}
\begin{array}{ll}
\displaystyle 
\mathcal{L}(\mathbf{w})=\frac{1}{2}\sum_{i=1}^N [d(\hat{y}_i, f(\hat{\mathbf{x}}_i;\mathbf{w}))]^2 \ +\ \text{regularization},
\end{array}
\end{eqnarray}
which measures the distance $d$ (e.g. the Euclidean distance), between the targets $\hat{y}_i$ furnished by the training set and the values $f(\hat{\mathbf{x}}_i;\mathbf{w})$ predicted by the model;
\item choose a \emph{metric} (i.e., an indicator) to monitor and measure the \emph{performance} of the model (which generally coincides with the loss function, although this is not mandatory, see Remark \ref{rem:loss-metric} and Sect. \ref{sec:performance});
\item \emph{train} the model by \emph{optimizing} the loss function ${\mathcal L}$, i.e., find
\begin{equation}\label{eq:minimizer}
\mathbf w^*=\argmin{\mathbf w \in \mathbb R^M}\mathcal L(\mathbf w)
\end{equation}
and define the ``optimal'' model $f(\cdot;\mathbf w^*)$;
\item measure the performance of $f(\cdot;\mathbf w^*)$ by evaluating the metric chosen in step 4 on different data than those used in step 5 (the so--called \textit{test data}). 
\end{enumerate}

The regularization term in (\ref{eq:loss-function}) can serve multiple purposes, as will be seen later in this paper.
Indeed, if $\mathcal L$ is not convex, it could show more than one (local) minimizer.}
\end{example}

\begin{remark}\label{rem:loss-metric} The primary purpose of the loss function is to guide the optimization process during training, while a metric is a measure used to evaluate the performance of a ML model. Unlike loss functions, metrics are used purely for evaluation purposes and do not influence the training process. 

The metric introduced in step 4 does not necessarily coincide with the loss function $\mathcal L$. Indeed, while the metric is not required to be differentiable with respect to the parameters $\mathbf w$, the loss function is. This is because gradient--like optimization algorithms that are invoked during the training phase require the loss function to be sufficiently regular, together with their gradients that are used to update the model parameters.
\end{remark}

In Example \ref{ex:classification}, we observe the typical characteristics of a ML approach to a classification problem. Rather than relying on prior knowledge about the context, such as how a dog looks, the method instead depends on a large dataset of images, referred to as the training set, and allows the program to autonomously identify the features that distinguish each class. A crucial element of this approach is the availability of sufficient and representative data, which enables the algorithm to extract meaningful information and achieve accuracy on unseen data. Equally important is the choice of candidate models, which must be expressive enough to capture the complexity of the input--output relationship being learned.

As a second example, let us examine an application of \emph{least--squares linear regression}. While this method is quite old -- its first application dates back to 1795 and is attributed to Gauss \cite{stigler1981gauss}, long before the term "Machine Learning" was coined -- it can still be interpreted as a learning process. This example, although deliberately simplified compared to the complexity of ML models commonly used in real--world applications, remains valuable from a pedagogical perspective as an introduction to ML concepts, particularly for those with a background in applied mathematics and scientific computing.
We will elaborate more on analogies and differences between
ML and least--squares algorithms in Sec. \ref{sec:conclusions}.

\begin{example} \label{ex:linear-regression}
Consider the following task: 
``Predict house prices $y$ based on the data $\mathbf x$ including the following information: square footage, number of bedrooms, number of bathrooms, the presence of brick construction, the number of offers received, and the neighbourhood in which each house is located.''
This is an example of \emph{regression tasks}.

In analogy with the procedure in Example 1, a ML strategy could be set up as follows:
\begin{enumerate}
\item set up the experience, by collecting a set of training data $\hat{\mathbf{x}}_i,\hat{y}_i$ (houses' information and the corresponding prices) for $i=1,\ldots,N$,
\item select a set of candidate models of the form
\begin{equation}\label{eq:model-regression}
y=f(\mathbf{x};\mathbf{w}) = \mathbf{m}\cdot\mathbf{x}+q,\quad \mbox{with}\quad 
\mathbf{w}=[\mathbf{m},q]
\end{equation}
which states a correlation between the input $\mathbf x$ and the output $y$;
\item define the loss function
\begin{eqnarray}\label{eq:loss-regression}
 \mathcal{L}(\mathbf{w})=
\frac{1}{2}\sum_{i=1}^N (\hat{y}_i-(\mathbf{m}\cdot \hat{\mathbf{x}}_i+q))^2;
\end{eqnarray} 
\item choose a metric to monitor and measure the performance of the model (see Remark \ref{rem:loss-metric} and Sect. \ref{sec:performance});
\item train the model by optimizing the loss function ${\mathcal L}$, i.e., solve the problem
\begin{equation}\label{eq:minimizer-regression}
\mathbf w^*=\argmin{\mathbf w \in \mathbb R^M}\mathcal L(\mathbf w)
\end{equation}
and define the ``optimal'' model $f(\cdot;\mathbf w^*)$;
\item measure the performance of $f(\cdot;\mathbf w^*)$ by evaluating the metric chosen in the step 4. 
\end{enumerate}


Since the loss function $\mathcal L$ (\ref{eq:loss-regression}) is convex, it features a unique minimizer $\mathbf w^*$ that is the unique solution of the first--order optimality condition $\nabla\mathcal L(\mathbf w)=0$.
In particular, in the simpler case when $\mathbf x\in {\mathbb R}$ (let us rewrite it $x$), so that the unknown parameters stored in $\mathbf w$ are the real values $m$ and $q$, the gradient of the loss function reads
\begin{eqnarray*}
\nabla\mathcal{L}=
\left[\begin{array}{c}
\frac{\partial\mathcal{L}}{\partial m}\\[4mm]
 \frac{\partial\mathcal{L}}{\partial q}\\
\end{array}\right] =
\left[\begin{array}{c}
\displaystyle -\sum_{i=1}^N(\hat{y}_i- (m\hat{x}_i+q))\hat x_i\\
\displaystyle -\sum_{i=1}^N(\hat{y}_i- (m\hat{x}_i+q))
\end{array}\right],
\end{eqnarray*}
and solving $\nabla\mathcal{L}=\mathbf{0}$ is equivalent to solving the so--called normal equations
\begin{eqnarray*}
\left[\begin{array}{cc}
\displaystyle \sum_{i=1}^N \hat{x}_i^2 & \displaystyle \sum_{i=1}^N \hat{x}_i\\
 \displaystyle \sum_{i=1}^N \hat{x}_i & N
\end{array}\right]
\left[\begin{array}{c}
m\\[4mm]
q
\end{array}\right] =
\left[\begin{array}{c}
\displaystyle \sum_{i=1}^N \hat{x}_i\hat{y}_i\\
 \displaystyle \sum_{i=1}^N \hat{y}_i 
\end{array}\right],
\end{eqnarray*}
i.e., a $2\times 2$ linear system whose unknowns are the two parameters $m$ and $q$ of the linear model $f( x; m,q)= m \, x+q$.
\end{example}

In Example \ref{ex:linear-regression}, the training process, i.e., the solution of the minimization problem \eqref{eq:minimizer-regression}, reduces to solving a linear system. This simplification stems from two key choices: first, the model \eqref{eq:model-regression} is linear in the parameters $\mathbf{w}$; second, the loss function \eqref{eq:loss-regression} is quadratic in the errors. These two characteristics are the hallmarks of the linear least--squares regression method, which is computationally efficient, as it requires only solving a linear system rather than engaging in an iterative optimization process.
However, choosing a polynomial of degree 1 as our model imposes a significant prior on the input--output relationship. In ML applications, much richer model classes are typically favoured over simple low--degree polynomial functions. While such an enrichment introduces substantial algorithmic and computational challenges, it enables the representation of far more complex input--output relationships than those achievable with basic models like straight lines or polynomials.

Although Example \ref{ex:linear-regression} can formally be considered a learning process, typical ML problems differ in the following ways:

\begin{enumerate}
\item the size $n$ of data $\mathbf x$ is often very large, implying severe memory space limitations; 
\item the number $M$ of free parameters (i.e., the size of the array $\mathbf w$) is typically very large. Notice that the larger the number of free parameters, the greater the complexity in minimizing the loss function $\mathcal L$; 
\item the model $f(\mathbf{x};\mathbf{w})$ is often non--linear with respect to the parameters ${\mathbf{w}}$  (this implies that manual differentiation is very tricky and should be replaced by automatic differentiation) so that $\nabla \mathcal L$ is not linear. Consequently, iterative optimization algorithms, like, e.g, gradient descent, stochastic gradient descent, Newton, Broyden--Fletcher--Goldfarb--Shanno (BFGS) methods \cite{Quarteroni-CS-2014} should be invoked;
\end{enumerate}

\subsubsection{Machine Learning Tasks}\label{sec:task}
ML encompasses various tasks that can be broadly categorized as follows based on the type of learning and the goals of the algorithm. Here are the major tasks in ML:

\begin{description}
\item[Regression:] predict a continuous--valued vector $\mathbf{y}\in\mathbb{R}^m$ starting with a training set composed of pairs of an input signal and a ground--truth value.  Examples include image recognition and medical diagnosis of continuous quantities of interest.
\item[Classification:] assign a label from a finite set $\{1,\ \ldots, C\}$ (indicating a class) to any given objects belonging to a predefined category, for instance, the label $y$ of an image $\mathbf x$.
 The training set is composed of pairs of an input signal and a ground--truth value (or a label). Although the output should be an integer or a label, typically, the standard approach consists of predicting, for each potential class, a score (which is a real value) or a probability of belonging to that class, 
 so that the correct class has the maximum score.
Examples include image classification, spam filters, and fraud detection.
\item[Clustering:] partition a set of items into subsets by finding similarities and differences. Examples include image segmentation and community detection in biology, social sciences, and economics.
\item[Association:] find the probability of the co--occurrence of items in a collection. Examples include recommendation systems and market basket analysis.
\item[Density estimation:]  unveal the common pattern of a population. Possible applications are in anomaly detection and the generation of realistic scenarios.
\item[Dimensionality reduction:] provide a compact (low--dimensional) description of complex (high--dimensional) data and it can be applied to tackle multi--query problems and data compression.
\end{description}

%
%
\subsubsection{Machine Learning Experience}\label{sec:experience}
The ML tasks listed above can be categorized based on the type of learning, according to this classification.
\begin{description}
    \item[Supervised Learning:] 
     both inputs $\hat{\mathbf{x}}_i$ and outputs   $\hat{\mathbf{y}}_i$ are provided in the training set. ML tasks that fall here are regression and classification.
     
\item[Unsupervised Learning:]
      only inputs   $\hat{\mathbf{x}}_i$   are provided in the training set and we do not know what the result could be. Instances of this type of learning are:
      clustering, associations, density estimation, and dimensionality reduction.

\item[Reinforcement Learning:]
     (special kind of supervised learning) the training set is not fixed a--priori, but it changes dynamically, and the learning machine itself explores the input space.  It involves an agent that learns to make decisions by performing certain actions and receiving rewards or penalties. The agent learns to maximize cumulative rewards over time. Examples include game playing and robotic control.
\end{description}

A further approach is sometimes called \emph{Semi--supervised Learning} and occurs when only a small subset of the original training set is actually labelled.

In all cases, ML algorithms aim at learning and predicting results starting from data, by unearthing the patterns and relationships that hide beneath them. In the next sections, we limit our essay to the supervised learning situation.

%
%

\subsubsection{Machine Learning Performance measurement}\label{sec:performance}
To evaluate the ability of an ML model to succeed, we can evaluate
a performance measurement $P$, expressed by a \emph{metric}. Typically, $P$ depends on the task the model has to achieve. A general form of $P$ is
\begin{eqnarray}\label{eq:metric}
P =
\frac{1}{N}\sum_{i=1}^N d_M(\hat{\mathbf{y}}_i, \mathbf{y}_i),
\end{eqnarray}
where $\hat{\mathbf{y}}_i\in \mathbb{R}^m$ are the target values  (or labels),
$\mathbf{y}_i=\mathbf f(\hat{\mathbf{x}}_i;\mathbf{w})\in\mathbb R^m$ are the predicted values computed by the model $\mathbf f:\mathbb R^n\to \mathbb R^m$ starting from the known inputs $\hat{\mathbf x}_i\in\mathbb R^n$, and $d_M:\mathbb R^m\times \mathbb R^m \to \mathbb R$ is a suitable distance.
The most used metrics are:
\begin{description}
    \item[Mean Square Error (MSE):] $d_M$ in (\ref{eq:metric}) is  the euclidian distance
$$d_M(\hat{\mathbf{y}}_i, \mathbf{y}_i)= \frac{1}{m}
\|\hat{\mathbf{y}}_i-\mathbf{y}_i\|_2^2 
= \frac{1}{m}\sum_{j=1}^m \left((\hat{\mathbf{y}}_i)_j- (\mathbf{y}_i)_j\right)^2,
$$
so that   
$$MSE=\frac{1}{N\cdot m}\sum_{i=1}^N 
 \|\hat{\mathbf{y}}_i-\mathbf{y}_i\|_2^2. $$
\item [Root Mean Square Error (RMSE):] it is the square root of \textit{MSE}:
$$RMSE=\sqrt{MSE}.$$
\item[Mean Absolute Error (MAE):] $d_M$ in (\ref{eq:metric}) is the distance induced by the $1-$norm:
$$d_M(\hat{\mathbf{y}}_i, \mathbf{y}_i)=\frac{1}{m} 
\|\hat{\mathbf{y}}_i-\mathbf{y}_i\|_1=\frac{1}{m}\sum_{j=1}^m \left|(\hat{\mathbf{y}}_i)_j- (\mathbf{y}_i)_j\right|,$$
so that
$$MAE=\frac{1}{N\cdot m}\sum_{i=1}^N 
\|\hat{\mathbf{y}}_i-\mathbf{y}_i\|_1.$$ 
\item[Cross--Entropy (CE):] $d_M$ in (\ref{eq:metric}) is defined by
\begin{eqnarray*}
d_M(\hat{\mathbf{y}}_i, \mathbf{y}_i)= -\sum_{j=1}^m (\hat{\mathbf{y}}_i)_j \log(\mathbf{y}_i)_j, \quad\text{with } 
 (\hat{\mathbf{y}}_i)_j =\left\{\begin{array}{ll}
1 & \mbox{if sample $i$}\\
 &   \mbox{belongs to class $j$}\\[2mm]
0 & \mbox{otherwise,}
\end{array}\right.
\end{eqnarray*}
so that 
\begin{equation}
    CE= -\frac{1}{N}\sum_{i=1}^N\sum_{j=1}^m
    (\hat{\mathbf{y}}_i)_j \log(\mathbf{y}_i)_j.
\end{equation}
This metric derives from the \emph{Maximum Likelihood Estimator} for categorical variables.
\end{description}
While \textit{MSE}, \textit{RMSE}, and \textit{MAE} are typically used to measure the performance of regression tasks, \textit{CE} is used for classification. 

\subsubsection{Machine Learning Models}\label{sec:ML-models}

In abstract terms, an ML model is a mathematical function $\mathbf{f}$ that maps an input $\mathbf x$ belonging to a space $\mathcal X\subseteq\mathbb R^n$ into an output $\mathbf y\in\mathcal Y\subseteq \mathbb R^m$ and depends on a set of parameters $\mathbf w\in \mathcal M \subseteq\mathbb R^M$: 
\begin{equation}
\mathbf{f}: \mathcal{X}\times \mathcal{M} \subseteq \mathbb{R}^{n+M} \to \mathcal{Y}\subseteq \mathbb{R}^m.
\end{equation}
$\mathcal X$ is named \emph{input space}, $\mathcal Y$ \emph{target space}, while $\mathcal X^{\mathcal Y}$ is the space of all functions $\mathbf f:{\mathcal X}\to {\mathcal Y}$.
There are infinitely many different ways to define a function $\mathbf f\in\mathcal X^{\mathcal Y}$, thus typically, the search of the model is limited to a set $\mathcal H\subset \mathcal X^{\mathcal Y}$ named \emph{hypothesis space} (in the seminal book \cite{Mitchell1997}, the \emph{hypothesis} was what now we call the model).  We list here some common choices used in ML for the hypothesis space.
\begin{description}
\item[Linear models] are largely used in regression tasks to predict the behaviour of systems, from biological to environmental and social sciences to business. They take the form 
$\mathbf{f}(\mathbf{x};\mathbf{w})=W\mathbf{x}+\mathbf{b},$ 
with $W\in\mathbb R^{m\times n}$ and $\mathbf b\in \mathbb R^m$; the parameters' array $\mathbf w$ collects all the entries of both $W$ and $\mathbf b$, i.e.,  $\mathbf{w}=\{W_{11},\ W_{12},\ldots, W_{mn},\ b_1, \ldots, b_m\}$;
\item[Polynomial models] are used in regression tasks when the data do not fit a linear model. They take the form (considering for simplicity the case when $n = 1$) $\mathbf{f}(x;\mathbf{w})=a_0+a_1x+a_2x^2+\ldots+a_q x^q$, with $a_i\in \mathbb R$ for $i=1,\ldots, q\in\mathbb N$, here  $\mathbf{w}=(a_0,\ a_1,\ldots, a_q)$. These models are also linear with respect to the parameters;
\item[Decision trees] are non--parametric supervised learning methods used for classification and regression. They provide an answer based on a series of binary decisions inferred from the data features. The picture {in Fig. \ref{fig:decision-tree}} provides a simplified example of a decision tree to predict the cardiovascular risk for a diabetic man between 40 and 69 years old.
\begin{figure}
\begin{center}
\scalebox{0.8}{
\begin{tikzpicture}
[level distance =1cm,sibling distance=6cm,-,thick]
\footnotesize
\node {Age}
    child {node {$40-49$}
        [sibling distance=1.3cm]
        child {node {smoker}[-]
            child {node {no}[-]
                child {node {\colorbox{lightgreen}{low}}[->]}
            }
            child {node {yes}[-]
                child{node{pressure}[-]
                    child{node{$<150$}[-]
                        child{node{\colorbox{lightgreen}{low}}[->]}
                        }
                    child{node{$>150$}[-]
                        child{node{cholesterol}[-]
                            child{node{$<220$}[-]
                                child{node{\colorbox{lightorange}{medium}}[->]}
                                }
                            child{node{$>220$}[-]
                            child{node{\colorbox{lightred}{high}}[->]}
                            }
                            }
                        }
                }    
            }
        }
    }   
    child {node {$50-59$}[-]
        [sibling distance=4.2cm,-]
        child {node {smoker}[-]
            child{node{no}[-]
                child{node{pressure} 
        [sibling distance=1.3cm,-]
                    child{node{$<150$}[-]
                        child{node{\colorbox{lightgreen}{low}}[->]}
                        }
                        child{node{$>150$}[-]
                            child{node{cholesterol}[-]
                                child{node{$<220$}[-]
                                    child{node{\colorbox{lightorange}{medium}}[->]}
                                    }
                                child{node{$>220$}[-]
                                    child{node{\colorbox{lightred}{high}}[->]}
                                    }
                                }
                            }
                    }
                }
            child{node{yes}[-]
                child{node{pressure}[sibling distance=1.3cm,-]
                    child{node{$<150$}[-]
                        child{node{cholesterol}[-]
                                child{node{$<220$}[-]
                                    child{node{\colorbox{lightorange}{medium}}[->]}
                                    }
                                child{node{$>220$}[-]
                                    child{node{\colorbox{lightred}{high}}[->]}
                                    }
                                }
                            }
                        child{node{$>150$}[-]
                            child{node{\colorbox{lightred}{high}}[->]
                                 }
                            }
                    }
                }
        }
    }
    child {node {$60-69$}
        [sibling distance=1.3cm,-]
        child {node {smoker}[-]
            child {node{no}[-]
                child {node{pressure}[-]
                    child{node{$<150$}[-]
                        child{node{cholesterol}[-]
                            child{node{$<220$}[-]
                                child{node{\colorbox{lightorange}{medium}}[->]}
                                }
                            child{node{$>220$}[-]
                                child{node{\colorbox{lightred}{high}}[->]}
                                }
                            }
                        }
                    child{node{$>150$}[-]
                        child{node{\colorbox{lightred}{high}}[->]}
                        }
                    }
            }
            child {node{yes}[-]
                child{node{\colorbox{lightred}{high}}[->]}
            } 
        }
    };
\end{tikzpicture}
}
\end{center}
\caption{A simplified example of a decision tree to predict the
cardiovascular risk for a diabetic man between 40 and 69 years old}
\label{fig:decision-tree}

\end{figure}
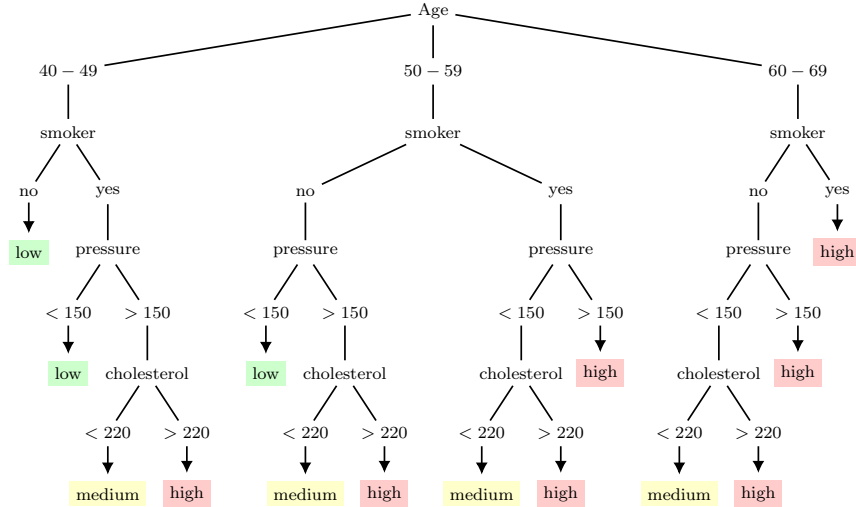

Unfortunately, these models are generally not robust with respect to different sets of data; 

\item[Random forests] represent an improvement of decision trees for both classification and regression tasks. They construct multiple decision trees by choosing subsets of data in the training set randomly and with repetition (a technique named \emph{bootstrapping} in statistics), then they average the output model (this strategy is named ``bootstrap aggregating'' or, more simply, \emph{bagging}). 
Despite their complexity, random forests provide robust and accurate models with reduced variance and reduced overfitting (the issue of overfitting models will be extensively analysed in the following);

\item[Maximal Margin Classifiers (MMC)] are ML methods for classification. In the simpler case of separating a set of data into two classes, they find a hyperplane $H=\{\mathbf x\in \mathbb R^n:\ \beta_0+\boldsymbol{\beta} \cdot \mathbf x\}$, with $\beta_0\in\mathbb R$ and $\boldsymbol\beta\in\mathbb R^n$, to classify the points $\mathbf x_i$ of a dataset by maximizing the margin, i.e., the minimum distance between the points of the two classes of the training set and the hyperplane itself. The points of the training set realizing the margin are called ``support vectors''. Introducing the array $\mathbf w=[\beta_0,\boldsymbol{\beta}]$, the model reads
$f(\mathbf x;\mathbf w)=\beta_0+\boldsymbol{\beta}\cdot\mathbf x$, if $f(\mathbf x;\mathbf w)>0$ the point $\mathbf x$ belongs to a class, otherwise to the other. This classifier cannot be applied to most data sets, since very often data belonging to different classes are not separable by a linear manifold;

\item[Support Vector Machines (SVM)] generalize the idea of MMC to the situation of classes not separable by a linear manifold. The idea consists of expanding the vector $\boldsymbol\beta$ (which is the normal vector to the hyperplane) with respect to $N$ training points $\mathbf x_i$, i.e., $\boldsymbol\beta=\sum_{i=1}^N\alpha_i\mathbf x_i$ so that $\boldsymbol\beta\cdot\mathbf x=\sum_{i=1}^N \alpha_i \mathbf x_i\cdot \mathbf x$, and in replacing 
the inner product $\mathbf x_i\cdot \mathbf x$ with a suitable non--linear kernel $\kappa(\mathbf x_i,\mathbf x)$. Then, the hyperplane is replaced by the surface $H=\{\mathbf x\in \mathbb R^n:\ \beta_0+\sum_{i=1}^N \alpha_i \kappa(\mathbf x_i,\mathbf x)=0\}$ and the model is
$f(\mathbf x;\mathbf w)=\beta_0+\sum_{i=1}^N \alpha_i \kappa(\mathbf x_i,\mathbf x)$. Notice that $H$ is the zero level set of $f(\mathbf x;\mathbf w)$. Now, the parameter array $\mathbf w$ includes $\beta_0$, the coefficients $\alpha_i$ and the parameters defining the kernel that would be polynomial, radial, or of other suitable types \cite{James-statistical-learning}. 

\item[Artificial Neural Networks (NN)] are non--linear models whose elementary core is the perceptron, the artificial neuron  introduced in 1958 by Rosenblatt, which in turn generalized the first artificial neuron invented in 1943 by Warren McCulloch and Walter Pitts. In both cases, artificial neurons mimic the logical behavior of a biological neuron.   
A biological neuron receives impulses from other connected neurons through dendrites, elaborates the information in its nucleus, and if the result of this elaboration exceeds a determined threshold, a new impulse is sent to subsequent neurons through the axon, otherwise, the signal is stopped (see Fig.\ref{fig:ANN} left).

A \emph{perceptron} (see Fig. \ref{fig:ANN}, right) is a mathematical model which receives 
one or more inputs $x_j$, computes the quantity
\begin{equation}
z=\sum_{j}w_{j}x_j+b,
\end{equation}
where $w_j$ are named \emph{weights} and $b$ is the so--called \emph{bias}, and returns the output 
\begin{equation}
y=\sigma(z),
\end{equation}
where $\sigma$, called \emph{activation function}, is typically a continuous and regular function approximating the Heavyside function
\begin{eqnarray} \label{eqn:heviside}
H(z)=\left\{\begin{array}{ll}
1 & \mbox{ if }z> 0\\
0 & \mbox{otherwise}.
\end{array}\right. 
\end{eqnarray}
In this context, the \emph{weights} $w_j$ represent the strength of the connections between neurons, modulating the influence of each input $x_j$ on the final output. The \emph{bias} $b$ can be interpreted as the activation threshold, but with the opposite sign: if the weighted sum of the inputs exceeds $-b$, the neuron fires and returns 1; if the weighted sum is less than $-b$, the neuron remains inactive and returns 0. The bias allows the neuron to shift its activation threshold. 

While the perceptron draws inspiration from the functioning of biological neurons, it is a highly simplified model compared to their actual behaviour. Advanced models of biological neurons account for the intricate ionic dynamics through differential equations (similar to those discussed in Sec.~\ref{sec:IHM}), which are compactly represented here by the activation function. They also consider diffusive processes using diffusion--reaction PDEs (as further detailed in Sec.~\ref{sec:IHM}), which in this case are approximated by simple linear combinations. However, the analogy with biological neurons serves only as an inspiration, and the simplicity of this model is seen as an advantage in the ML context. This simplicity enables indeed the training of models with very large dimensions, given the low computational cost of each individual neuron. As we will see, the preferred approach is often to design architectures with a very high number of simple units, rather than a small number of complex units.

\begin{figure}[h!]
    \begin{center}
    \includegraphics[width=0.4\textwidth]{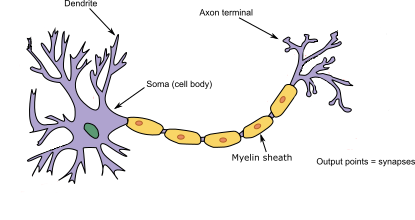}\quad 
    \includegraphics[width=0.3\textwidth]{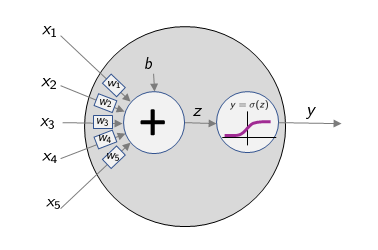}
 \end{center}
    \caption{A biological neuron (left) and the perceptron (right)}
    \label{fig:ANN}
\end{figure}
In the original model of McCulloch and Pitts (1943) named \emph{threshold logic unit}, the input and output were binary values $0/1$, the weights were all equal to 1 and the activation function was the  Heaviside function.

Different types of activation functions are shown in Fig. \ref{fig:activation}, the most common are:
\begin{itemize}[noitemsep]
\item Rectified Linear Unit (ReLU): $\sigma(z)=\text{ReLU}(z)=\max\{0,z\}$;
\item Sigmoid (or Logistic): $\sigma(z)=\displaystyle\frac{1}{1+e^{-z}}$. Its derivative $\sigma'(z)=\sigma(z)(1-\sigma(z))$ is very easy to evaluate and this makes the sigmoid function very attractive for the loss function minimization process;
\item Hyperbolic tangent: $\sigma(z)=\tanh(z)$,
\end{itemize}
and their relative performance is application--specific.

\begin{figure}[h!]
    \centering
    \begin{subfigure}[t]{0.25\textwidth}
        \centering
        \includegraphics[trim=0 100 0 0,width=\textwidth]{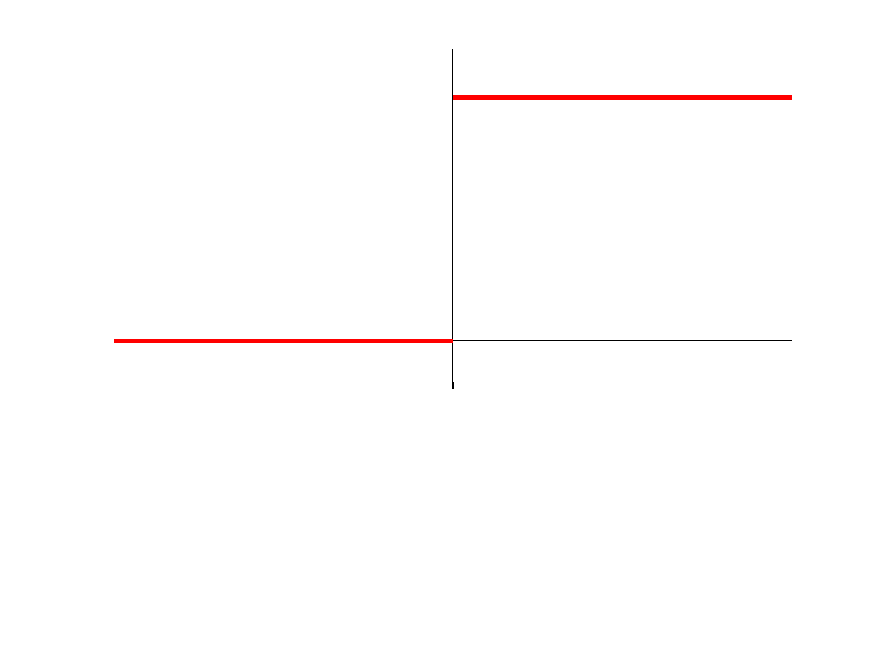}
        \caption{}
    \end{subfigure}%
    \begin{subfigure}[t]{0.25\textwidth}
        \centering
        \includegraphics[trim=0 100 0 0,width=\textwidth]{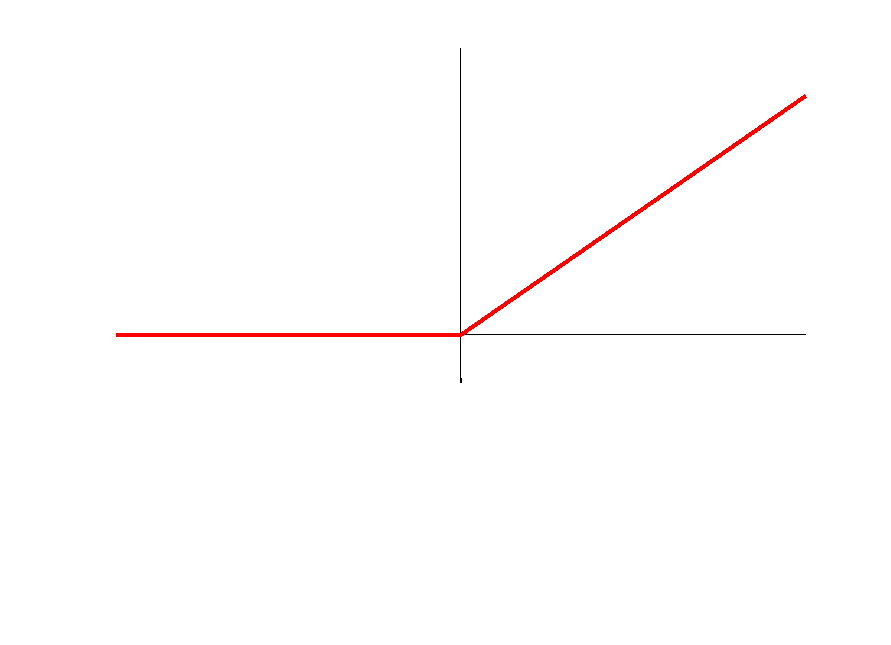}
        \caption{}
    \end{subfigure}%
     \begin{subfigure}[t]{0.25\textwidth}
        \centering
        \includegraphics[trim=0 100 0 0,width=\textwidth]{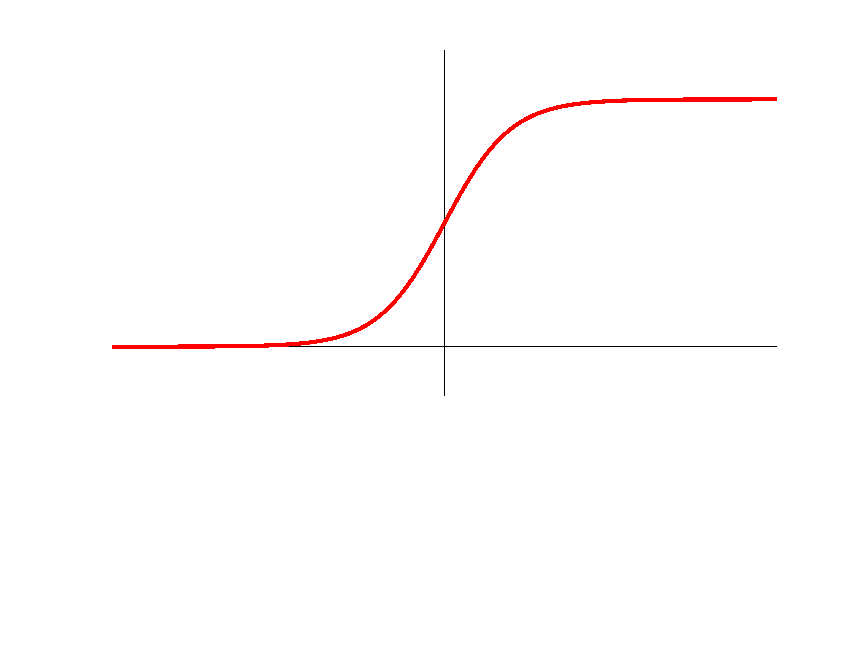}
        \caption{}
    \end{subfigure}%
    \begin{subfigure}[t]{0.25\textwidth}
        \centering
        \includegraphics[trim=0 50 0 50,width=\textwidth]{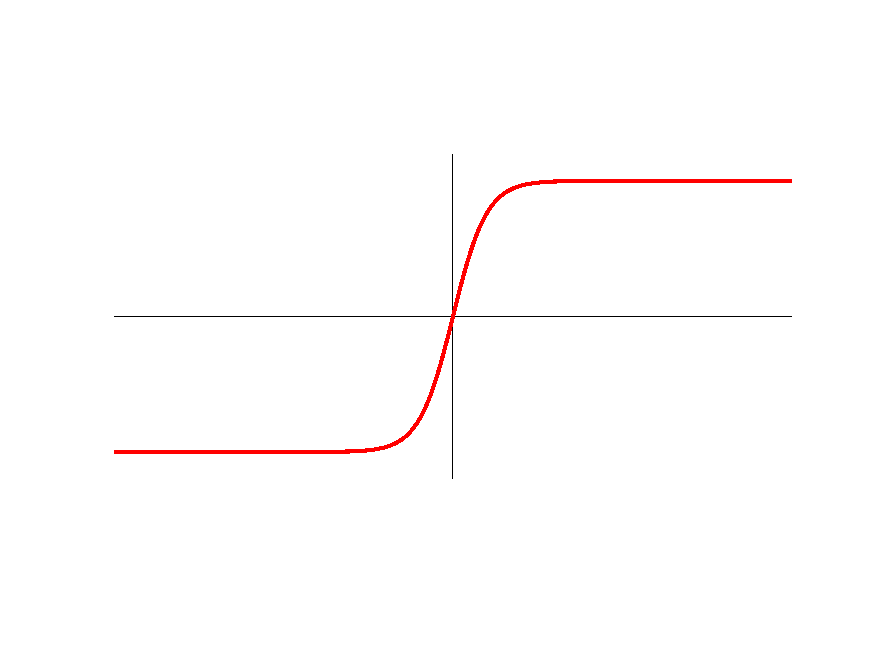}
        \caption{}
    \end{subfigure}%
    \caption{Activation functions: the Heaviside function (a), ReLU (b), sigmoid (c), and hyperbolic tangent (d)}
    \label{fig:activation}
\end{figure}

We may alter the transition's steepness and location in the sigmoid function by scaling and shifting the argument or, in the language of neural networks, by weighting and biasing the input, see Fig. \ref{fig:sigmoid}

\begin{figure}[h!]
    \centering
    \includegraphics[trim=0 100 0 0, width=0.4\textwidth]{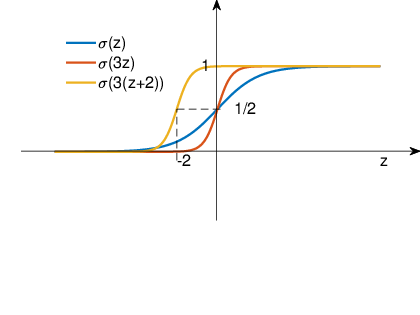}
    \caption{The sigmoid function with different transitions' steepness and location}
    \label{fig:sigmoid}
\end{figure}

The mathematical model associated with the perceptron reads
\begin{equation}\label{eq:perceptron}
y=f(\mathbf{x};\mathbf{w})=\sigma\left(\sum_j w_j x_j +b\right)
\end{equation}
where 
$\mathbf{w}=[w_1,\ w_2,\ldots,w_n,\ b]$ is the parameters' array.

A \emph{Feed Forward Neural Network} (FFNN) (or \emph{multilayer perceptron}) is a network composed of more layers each one including one or more perceptrons (which are usually named neurons) so that the outputs of the neurons of one layer become the inputs for the neurons of the next layer.
The first and the last layers are the \emph{input layer} and \emph{output layer}, respectively,  while the intermediate ones are named \emph{hidden layers} (see Fig. \ref{fig:FFNN}).
\begin{figure}
    \centering
    \includegraphics[width=0.9\linewidth]{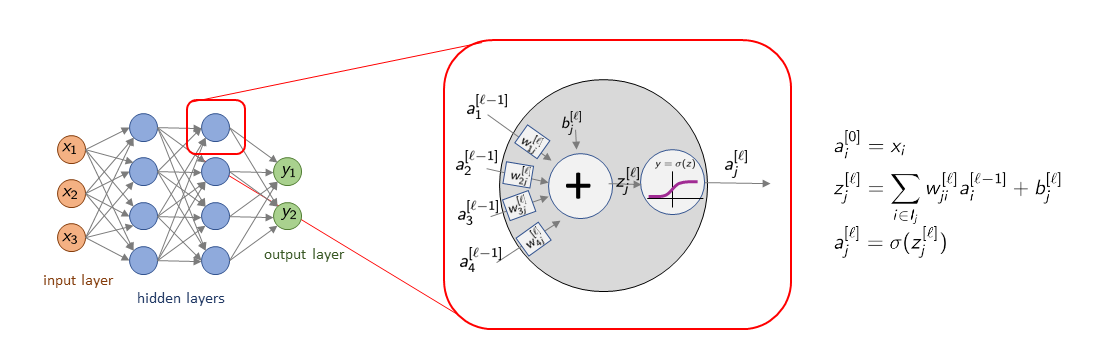}
    \caption{A Feed Forward Neural Network with $L=3$ and the zoom on a neuron at the layer $\ell=2$. } 
    \label{fig:FFNN}
\end{figure}
Let $L$ denote the total number of layers excluding the input one and $N_\ell$ the number of neurons of the generic layer $\ell\in\{0,\ldots, L\}$. We denote by $W^{[\ell]}\in \mathbb{R}^{N_\ell\times N_{\ell-1}}$,\ 
$\mathbf{b}^{[\ell]}\in
\mathbb{R}^{N_\ell}$ the matrix of the weights and biases of the neurons in layer $\ell$ and by $\mathbf{w}=\{(W^{[1]},\mathbf{b}^{[1]}),\ldots,
(W^{[L]},\mathbf{b}^{[L]})\}$ the list of all matrices and vectors containing weights and biases of all neurons. Finally, $\sigma$ is the activation function that, when applied to an array, works simultaneously on each component.
The mathematical model of a FFNN is a function 
\begin{equation}
    \mathbf{f}(\cdot;\mathbf{w}):\mathbb{R}^n\to \mathbb{R}^m
\end{equation}
that, given an instance of the parameters $\mathbf w$, takes the input  $\hat{\mathbf x}\in \mathbb R^n$ and computes the output $\mathbf y$ by 
following the instructions of algorithm \ref{alg:FFNN}. $\mathbf{f}$ is also called \emph{realization function} of the NN.

\begin{algorithm}[tb!]
\caption{Feed Forward Neural Network} \label{alg:FFNN}
\begin{algorithmic}[0]
\Procedure{FFNN\,}{$\hat{\mathbf x}$, $\mathbf w$}
\State $\mathbf{a}^{[0]}=\hat{\mathbf{x}}$
\For{$\ell=1,\ldots,L$} 
\State $\mathbf{z}^{[\ell]}=
W^{[\ell]}\mathbf{a}^{[\ell-1]}+\mathbf{b}^{[\ell]}$
\State $\mathbf{a}^{[\ell]}=\sigma(\mathbf{z}^{[\ell]})$
\EndFor
\State \Return{$\mathbf{y} =\mathbf{f} (\hat{\mathbf{x}};\mathbf{w})=\mathbf{a}^{[L]}$}
\EndProcedure
\end{algorithmic}
\end{algorithm}

The function $\mathbf f$ can be written as a composition of functions
\begin{equation}\label{eq:f_NN}
\mathbf f= \mathbf T^{[L]} \circ \ldots \circ  \sigma \circ \mathbf T^{[\ell]} \circ \ldots
\circ\sigma\circ \mathbf T^{[1]},
\end{equation}
where, for any $\ell=1,\ldots,L$, $\mathbf T^{[\ell]}: \mathbb R^{N_{\ell-1}}\to \mathbb R^{N_\ell}$ is defined by $\mathbf T^{[\ell]}(\mathbf x)=W^{[\ell]}\mathbf x+\mathbf b^{[\ell]}.$

FFNNs are typically characterized by the following indices:
\begin{itemize}[noitemsep]
    \item the total number of neurons
$\displaystyle N(\mathbf{w})=\sum_{\ell=0}^L N_\ell$, 
\item the width of the network $\displaystyle N_{\max}(\mathbf{w})=\max_{\ell=0,\ldots,L} N_\ell$,
\item  the total number of non--null weights and biases $\displaystyle M(\mathbf{w})=\sum_{\ell=0}^L M_j(\mathbf{w})$. 
\end{itemize}
In the example of Fig. \ref{fig:FFNN} we have $N(\mathbf w)=13$ and a total of $M(\mathbf w)=46$ weights and biases, 16 in the layer 1, 20 in the layer 2 and 10 in the layer 3.

FFNNs with $L=2$ (only one hidden layer) are named \emph{shallow NNs}, while when $L>2$ (at least two hidden layers) we speak about \emph{Deep Neural Networks (DNNs)}.

The number $L$ of layers and the numbers $N_\ell$ of neurons per layer (jointly with the activation function $\sigma$) characterize the \emph{architecture} of a NN and are usually named \emph{hyperparameters} of the NN.

On the contrary, weights and biases, usually named \emph{parameters} of the NN, are not given a--priori, and the challenge for NNs is to determine the parameters values that minimize the loss function and guarantee the best performance. Especially in deep NNs, this optimization problem is highly non--linear and non--convex, making classical global optimization algorithms cumbersome and inefficient. Starting from the straightforward but winning observation that partial derivatives of the loss function with respect to weights and biases can be calculated by using the chain rule and automatic differentiation \cite{Griewank2003, Baydin2018} (such a differentiation technique is known as \emph{backpropagation} in the ML community), more efficient gradient--based methods have superseded classical optimization algorithms, giving a tremendous boost to NNs (see Sect. \ref{sec:optimization}).

NNs based on the threshold logic unit of McCulloch and Pitts represent the first generation of NNs. They are not learnable because there are no weights or biases.
NNs based on Rosenblatt's sigmoid neuron represent the second generation of NNs. With respect to the first generation, they are learnable, but they are again fundamentally different from biological neurons. In fact, they model input and output in a continuous setting instead of spikes, and they process signals synchronously instead of taking care of spike timing and frequency. A third generation of NNs, based on spiking neurons was introduced in 1997 \cite{Maass1997}. These NNs can process spikes coming from different neurons at different times
and the spike timing itself plays a fundamental role in overcoming the threshold function and generating the action potential for successive neurons. Different spiking neuron models have been proposed so far. The most biologically plausible one, but at the same time very expensive and infeasible in
large--scale simulations, is based on the Hodgkin–-Huxley model \cite{Hodgkin1952}. Cheaper models, but poorer from a biological viewpoint, are the Integrate and Fire ones \cite{Maass1997,Yamazaki2022}. Despite 2nd generation NNS, spiking NNs still do not have solid training methods \cite{Yamazaki2022}.
\item[Convolutional Neural Networks (CNN)] are very effective NNs typically employed for signal processing and computer vision tasks. In CNNs, the weight matrices are circulating matrices whose kernels act as filters on the inputs and allow the extraction of typical features from the signal.
\end{description}

In addition to those contained in this list, far from being exhaustive, numerous other ML models are available depending on the specific application -- such as signal analysis, image analysis, and text analysis. The range of models in the literature continues to expand rapidly.
Finding a model means, first of all, choosing the hypothesis space (see Sect. \ref{sec:supervised_learning} and \ref{sec:hyperparameters}) and then characterizing the model itself by determining its parameters $\mathbf w$ by invoking optimization methods (see Sect. \ref{sec:optimization}). 

\subsubsection{Setting of supervised learning and error analysis}\label{sec:supervised_learning}

Starting from a finite set of samples representing inputs and targets (the training set), we are interested in finding a function that models the unknown relation between them, with the ultimate goal of predicting the output once a new input is provided. To move on we need the following ingredients:
\begin{enumerate}
\item an input space $\mathcal{X}\subseteq \mathbb{R}^n$, a target space $\mathcal Y\subseteq \mathbb{R}^m$, and the space $\mathcal{Y}^{\mathcal{X}}$ of all the functions $\mathbf{f}:\mathcal{X}\to\mathcal{Y}$,
\item a hypothesis space
$\mathcal{H}\subset \mathcal{Y}^{\mathcal{X}}$ where the model is searched for,
\item 
a joint probability distribution function $P_{\mathcal{X,Y}}$ defined on some $\sigma$-algebra on $\mathcal{X}\times\mathcal{Y}$ that simultaneously represents the distribution of inputs as well as the conditional probability of the target $\mathbf y$ being appropriate for an input $\mathbf x$, 
\item a loss metric\footnote{the term ``loss'' is used for the metric as well as to define the loss function (\ref{eq:loss-function}), however, the loss metric $d_M$ can be different from the metric $d$ in the loss function.} $d_M:\mathcal{Y}\times\mathcal{Y} \to [0,+\infty]$
measuring the mismatch between the output provided by the function $\mathbf f$ and the target,
\item the \emph{expected risk}
\begin{equation}\label{eq:expected_risk} 
R(\mathbf{f})=\mathbb{E}[d_M(\mathbf{y},\mathbf{f}(\mathbf{x}))]=\int_{\mathcal{X}\times\mathcal{Y}}d_M(\mathbf{y},\mathbf{f}(\mathbf{x})) dP_{\mathcal{X},\mathcal{Y}}(\mathbf{x},\mathbf{y}),
\end{equation}
i.e., the expected value of the loss metric when applied to the whole space  
$\mathcal{X}\times\mathcal{Y}$. It is a measure of how well a model performs on average when making predictions,

\item a training set
$S=\{(\hat{\mathbf{x}}_i,\hat{\mathbf{y}}_i),\ i=1,\ldots,N\}$ with $\hat{\mathbf x}_i\in \mathcal X$ and $\hat{\mathbf y}_i\in\mathcal Y$,
\item the \emph{empirical risk}
\begin{equation}\label{eq:empirical_risk}
R_{S,N}(\mathbf{f})=\frac{1}{N}\sum_{i=1}^N d_M(\hat{\mathbf{y}}_i,
\mathbf{f}(\hat{\mathbf{x}}_i))
\end{equation}
that represents an approximation of the expected risk $R(\mathbf f)$ evaluated on the training set $S$.
\end{enumerate}

\begin{figure}
    \centering
    \includegraphics[width=0.3\linewidth]{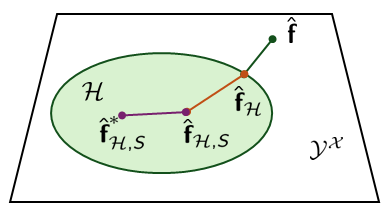}
    \caption{The ideal model $\hat{\mathbf f}$ and the really computed model $\hat{\mathbf{f}}^*_{\mathcal{H},S}$ }
    \label{fig:goals}
\end{figure}

The ideal goal would consist of finding a function $\hat{\mathbf f}\in \mathcal Y^{\mathcal X}$ that minimizes the expected risk, i.e.,
\begin{equation}\label{eq:goal1}
\hat{\mathbf{f}}=\argmin{\mathbf{f}\in\mathcal{Y}^{\mathcal{X}}} R(\mathbf{f}).
\end{equation}
Such function $\hat{\mathbf f}$ is the \textit{optimal} (in the sense of \eqref{eq:goal1}) input--output relationship corresponding to the joint probability distribution $P_{\mathcal{X,Y}}$.
However, looking for the unknown model $\hat{\mathbf{f}}$ in the space $\mathcal{Y}^{\mathcal{X}}$ is unfeasible because of the infinite dimension of the space itself, thus a first approximation of (\ref{eq:goal1}) consists in looking for a model in a suitable hypothesis space $\mathcal H$ subset of $\mathcal{Y}^{\mathcal{X}}$, that is, computing
\begin{equation}\label{eq:goal2}
\hat{\mathbf{f}}_{\mathcal{H}}=\argmin{\mathbf{f}\in\mathcal{H}} R(\mathbf{f}).
\end{equation}
Unfortunately, (\ref{eq:goal2}) is unfeasible too because the expected risk would require both sampling the whole space $\mathcal X\times \mathcal Y$ and explicitly knowing the probability distribution function $P_{\mathcal X\times \mathcal Y}$. So a new model $\hat{\mathbf{f}}_{\mathcal{H}, S}$ that relies solely on the training set $S$ is defined by solving the surrogate minimization problem
\begin{equation}\label{eq:goal3}
\hat{\mathbf{f}}_{\mathcal{H},S}=\argmin{\mathbf{f}\in\mathcal{H}} R_{S,N}(\mathbf{f}).
\end{equation}
Because $\hat{\mathbf{f}}_{\mathcal{H},S}$ is typically searched for by iterative optimization procedures, we finally find an approximation  $\hat{\mathbf{f}}^*_{\mathcal{H},S}$ of it (see Fig. \ref{fig:goals}).

We can measure the \emph{overall error} between the ideal model $\hat{\mathbf f}$ and the computed model $\hat{\mathbf{f}}^*_{\mathcal{H},S}$ by evaluating the difference between the corresponding expected risks $|R(\hat{\mathbf{f}})-R(\hat{\mathbf{f}}^*_{\mathcal{H},S})|$, that results in the sum of three different errors \cite{Berner_Grohs_Kutyniok_Petersen_2022, Guhring_Raslan_Kutyniok_2022}: 
\begin{eqnarray}\label{eq:total_error}
\begin{array}{llll}
|R(\hat{\mathbf{f}})-R(\hat{\mathbf{f}}^*_{\mathcal{H},S})|  
& \leq & |R(\hat{\mathbf{f}})-R(\hat{\mathbf{f}}_{\mathcal{H}})| &
\mbox{\emph{(approximation error)}}\\[4mm]
&+ &|R(\hat{\mathbf{f}}_{\mathcal{H}})-R_{S,N}(\hat{\mathbf{f}}_{\mathcal H})| &
\mbox{\emph{(generalization error)}}\\[4mm]
&+ &|R_{S,N}(\hat{\mathbf{f}}_{\mathcal{H}})-R_{S,N}(\hat{\mathbf{f}}^*_{\mathcal{H},S})| &
\mbox{\emph{(optimization error)}}\\[4mm]
&+ &|R_{S,N}(\hat{\mathbf{f}}^*_{\mathcal{H},S})-
R(\hat{\mathbf{f}}^*_{\mathcal{H},S})|&
\mbox{\emph{(bounded by the uniform}}\\
&  & & \hfill{\mbox{\emph{generalization error)}}}.
\end{array}
\end{eqnarray}
The \emph{approximation error} measures how well the hypothesis space $\mathcal H$ approximates the functional space $\mathcal Y^{\mathcal X}$, the \emph{generalization error} measures the gap between the performance of a neural network on the training set and unseen data, indicating how accurately the network can predict outcome values for previously unobserved inputs, while the \emph{optimization error} is related to the minimization algorithm. The fourth term in (\ref{eq:total_error}) can be bounded by the \emph{uniform generalization error}\cite{Berner_Grohs_Kutyniok_Petersen_2022}
\begin{equation}\label{eq:uniform_gen_err}
\sup_{\mathbf f\in\mathcal H}|R(\mathbf f)-R_{S,N}(\mathbf f)|.
\end{equation}
Thus, the design of a successful ML model is based on three ingredients: the hypothesis space, the training set, and the optimization algorithm. 

%
%
\null\textbf{The approximation error.} It measures the \emph{expressivity} of the ML model, i.e. how well a given input--output function can be approximated by functions belonging to the hypothesis space.
We recall here a few fundamental results about approximation error in FFNNs, distinguishing between shallow and deep FFNNs, while we refer to \cite{GrohsKutyniok2023} and the references therein for an in--depth discussion of this topic.

Shallow FFNNs (those with $L=2$, one hidden layer) are universal approximators thanks to the ``Universal Approximation Theorem'' of Cybenko \cite{Cybenko1989}. This means that, under suitable regularity assumptions on the activation function, any continuous function on a compact set can be approximated by a suitable shallow FFNN up to an arbitrary precision. 
The following definition serves in the statement of Cybenko's theorem.

\begin{definition}\cite{Cybenko1989}
Let $K\subset\mathbb{R}^n$ be a compact set. A continuous function $\sigma:\mathbb{R}\to{\mathbb R}$ is said \emph{discriminatory} with respect to $K$ if, for every finite, signed, regular Borel measure $\mu$ on $K$, we have that 
$$\left(\int_K \sigma(\mathbf{w}\cdot\mathbf{x}+b)d\mu(\mathbf{x})=0,\ \mbox{for all } \mathbf{w}\in\mathbb{R}^n \mbox{ and } b\in \mathbb{R}\right) \Rightarrow
\mu=0.$$
\end{definition}
Examples of discriminatory functions are the sigmoid and ReLU, while polynomials are not discriminatory.


\begin{theorem}[Cybenko (1989)]\cite[Thm 1]{Cybenko1989}\label{thm:Cybenko}
Let $K\subset\mathbb{R}^n$ be a compact set and $\sigma\in{\mathcal C}({\mathbb R})$ a discriminatory activation function on $K$. Then, given any real--valued function $\hat f$ continuous on $K$ and tolerance $\varepsilon>0$, there exists a shallow FFNN with weights and biases $\mathbf{w}=
\{(W^{[1]},\mathbf{b}^{[1]}),(W^{[2]},\mathbf{b}^{[2]})\}$ and real output $y=f (\cdot;\mathbf{w})$ such that $\|\hat{f}-f (\cdot;\mathbf{w})\|_\infty \leq \varepsilon.$
\end{theorem}

Denoting by $N_1$ the number of hidden neurons, Cybenko's theorem ensures that the larger the number $N_1$ of parameters, the lower the tolerance $\varepsilon$. But, how does $N_1$ depend on $\varepsilon$?
The following theorem, which summarizes the results of \cite{Maiorov1999, Maiorov-Pinkus, Mhaskar1996}
(see also \cite[Thm 3.6 and 3.8]{Guhring_Raslan_Kutyniok_2022}) provides lower and upper bounds for $N_1$ versus the tolerance $\varepsilon$. 
\begin{theorem}[Lower and upper complexity bound for shallow NNs] 
Let $n\geq 2$, $s\in\mathbb N$, $I$ an open real interval, $\sigma:\mathbb{R}\to\mathbb{R}$ such that $\sigma|_I\in\mathcal{C}^\infty(I)$ and $\sigma^{(k)}(x_0)\not = 0$ for some $x_0\in I$ and all integers $k\ge 1$. Then, for any
$\varepsilon>0$ and for any $\hat{f}\in W^{s,2}((0,1)^n)$ with $\|\hat f\|_{W^{s,2}((0,1)^n)}\leq 1$, there exists a shallow FFNN with weights
$\mathbf{w}=\{(W^{[1]},\mathbf{b}^{[1]}),(W^{[2]},\mathbf{b}^{[2]})\}$ depending on $\hat f$ and $\varepsilon$, real output $y=f (\cdot;\mathbf{w})$, and $N_1$ hidden neurons such that
$\|\hat{f}-f (\cdot;\mathbf{w})\|_2\leq \varepsilon$
and
\begin{equation}\label{eq:lower-upper-bounds}
    \left(\frac{1}{\varepsilon}\right)^{(n-1)/s} \lesssim N_1\lesssim \left(\frac{1}{\varepsilon}\right)^{n/s}.
\end{equation}
\end{theorem}

The estimate (\ref{eq:lower-upper-bounds}) is sharp for functions in $W^{s,2}((0,1)^n)$, but it can be improved if $\hat f$ is analytic \cite{Mhaskar1996}. We conclude that the smaller the tolerance $\varepsilon$, the larger the number $N_1$ of hidden neurons. Moreover,  
for fixed $\varepsilon$: 
\begin{itemize}[noitemsep]
    \item the higher the regularity $s$ of $\hat f$, the lower the number $N_1$ of hidden neurons,
    \item the larger the input dimension $n$, the larger the number $N_1$ of hidden neurons that guarantee the approximation error is lower than $\varepsilon$.
\end{itemize} 
The latter observation highlights that shallow FFNNs suﬀer from the \emph{curse of dimensionality}, i.e. the number of parameters of a shallow FFNN grows exponentially in the input dimension $n$.

In the case of deep FFNNs ($L>2$) with ReLU activation function, the following theorem by Yarotsky states that any function in  $W^{s,\infty}([0,1]^n)$ with bounded norm can be approximated by a suitable deep ReLU network.
\begin{theorem}[Yarotsky (2017)]\cite[Thm 1]{Yarotsky2017}
    For any $n,\ s \in{\mathbb N}_{s\geq 1}$, and for any $\varepsilon\in(0,1)$, there exists a network architecture with ReLU activation function $\sigma$, $L \lesssim 1+\log_e\frac{1}{\varepsilon}$ layers, $M \lesssim \varepsilon^{-n/s}\left(1+\log_e\frac{1}{\varepsilon}\right)$ non--null parameters, and real output such that, for any $\hat{f}\in W^{s,\infty}([0,1]^n)$ with $\|\hat f\|_{W^{s,\infty}([0,1]^n)}\leq 1$, there exist parameters
$\mathbf{w}=\{(W^{[1]},\mathbf{b}^{[1]}),\ldots, (W^{[L]},\mathbf{b}^{[L]})\}$ depending on $\hat f$ and $\varepsilon$, such that, denoting with $y=f (\cdot;\mathbf{w})$ the output of the NN, it holds
$\|\hat{f}-f (\cdot;\mathbf{w})\|_\infty \leq \varepsilon$.
\end{theorem}

The next theorem extends Yarotsky's result to functions belonging to more general Sobolev spaces and provides an upper complexity bound of the NN.
\begin{theorem}[G\"{u}ring, Kutyniok, Petersen (2020)]\cite[Thm. 4.1]{Guhring2020}
Let $n\in{\mathbb N}$, $k\in \mathbb N_{k\geq 2}$, $p\in[1,\infty]$, $B>0$, and $s\in[0,1]$. For any $\varepsilon\in(0,1/2)$, there exists a network architecture with ReLU activation function $\sigma$, real output,
$L\lesssim \log_2\varepsilon^{-k/(k-s)}$ layers, $M\lesssim \varepsilon^{-n/(k-s)}\log_2(\varepsilon^{-k/(k-s)})$ non--null parameters, 
$N\lesssim \varepsilon^{-n/(k-s)}\log_2(\varepsilon^{-k/(k-s)})$ neurons such that, for any $\hat{f}\in W^{k,p}((0,1)^n)$ with $\|\hat f\|_{W^{k,p}((0,1)^n)}\leq B$, there exist parameters
$\mathbf{w}=\{(W^{[1]},\mathbf{b}^{[1]}),\ldots, (W^{[L]},\mathbf{b}^{[L]})\}$, depending on $\hat f$ and $\varepsilon$ and it holds
$\|\hat{f}-f (\cdot;\mathbf{w})\|_{W^{s,p}((0,1)^n)} \leq \varepsilon$.
\end{theorem}
A lower complexity bound for deep NNs is stated by the following theorem, which is a generalization of \cite[Thm. 4a]{Yarotsky2017}. 
\begin{theorem}[G\"{u}ring, Kutyniok, Petersen (2020)]\cite[Thm. 4.3]{Guhring2020}
Let $n\in{\mathbb N}$, $k\in \mathbb N_{k\geq 2}$, $B>0$, and $s\in\{0,1\}$. If $\varepsilon\in(0,1/2)$, and there exists a ReLU network architecture with real output, $L$ layers, $M$ non--null parameters such that for any $\hat{f}\in W^{k,\infty}((0,1)^n)$ with $\|\hat f\|_{W^{k,\infty}((0,1)^n)}\leq B$ it holds $\|\hat{f}-f (\cdot;\mathbf{w})\|_{W^{s,p}((0,1)^n)} \leq \varepsilon,$ then the architecture
must have at least $M=c\, \varepsilon^{-n/(2k-2s)}$ non--null parameters, where $c$ is a positive constant depending on $n,\ k,\ B,$ and $s$.
\end{theorem}

Deep FFNNs do not break the curse of dimensionality, however, they can represent some functions that shallow FFNNs cannot (e.g. compactly supported functions \cite{Peterson2022}).
Deep FFNNs are typically more efficient in approximating functions, i.e., for a given tolerance, the required total number $M$ of non--null parameters is lower than in shallow FFNNs, and with the same number $M$ of parameters, deep FFNNs can represent more complex functions than shallow FFNNs.

\null\textbf{The generalization error.}
This is the central and crucial aspect of ML algorithms: we want our model to be able to predict the outcome value for unobserved inputs and not only for training data. The ability of an ML algorithm to perform well on previously unobserved inputs is called \emph{generalization} \cite{Goodfellow-et-al-2016}. In designing an ML algorithm, we have two goals to bear in mind which, unfortunately, are often in competition with each other: 
\begin{enumerate}[noitemsep]
\item minimizing the empirical risk $R_{S,N}(\hat{\mathbf f}_{{\mathcal H},S})$ (see (\ref{eq:empirical_risk})) (also known as training error),  
\item minimizing the gap $|R(\hat{\mathbf f}_{\mathcal H})-R_{S,N}(\hat{\mathbf f}_{\mathcal H})|$, i.e. the \emph{generalization error}, between the expected risk and the empirical risk, for any $\hat{\mathbf f}_{\mathcal H}$ in the hypothesis space $\mathcal H$.
\end{enumerate}

\begin{figure}
    \centering
    \includegraphics[width=0.7\linewidth]{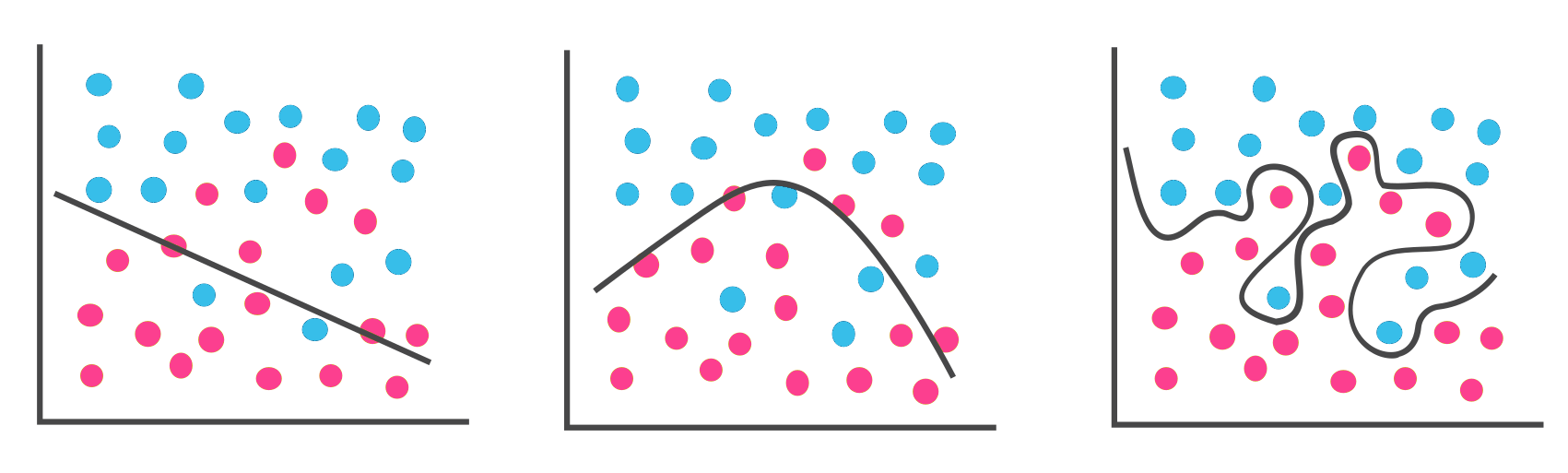}
    \caption{Underfitting (left), optimal fitting (centre), and overfitting (right) for a classification task}
    \label{fig:overfitting}
\end{figure}

The first goal can be accomplished by enriching the hypothesis space $\mathcal H$, but if $\mathcal H$ fits too much the training set $S$ we can incur in what is named \emph{overfitting}. This means that the network performs very accurately on the training data but cannot generalize well to new data, because the fitting process has focused too heavily on the unimportant and
unrepresentative ``noise'' in the given training set. Overfitting would generate a large gap between empirical risk and expected risk, i.e., a large generalization error.
On the contrary, when the hypothesis space is too coarse and the model is unable to represent the relation between inputs and targets accurately, we speak about \emph{underfitting}, here the empirical risk is large. A graphical interpretation of overfitting and underfitting is reported in Fig. \ref{fig:overfitting} relatively to a classification task.

A possible indicator of the flexibility and ability of a model to fit diverse data is called \emph{capacity}. It suggests how rich the hypothesis space $\mathcal{H}$ must be in the whole space $\mathcal Y^{\mathcal X}$ to model accurately the relation between inputs and targets. Small capacity corresponds to inaccurate fitting of the training data (underfitting), while high capacity may induce overfitting.

\begin{figure}[h!]
    \centering
    \includegraphics[width=0.6\linewidth]{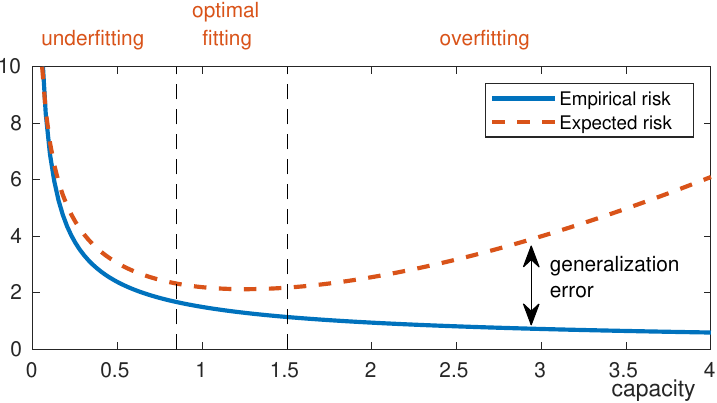}
    \caption{The empirical risk (training error), the expected risk, and the generalization error versus the capacity of the hypothesis space $\mathcal H$}
    \label{fig:capacity}
\end{figure}

The capacity is determined by the choice of the so--called 
\emph{hyperparameters}, like for example, the polynomial degree for polynomial models, or the number of neurons and layers for NNs (notice that hyperparameters are different from the parameters $\mathbf w$ that characterize the NN and are not optimized during the training).
Thus, the capacity increases with the polynomial degree in polynomial models, or with the number of neurons and layers in NNs.
The larger the capacity of the model, the smaller the empirical risk $R_{S,N}(\hat{\mathbf{f}}_{{\mathcal H},S})$ (see the blue curve in Fig. \ref{fig:capacity}). 
The expected risk instead increases with the capacity when the latter is too large (red curve in Fig. \ref{fig:capacity}) and features a minimum, in correspondence with which the fitting (and consequently the capacity) is optimal.

\emph{Cross--validation}.
A useful technique to combat overfitting, avoiding the (very cumbersome and often impossible) evaluation of the generalization error, consists of
splitting a given data set $\widehat S=\{(\hat{\mathbf{x}}_i,\hat{\mathbf{y}}_i),\ i=1,\ldots,\widehat N\}$ into two distinct groups:
\begin{itemize}
    \item the \emph{training set} $S=\{(\hat{\mathbf{x}}_i,\hat{\mathbf{y}}_i),\ i=1,\ldots,N\}$, used to drive the process that iteratively updates weights and biases of the NN by minimizing the empirical risk on it,
    \item the \emph{validation set} $S_{valid}=\{(\hat{\mathbf{x}}_i,\hat{\mathbf{y}}_i),\ i=N+1,\ldots,N+N_{valid}(=\widehat N)\}$, used to judge the performance of the current network. 
\end{itemize}

We define the \emph{training error}
\begin{equation}\label{eq:training_error}
 \mathcal E_{train}=R_{S,N}(\mathbf{f}) =\frac{1}{N}\sum_{i=1}^{N} d_M(\hat{\mathbf{y}}_i, {\mathbf f} (\hat{\mathbf{x}}_i;\mathbf w))
\end{equation}
and the \emph{validation error}
\begin{equation}\label{eq:test_error}
    \mathcal{E}_{valid}=R_{S_{valid},N_{valid}}(\mathbf{f}) =\frac{1}{N_{valid}}\sum_{i=N+1}^{N+N_{valid}} d_M(\hat{\mathbf{y}}_i, \mathbf f(\hat{\mathbf{x}}_i;\mathbf w)).
\end{equation}

Underfitting corresponds to large training errors $\mathcal{E}_{train}$ and can be mitigated by increasing the model capacity, i.e., by augmenting the number of parameters. 

\begin{figure}
    \centering
    \includegraphics[width=0.7\linewidth]{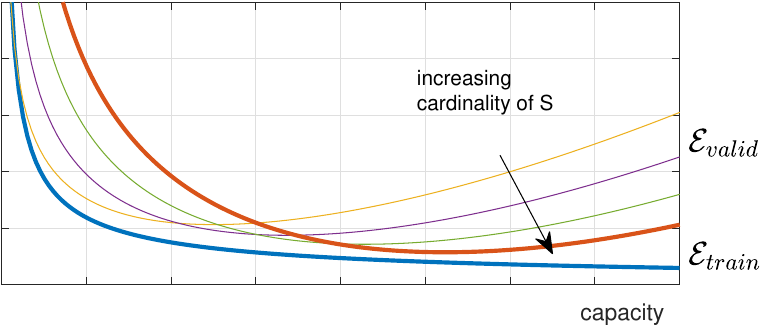}
    \caption{Train and validation errors versus the capacity of the hypothesis space $\mathcal H$}
    \label{fig:etest_evalid}
\end{figure}

Overfitting corresponds to the situation where the optimization process is driving down the training error, but the validation error  $\mathcal E_{valid}$ is no longer decreasing (so the
performance on unseen data does not improve), see Fig. \ref{fig:etest_evalid}. It is reasonable to
stop the optimization process at a stage where no improvement is
seen in the validation error.
Overfitting corresponds to large values for $\mathcal{E}_{test}-\mathcal{E}_{valid}$ and can be treated in different ways: by increasing the cardinality of the set $S$, or decreasing the capacity of the model (we speak about regularization by parsimony), or again by employing regularization techniques \cite[Ch. 5] {Goodfellow-et-al-2016}.

\medskip 
In the framework of binary classification, the generalization error can be estimated by exploiting the \emph{Vapnik--Chervonenkis (VC) dimension} of the hypothesis space \cite{Vapnik1974}. Following \cite{Mitchell1997}, the VC dimension measures the complexity of the hypothesis space $\mathcal H$ over an instance space $S$, not by the number of distinct models in $\mathcal H$, but instead by the number of distinct points from $S$ that can be completely discriminated using $\mathcal H$.

We say that a set of points $S$ is \emph{shattered} by a hypothesis space $\mathcal H$ if there exists a model in $\mathcal H$ that correctly classifies all the points of $S$, whatever their labels are. The VC dimension $d_{\mathcal H}$ of a hypothesis space $\mathcal H$ is the largest cardinality of a set $S$ that is shattered by $\mathcal H$ \cite{Peterson2022,GrohsKutyniok2023,Mitchell1997}.
For instance, if $S\subset \mathbb{R}^2$ and $\mathcal H$ is the set of linear binary classifiers, then $d_{\mathcal H}=3$, because any set of three points labelled by 0 or 1 can be correctly discriminated by a line, while not all the sets of four points labelled by 0 or 1 can.

The following results bound the generalization error in terms of the cardinality of the set $S$ and the VC dimension $d_{\mathcal H}$.
With probability at least $(1-\eta)$ one has  \cite{Bottou2018}
\begin{equation}\label{eq:bound_generalization_error}
\sup_{\mathbf{f}\in\mathcal{H}}|R(\mathbf{f})-R_{S,N}(\mathbf{f})|\lesssim \sqrt{\frac{1}{2N}\log\frac{2}{\eta}+\frac{d_{\mathcal H}}{N} \log\frac{N}{d_{\mathcal{H}}}}.
\end{equation}
This important estimate shows that:
\begin{itemize}[noitemsep]
\item if $d_{\mathcal H}$ is fixed the generalization error vanishes when $N\to +\infty$, i.e., when the cardinality of the training set increases;
\item if $N$ is fixed, the generalization error increases when $d_{\mathcal{H}}\to +\infty$, i.e., 
  too large capacity is bad if the training set is not large enough;
\item the cardinality $N$ of $S$ and the VC dimension $d_{\mathcal{H}}$ should increase at the same rate to maintain the same gap and guarantee the best performance.
\end{itemize}

\subsubsection{Optimization methods for training}\label{sec:optimization}

It is worth noticing that ML is strongly rooted in optimization algorithms, in a way that is far more complex than for classical optimization.
Indeed, the ideal goal of the training step would be to reach the best performance measure $P$ of the model. This would imply minimizing the expected risk $R$ among all the functions $\hat{\mathbf f}$ belonging to the hypothesis space $\mathcal H$. However, this would involve an intractable integral on the input--output space $\mathcal X\times \mathcal Y$, and neither it is possible to restrict the analysis to the test set because the only set we can use in this phase is the training one.
Thus, what we can actually do during the training is to minimize the loss function, that is the empirical risk $R_{S,N}$ by exploiting the samples of the training set. We can thus optimize the performance only indirectly.

Given the training set
$S=\{(\hat{\mathbf x}_i,\hat{\mathbf y}_i),\ i=1,\ldots, N\}\subset \mathbb R^n\times \mathbb R^m,$  
let us consider the loss function
\begin{equation}\label{eq:loss}
\mathcal{L}(\mathbf{w})=\frac{1}{N}\sum_{i=1}^N\frac{1}{2}
[d(\hat{\mathbf{y}}_i,\mathbf{f}(\hat{\mathbf{x}}_i;\mathbf{w}))]^2,
\end{equation}
where $d:\mathbb R^m \times \mathbb R^m\to \mathbb R^+ $ is a suitable distance, and its optimal solution
\begin{equation}\label{eq:optimalw}
\mathbf{w}^*=\argmin{\mathbf{w}\in\mathbb{R}^M} \mathcal{L}(\mathbf{w}).
\end{equation}
To find $\mathbf w^*$, typically, descent methods are invoked (see, e.g., \cite{Nocedal, Quarteroni-CS-2014, Manzoni-Quarteroni-Salsa-book}). For a given initial guess $\mathbf w^{(0)}$, construct recursively a sequence $\{\mathbf w^{(k)}\}_{k\geq 0}$ such that
\begin{equation}\label{eq:descent-update}
\mathbf w^{(k+1)}=\mathbf w^{(k)}+\eta_k \mathbf d^{(k)}, \hskip 1.cm k\geq 1,
\end{equation}
where $\mathbf d^{(k)}$ is a \emph{descent direction} that leads towards the minimum of $\mathcal L$ and
$\eta_k\in{\mathbb R}$ is the step--size (named \emph{learning rate} in the ML community). The iterations are stopped when a suitable error estimator (typically the norm of the residual or the norm of the increment between two iterates) is less than a given tolerance $\varepsilon$, or when a maximum number of iterations has been reached. 

\begin{figure}[b!]
    \centering
    \includegraphics[width=0.5\linewidth]{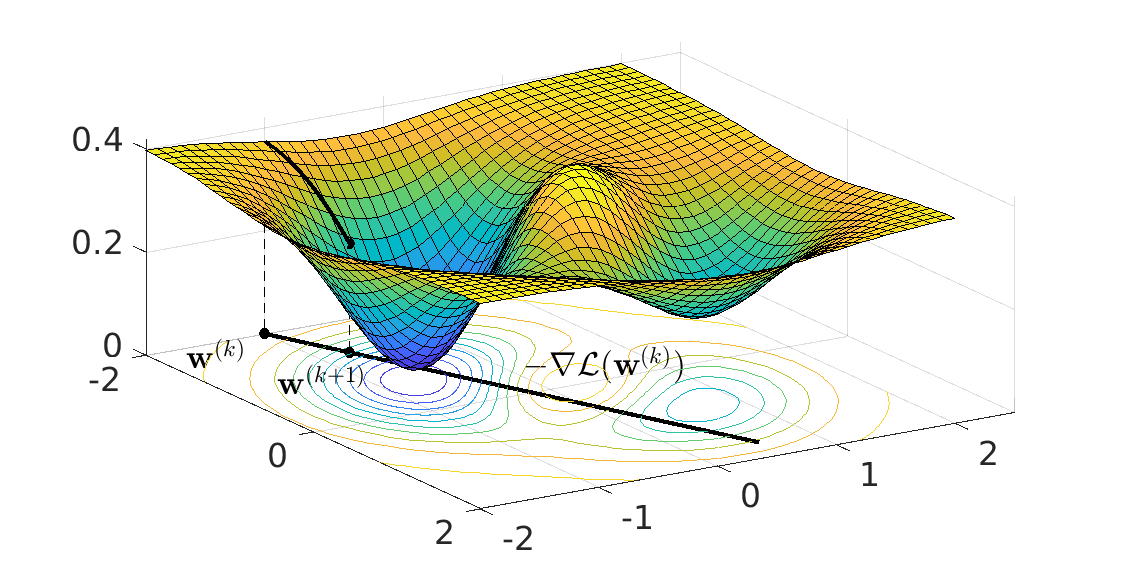}
    \quad
    \includegraphics[width=0.4\linewidth]{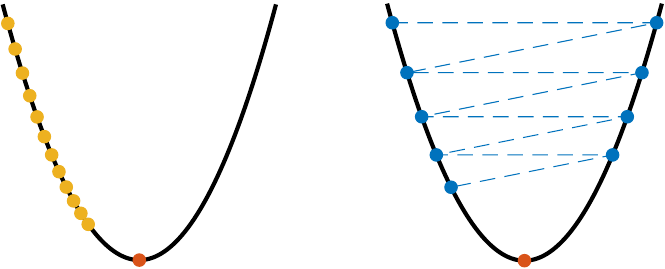}
    \caption{(Left) One update of the descent method with the descent direction given by the opposite of the gradient of the loss function. This function shows two local minima, the one close to $(-1,-1)$ is global while the other is local. (Right) Two non--adequate choices for the learning rate: one too small providing the yellow sequence and the other too large providing the blue one. Both of them make the convergence rate of the descent method very slow}
    \label{fig:descent-methods}
\end{figure}

If $\mathcal L$ is strictly convex on its domain (then $\mathcal L$ admits a unique minimizer), the sequence $\{\mathbf w^{(k)}\}$ converges to $\mathbf{w}^*$ for $k\to \infty$ for any $\mathbf w^{(0)}\in \mathbb R^M$. Notice that, whenever $\mathcal L$ is only convex in a neighbourhood of $\mathbf w^*$, the convergence is ensured if $\mathbf w^{(0)}$ is chosen sufficiently close to the solution. 
A loss function which is not globally strictly convex can future many local minimizers, see Fig. \ref{fig:descent-methods}. In such a case, descent methods can converge to any such minimizer, not necessarily to the one realizing $\min_{\mathbf w\in\mathbb R^M}\mathcal L(\mathbf w)$. This would depend on the choice of the initial guess.

Different criteria for choosing $\mathbf d^{(k)}$ are known in the literature. The most commonly used in classical optimization methods are reported in Tab \ref{tab:descent_directions} \cite{Nocedal, Quarteroni-CS-2014}.
Because $\nabla \mathcal L(\mathbf w^{(k)})$ points out the direction along which the function exhibits the steepest positive variation starting from $\mathbf w^{(k)}$, $-\nabla \mathcal L(\mathbf w^{(k)})$ leads towards a local minimizer of $\mathcal L$.
Thus, all the directions listed in Tab. \ref{tab:descent_directions} are descent directions, provided that $\mathcal L$ is strictly convex.

\begin{table} 
    \centering
    \begin{tabular}{lp{0.7\textwidth}}
    \toprule
    descent direction & $\mathbf d^{(k)}$\\
    \midrule
\emph{gradient:} & $-\nabla \mathcal L(\mathbf w^{(k)})$, \\
\emph{conjugate--gradient:} & $-\nabla \mathcal L(\mathbf w^{(k)})+\beta_k \mathbf d^{(k-1)}$, with $\beta_k$ depending on $\nabla \mathcal L(\mathbf w^{(k)})$ and $\nabla \mathcal L(\mathbf w^{(k-1)})$,\\
\emph{Newton:} & $-(H_{\mathcal L}(\mathbf w^{(k)}))^{-1}\nabla \mathcal L(\mathbf w^{(k)})$, with $H_{\mathcal L}$ the Hessian of $\mathcal L$,\\ 
\emph{quasi--Newton:} & $-H_k^{-1}\nabla \mathcal L(\mathbf w^{(k)})$, with $H_k$ a suitable approximation of  $H_{\mathcal L}$ (an example is the Broyden, Fletcher, Goldfarb, and Shanno (BFGS) method).\\
\bottomrule
    \end{tabular}
    \caption{Different descent directions in classical optimization methods}
    \label{tab:descent_directions}
\end{table}

The choice that requires the lowest computational cost at each iteration is the gradient one, although it does not necessarily provide the fastest convergence (i.e. the smallest number of iterations to converge up to a given tolerance $\varepsilon$). Conjugate--gradient directions require a slightly greater effort, but typically outperform gradient ones, while both Newton and quasi--Newton directions are more expensive, although they provide quadratic and superlinear convergence, respectively.

In classical optimization algorithms, the learning rate can be chosen dynamically by a line--search strategy \cite{Nocedal, Dennis-Schnabel} or taken constant during the iterations, i.e., $\eta_k\equiv \eta$ for any $k\geq 0$. The high cost of determining accurately a dynamic $\eta_k$ by, e.g., a line--search method or the trust method \cite{Dennis-Schnabel,Nocedal}, can be paid off by a faster convergence rate of the descent method. However, when the size of the problem is very large it could be more convenient to use a constant learning rate that, however, must be calibrated accurately, as we can see in Fig. \ref{fig:descent-methods}, right: too small or too large learning rates can drastically slow down the convergence rate of the descent method.

We further remark that when the loss function is globally strictly convex, finding its global minimizer is an achievable goal, although it could be costly. When instead, the loss function is not globally strictly convex, like in the left picture of Fig. \ref{fig:descent-methods}, finding the global minimizer is typically unfeasible or too computationally expensive and it is not worth looking for it (always bear in mind that we are minimizing the loss function, but the real goal should be optimizing the performance measure). 

Thus, instead of searching for the global minimizer of the loss function, we settle for looking for a ``good enough'' solution, i.e., a local minimizer featuring a low value of $\mathcal L$, or, to facilitate the solution process, it is often common to regularize the non--convex loss function by a strongly convex function (see Sect. \ref{sec:regularization}).

In the next sections, we recall the descent methods most widely used in ML optimization algorithms.
The following notation will be useful for all methods: let
\begin{equation}\label{eq:individual_loss}
\mathcal{L}_i(\mathbf{w})=\frac{1}{2}
[d(\hat{\mathbf{y}}_i,\mathbf{f}(\hat{\mathbf{x}}_i;\mathbf{w}))]^2,
\end{equation}
be the \emph{individual loss}, i.e., the contribution to the loss function of the $i-$th sample of $S$, so that (\ref{eq:loss}) can be interpreted as the arithmetic mean of the individual gradients, i.e., 
\begin{equation}\label{eq:loss2}
\mathcal{L}(\mathbf{w})=\frac{1}{N} \sum_{i=1}^N\mathcal{L}_i(\mathbf{w}).
\end{equation}

\null\textbf{Gradient Descent.}
The Gradient Descent (GD) method, also known as Steepest Descent, Full Gradient, or, again, Batch Gradient, is a descent method of the form (\ref{eq:descent-update}), where the descent directions are chosen equal to the opposite of the gradient, i.e., $\mathbf d^{(k)}=-\nabla \mathcal L(\mathbf w^{(k)})$. The learning rate can be chosen constant, i.e.  $\eta_k=\eta$ for any $k\geq 0$, or dynamically. In the latter case, a convenient strategy consists of choosing a monotone decreasing sequence $\eta_k\to 0^+$ when $k \to \infty$.
When chosen constant, $\eta$ represents a hyperparameter of the NN.
 The GD algorithm is sketched in Algorithm \ref{alg:GD}.
 
\begin{algorithm}[tb!]
\caption{Gradient Descent} \label{alg:GD}
\begin{algorithmic}[0]
\Procedure{GradientDescent\,}{training set $S$, $\mathbf{w}^{(0)}\in\mathbb{R}^M$}
\For{$k=0,\ldots,$ until convergence}
\State compute the learning rate $\eta_k$
\State $\mathbf{w}^{(k+1)}=\mathbf{w}^{(k)}-\eta_k 
\nabla\mathcal{L}(\mathbf{w}^{(k)})$
\EndFor
\EndProcedure
\end{algorithmic}
\end{algorithm}

\null\textbf{Stochastic Gradient Descent (SGD).}
When $N\gg 1$ the computation of the full gradient $\nabla\mathcal L(\mathbf w^{(k)})$ might be costly, thus a much cheaper alternative is to replace the mean of the individual gradients computed over all training points by the individual gradient at a single, randomly chosen, training point, 
i.e., choose randomly $i\in\{1,\ldots,N\}$ and consider the approximation
\begin{eqnarray}\label{eq:gradient_app1}
\nabla \mathcal{L}(\mathbf{w}^{(k)}) & \displaystyle =\frac{1}{N}\sum_{j=1}^N \nabla \mathcal{L}_j(\mathbf{w}^{(k)}) 
 \sim \nabla \mathcal{L}_i(\mathbf{w}^{(k)}).
\end{eqnarray}
This strategy can be intuitively justified thinking that it corresponds to restricting the training set $S$ to a single sample at each iteration. Such an approach can be interpreted as the limiting case of splitting the original sampling set into smaller and smaller training 
sets \cite{Bottou2018}.

Moreover, we remark once more that even the full gradient of the loss function $\mathcal L$ used by GD is an approximation of the gradient of the expected risk $R$ whose minimization should be our goal, thus, both $\nabla\mathcal L$ and $\nabla\mathcal L_i$ are a stochastic sampling of the gradient of the expected risk $R$. 

SGD was first proposed in \cite{Robbins1951} and, since then, many different variants of it have been derived. The simplest implementation consists of making a
``sampling with replacement'', i.e., the same ``$i$'' can be chosen at two successive steps, as done in Algorithm \ref{alg:SGD1}. Another implementation consists of making 
``sampling without replacement'' up to $N$ iterations, where $N$ is exactly the number of training samples. The cycle of $N$ iterations is named \emph{epoch}, see Algorithm \ref{alg:SGD2}. 

\begin{algorithm}[tb!]
\caption{Stochastic Gradient Descent (sampling with replacement)} \label{alg:SGD1}
\begin{algorithmic}[0]
\Procedure{StochasticGradientDescent\,}{training set $S$, $\mathbf{w}^{(0)}\in\mathbb{R}^M$}
\For{$k=0,\ldots,$ until convergence}
\State choose $i$ uniformly random from $\{1,2,\ldots,N\}$
\State compute the learning rate $\eta_k$
\State $\mathbf{w}^{(k+1)}=\mathbf{w}^{(k)}-\eta_k
\nabla\mathcal{L}_i(\mathbf{w}^{(k)})$
\EndFor
\EndProcedure
\end{algorithmic}
\end{algorithm}

\begin{algorithm}[tb!]
\caption{Stochastic Gradient Descent (sampling without replacement)} \label{alg:SGD2}
\begin{algorithmic}[0]
\Procedure{StochasticGradientDescent\,}{training set $S$, $\mathbf{w}^{(0)}\in\mathbb{R}^M$}
\State $k=0$
\For{$\overline{k}=0,\ldots,n_{epochs}-1$}
\State shuffle the integers $\{1,2,\ldots,N\}$ into a new order $\{k_1,\ldots, k_N\}$
\State \emph{\# next loop is an epoch}
\For{$i=1,\ldots,N$} 
\State compute the learning rate $\eta_k$
\State $\mathbf{w}^{(k+1)}=\mathbf{w}^{(k)}-\eta_k
\nabla\mathcal{L}_{k_i}(\mathbf{w}^{(k)})$
\State $k=k+1$
\EndFor 
\EndFor
\EndProcedure
\end{algorithmic}
\end{algorithm}

\begin{figure}
    \centering
    \includegraphics[width=0.5\linewidth]{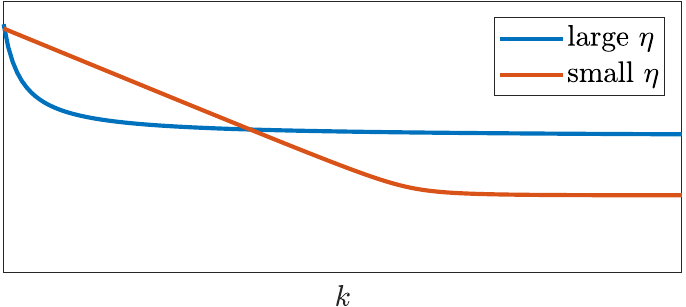}
    \caption{$\mathbb{E}[\mathcal{L}(\mathbf{w}^{(k)})-
\mathcal{L}(\mathbf{w}^*)]$ for the SDG method with two different constant choices of the learning rate $\eta$.}
    \label{fig:convergence_eta_const}
\end{figure}

%
%
\null\textbf{Choice of the learning rate.} 
The most popular and simple setting for the learning rate in SGD consists of taking a constant value for it: $\eta_k=\eta\ \forall k$. As already hinted before, choosing a suitable value for the learning rate is a major issue: too small or too large learning rates could be ineffective. In \cite[Thm 4.6]{Bottou2018}, the authors proved that, if $\mathcal L$ is strictly convex, if $\eta$ is not too large, then, in expectation, the sequence of function values $\{\mathcal L(\mathbf w^{(k)})\}$ converges near the optimal value, i.e., there exists a positive constant $C$ depending on $\mathcal L$ and its gradient, such that (see Fig. \ref{fig:convergence_eta_const})
\begin{equation}\label{eq:convergence_eta_const}
\mathbb E[\mathcal L(\mathbf w^{(k)})-\mathcal L(\mathbf w^*)]\xrightarrow{k\to\infty} \eta\, C,
\end{equation}
where $\mathbb E[\mathcal L(\mathbf w^{(k)})-\mathcal L(\mathbf w^*)]$ denotes the total expectation of $\mathcal L(\mathbf w^{(k)})-\mathcal L(\mathbf w^*)$ taken with respect to the joint distribution of all random variables (i.e. the stochastic choices of the samples) involved during the iterations from 1 to $k$.
From \eqref{eq:convergence_eta_const}, we can conclude that the smaller the learning rate, the higher the accuracy. However, decreasing the learning rate reduces the convergence speed, as illustrated in Fig.~\ref{fig:convergence_eta_const}.

To avoid very slow convergence, a strategy often employed in practice in SGD methods consists of starting with a large $\eta_0$ and reducing it progressively, for instance by halving $\eta_k$ when the loss function features a plateau, or reducing it with an inverse decay or an exponential rule:
   \begin{eqnarray}\label{eq:diminishing_eta}
       \eta_k=\frac{\beta}{\gamma+k} & \mbox{(inverse decay)}\\ 
       \eta_k=\eta_0 e^{-\gamma\, k} & \mbox{(exponential decay)}
   \end{eqnarray} 
   where $\gamma$ is a positive hyperparameter \cite{Bottou2018, Aggarwal2023}.

\null\textbf{GD vs SGD: costs and convergence rates.} 
Under the assumption that $\mathcal L$ is strongly convex, it is well known that GD features a \emph{linear convergence}, i.e., there exists $\rho\in(0,1)$ such that, for all $k\in{\mathbb N}$
\begin{equation}
    |\mathcal L(\mathbf w^{(k)})-\mathcal L(\mathbf w^*)| \lesssim \rho^k,
\end{equation}
meaning that at each iteration the error is reduced by a factor $\rho$ and, the lower $\rho$, the greater the reduction of the error. In the ML community, this type of convergence is known as ``geometric convergence''. If the iterations are stopped when the estimator of the error is less than a given tolerance $\varepsilon$, the total number of iterations required to guarantee $\rho^k\leq \varepsilon$ is proportional to $\log(1/\varepsilon)$.
Denoting by $C$ the cost for computing $\nabla \mathcal L_i(\mathbf w^{(k)})$, the cost of each iteration is about $N\cdot C$ and the total cost of GD is
\begin{equation}
    \mbox{cost of GD} \propto N\cdot C\cdot\log\frac{1}{\varepsilon}.
\end{equation}

For the SGD method, in \cite[Thm 4.7]{Bottou2018} the authors proved that, if $\mathcal L$ is strictly convex, for sampling without replacement and a learning rate is chosen as in (\ref{eq:diminishing_eta}), then
\begin{equation}\label{eq:SGD_conv}
     \mathbb{E}[\mathcal{L}(\mathbf{w}^{(k)})-\mathcal{L}(\mathbf{w}^*)]\lesssim \frac{1}{k},
\end{equation}
i.e., SGD features \emph{sublinear convergence}. 

The minimum number of iterations needed to satisfy convergence up to tolerance $\varepsilon$ is then proportional to $1/\varepsilon$ (we iterate until $1/k \leq \varepsilon$).
However, the cost for each iteration now is independent of the training set size $N$, thus, if $C$ always denotes the cost of computing $\nabla \mathcal L_i$, the total cost of SGD is
\begin{equation}
\mbox{cost of SGD}\propto \frac{C}{\varepsilon}.
\end{equation}
We can conclude that GD is more convenient when $N$ is moderate, SGD when $N$ is large.

Denoting by $\mathbf f_k=\mathbf f(\cdot; \mathbf w^{(k)})$ the model obtained after $k$ iterations of SGD, it is also possible to prove \cite{Bottou2018} that, in the limit $N \to + \infty$,
\begin{equation}\label{eq:SGD_conv_R}
\mathbb{E}[R(\mathbf{f}_k)-R(\hat{\mathbf{f}}_{\mathcal H})]\lesssim \frac{1}{k},
\end{equation}
i.e., SGD provides on the expected risk $R$ the same sublinear convergence  of (\ref{eq:SGD_conv}), meaning that, if $N$ is large enough with respect to $k$, SGD behaves similarly on both the expected and empirical risks. 
This suggests that SGD is less susceptible to overfitting than GD.

\null\textbf{Minibatch SGD and other noise reduction methods.}
SGD suffers from noisy gradient estimates that prevent the method from converging to the desired solution when the learning rates $\eta_k$ are taken constant, or provide sublinear convergence for a decreasing sequence $\eta_k$, for $k\geq 1$. One possible strategy to overcome this drawback consists of approximating $\nabla \mathcal L(\mathbf w^{(k)})$ with the linear combination of a small number $N_b$ of individual gradients (the so--called \emph{minibatch}):
\begin{equation}
\nabla \mathcal{L}(\mathbf{w}^{(k)}) \sim \frac{1}{N_b}\sum_{j=1}^{N_b} \nabla \mathcal{L}_{k_j}(\mathbf{w}^{(k)}),
\end{equation}
whose indices $k_j$ are randomly chosen in the set $\{1,\ldots,N\}$, as in Algorithm \ref{alg:SGD_minibatch}. This is a compromise between the full gradient method and the stochastic one to incorporate new gradient information and construct a more trustworthy step.

The convergence analysis provided in \cite{Bottou2018}  with a constant learning rate and $N_b \ll N$ shows that SDG and Minibatch SGD have comparable costs. However, Minibatch SGD can better benefit from GPU processors by adjusting the minibatch size to match the capacity of the GPU registers and optimizing data loading from the GPU memory to individual GPU cores.

\begin{algorithm}[tb!]
\caption{Minibatch SGD} \label{alg:SGD_minibatch}
\begin{algorithmic}[0]
\Procedure{MinibatchSGD\,}{training set $S$, $\mathbf{w}^{(0)}\in\mathbb{R}^M$}
\For{$k=0,\ldots,$ until convergence}
\State choose $N_b$ integers $k_1,\ldots, k_{N_b}$ uniformly random in $\{1,2,\ldots,N\}$
\State compute the learning rate $\eta_k$
\State $\displaystyle \mathbf{w}^{(k+1)}=\mathbf{w}^{(k)}-\eta_k 
\frac{1}{N_b}\sum_{j=1}^{N_b}\nabla\mathcal{L}_{k_j}(\mathbf{w}^{(k)})$
\EndFor
\EndProcedure
\end{algorithmic}
\end{algorithm}

Other methods aimed at reducing noisy gradient estimates and improving the rate of convergence from sublinear to linear regimes are the Dynamic Sample Size and the Gradient Aggregation methods.

Like Minibatch SDG, also \emph{Dynamic Sample Size} (DSS) methods use minibatches to approximate the gradient, but while in Minibatch SGD the minibatch size $N_b$ is constant during the iterations, in DSS it gradually increases (geometrically as a function of the iteration counter $k$) thus employing more and more accurate gradient information as the iterations proceed: given $\tau>1$, at the iteration $k$, the minibatch size is defined by $N_{b,k}\approx \tau^{k-1}$. Under the assumption that $\mathcal L$ is strictly convex, the convergence order of the resulting method is linear by maintaining the same computational complexity of SGD \cite[Cor. 5.2, Thm 5.3]{Bottou2018}.

\emph{Gradient Aggregation} (GA) methods reuse and/or revise gradient information computed and stored during previous iterations: the new descent direction is computed as a weighted average of the new direction with the old ones. An instance of these methods is the \emph{Stochastic Variance Reduced Gradient (SVRG)} \cite{Johnson2013,Bottou2018}.

\null
\textbf{SGD methods with momentum.}
Typically, (stochastic) gradient directions result in a zigzag behaviour thus cancelling the progress made in the previous steps.
This is not surprising because, even in the simpler case of solving symmetric positive definite linear systems, full gradient directions can repeat after a few steps  (they are linearly dependent), making the convergence of the gradient method very slow. On the contrary, conjugate--gradient directions (see Tab. \ref{tab:descent_directions}), with $\beta_k$ computed optimally, result all linearly independent and, if the matrix is symmetric and positive definite, the convergence of the Conjugate Gradient method is guaranteed at most in $N$ iterations if $N$ is the size of the system \cite{Quarteroni-NumerMath-2007}. 

While in SGD the descent direction is simply the stochastic gradient, in SGD with momentum the descent direction is a linear combination of the stochastic gradient and the previous descent direction, which in the ML community is named \emph{momentum}, resulting in a method very close to the Conjugate Gradient method.

The additional cost is due to the introduction of another hyperparameter $\rho\in[0,1)$, named smoothing parameter (or momentum parameter), which prevents GD from slowing down in flat regions of the loss function and getting trapped around local minima. SGD with momentum is reported in Algorithm \ref{alg:SGD-momentum}.

\begin{algorithm}[tb!]
\caption{SGD with momentum} \label{alg:SGD-momentum}
\begin{algorithmic}[0]
\Procedure{SGDwithMomentum\,}{training set $S$, 
$\mathbf{w}^{(0)}\in\mathbb{R}^M$, momentum parameter $\rho\in(0,1)$}
\State $\mathbf{v}^{(0)}=\textbf 0$
\For{$k=0,\ldots,n_{epochs}-1$}
\State choose $N_b$ integers $k_1,\ldots, k_{N_b}$ uniformly random in $\{1,2,\ldots,N\}$
\State $\displaystyle \mathbf{g}^{(k+1)}=
\frac{1}{N_b}\sum_{i=1}^{N_b}\nabla\mathcal{L}_{k_i}(\mathbf{w}^{(k)})$
\State compute the learning rate $\eta_k$
\State $\mathbf v^{(k+1)}=-\eta_k \mathbf g^{(k)}+\rho \mathbf v^{(k)}$
\State $\mathbf w^{(k+1)}=\mathbf w^{(k)}+\mathbf v^{(k+1)}$
\EndFor
\EndProcedure
\end{algorithmic}
\end{algorithm}
 Notice that $\mathbf v^{(k)}$ is aligned with the descent conjugate--gradient direction
of Tab. \ref{tab:descent_directions} and, when $\rho=0$, we recover Minibatch SGD. Moreover, the updating of $\mathbf w^{(k+1)}$ in Algorithm \ref{alg:SGD-momentum} can be equivalently written as
\begin{equation}
    \mathbf w^{(k+1)}=\mathbf w^{(k)}-\eta_k \mathbf g^{(k)}+\rho_k(\mathbf w^{(k)}-\mathbf w^{(k-1)}).
\end{equation}
When  $\eta_k=\eta$, $\rho_k=\rho$ for each $k\geq 0$, this method reduces to the \emph{heavy ball} method \cite{Polyak1964, Bottou2018} that outperforms the steepest descent for certain functions. Moreover, if $\mathcal L$ is strictly convex quadratic and $\eta_k$ and $\rho_k$ are chosen optimally at each iteration, i.e., they are the values that minimize $\mathcal L(\mathbf w^{(k+1)})$, then SGD with momentum coincides with the Conjugate Gradient method \cite{Bottou2018}. 

\null\textbf{Nesterov momentum method.}
It is a variant of SGD with momentum and it was proposed in \cite{Nesterov1983}. The gradient of the loss function is computed not on the last iteration $\mathbf w^{(k)}$, but on a linear combination of the two last iterations, as shown in Algorithm \ref{alg:Nesterov}. This method performs well when it is applied to the full gradient descent method (i.e. $N_b=N$) as it improves the convergence rate from $\mathcal O(1/k)$ to $\mathcal O (1/k^2)$, while when $N_b \ll N$ the Nesterov momentum method behaves like SGD. 
\begin{algorithm}[tb!]
\caption{Nesterov momentum} \label{alg:Nesterov}
\begin{algorithmic}[0]
\Procedure{Nesterov\,}{training set $S$, 
$\mathbf{w}^{(0)}\in\mathbb{R}^M$}
\For{$k=0,\ldots,n_{epochs}-1$}
\State $\widetilde{\mathbf w}^{(k)}=\mathbf w^{(k)}+\rho_k(\mathbf w^{(k)}-\mathbf w^{(k-1)})$
\State choose $N_b$ integers $k_1,\ldots, k_{N_b}$ uniformly random in $\{1,2,\ldots,N\}$
\State $\displaystyle \mathbf{g}^{(k+1)}=
\frac{1}{N_b}\sum_{i=1}^{N_b}\nabla\mathcal{L}_{k_i}(\widetilde{\mathbf{w}}^{(k)})$ 
\State compute the learning rate $\eta_k$
\State $\mathbf v^{(k+1)}=-\eta_k \mathbf g^{(k)}+\rho_k \mathbf v^{(k)}$ \quad with   $\rho_k\in[0,1)$
\State $\mathbf w^{(k+1)}=\widetilde{\mathbf {w}}^{(k)}-\eta_k \mathbf g^{(k)}$
\EndFor
\EndProcedure
\end{algorithmic}
\end{algorithm}

\null\textbf{Adaptive learning rate.}
The convergence rate of descent methods, and consequently the model's performance, is strongly sensitive to the learning rate $\eta$. Descent methods with momentum mitigate this drawback but at the price of introducing a second hyperparameter $\rho$.

A reliable alternative consists of considering diagonal scaling matrices that calibrate the learning rate to the descent direction component--wise. The purpose is to weigh less the larger components of the gradient (if along a direction the loss function decreases very rapidly we want to take a small step) and more the smaller ones (along the direction where the loss function features a weaker variation we are allowed to take a longer step). 

The first method designed with this strategy is the \emph{AdaGrad} method, introduced in \cite{Duchi2011}. It defines the scaling matrix by accumulating, element--by--element, the square of each component of the gradient $\mathbf g^{(\ell)}$ at all previous iterations $\ell\leq k$. However, the accumulation of gradient information from the beginning of the iterations reduces the learning step too quickly, thus AdaGrad works well when the loss function is convex, but its performance is not so good for non--convex problems.

\medskip
The \emph{Root Mean Square Propagation} (RMSProp) method, proposed in \cite{Hinton2012}, unlike AdaGrad, discards the oldest part of the gradient history so that when it falls into a convex bowl, the convergence gets faster. It results that RSMProp works better than AdaGrad, especially for non--convex functions and it is effective for many deep NNs. RMSProp is reported in Algorithm \ref{alg:RMSProp} in its original form. We notice that divisions and square roots in the expression $\eta/(\delta+\sqrt{\mathbf r^{(k+1)}})$ are applied element by element so that the result of this operation is an array of the same dimension of $\mathbf g^{(k)}$. Typical choices for the hyperparameters are $\rho=0.9$ and $\delta=10^{-7}$, while $(\mathbf{u}\odot\mathbf{v})_j=u_jv_j$ denotes the Hadamard product between the two vectors. 

\begin{algorithm}[tb!]
\caption{RMSProp} \label{alg:RMSProp}
\begin{algorithmic}[0]
\Procedure{RMSProp\,}{training set $S$, 
$\mathbf{w}^{(0)}\in\mathbb{R}^M$, hyperparameters $\eta,\ \rho,\ \delta\in \mathbb R$}
\State $\mathbf{r}^{(0)}=\textbf 0$
\For{$k=0,\ldots,n_{epochs}-1$}
\State choose $N_b$ integers $k_1,\ldots, k_{N_b}$ uniformly random in $\{1,2,\ldots,N\}$
\State $\displaystyle \mathbf{g}^{(k+1)}=
\frac{1}{N_b}\sum_{i=1}^{N_b}\nabla\mathcal{L}_{k_i}(\mathbf{w}^{(k)})$ 
\State $\mathbf r^{(k+1)}= \rho \mathbf r^{(k)}+(1-\rho) \mathbf g^{(k)}\odot \mathbf g^{(k)}$
\State $\displaystyle \mathbf w^{(k+1)}=\mathbf w^{(k)}-\frac{\eta}{\delta+\sqrt{\mathbf r^{(k+1)}}}\odot \mathbf g^{(k)}$
\EndFor
\EndProcedure
\end{algorithmic}
\end{algorithm}

\medskip\emph{Adam (Adaptive Moments) method.} This method was proposed in \cite{Kingma2015} it uses moments and shares some similarities with RMSProp. Adam considers both the first and second--order moments of the gradient and incorporates momentum with exponential weighting in both moments, see Algorithm \ref{alg:adam}. It requires four hyperparameters whose typical values are $\beta_1=0.9$, $\beta_2=0.999$, $\varepsilon=10^{-8}$, $\eta=10^{-3}$.

It is computationally efficient and features little memory requirement and it is appropriate for problems with very noisy and/or sparse gradients.

\begin{algorithm}[tb!]
\caption{Adam} \label{alg:adam}
\begin{algorithmic}[0]
\Procedure{Adam\,}{training set $S$, $\mathbf{w}^{(0)}\in\mathbb{R}^M$, hyperparameters $\beta_1,\ \beta_2\in(0,1)$, 
$\varepsilon,\ \eta>0$}
\State $\mathbf{m}_1^{(0)}= 
\mathbf{m}_2^{(0)}=\mathbf{0}\in\mathbb{R}^M$,
\For{$k=0,\ldots,$ until convergence}
\State choose $N_b$ integers $k_1,\ldots, k_{N_b}$ uniformly at random in $\{1,2,\ldots,N\}$
\State $\mathbf{g}^{(k)}=\frac{1}{N_b}\sum_{i=1}^{N_b}\nabla\mathcal{L}_{k_i}(\mathbf{w}^{(k)})$
\State $\mathbf{m}_1^{(k+1)}=\beta_1\mathbf{m}_1^{(k)}+
(1-\beta_1)\mathbf{g}^{(k)}$ ($1^{st}$ order moment)
\State $\mathbf{m}_2^{(k+1)}=\beta_2\mathbf{m}_2^{(k)}+
(1-\beta_2)\mathbf{g}^{(k)}\odot\mathbf{g}^{(k)}$ 
($2^{nd}$ order moment)
\State $ \hat{\mathbf{m}}_1^{(k+1)}=\mathbf{m}_1^{(k+1)}/(1-\beta_1^{k+1})$\;
\State $\displaystyle \hat{\mathbf{m}}_2^{(k+1)}=\mathbf{m}_2^{(k+1)}/(1-\beta_2^{k+1})$
\State $\mathbf{w}^{(k+1)}=\mathbf{w}^{(k)}-\frac{\eta}{\varepsilon+\sqrt{ \hat{\mathbf{m}}_2^{(k+1)}}}\odot  \hat{\mathbf{m}}_1^{(k+1) }$
\EndFor
\EndProcedure
\end{algorithmic}
\end{algorithm}

\null\textbf{Second--order and other methods.} Loss functions in  deep NNs often suffer from high non--linearity and ill--conditioning, the latter to be understood as instability with respect to data: in correspondence with small variations of the parameters the loss function might feature strong variations. We can think that the loss function can show cliffs preceded by gentle hills.
This makes the convergence of first--order methods very slow because the gradient alone does not give very accurate information to achieve the minimum. 

Methods using second--order derivatives of the loss function overcome such a drawback by capturing essential information about the curvature of the loss function. Examples include: Newton, Hessian--free Inexact Newton, 
Stochastic Quasi--Newton, and
Gauss--Newton \cite{Bottou2018}. However, these methods are far less successful than gradient descent methods in deep NNs training for many reasons, above all the difficulty of scaling to large NNs and the presence of many saddle points. We refer to \cite{Goodfellow-et-al-2016, Bottou2018} and the references therein for their description.

Stochastic Trust Region Methods (STRM) \cite{Blanchet2019, Curtis2019, bellavia0222} are alternatives to descent methods. While the latter compute the descent direction and the learning rate separately, STRMs define a quadratic model that locally approximates the loss function and directly provide the displacement from the last iterate, thus including both direction and step--length. Under suitable probabilistic assumptions, these methods find (in expectation) an $\varepsilon-$approximate minimizer at most in $\mathcal O(\varepsilon^{-2})$ inexact evaluations of the function and its derivatives \cite{bellavia0222}.

\subsubsection{Backpropagation}\label{sec:backpropagation}

As already pointed out, given a generic input $\hat{\mathbf x}\in{\mathbb R}^n$ and a corresponding target $\hat{\mathbf y}$, provided that the weights and biases stored in $\mathbf{w}=\{(W^{[1]}, \mathbf{b}^{[1]}), \ldots, (W^{[L]}, \mathbf{b}^{[L]})\}$ are known,  we can compute the output $\mathbf y=\mathbf f(\hat{\mathbf x};\mathbf w)=\mathbf a^{[L]}$ of the Feed Forward Neural Network (FFNN) with Algorithm \ref{alg:FFNN}. For the sake of simplicity, we assume that the distance $d$ used to define the loss function (\ref{eq:loss}) is the Euclidean norm, i.e.,
$d( \mathbf{y}_1,\mathbf y_2)=\|\mathbf{y}_1- \mathbf y_2\|_2$.

The backpropagation technique is a method to compute efficiently the gradients of the individual losses (\ref{eq:individual_loss}) during the training of a NN in order to determine the unknown weights and biases that minimize the loss function. For any $\hat{\mathbf x}$ in the input space, define  
\begin{equation}
    \mathcal L_{\hat{\mathbf x}}(\mathbf w)=\frac{1}{2}
    \|\hat{\mathbf{y}}-\mathbf{f} (\hat{\mathbf{x}};\mathbf{w})\|^2_2=
    \frac{1}{2}
    \|\hat{\mathbf{y}}- \mathbf a^{[L]}\|^2_2, 
\end{equation}
which can be seen as a particular instance of the individual loss (\ref{eq:individual_loss}).

The idea of backpropagation consists of applying the chain rule for computing the partial derivatives of $\mathcal L_{\hat{\mathbf x}}$ with respect to weights and biases, starting from the output $\mathbf a^{[L]}$ of the last layer $L$ and proceeding backwards, layer by layer with $\ell=L-1, L-2,\ldots, 1$, and leveraging both the variables $\mathbf z^{[\ell]}=W^{[\ell]}\mathbf a^{[\ell-1]}+\mathbf b^{[\ell]}$ introduced in Algorithm \ref{alg:FFNN} and the activation function $\sigma$.  The superindex $\ell$ stands for the number of the layer, more precisely $w_{j k}^{[\ell]}$ is the coefficient that weighs the output $a_k^{[\ell-1]}$ (of the neuron $k$ of the layer $\ell-1$) playing the role of input of the neuron $j$ of the layer $\ell$, while $b_j^{[\ell]}$ is the bias of the neuron $j$ at layer $\ell.$

The use of backpropagation to compute gradients fits into a very general framework of techniques known as \emph{automatic differentiation} or algorithmic differentiation \cite{Baydin2018, Griewank1989, Wengert1964}.

Following \cite{Higham2019} and recalling that $W^{[\ell]}\in \mathbb R^{N_\ell\times N_{\ell-1}}$, the component $j$ of the array $\mathbf z^{[\ell]}$, i.e.,
\begin{equation}\label{eq:weighted-input}
z_j^{[\ell]}=\sum_{k=1}^{N_\ell-1}w_{jk}^{[\ell]}a_k^{[\ell-1]} +b_j^{[\ell]} \hskip 1.cm \mbox{with } j=1,\ldots, N_\ell 
\end{equation}
is referred as the \emph{weighted input} for neuron $j$ at layer $\ell$ (see Fig. \ref{fig:FFNN}), while the quantity 
\begin{equation}\label{eq:delta_jell}
    \delta^{[\ell]}_j=\frac{\partial\mathcal L_{\hat{\mathbf x}}}{\partial z_j^{[\ell]}} 
\end{equation}
measures the sensitivity of the loss function $\mathcal L_{\hat{\mathbf x}}$ to the weighted input for the neuron $j$ at layer $\ell$. Then, let $\boldsymbol\delta^{[\ell]}\in \mathbb R^{N_\ell}$ be the array with entries given by (\ref{eq:delta_jell}).
As a consequence of the chain rule, the following relations hold \cite[Lemma 5.1]{Higham2019}:
\begin{eqnarray}\label{eq:backprop}
\begin{array}{rcll}
    \boldsymbol\delta^{[L]}&=&\sigma'(\mathbf{z}^{[L]})\odot
(\mathbf{a}^{[L]}-\hat{\mathbf{y}})\\[3mm]
\boldsymbol\delta^{[\ell]}&=&\sigma'(\mathbf{z}^{[\ell]})\odot
(W^{[\ell+1]})^T\boldsymbol\delta^{[\ell+1]} & \mbox{ for }1\leq \ell\leq L-1\\[3mm]
\displaystyle \frac{\partial\mathcal{L}_{\hat{\mathbf x}}}{\partial b_j^{[\ell]}}
&=&\delta_j^{[\ell]} & \mbox{ for }1\leq \ell\leq L\\[4mm]
\displaystyle \frac{\partial\mathcal{L}_{\hat{\mathbf x}}}{\partial w_{jk}^{[\ell]}}
&=&
\delta^{[\ell]}_j a_k^{[\ell-1]}& \mbox{ for }1\leq \ell\leq L,
\end{array}
\end{eqnarray}
where the first derivative $\sigma'$ of the activation function  is evaluated component--wise on $\mathbf z^{[\ell]}$ and $\mathbf u\odot\mathbf v$ is the Hadamard product defined by
$(\mathbf u\odot\mathbf v)_i=u_i\, v_i.$

The derivative $\partial\mathcal{L}_{\hat{\mathbf x}}/\partial w_{jk}^{[\ell]}$ measures how much $\mathcal{L}_{\hat{\mathbf x}}$ changes when we make a small perturbation on $w_{jk}^{[\ell]}$ and, similarly, $\partial\mathcal{L}_{\hat{\mathbf x}}/\partial b_{j}^{[\ell]}$.

Formulas (\ref{eq:backprop}) can be rewritten in a backward loop as described in Algorithm \ref{alg:backprop}, while in Algorithm \ref{alg:sdg_backprop} we train a NN by combining the basic version of Stochastic Gradient Descent (sampling with replacement, see Algorithm \ref{alg:SGD1}) with constant learning rate and backpropagation.
For safe of clarity, we omitted the iteration index $k$ on the weight matrices $W^{[\ell]}$ and biases arrays $\mathbf b^{[\ell]}$.


\begin{algorithm}[tb!]
\caption{BackPropagation} \label{alg:backprop}
\begin{algorithmic}[0]
\Procedure{BackPropagation\,}{$\hat{\mathbf y}$, $\mathbf z^{[\ell]}, \mathbf a^{[\ell]}$, $W^{[\ell]}$ for $\ell=1,\ldots, L$}
\State $\boldsymbol\delta^{[L]}=\sigma'(\mathbf{z}^{[L]})\odot
(\mathbf{a}^{[L]}-\hat{\mathbf{y}})$
\State $\displaystyle \frac{\partial\mathcal{L}_{\hat{\mathbf x}}}{\partial b_j^{[L]}} =\delta_j^{[L]}$, \qquad 
$\displaystyle \frac{\partial\mathcal{L}_{\hat{\mathbf x}}}{\partial w_{jk}^{[L]}}= \delta^{[L]}_j a_k^{[L-1]}$

\For{$\ell=L-1,\ldots,1$}
\State $\boldsymbol\delta^{[\ell]}=\sigma'(\mathbf{z}^{[\ell]})\odot
(W^{[\ell+1]})^T\boldsymbol\delta^{[\ell+1]}$
\State $\displaystyle \frac{\partial\mathcal{L}_{\hat{\mathbf x}}}{\partial b_j^{[\ell]}}
=\delta_j^{[\ell]}$, \qquad 
$\displaystyle \frac{\partial\mathcal{L}_{\hat{\mathbf x}}}{\partial w_{jk}^{[\ell]}}=
\delta^{[\ell]}_j a_k^{[\ell-1]}$
\EndFor
\EndProcedure
\end{algorithmic}
\end{algorithm}


\begin{algorithm}[tb!]
\caption{Stochastic Gradient Descent with Backpropagation} \label{alg:sdg_backprop}
\begin{algorithmic}[0]
\Procedure{SgdBackProp\,}{training set $S$, initial guesses for $W^{[\ell]}, \ \mathbf b^{[\ell]}$ for $\ell=1,\ldots, L$, constant learning rate $\eta$}
\For{$k=1, \ldots, n_{iter}$}
\State choose $i$ uniformly random in $\{1,\ldots,N\}$
\State \emph{\# $\mathbf{x}_i$ is the current training data point}
\State
\State \emph{\# feed forward step}
\State $\mathbf{a}^{[0]}=\mathbf{x}_i$
\For{$\ell=1,\ldots L$}
\State $\mathbf{z}^{[\ell]}=W^{[\ell]}\mathbf{a}^{[\ell-1]}+\mathbf{b}^{[\ell]}$
\State $\mathbf{a}^{[\ell]}=\sigma(\mathbf{z}^{[\ell]})$
\State $\mathbf s^{[\ell]}=\sigma'(\mathbf{z}^{[\ell]})$
\EndFor
\State
\State \emph{\# backpropagation step}
\State $\boldsymbol\delta^{[\ell]}=\mathbf s^{[\ell]} \odot(\mathbf{a}^{[\ell]}-\mathbf{y}_i)$
\For{$\ell=L-1,\ldots,1$}
\State $\boldsymbol\delta^{[\ell]}=\mathbf s^{[\ell]} \odot((W^{[\ell+1]})^T\boldsymbol\delta^{[\ell+1]})$
\EndFor
\State
\State \emph{\# SGD update}
\For{$\ell=L,\ldots,1$}
\State  $W^{[\ell]}\leftarrow W^{[\ell]}-\eta \boldsymbol\delta^{[\ell]}(\mathbf{a}^{[\ell-1]})^T$
\State $\mathbf{b}^{[\ell]}\leftarrow\mathbf{b}^{[\ell]}-\eta\boldsymbol\delta^{[\ell]}$
\EndFor
\EndFor
\EndProcedure
\end{algorithmic}
\end{algorithm}

\subsubsection{Penalty--based regularization}\label{sec:regularization}

Regularizing a loss function aims to reduce overfitting, the typical approach consists of adding a strictly convex term to the loss as follows
\begin{equation}\label{eq:loss_regularized}
\mathcal{L}(\mathbf{w})=\frac{1}{N}\sum_{i=1}^N\frac{1}{2}
\|\hat{\mathbf{y}}_i-\mathbf{f} (\hat{\mathbf{x}}_i;\mathbf{w})\|^2_2
+\lambda H(\mathbf{w}),
\end{equation}
where
$\lambda>0$ is the regularization parameter and $H$ is the so--called \emph{regularization term}. 

The most common choices for $H$ are the \emph{Tichonov regularization} (also known as $L^2-$re\-gularization) 
\begin{equation}\label{eq:tichonov}
H(\mathbf{w})=\|\mathbf{w}\|^2_2=\sum_{j=1}^M w_j^2
\end{equation}
and the \emph{Least Absolute Shrinkage and Selection Operator (LASSO)} ($L^1-$regularization) 
\begin{equation}\label{eq:sparsisty_reg}
H(\mathbf{w})=\|\mathbf{w}\|_1=\sum_{j=1}^M |w_j|.
\end{equation}

Because large weights may lead to neurons that are too sensitive to their inputs and, hence, less reliable when unseen data are presented, the idea is to penalize them by introducing a regularization term.
Typically the regularization is only applied to weights, but not to biases.

From the accuracy point of view, the $L^2-$regularization usually outperforms the $L^1-$re\-gularization, nevertheless, the latter produces sparse solutions, i.e., most of the weights found are null or very small (compared to machine precision), meaning that the corresponding input or neuron is meaningless and can be dropped from the network. Thus, $L^1-$regularization can be used to estimate the features (corresponding to non--null weights) which characterize the application \cite{Aggarwal2023}.

From the computational point of view, we notice that regularization makes a very minor and inexpensive change to the backpropagation algorithm.

\subsubsection{Tuning of hyperparameters}\label{sec:hyperparameters}

Hyperparameters in ML are parameters used to control the behaviour of NNs and optimization algorithms and must be distinguished from the ``parameters'' (weights and biases) of the NN. 
We can distinguish between hyperparameters that control the capacity of the hypothesis space $\mathcal H$, such as: 
\begin{itemize}[noitemsep]
\item number of layers of a NN, 
\item number of neurons per layer,
\end{itemize}
and hyperparameters that control the optimization process, such as: 
\begin{itemize}[noitemsep]
\item number of iterations,
\item number of epochs,
\item size of the minibatches (if used), 
\item learning rate in SGD if it is taken constant, 
\item momentum constant $\rho$, 
\item type of regularization,
\item regularization coefficient.
\end{itemize}

Hyperparameters are not subject to the same optimization process of weights and biases, but their tuning is performed either manually by trial and error analysis, by exploiting experience or literature data, or automatically by exploiting further (computationally intensive) optimization algorithms. In this latter case, we speak about encapsulating tuning. 

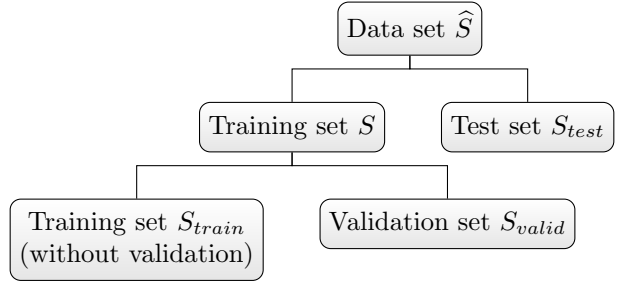
\begin{figure}
    \centering
    \begin{forest}
  for tree={
    l sep'+=.2cm,
    s sep'+=.5cm,
    rounded corners,
    draw, 
    align=center,
    top color=white, 
    bottom color=gray!20,
  },
  forked edges,
  my label/.style={
    edge label={
      node [midway, below right, align=center] {#1}
    },
  }
   [Data set $\widehat S$,edge label={node[midway,left]{z=0}} 
    [Training set $S$,  
        [Training set $S_{train}$\\ (without validation), ]
        [Validation set $S_{valid}$, ]
    ] 
    [Test set $S_{test}$, for children={font=\bfseries} ] 
   ]
  \end{forest}
    \caption{Subdivision of the dataset in training, test and validation sets}
    \label{fig:cross-validation}
\end{figure}

Let us consider the manual tuning of the hyperparameters controlling the hypothesis space. The idea consists of extracting from the training set $S$ a subset $S_{valid}$, named \emph{validation set} of $N_{valid}$ samples to be used to monitor the performance of the models depending on different instances of the hyperparameters, and then in reducing the real training set to $S_{train}=S\setminus S_{valid}$ (see Fig. \ref{fig:cross-validation})

\begin{figure}
    \centering
    \includegraphics[width=0.5\linewidth]{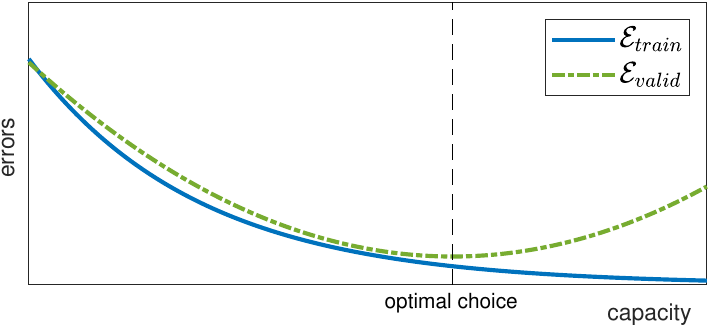}
    \caption{Comparison of training and validation empirical risks for hyperparameters tuning}
    \label{fig:hypertuning}
\end{figure}

By comparing the empirical risks 
\begin{eqnarray}
\mathcal E_{train}=R_{S_{train},N_{train}}(\mathbf f)=\frac{1}{N_{train}} \sum_{i=1}^{N_{train}}
d_M(\hat{\mathbf y}_i,\mathbf f(\hat{\mathbf x}_i))\\
\mathcal E_{valid}=R_{S_{valid},N_{valid}}(\mathbf f)=\frac{1}{N_{valid}}\sum_{i=N_{train}+1}^{N_{train}+N_{valid}}
d_M(\hat{\mathbf y}_i,\mathbf f(\hat{\mathbf x}_i))
\end{eqnarray}
of training and validation sets versus the capacity of the hypothesis space (see Fig. \ref{fig:hypertuning}), typically it happens that, while the training error is almost decreasing versus the capacity of $\mathcal H$, the validation error initially decreases and then grows up. We can identify the optimal values of the hyperparameters as those providing the minimum validation error. For larger capacity, we fall into the overfitting regime.

Once the hyperparameters have been chosen, the test set $S_{test}$ (which contains samples unseen from the set $S$) is used to estimate the expected risk for the corresponding model.

It is not rare that the dataset, from which we extract the training, the validation, and the test sets, is small. In this case, splitting the dataset into fixed subsets can be problematic because of the statistical uncertainty around the estimated average errors. 

\emph{$k-$fold cross--validation} is a way to overcome this problem. In the case of hyperparameters tuning, for a fixed instance of the hyperparameters, it consists of repeating the training and validation on different disjoint couples of training and validation subsets chosen within the dataset, and then computing the average of the validation errors. Repeating the procedure on different instances of the hyperparameters, the optimal choice is the one providing the smallest averaged validation error.
A similar procedure can be used to optimize the hyperparameters controlling the optimization solver.

Let us denote by $\boldsymbol\vartheta$ the array collecting all the unknown hyperparameters and $\boldsymbol\Theta=\{\boldsymbol\vartheta\}$ the set of all the instances $\boldsymbol\vartheta$ we want to consider. 
Algorithm \ref{alg:kcross} reports a version of $k-$fold cross--validation with given instances of the hyperparameters. Alternatively, starting from an instance of $\boldsymbol\vartheta$, after the loop on $j$, we can compute a new value of $\boldsymbol\vartheta$ with some heuristics or optimization algorithms and repeat lines from 4 to 10 with the new hyperparameters.
Once $\boldsymbol\vartheta^*$ has been found, the model with this hyperparameters setting is trained on the original training set $S$.

\begin{algorithm}[tb!]
\caption{$k-$fold cross validation} \label{alg:kcross}
\begin{algorithmic}[1]
\Procedure{kFoldCrossValidation\,}{dataset $S$, hyperparameters set $\boldsymbol\Theta$, $k\in\mathbb N$}
\State Split $S$ in $k$ non--overlapping subsets, such that  $S=\cup_j S_j$
\For{$\boldsymbol\vartheta \in \boldsymbol\Theta$}
\For{$j=1, \ldots, k$}
\State $S_{valid}=S_j$  (validation set)
\State $S_{train}=S\setminus S_{valid}$ (training set)
\State train the model on $S_{train}$
\State validate the model on $S_{valid}$ and compute $\mathcal E_{valid}(j,\boldsymbol\vartheta)$
\EndFor
\State $\overline{\mathcal E}_{valid}(\boldsymbol\vartheta)=\frac{1}{k}\sum_{j=1}^k \mathcal E_{valid}(j,\boldsymbol\vartheta)$ (averaged validation error)
\EndFor
\State \Return{$\boldsymbol\vartheta^*=\argmin{\boldsymbol\vartheta\in \boldsymbol\Theta}\overline{\mathcal E}_{valid}(\boldsymbol\vartheta)$}
\EndProcedure
\end{algorithmic}
\end{algorithm}

Common techniques to choose the hyperparameters to test include \emph{grid search}, which exhaustively tries combinations of hyperparameters, and random search, where random configurations are sampled. More advanced methods, such as \emph{Bayesian optimization}, adaptively focus the search based on previous results, often leading to faster convergence to optimal hyperparameters.

It is important to note that different metrics can be used for evaluating model performance during training and validation phases, depending on the specific task or the goal of the model. For instance, during training, a loss function such as mean squared error (MSE) or cross--entropy might be minimized, while during validation, one might monitor metrics such as accuracy, precision, recall, or F1--score \cite{Goodfellow-et-al-2016}.
These validation metrics offer more meaningful insights, especially for imbalanced datasets or tasks with varying class distributions.

\subsection{A quick glance at Deep Learning models}\label{sec:deep-models}


Deep Neural Networks (DNNs) can be designed with layers of several types, each with a specific task.

\subsubsection{Model components}\label{sec:models-components}

In this Section, we briefly resume the most common types of layers that constitute a deep learning model. 
From a mathematical perspective, a layer is a linear or non--linear operator that maps inputs (the output of the previous layer) to outputs (the input of the next layer). A layer may or may not depend on parameters and/or hyperparameters. In general, a model is the composition of several layers.
Far from being exhaustive, we refer e.g. to \cite{Goodfellow-et-al-2016, Aggarwal2023, Fleuret2024} for an in--depth description of these topics.

\begin{description}
\item[Fully connected layers.] Its neurons receive the output indiscriminately from all the neurons of the previous layer and completely ignore a possible signal structure, without putting any restrictions on which connections are possible and which are not.
At the algebraic level, they implement the basic affine transformation (\ref{eq:weighted-input}). Typically, they are followed by a non--linear transformation (activation function), acting component--wise. See Fig.~\ref{fig:layers}~(a).

\item[Convolutional layers.] They are designed to extract features from grid--structured inputs (e.g., 2D image, but also text, time--series and sequences) which have strong spatial dependences in local regions of the grid. They extract similar feature values from local regions with similar patterns, by applying repeatedly a kernel matrix to the input signal. For instance, let us consider a digital image of $32\times 32$ pixels and denote by $X\in \mathbb R^{32\times 32}$ the array whose entries are the colour intensity (for simplicity, we suppose that the image has only one channel, i.e., it is a grey-scale image). The kernel is a small matrix $K\in\mathbb R^{k\times k}$, with $k$ much smaller than both the width and height of $X$ and whose entries are optimized during the training of the NN, i.e., they are the trainable parameters. Let us take $k=3$, add a ring of null elements to the matrix $X$ (this operation is named \emph{padding}), and call with $\tilde X$ the extended matrix whose indices go from 0 to 33. 
The kernel $K$ interacts with each sub--matrix $\tilde X^{ij}$ of size $3\times 3$ of $\tilde X$, centred at the pixel $(i,j)$ with $i,j=1,\ldots, 32$ and stride $s$, as follows:
\begin{equation}
    c_{\hat\imath\hat\jmath}=\sum_{m=1}^k\sum_{n=1}^k K_{mn} \tilde X^{ij}_{mn}, \qquad\mbox{ with }  \hat\imath=\frac{i-1}{s}+1,\ \hat\jmath=\frac{j-1}{s}+1.
\end{equation}
Then the values $c_{\hat\imath\hat\jmath}$ are stored in a new array containing the output of the layer. In Fig. \ref{fig:layers} (b) the blue sub--matrix is $\tilde X^{ij}$, while the red pixel is the corresponding $c_{\hat\imath\hat\jmath}$. In more general situations, the size $k$ of the kernel, the number of added rings during the padding operation, and the stride $s$ used to go through the matrix $X$ can be arbitrary values and are hyperparameters of the network.

\item[Pooling.] It is a strategy to reduce the signal size (and the number of neurons of the successive layer). Ideally, it summarizes the information from more neurons, for example by extracting the maximum from a set or computing the average. See Fig. \ref{fig:layers} (c). 

\item[Dropout.] It is a layer that provides regularization through random and independent removal of some neurons from the network \cite{srivastava2014dropout}. The goal is to reduce overfitting and prevent neurons from being too specialized. This layer is only applied during the training phase and can also be interpreted as a noise injection, making the training step more robust. It requires a hyperparameter $p$ expressing the percentage of neurons that are turned off. See Fig. \ref{fig:layers} (d).

\item[Normalizing layers.] They facilitate the training by forcing the empirical mean and the variance of the output to 0 and 1, respectively. In particular, we refer to batch normalization layers introduced in \cite{ioffe2015batch}.
\begin{figure}
    \begin{center}
    \begin{subfigure}[b]{0.4\textwidth}
    \includegraphics[trim=0 50 0 0, width=\linewidth]{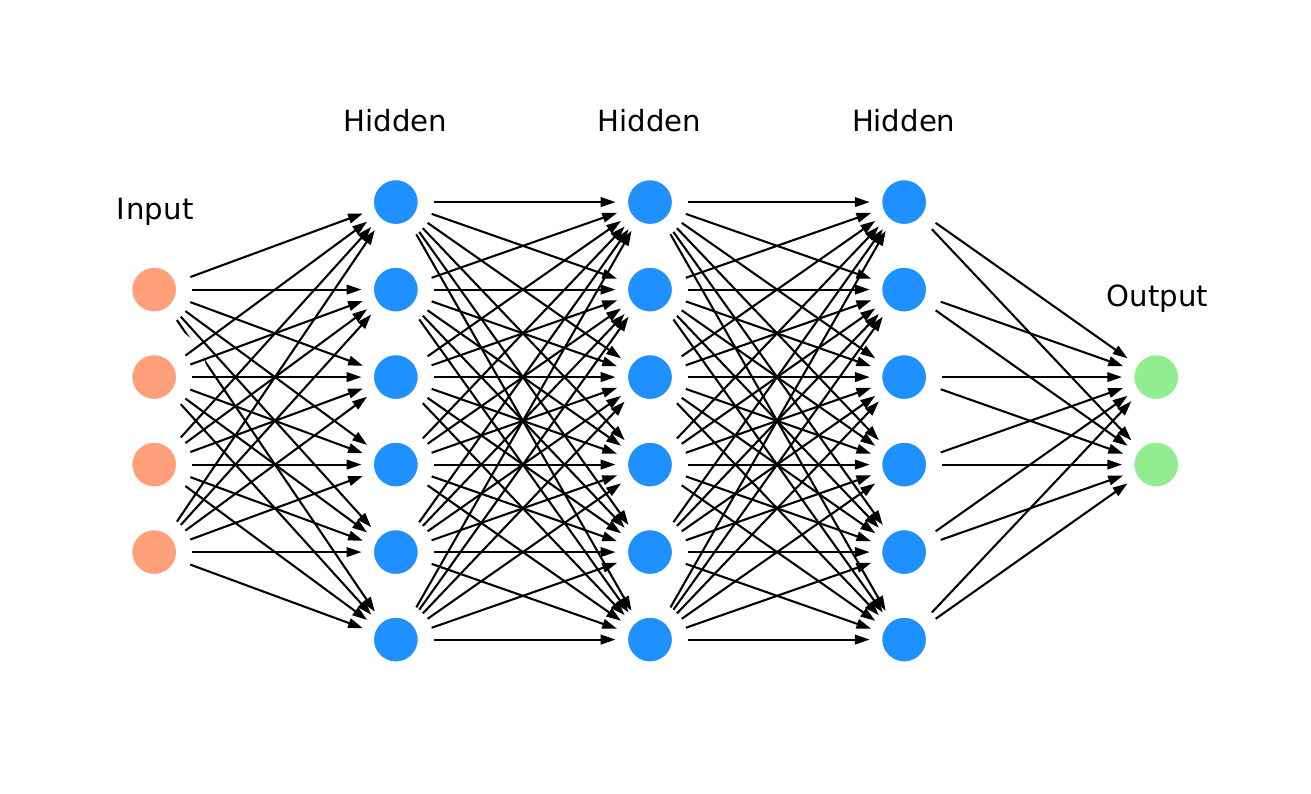}
    \caption{}
    \end{subfigure}\quad
    \begin{subfigure}[b]{0.3\textwidth}\includegraphics[width=\linewidth]{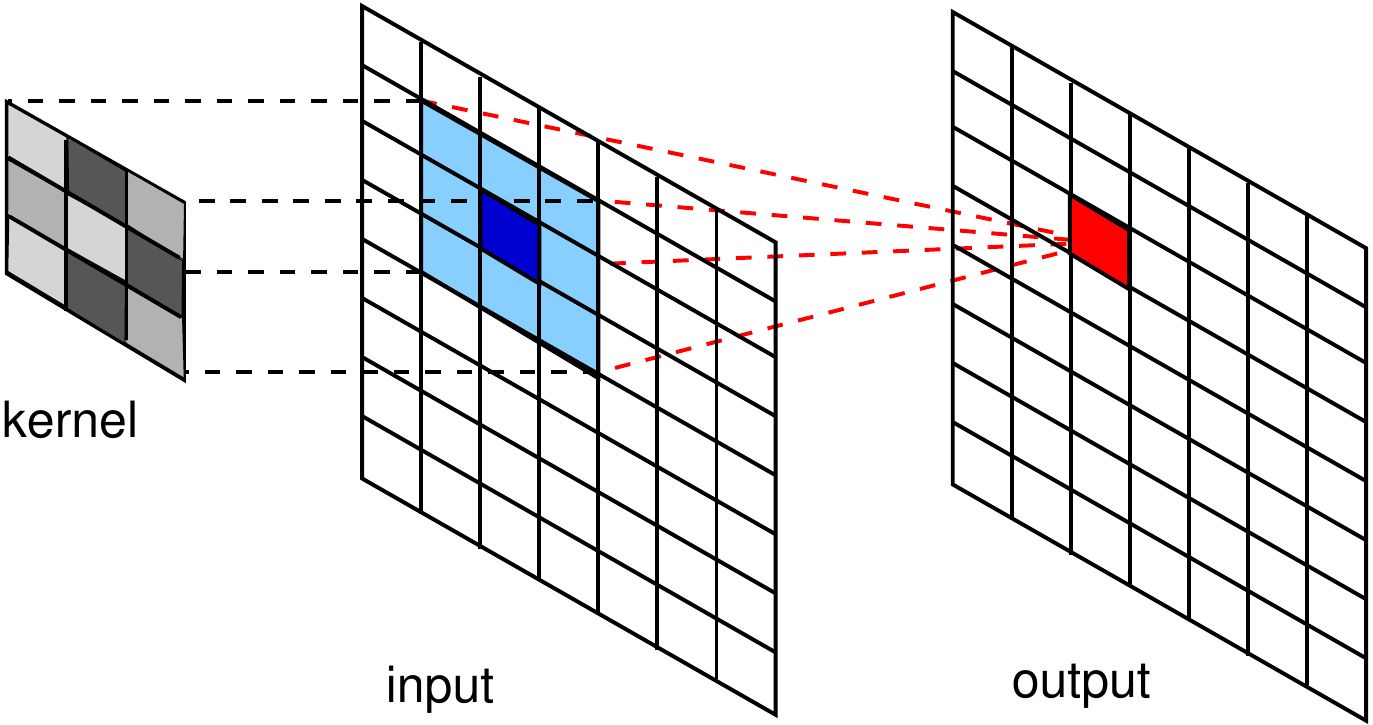}
    \caption{}
    \end{subfigure}\\
    \begin{subfigure}[b]{0.25\textwidth}\includegraphics[width=\linewidth]{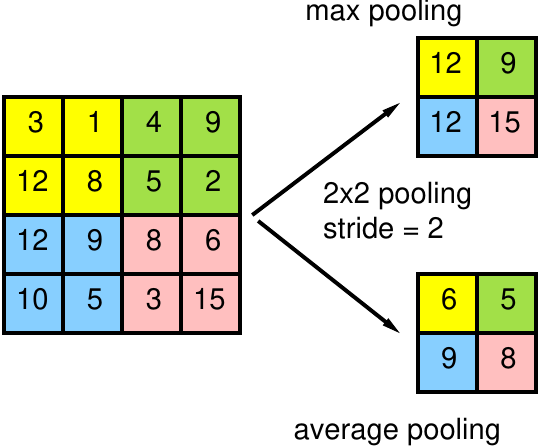}
    \caption{}
    \end{subfigure}\hskip 1.cm
    \begin{subfigure}[b]{0.4\textwidth}\includegraphics[trim=0 50 0 0, width=\linewidth]{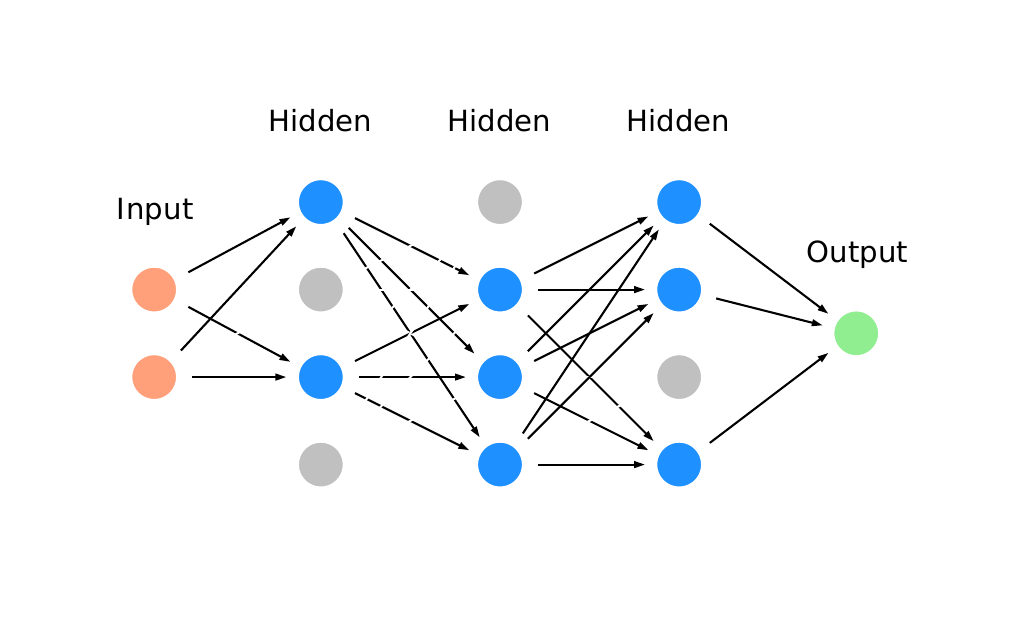} 
    \caption{}
    \end{subfigure}
    \end{center}
    \caption{Some of the most common model layers in NNs. Fully connected layers (a), a convolutional layer (b), pooling (c) and dropout (d).}
    \label{fig:layers}
\end{figure}

\item[Attention layers.] They model dependencies between tokens in a sequence. A token could be a unit of text (such as a word, punctuation mark, or even a subword) or a numeric array.
Each token interacts with every other token, allowing the model to compute the importance (or ``attention weight'') of each one relative to others, regardless of their position in the sequence. They are mainly used in Large Language Models (to make grammatical or semantic decisions in Natural Language Processing), but also in Graph Neural Networks and image processing \cite{Liu2021}. The attention mechanisms have been introduced in \cite{Bahdanau2015} in the context of Recurrent Neural Networks (see next paragraph) for translation tasks. Since then, the idea of attention has evolved and has been generalized to multi--head attention mechanisms for the design of the Transformer in \cite{Vaswani2017} and later to image processing as in the Swin Transformer \cite{Liu2021}. Among attention mechanisms, there are \emph{self--attention} ones, whose goal is to compute a representation of the input sentence.

An attention layer acts as follows. Given a sentence made of $n$ tokens, the first step consists of converting each token into its corresponding \emph{word embedding}, which is a real array of fixed length $n_{model}$. 
During this step, called \emph{embedding}, a matrix $X\in \mathbb R^{n\times d_{model}}$ is created, by using ad--hoc NNs (like Word2Vec \cite{mikolov2013word2vec} and GloVe \cite{pennington2014glove}) or precomputed vocabularies available on the web. Word embeddings are generated with the goal that similarity between words is measured by a suitable distance between the corresponding word embedding vectors.

Unfortunately, the word embedding matrix $X$ is oblivious to the absolute position of the tokens in the sentence, which is instead crucial information to analyse the input correctly. This drawback is overcome by adding to the $i-$th row  of $X$ the \emph{positional encoding} of the $i-$th token of the sentence: given $n$ vectors of length $d$, the associated positional encoding $PE_i$ of the token $i$ is a row vector of length $d$ that can be defined by \cite{Vaswani2017}:
\begin{eqnarray}\label{eq:positional-encoding}
PE_{i,j}=\left\{\begin{array}{ll}
\sin\left(\frac{i}{T^{j/d}}\right) & \mbox{for } 0\leq j\leq d-1\mbox { even}\\
\cos\left(\frac{i}{T^{(j-1)/d}}\right) & \mbox{for }  1\leq j\leq d-1\mbox{ odd},
\end{array}\right.
\end{eqnarray}
where $T=10^4$. Alternatively, positional encodings could be learned during the training, however, in \cite{Vaswani2017} it is shown that the latter approach produces nearly identical results to those obtained with (\ref{eq:positional-encoding}).

The second step consists of building three matrices: the queries $Q$, the keys $K$ and the values $V$. In the case of self--attention mechanisms, they are generated by projecting the word embeddings matrix $X$ onto three linear spaces using three trainable matrices $\widetilde W_Q,\ \widetilde W_K\in \mathbb R^{d_{model}\times d_k}$ and $\widetilde W_V\in \mathbb R^{d_{model}\times d_v}$ so that
\begin{equation}\label{eq:QKV-selfattention}
    Q=X\, \widetilde W_Q, \qquad K=X\, \widetilde W_K,\qquad V=X\, \widetilde W_V.
\end{equation}
$Q\in \mathbb R^{n\times d_k}$ is the query matrix, $K\in \mathbb R^{n\times d_k}$ the key matrix, and $V\in \mathbb R^{n\times d_v}$ the value matrix.
In other attention mechanisms, $Q$, $K$, and $V$ can be computed starting from different sources.

Let us consider the attention layer introduced in \cite{Vaswani2017}: given the matrices $Q,\ K,$ and $V$, the output of the attention layer is the matrix $A\in \mathbb R^{n\times d_v}$
\begin{equation}\label{eq:attention}
    A=  Attention(Q,K,V)=\text{softmax}\left(\frac{Q \,K^T}{\sqrt{d_k}}\right)V.
\end{equation}
Given an array $\mathbf x\in \mathbb R^n$, the softmax function is defined by
\begin{equation}\label{eq:softmax_def}
    (\text{softmax}(\mathbf x))_i=\frac{e^{\mathbf x_i}}{\sum_j e^{\mathbf x}_j}, \qquad i=1,\ldots, n.
\end{equation} 
In (\ref{eq:attention}), the softmax acts row by row on the matrix argument.
A representation of the attention layer is given in Fig. \ref{fig:attention}, left. By computing the matrix product $Q\,K^T$ (i.e., the dot products between any possible pairs of vectors of $Q$ and $K$), we are evaluating the (dis)similarity between queries and keys, then the product between the softmax result and the matrix $V$ provides the projection on the values set. The attention layer (\ref{eq:attention}) is also named \emph{Scaled Dot--Product Attention}.

In \emph{multi--head attention} layers, $n_h$ attention outputs $A_i$ are computed as in (\ref{eq:attention}) after projecting $Q$, $K$ and $V$ on $n_h$ different corresponding spaces:
\begin{eqnarray}\label{eq:multihead-attention1}
\begin{array}{l}
\mbox{for } i=1,\ldots, n_h,\\[2mm]
\qquad     A_i= Attention(Q W_i^Q,K W_i^K, V W_i^W),
\end{array}
\end{eqnarray}
where $W_i^Q,\ W_i^K\in \mathbb R^{d_{model}\times d_k}$,
 and $W_i^V\in \mathbb R^{d_{model}\times d_v}$ are trainable arrays. Then, the outputs $A_i$ are concatenated and projected on the ultimate space as follows:
\begin{eqnarray}\label{eq:multihead-attention2}    
H&=&MultiHeadAttention(Q,K,V)=\text{concat}(A_1,\ldots,A_{n_h})W^O,
\end{eqnarray}
where $W^O\in \mathbb R^{n_h d_v\times d_{model}}$ is another trainable array, see Fig. \ref{fig:attention}, right.
The idea of multi--head attention mechanisms is that each set $\{W_i^Q,\ W_i^K,\ W_i^V\}$ might incorporate different representations of the input (as different heads elaborate them) allowing a deeper and complete analysis of the data.

\end{description}

\begin{figure}
    \begin{center}
    \includegraphics[width=0.8\linewidth]{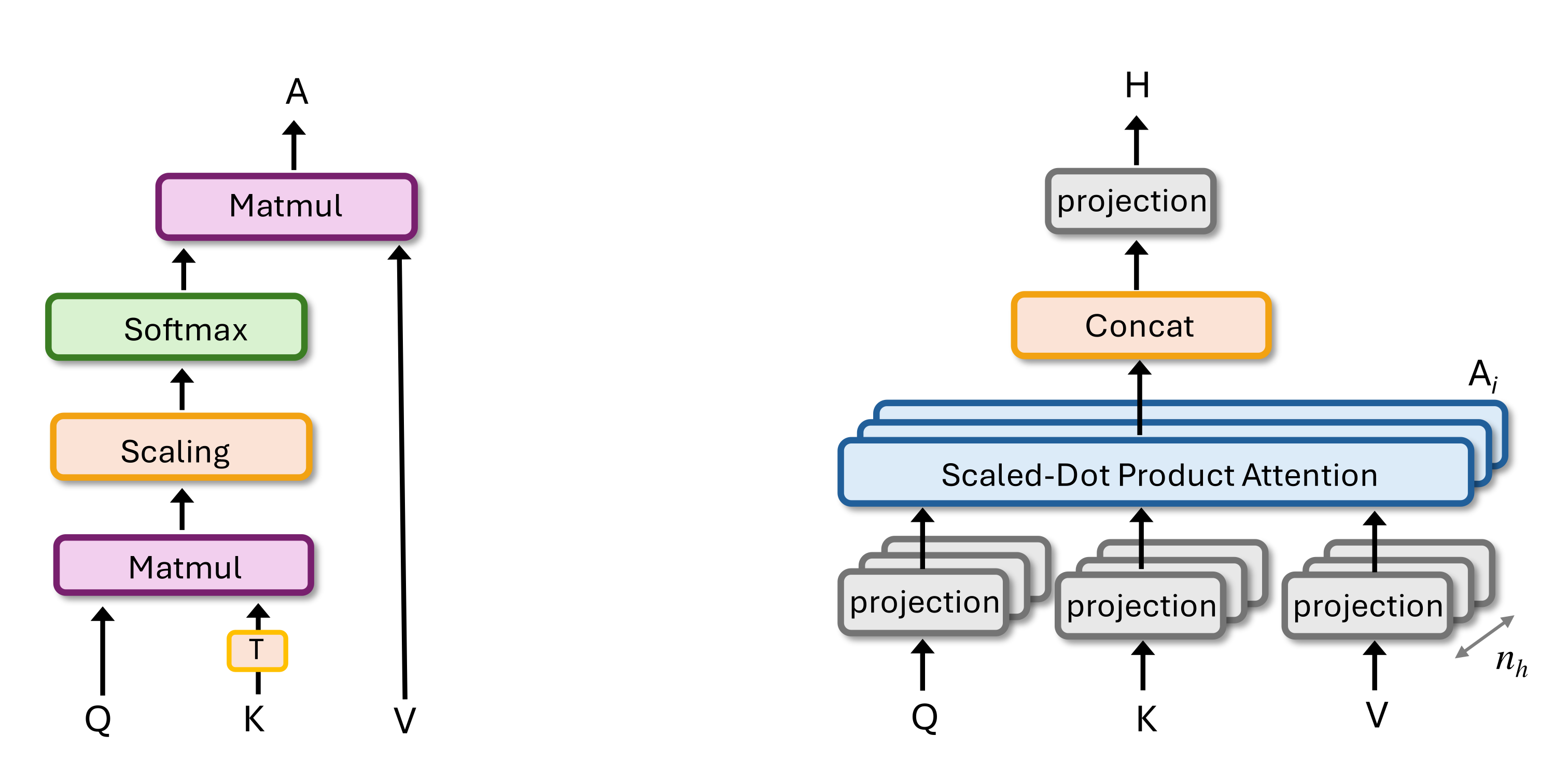}
    \end{center}
    \caption{Scaled Dot--Product Attention layer (\ref{eq:attention}) on the left, and Multi--Head Attention layer (\ref{eq:multihead-attention1})--(\ref{eq:multihead-attention2}) on the right}
    \label{fig:attention}
\end{figure}

\subsubsection{Architectures}\label{sec:architectures}

In this Section, we review the most acclaimed deep learning architectures nowadays. We classify them starting from the learning type: supervised and unsupervised. The major domains of applications are Natural Language Processing (NLP) and Computer Vision.

A remark is appropriate before presenting the list: the vast diversity of neural network architectures makes it challenging to apply rigid classifications, as many models defy straightforward categorization. Architectures traditionally associated with supervised learning, such as fully connected or convolutional networks, are also frequently found in unsupervised or semi-supervised contexts. Similarly, layers initially designed for one specific task can often be repurposed or combined within different learning paradigms, as seen with autoencoders, which can integrate a variety of layers--fully connected, convolutional, or graph neural layers--despite their unsupervised training approach.
This fluidity illustrates that any classification scheme for neural networks should be approached with flexibility and a recognition of the inherent overlaps across architectures. Rather than a rigid taxonomy, it may be more accurate to view these categories as guiding frameworks, acknowledging that combinations and adaptations will often transcend traditional boundaries.

\null\textbf{Supervised learning architectures.}
\begin{description}
\item [Feed Forward Neural Networks (FFNN)] are the simplest NNs. They are formed by fully connected layers in which connections between neurons do not form cycles, see Sect. \ref{sec:ML-models}. Applications are classification, regression, and pattern recognition tasks (see Fig. \ref{fig:layers} (a)). 
\item [Convolutional Neural Networks (CNN)] are designed to work with grid-structured inputs which have strong spatial dependences in local regions of the grid (e.g. 2D image, but also text, time--series and sequences). The goal is to extract features from inputs and CNNs tend to create similar feature values from local regions with similar patterns. CNNs are mainly composed of convolutional layers, but also of fully connected layers, pooling and dropout, see Fig. \ref{fig:cnn}. Applications are image classification, object detection, and facial recognition. Instances are: LeNet (1998, used for digit recognition),
AlexNet (2012), VGG16 (2014), ResNet (2015), GoogLeNet/InceptionV1 (2014). 
\begin{figure}
    \centering
    \includegraphics[width=0.7\linewidth]{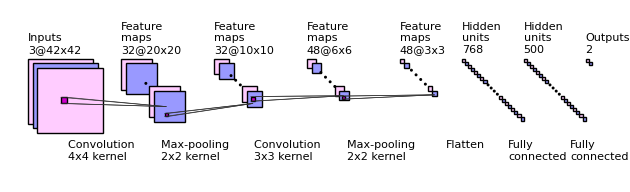}
    \caption{An example of Convolutional Neural Network}
    \label{fig:cnn}
\end{figure}

\item [Recurrent Neural Networks (RNN)] are designed to process sequential data of any length. They factor the computation in time--steps at which a hidden state $\mathbf h_t$ is generated as a function of both the previous state $\mathbf h_{t-1}$ and input $\mathbf x_t$. The hidden state $\mathbf h_t$ implicitly keeps the memory of the whole history of the process and contributes to generating the output $\mathbf y_t$. More precisely, at the generic time $t$, a traditional RNN computes
\begin{eqnarray*}
    \mathbf h_t&=&\sigma(W\mathbf h_{t-1}+U\mathbf x_t +\mathbf b)\\
    \mathbf o_t&=& V\mathbf h_t+\mathbf c\\
    {\mathbf y}_t&=&\text{softmax}(\mathbf o_t),
\end{eqnarray*}
where the softmax function has been defined in (\ref{eq:softmax_def}) and is used to interpret the entries of ${\mathbf y}_t$ as probabilities.
The trainable weights matrices $W$, $U$, and $V$, as well as the biases $\mathbf b$ and $\mathbf c$, are the same for any $t$, thus the total number of trainable parameters remains fixed even for longer inputs, see Fig. \ref{fig:rnn}. 

Let us suppose that the RNN has to foresee the word of a sentence in position $t$, this means that we have $t$ inputs $\mathbf x_1,\ldots,\mathbf x_t$ that are word embeddings 
and we want to foresee the next word $\mathbf x_{t+1}$. The output $\mathbf y_t$ is a vector whose generic entry $(\mathbf y_t)_i$ is the probability of the $i-$th word of the vocabulary to be the next word. This is an example of \emph{sequence--aligned} RNN (at each time--step the input is required and the output is produced), while other examples of RNNs can contemplate input only at the first time--step (a possible application is in music generation) or output only at the last (a possible application is in sentiment classification).

Unfortunately, training traditional RNNs using backpropagation is cumbersome because the loss function suffers from exploding and vanishing gradients, especially when dealing with long sentences. Gated RNNs, Long Short--Term Memory (LSTM, 1997), and Gated Recurrent Units (GRU, 2014) bypass such drawbacks, however, the sequential nature of RNNs precludes parallelization within training. Applications of RNNs are time series prediction, language modelling, machine translation, music generation, and speech recognition. They can also be considered to simulate dynamical systems and solve time--dependent differential equations (see Sect.~\ref{sec:operator-learning-time-dependent}). 
\begin{figure}
    \centering
    \scalebox{0.6}{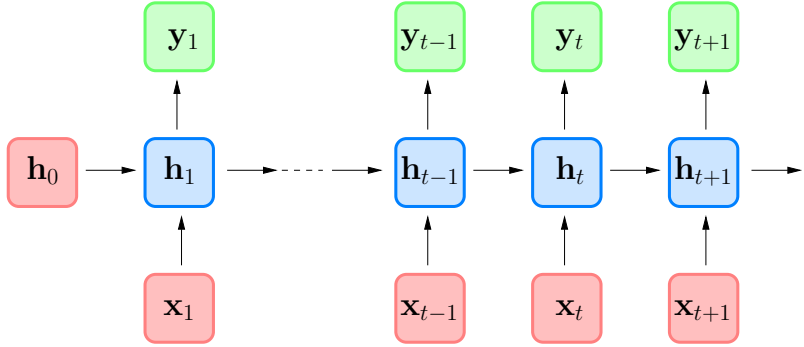}
    \caption{The graph of a traditional Recurrent Neural Network}
    \label{fig:rnn}
\end{figure}
\item [Autoregressive models] are architectures with a structure similar to that of RNNs. However, while RNNs feature an internal variable $\mathbf h_t$ which keeps the memory of the system, i.e. for any $t$, $\mathbf h_t$ feeds $\mathbf h_{t+1}$, in autoregressive models the system's memory is explicitly transferred by the output variable. This means that the output $\mathbf y_{t+1}$ explicitly depends on previous outputs typically belonging to a temporal window of fixed length $k$, i.e., $\mathbf y_{t+1}=f(\mathbf y_t,\ \mathbf y_{t-1},\ldots,\ \mathbf y_{t_k})$. 

\item [Graph Neural Networks (GNN)] 
are designed to work directly on graph--structured data, where nodes represent objects or concepts and edges represent relationships between them. We refer to \cite{Wu20214} for a survey on GNNs. They are applied to detect objects in images, analyse social networks, develop recommendation systems, and study and design compounds in chemistry and biology.
Examples of GNNs are Graph Convolutional Networks (2016), Graph Attention Networks (2017), and Message Passing Neural Networks (2017).

To briefly present GNNs, we follow the setting of \cite{Gori2005, Scarselli2009}.
GNNs work on data represented in graph domains $G=(N,E)$, where $N$ and $E$ are the sets of nodes and edges, respectively. 
To give an example, in object detection tasks, a node represents a homogeneous region of the image, while an edge is established between adjacent regions. Typically, some labels expressing features are associated with nodes and edges, like, e.g. the area or the colour of a homogeneous region. Labels are stored in an array $\boldsymbol\ell$ and $\boldsymbol\ell_S$ denotes the labels of all the entities of the subset $S$ of the graph. Finally, each node $n$ is characterized by a state $\mathbf{x}_n$, which is an array that contains a representation of the associated object or concept.

The state $\mathbf x_n$ depends on the states and labels of the neighbouring nodes and edges, thus all states are strictly related one each other. We denote by $\text{ne}[n]$ the set of nodes in the neighbourhood of the node $n$, and by $\text{co}[n]$ the set of the edges having the node $n$ as an endpoint. Then, we express relations between the entities of the graph by 
\begin{equation}\label{eq:gnn-state}
    \mathbf x_n=f_{\mathbf w}(\boldsymbol \ell_n, \boldsymbol\ell_{\text{co}[n]}, \mathbf x_{\text{ne}[n]}, \boldsymbol \ell_{\text{co}[n]}).
\end{equation}
$f_{\mathbf w}$ is the so--called \emph{local transition function} depending on a set of learnable weights $\mathbf w$, typically it is implemented by a feed--forward NN.
Then, a learnable local output function $g_{\mathbf w}$, again implemented by a feed--forward NN, provides the output $\mathbf o_n$ of the node $n$:
\begin{equation}\label{eq:gnn-output}
    \mathbf o_n=g_{\mathbf w}(\mathbf x_n, \boldsymbol\ell_n).
\end{equation}
For instance, in object detection tasks, we have to recognize which nodes of the graph $G$ belong to a sub--graph $S$ of $G$, so the output $\mathbf o_n$ will be equal to 1 if the node $n$ belongs to $S$, 0 otherwise.

By stacking all the states, output, labels and labels of the nodes into the arrays $\mathbf x$, $\mathbf o$, $\boldsymbol\ell$, and $\boldsymbol\ell_N$, respectively, the equations
\begin{eqnarray}\label{eq:GNN-system}
\begin{array}{l}
\mathbf x=F_{\mathbf w}(\mathbf x,\boldsymbol\ell)\\
\mathbf o=G_{\mathbf w}(\mathbf x, \boldsymbol\ell_N)
\end{array}
\end{eqnarray}
provide the state and output of each node of the graph. $F_{\mathbf w}$ and $G_{\mathbf w}$ are the \emph{global transition function} and \emph{global output function}, respectively. Notice that the output $\mathbf o_n$ of the node $n$ only depends on its states and labels. See Fig. \ref{fig:GNN}.

Provided that $F_\mathbf w$ is a contraction, (\ref{eq:GNN-system})$_1$ admits a unique solution $\mathbf x$ that is the unique fixed point of $F_\mathbf w$ which is the limit for $t\to\infty$ of the fixed point iterations 
\begin{eqnarray}\label{eq:GNN-fixedpoint}
    \left\{\begin{array}{ll}
    \mathbf x^{(0)} \mbox{ given}\\
    \mathbf x^{(t+1)}=F_{\mathbf w}(\mathbf x^{(t)},\boldsymbol\ell)& t\geq 0.
    \end{array}\right.
\end{eqnarray}
Solving these iterations up to time $T$ (each step $t$ plays the role of a pseudo time instant), the solution $\mathbf x^{(T)}$ is taken as an approximation of $\mathbf x$.

When $f_\mathbf w$ and $g_{\mathbf w}$ are implemented by a feed--forward NN, the connections between neurons of the global network can be divided into internal and external. The internal ones are determined by the functions $f_\mathbf w$ and $g_{\mathbf w}$ used to implement each single unit  (\ref{eq:gnn-state}) and (\ref{eq:gnn-output}), while the external ones are determined by the edges of the graph. 

Parameters learning is achieved by minimizing the Mean Square Error $\mathbf e_{\mathbf w}$ between target and computed values. A gradient--descent method with backpropagation (as seen in Sects. \ref{sec:optimization} and \ref{sec:backpropagation}) can be applied. The fixed point iterations (\ref{eq:GNN-fixedpoint}) constitute the forward step inside each iteration of the gradient method as described in Algorithm \ref{alg:GNN}.
A \emph{Graph Neural layer} corresponds to a time instant and contains a copy of all the entities of the network.
    \begin{algorithm}[tb!]
    \begin{algorithmic}
    \Function{\textsc{Main}}{}
    \State{\emph{\# $\eta$ hyperparameter}}
    \State{initialize $\mathbf w$}
    \State{compute $\mathbf x$=\textsc{Forward}($\mathbf w$)}
    \Repeat
    \State{compute $\nabla \mathbf e_{\mathbf w}$=\textsc{Backpropagation}($\mathbf x, \mathbf w$)}
    \State{update the weights $\mathbf w=\mathbf w-\eta \nabla \mathbf e_{\mathbf w}$}
    \State{compute $\mathbf x$=\textsc{Forward}($\mathbf w$)}
    \Until{a stopping criterion}
    \EndFunction
    \Function{\textsc{Forward}($\mathbf w$)}{}
    \State{\emph{\# $\varepsilon_F$ hyperparameter}}
    \State{initialize $\mathbf x^{(0)}$, $t=0$, $err=\varepsilon_F+1$}
    \While{$err\geq \varepsilon_F$}
    \State{$\mathbf x^{(t+1)}=F_{\mathbf w}(\mathbf x^{(t)},\boldsymbol\ell)$}
    \State{$err=\|\mathbf x^{(t+1)}-\mathbf x^{(t)}\|$}
    \State{$t=t+1$}
    \EndWhile\\
    \Return{$\mathbf x^{(t)}$}
    \EndFunction
    \end{algorithmic}
    \caption{The GNN model}\label{alg:GNN}
    \end{algorithm}

\begin{figure}
    \centering
    \includegraphics[width=0.9\linewidth]{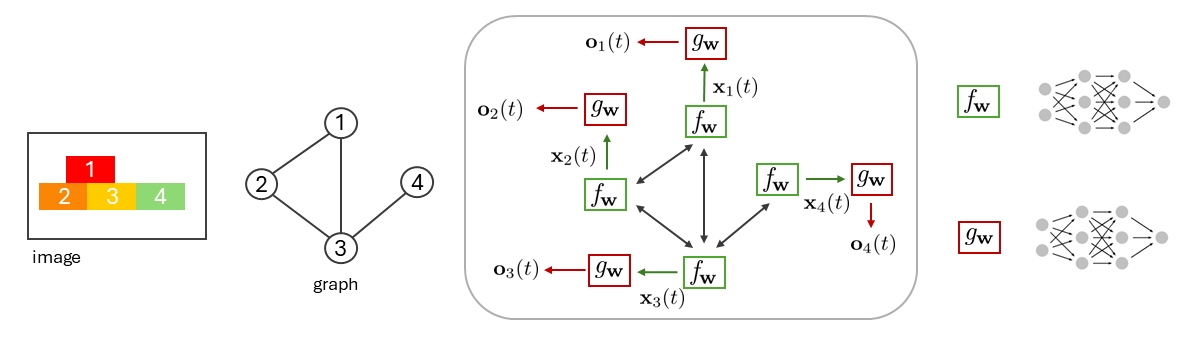}
    \caption{Graph Neural Network. Different coloured regions of an image are associated with the nodes of a graph, the edges express contact between two different regions. The grey box represents the equation (\ref{eq:GNN-system}), where $f_{\mathbf w}$ and $g_{\mathbf w}$ are feed--forward NNs.}
    \label{fig:GNN}
\end{figure}
\item[Encoder--decoder] is a neural network architecture used to transform an input sequence into an output sequence. The encoder compresses the input into an intermediate representation (typically named code), and the decoder generates the output from this encoded representation. It is used in tasks such as automatic translation and text summarization, but also in more complex architectures with different purposes. See Fig. \ref{fig:encoder-decoder}.
\begin{figure}
    \centering
    \includegraphics[width=0.5\linewidth]{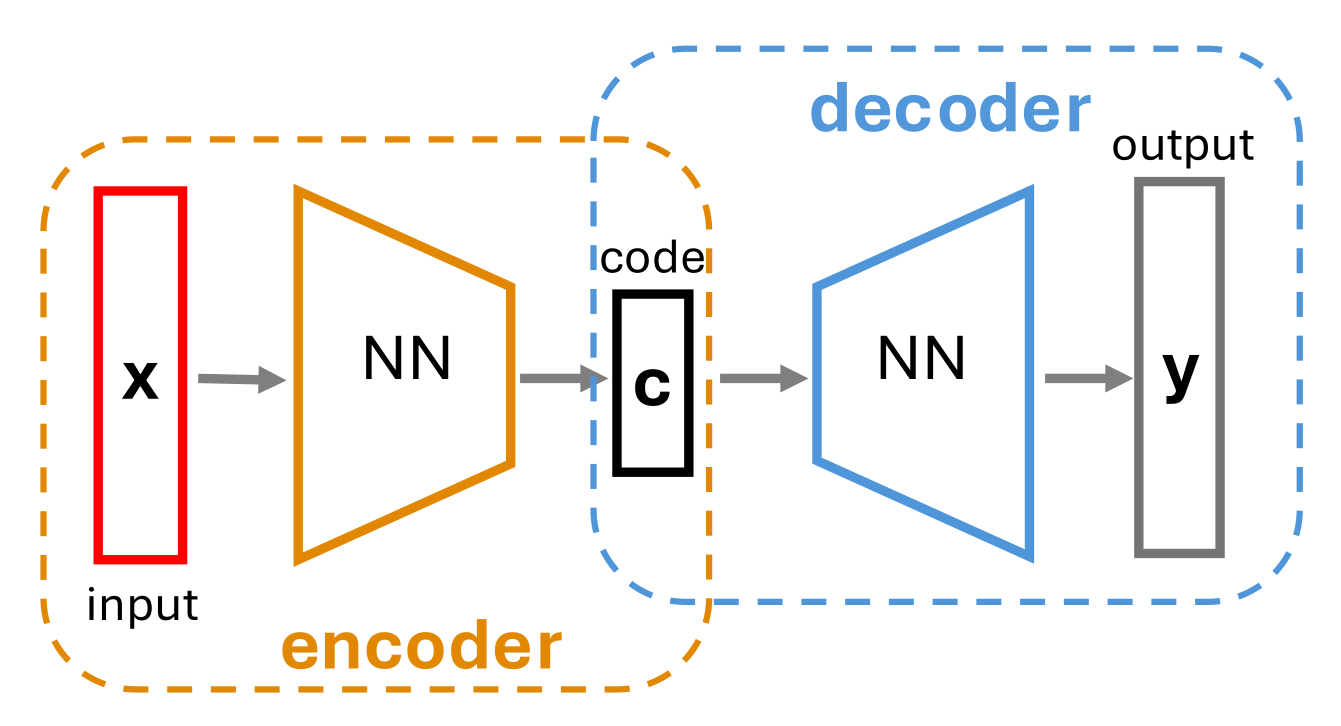}
    \caption{Encoder--decoder architecture}
    \label{fig:encoder-decoder}
\end{figure}
\item [Attention models] use (self--)attention mechanisms to process input data,
they can be composed of multi--head attention layers, FFNNs, and, typically, feature an encoder--decoder structure. The most famous attention model is the Transformer \cite{Vaswani2017} that, with its multi--head attention layers, can exploit parallelism and be highly efficient and scalable. Attention models outperform RNNs, which instead have a sequential structure.

Applications are Natural Language Processing (NLP), machine translation, text generation, and image classification.
Examples are Transformer (2017), BERT (Bidirectional Encoder Representations from Transformers, 2018), GPT (Generative Pre-trained Transformer, 2018), ViT (Vision Transformer, 2020), GPT--3 (2020), PaLM (2022), LaMDA (2022), GPT--4 (2023). 
BERT, GPT--x, LaMDA are instances of Large Language Models (LMM), i.e., deep NNs able to achieve general--purpose language generation and other natural language processing tasks such as classification. Recently they have been used in the context of Scientific Machine Learning \cite{herde2024poseidon}.
\end{description}

Among the attention models, we briefly describe the Transformer, designed for natural language processing, and the Swin Transformer, for image processing.

\begin{description}
    \item[The Transformer] is a transduction model introduced in \cite{Vaswani2017} relying entirely on self--attention mechanisms to compute a representation of its input and output, without using RNNs or convolution. It has an encoder--decoder structure and it is depicted in Fig. \ref{fig:transformer}. These architectures have been introduced in the context of generative AI, but recently they have also been employed in Scientific Machine Learning (SciML). An example of a transformer in SciML is the architecture Poseidon \cite{herde2024poseidon} described in Sect. \ref{sec:operator-learning}. 
    
    A transformer is a complex architecture which includes many layers of different nature, like, e.g., fully connected, normalization, projection, and multi--head attention layers. It acts as follows.
    First, in the \emph{input embedding} phase, the input sentence is decomposed into $n$ tokens which are transformed into the corresponding word embeddings of fixed length $d_{model}$. The matrix $X$ containing the word embeddings is added to the positional encoding matrix $PE$ (see \ref{eq:positional-encoding}) of the tokens and the sum becomes the input of the encoder. The encoder is made of $N$ blocks, each one composed of two sub--layers: one multi--head attention layer (\ref{eq:multihead-attention1})--(\ref{eq:multihead-attention2}) and one feed--forward layer, both of them combined with a residual connection and followed by a normalization layer. If the output of the sub--layer is $y=S(x)$, applying a \emph{residual connection} to the sub--layer itself means replacing $y=S(x)$ with $y=x+S(x)$. Residual connections were introduced in \cite{he2015resnet} to accelerate the backpropagation phase during the optimization process and to address the vanishing gradient problem.

    The decoding phase is an iterative process that stops when the last output token is produced. At each iteration, the input of the decoder is composed of the full output of the encoder and the (partial) output of the previous iteration of the decoder phase (at each iteration, one new token is produced). More precisely, the previous iteration's output is transformed into the corresponding word embeddings and added to its positional encoding before entering the decoder phase. The decoder is composed of $N$ equal sub--blocks, each one composed of three sub--layers. The first one is a masked multi--head attention layer with a residual connection followed by a normalization layer, its input is the output of the decoder at the previous iteration and the masked attribute means that each input token can only interact with the ones on its left. The second sub--layer is another multi--head attention layer with residual connection, followed by a normalization layer. The output of the first sub--layer is used to generate the query matrix, while the full output of the encoder is used to generate the key and value matrices. The third and last sub--layer of each block is a feed--forward layer with a residual connection followed by a normalization layer. At the end of the $N$ blocks, a linear transformation (projection layer) maps the decoder's output dimension to the size of the vocabulary. Successively, a softmax layer provides the probabilities associated with each token of the vocabulary and the next token is selected for the translation.

    The base model of the Transformer has $N=6$ blocks in both the encoder and the decoder, and has been trained on a dataset of about 4.5 million sentence pairs which required about $3.3\cdot 10^{18}$ floating point operations. The Adam method has been used to minimize the loss function.
\begin{figure}
    \centering
    \includegraphics[width=0.6\textwidth]{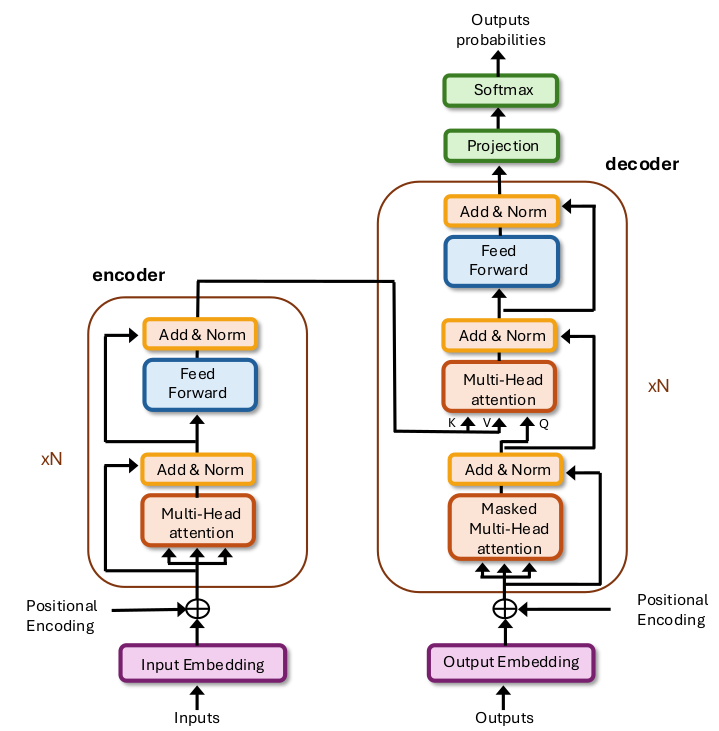}
    \caption{The Transformer\cite{Vaswani2017} }  \label{fig:transformer}
\end{figure}

    \item[Vision Transformer (ViT)]
    is a transformer applied to 2D images to answer tasks like image classification, object detection, and semantic segmentation \cite{dosovitskiy2021}. The idea is as follows: a 2D image $\mathbf x$ of $H\times W$ pixels and $C$ channels is first reshaped into a sequence of $N=HW/P^2$ patches of size $P\times P$ (for instance $P=4$) and $C$ channels, then the patches are flattened to 1D arrays $\mathbf x_p^{(i)}\in\mathbb R^{1\times (P^2\cdot C)}$ which play the role of tokens that feed a traditional transformer. 
    A trainable vector $\mathbf x_{class}\in\mathbb R^{1\times (P^2\cdot C)}$ is usually prepend to the patches $\mathbf x_p^{(i)}$, to achieve the image classification task (like in the BERT model). 
    Then, tokens are projected onto a linear space of dimension $D$ by a matrix $E\in \mathbb R^{D\times (P^2\cdot C)}$
    and the output is added to the positional encodings $E_{pos}\in \mathbb R^{(N+1)\times D}$ of the tokens, like in the Transformer. Denoting by $L$ the total number of layers, by MSA a Multi--Head Self--Attention layer (see Algorithm \ref{alg:MSA}), by LN a normalization layer, and by MLP a MultiLayerPerceptron with one hidden layer, the ViT model reads as in Algorithm \ref{alg:vit}.
    \begin{algorithm}[tb!]
    \begin{algorithmic}
    \State{\emph{Self--Attention Projection}}
    \State{$Q=X\ \widetilde W_Q,\quad K=X\, \widetilde W_K,\quad V=X\, \widetilde W_V$}
    \State{\emph{Multi--Head Attention}}
    \For{$i=1,\ldots, n_h$}
        \State{$Q_i=Q W^Q_i,\ K_i=K W^K_i,\ V_i=V W_i^V$}
        \State{$A_i=Attention(Q_i,K_i,V_i)=
        \text{softmax}\left(\frac{Q_i K_i^T}{\sqrt{d_k}}\right)V_i$}
    \EndFor
    \State{$H=\text{concat}(A_1,\ldots,A_{n_h})W^O$}
    \end{algorithmic}
    \caption{The Multi--Head Self--Attention layer: $H=\text{MSA}(X)$}\label{alg:MSA}
    \end{algorithm}
    \begin{algorithm}[tb!]
    \begin{algorithmic}
    \State{$Z_0=[\mathbf x_{class};\ \mathbf x_p^{(1)}E;\ldots;\mathbf x_p^{(N)} E]+ E_{pos}$}
    \For{$\ell=1,\ldots, L$}
        \State{$\hat Z_\ell=\text{MSA}(\text{LN}(Z_{\ell-1}))+Z_{\ell-1}$}
        \State{$Z_\ell=\text{MLP}(\text{LN}(\hat Z_\ell))+\hat Z_\ell$}
    \EndFor
    \State{$Y=\text{LN}((Z_L)_{1j})$}    
    \end{algorithmic}
    \caption{The ViT model: $Y=ViT(\mathbf x_p^{1},\ldots,\mathbf x_p^{(N)})$}\label{alg:vit}
    \end{algorithm}
    
    \item[Swin Transformer] (Shifted Windows Transformer) is a vision transformer, proposed in \cite{Liu2021}, designed to be a general--purpose backbone network for both image classification and dense recognition tasks. It introduces two new ideas compared to ViT: 
    \emph{(i)} it considers two different partitions of the input image in non--overlapping windows (one, named regular, and the other one shifted compared to the regular) containing $M\times M$ patches (see Fig. \ref{fig:swin1} (b) and (c)) and compute self--attention locally within the non--overlapping windows. These two window partitions are used alternately in consecutive blocks so that the shifted partition allows the detection of connections between the neighbouring non--overlapping windows of the regular one;
    \emph{(ii)} it uses a hierarchical representation of the image, to capture features at different scales, see Fig. \ref{fig:swin1} (a). Starting from an initial partition of the original image in $P\times P$ patches of $4\times 4$ pixels, they are grouped in non--overlapping windows of $M\times M$ patches (for both the regular and shifted partitions). At the next level, patches are merged by concatenating features of each group of $2\times 2$ neighbouring patches, while the number of patches inside each window is kept fixed and equal to $M\times M$. The Swin Transformer architecture is composed of a fixed number of stages, each one corresponding to a level of the hierarchical partition. The larger the number of hierarchical levels, deeper the network. The first stage contains a linear embedding layer followed by an even number of Swin Transformer Blocks (STB), while the successive stages contain a patch merging layer followed by an even number of Swin Transformer Blocks (STB). In Algorithm \ref{alg:swin} a pair of two consecutive STBs, the first on the regular partition and the second on the shifted one, is formulated. Then let LN be a normalization layer, W-MSA a Window Multi--Head Self--Attention layer applied to the regular partition, MLP (Multi Layer Perceptron) a 2--layer feed--forward NN with GELU activation function\footnote{$GELU(x)=x\int_{-\infty}^x f(t)\, dt$, where $f(t)$ is the normal probability density function. It is approximated by the formula  $GELU(x)=0.5\,x\,(1+\tanh[\sqrt{2\pi}\ (x+0.044715\,x^3)])$}, and SW-MSA a Shifted Window Multi--Head Self--Attention layer similar to W-MSA, but applied to the shifted partition. Notice that each window of the two partitions, as well as each head of the multi--head attention layer, has independent trainable matrices and biases.
    %
    \begin{algorithm}[tb!]
    \begin{algorithmic}
    \For{$w=1,\ldots, n_W$} \emph{\quad ($w$ = index of the window)}
    \State{$Q=X^w\ (\widetilde W_Q)^w,\quad K=X^w\, (\widetilde W_K)^w,\quad V=X^w\, (\widetilde W_V)^w$}
    \State{\emph{Multi--Head Attention}}
    \For{$i=1,\ldots, n_h$}\emph{\quad ($i$ = index of the head)}
        \State{$Q_i^w=Q (W^Q)^w_i,\ K_i^w=K (W^K)^w_i,\ V_i=V (W^V)^w_i$}
        \State{$A_i^w=Attention(Q_i^w,K_i^w,V_i^w)=
        \text{softmax}\left(\frac{Q_i^w(K^w_i)^T}{\sqrt{d_k}}+B^w_i\right)V^w_i$}
    \EndFor
    \State{$H^w=\text{concat}(A_1^w,\ldots,A^w_{n_h})(W^O)^w$}
    \EndFor
    \State{$H=\text{concat}(H^1,\ldots,H^{n_W})$}
    \end{algorithmic}
    \caption{The Window Multi--Head Self--Attention layer: $H=\text{W-MSA}(X)$, with $X=\text{concat}(X^1,\ldots, X^{n_W})$ and $X^w$ features in the $w-$th window}\label{alg:W-MSA}
    \end{algorithm}
    %
    %
    \begin{algorithm}[tb!]
        \begin{algorithmic}
            \State{\emph{regular windows partition}}
            \State{$\hat{Z}=\text{W-MSA}(\text{LN}(X))+X$}
            \State{$Z=\text{MLP}(\text{LN}(\hat Z))+\hat Z$}
           \State{\emph{shifted windows partition}}
            \State{$\hat{Y}=\text{SW-MSA}(\text{LN}(Z))+Z$}
            \State{$Y=\text{MLP}(\text{LN}(\hat Y))+\hat{Y} $}
        \end{algorithmic}
    \caption{Two successive Swin Transformer Blocks at a fixed stage. $n_W$ denotes the number of non--overlapping windows of the regular partition, $n_{SW}$ that of the shifted partition}
    \label{alg:swin}
    \end{algorithm}
 Different to ViT, no absolute positional encoding $E_{pos}$ is applied. However, both W-MSA and SW-MSA include trainable relative position bias matrices $B^w_i$ in the attention layer.

    \begin{figure}
    \begin{subfigure}[t]{0.3\textwidth}
        \centering
     \includegraphics[width=\linewidth]{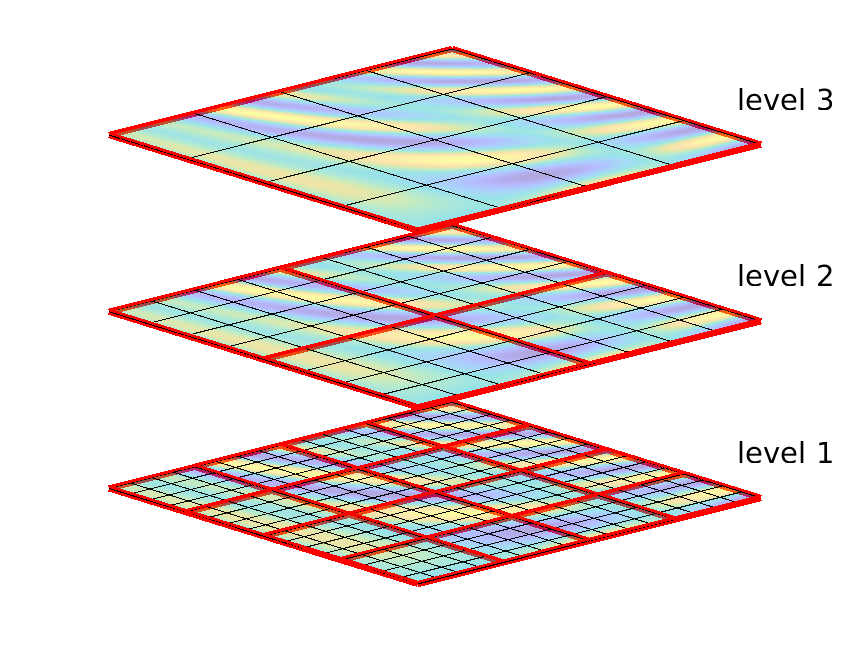}
     \caption{}
     \end{subfigure}
    \begin{subfigure}[t]{0.3\textwidth}
        \centering
     \includegraphics[width=\linewidth]{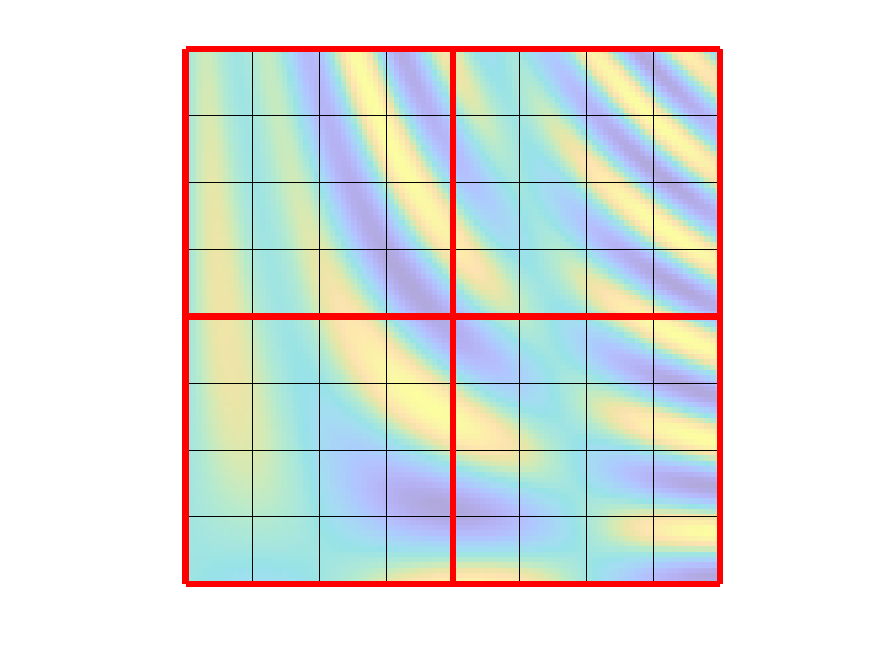}
     \caption{}
     \end{subfigure}
    \begin{subfigure}[t]{0.3\textwidth}
        \centering
     \includegraphics[width=\linewidth]{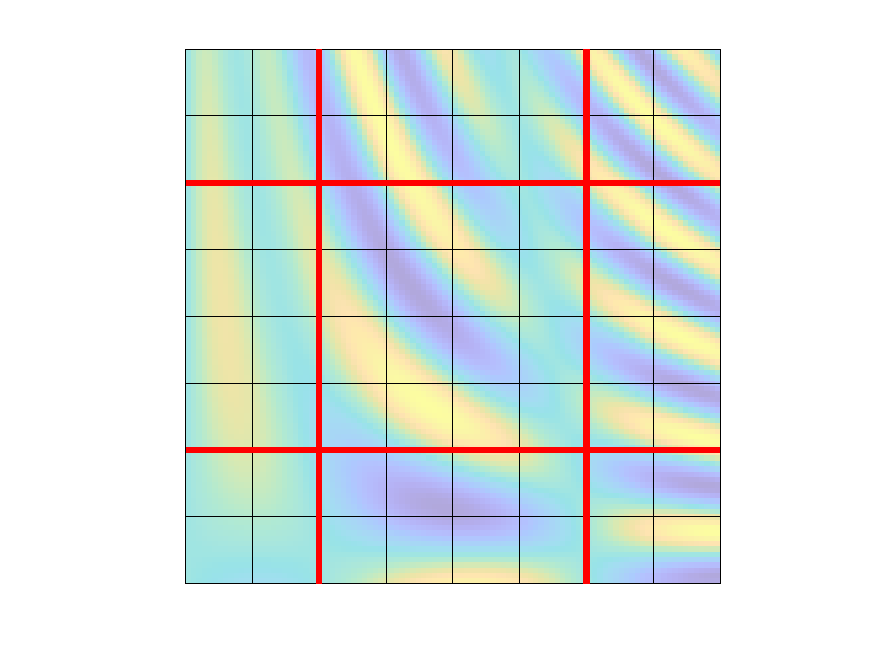}
     \caption{}
     \end{subfigure}
     \caption{Swin Transformer. (a) The hierarchical partition in the Swin Transformer. Red lines bound the windows, while black lines the patches. The first level is the lower one, the number of windows is reduced by $1/4$ moving from one level to the next and the number of features of each patch inside a window is doubled compared to the previous level. (b) The regular partition in windows of $4\times 4$ patches and (c) the corresponding shifted partition.}
     \label{fig:swin1}
      \end{figure}
      
\item[SwinV2 Transformer] has been proposed in \cite{liu2022swinV2} as an improvement of the Swin Transformer, to overcome training instabilities due to exploding activation values and gradients in the presence of deep networks and large images. More precisely, in SwinV2: \emph{(i)} the normalization layers LN are postponed to both the MLP layer and the attention layers W-MHA, SW-MHA; \emph{(ii)} the scaled dot product $Q K^T/\sqrt{d_k}$ in the attention layer (see row 6 of Algorithm \ref{alg:W-MSA}) is replaced by a scaled cosine attention approach:
\begin{eqnarray}\label{eq:W-MHA-cos}
    A_i^w=Attention(Q_i^w,K_i^w,V_i^w)=
        \text{softmax}\left(\frac{\cos(Q_i^w(K^w_i)^T)}{\tau}+B^w_i\right)V^w_i
    \end{eqnarray}
    where given two matrices $A$ and $B$ of the same size and denoting by $\mathbf a_i$ and $\mathbf b_j$ their row vectors, $\cos(A,B)$ is a matrix whose entries are $(\cos(A,B))_{jk}= 
    \frac{\mathbf a_j\cdot \mathbf b_k}{\|\mathbf a_j\|\, \|\mathbf b_k\|}$; $\tau$ is a learnable parameter, and $B^w_i$ is the relative position bias between pixels inside the window; \emph{(iii)} a log--spaced position bias approach is adopted, so that the relative position biases can be smoothly transferred across windows resolutions. SwinV2 also uses a self--supervised pretraining method (SimMIM) to reduce the need for vast labelled images.
    %
    %
    \begin{algorithm}[tb!]
        \begin{algorithmic}
            \State{\emph{regular windows partition}}
            \State{$\hat{Z}=\text{LN}(\text{W-MSA}^{cos}(X))+X$}
            \State{$Z=\text{LN}(\text{MLP}(\hat Z))+\hat Z$}
           \State{\emph{shifted windows partition}}
            \State{$\hat{Y}=\text{LN}(\text{SW-MSA}^{cos}(Z))+Z$}
            \State{$Y=\text{LN}(\text{MLP}(\hat Y))+\hat{Y} $}
        \end{algorithmic}
    \caption{Two successive SwinV2 Transformer Blocks at a fixed stage. $n_W$ denotes the number of non--overlapping windows of the regular partition, $n_{SW}$ that of the shifted partition. W-MSA$^{cos}$ is analogous to Algorithm \ref{alg:W-MSA}, but with row 6 replaced by (\ref{eq:W-MHA-cos}). SW-MSA$^{cos}$ is analogous to  W-MSA$^{cos}$, but on the shifted windows partition.}
    \label{alg:swinv2}
    \end{algorithm}

\end{description}

\null\textbf{Unsupervised learning architectures.}
\begin{description}
\item [Generative Adversarial Networks (GANs)] \cite{Goodfellow-2014-GAN}  consist of two networks, a generator and a discriminator, that are trained together. Given a training set, the generator creates fake samples that mimic the distribution of the training set, starting from a latent space.
The discriminator tries to distinguish between real and generated data.
They compete with each other in the form of a zero--sum
game: one agent's gain corresponds to the other agent's loss. See Fig. \ref{fig:GAN}. 
GANs can generate new content that seems original (it preserves statistical similarities). 
Applications are image generation, style transfer, data augmentation, and image super--resolution. Examples are Vanilla GAN (2014), DCGAN (2015), CycleGAN (2017), and StyleGAN (2018).
\begin{figure}
    \centering
    \includegraphics[width=0.7\linewidth]{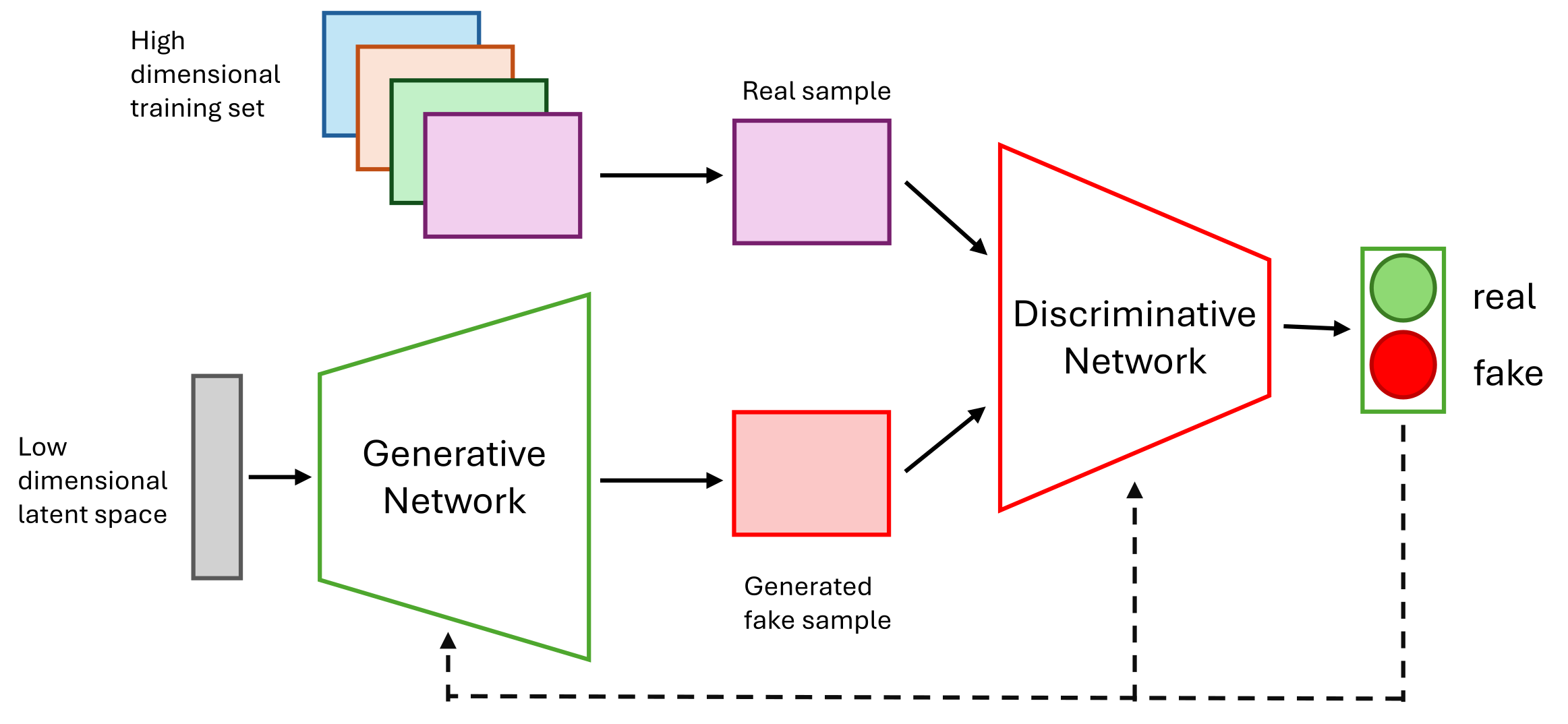}
    \caption{Generative Adversarial Network}
    \label{fig:GAN}
\end{figure}

\item [Autoencoders] are special instances of encoder--decoder architectures trained to attend to reconstruct their input after having reduced it to a lower dimension. They learn to compress/reduce the data (encoding) and then reconstruct them (decoding), see Fig. \ref{fig:encoder-decoder}. Applications are dimensionality reduction, anomaly detection, data compression, and generative modelling.
Examples are Vanilla Autoencoder (developed in the 80s), Sparse Autoencoder (2006), Denoising Autoencoder (2008), and Variational Autoencoder (VAE, 2013).
\end{description}


\null\textbf{Reinforcement learning architectures.}
They focus on learning policies for decision--making by interacting with an environment. Applications are game playing (e.g., AlphaGo (2015)), robotics, autonomous vehicles, and resource management. Examples of these architectures are Deep Q--Networks (DQN, 2015), Policy Gradient Methods, and Actor--Critic Methods.

\null\textbf{Hybrid architectures.} 
Finally, hybrid architectures (supervised and unsupervised) combine elements from different types of networks to leverage each one's strengths.
Applications are complex tasks requiring a mix of techniques, such as video analysis, multimodal learning, and more.
Examples are attention--based CNNs, which combine CNNs with attention mechanisms for better feature extraction, and Transformer--CNN Hybrids, which are used in vision tasks and combine local feature extraction with global context awareness.

Needless to say, good deep architectures should exhibit good trade--offs among multiple criteria\cite{Fleuret2024}:
\begin{enumerate}[noitemsep]
\item easiness of training,
\item accuracy of prediction,
\item memory footprint,
\item computational cost,
\item scalability.
\end{enumerate}

\begin{figure}
    \centering
    \includegraphics[width=\linewidth]{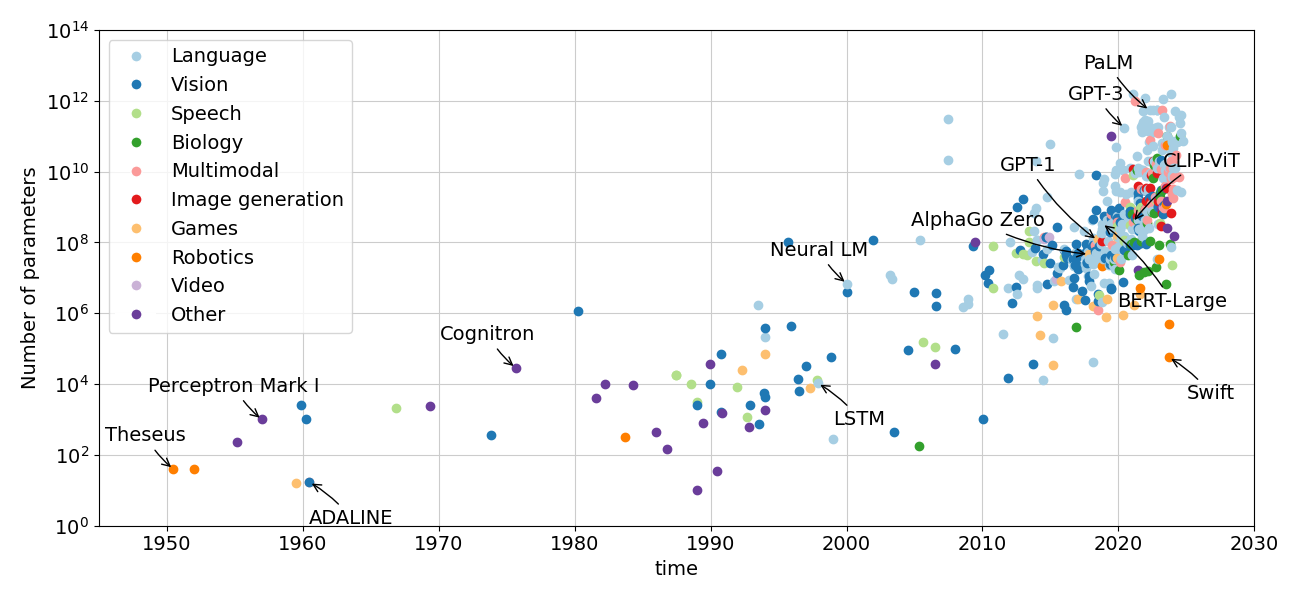}
    \caption{Number of trainable parameters for notable AI models. GPT-4 (2023) with $10^{20}$ trainable parameters is out of range.  Data from  \texttt{epochai.org} \cite{EpochNotableModels2024}}
    \label{fig:parameters_size}
\end{figure}

To give an idea of the complexity of some deep NNs, 
AlphaGo, developed by the Google DeepMind team in 2016 to play the board game Go, is composed by
2 NNs with 13 and 15 layers each, convolutional and fully connected, involving millions of parameters. 
The Large Language Model GPT--3 (produced by Microsoft's Turing Natural Language Generation team in 2020) is composed of 96 layers and has 175 billion parameters, while GPT-4 (the next generation of GPT--3, published in 2023) has 120 layers and 1500 billion parameters. In Fig. \ref{fig:parameters_size}, the number of unknown parameters for some of the most known NNs is shown.

\begin{figure}
    \centering
    \includegraphics[width=\linewidth]{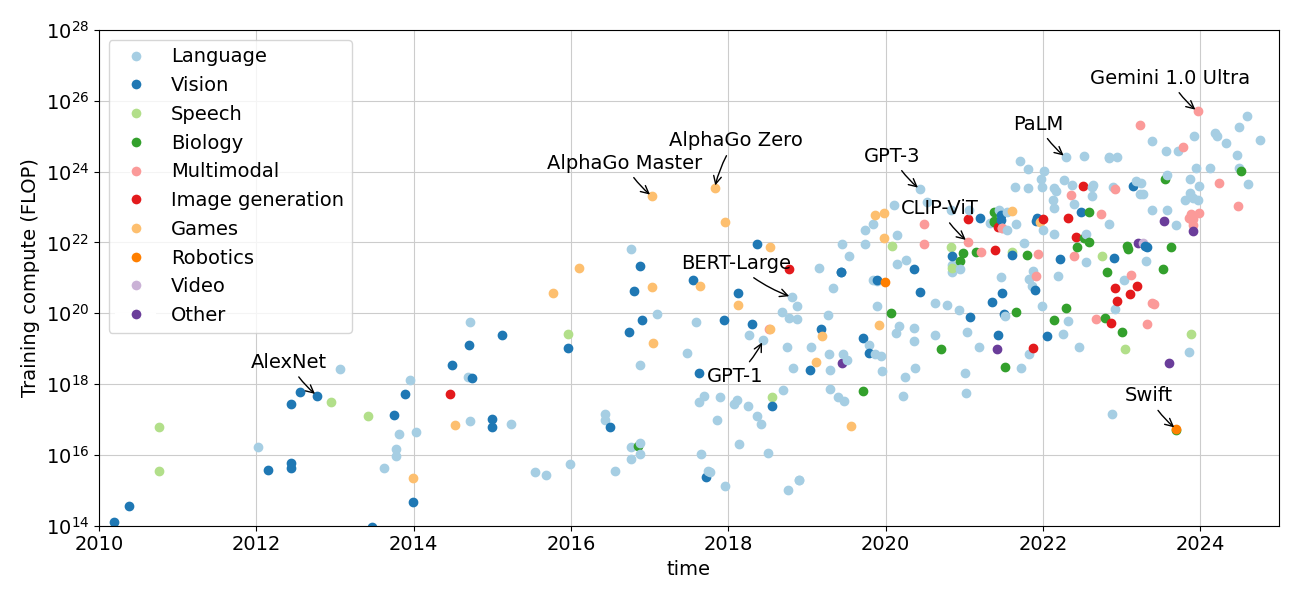}
    \caption{Training compute of notable models in the number of FLoating Point Operations (FLOP). Data from \texttt{epochai.org} \cite{EpochNotableModels2024}}
    \label{fig:flops_training}
\end{figure}

Clearly, the greater the complexity of the algorithm, the more powerful the hardware to train the NN. In 
Fig. \ref{fig:flops_training} we show the training costs in the number of FLoating Point Operations (FLOP) for notable deep learning models (data from \cite{EpochNotableModels2024}).

\subsubsection{Ultimate generation hardware: GPUs and TPUs}

The most suitable hardware to train neural networks are those equipped with GPUs (Graphics Processing Units) with tensor cores and TPUs (Tensor Processing Units). 

\emph{GPUs} are specialized and very efficient hardware processors originally designed to handle complex graphics, but then became crucial in artificial intelligence too, for training NNs. They are typically used when repeated matrix--vector operations are required, like, e.g., in the computation of weighted inputs (\ref{eq:weighted-input}), the realization of convolution layers, or the backpropagation phase and, more in general for any specific tasks requiring massive data processing. Multiple GPUs can be used to scale up processing power, allowing for the fastest computations.
The most powerful architectures, like e.g. those used to train GPT--4, feature several thousand parallel units and fast small local memories. Nevertheless, the huge number of GPUs is not sufficient to guarantee that the hardware is highly efficient. Indeed, it is essential that communication between GPU and CPU memories, often referred to as ``data transfer'', be as small as possible for several reasons, e.g., latency, bandwidth limitation, synchronization overhead, memory coherency issues, and, last but not least, power consumption. 

To move towards an energy--efficient computing (especially for training and inference of large AI models), recently NVIDIA GPUs have been equipped with dedicated tensor cores, while Google Cloud has designed custom AI accelerators (Tensor Processing Units -- TPUs), specifically designed to accelerate the computation of tensor operations\footnote{Tensors are used to represent signals, trainable parameters, and intermediate quantities.}.

\subsection{Topics related to Machine Learning not covered in this paper}
\label{sec:other_ML_topics}

For the sake of space, we had to select specific topics to discuss in this paper. However, we believe it is useful to bring to the reader's attention some other relevant topics in the context of ML. Interested readers are encouraged to explore these references for a deeper understanding.

Proper initialization of neural network parameters is crucial for effective training. Techniques such as Xavier initialization \cite{glorot2010understanding} and He initialization \cite{he2015delving} help mitigate {the shortcoming of} vanishing/exploding gradients, a common issue in training deep neural networks where gradients can become too small or too large, hindering effective learning \cite{bengio1994learning}. 

Transfer learning involves transferring knowledge from one model trained on a specific task to another model for a different but related task, significantly reducing the training time and improving performance \cite{pan2009survey}. Federated learning is a decentralized approach to ML where models are trained across multiple devices or servers holding local data samples, without exchanging them \cite{kairouz2021advances}.

Explainable AI encompasses methods and techniques developed to make the outputs of ML models understandable to humans, enhancing trust and transparency \cite{arrieta2020explainable,xu2019explainable, montavon2018methods}. 

Neural Tangent Kernel (NTK) theory provides a theoretical framework that helps in understanding the behaviour of infinitely wide neural networks during training \cite{jacot2018neural}. Overparameterization and double--descent describe the phenomenon where increasing the number of parameters in a model beyond a certain point can {improve} performance, contrary to traditional bias--variance trade--off expectations \cite{belkin2019reconciling}.


\section{Scientific Machine Learning}\label{sec:SML}


Many definitions of Scientific Machine Learning (SciML) can be found, more or less definite. We report some of them expressing different points of view. The one given in a report of 2021 by the US Department of Energy \cite{Doe2019} envisions scientific computing as fully supporting AI and states
\begin{quoteit}
``SciML is a core component of Artificial Intelligence (AI) and a computational technology that can be trained, with scientific data, to augment or automate human skills.''
\end{quoteit}
The next two definitions put scientific computing and AI on the same level and perceive advantages for both disciplines. The first one reads
\begin{quoteit}
``SciML is a merge of computational sciences and data--driven machine learning, implemented in software as a set
of abstractions to leverage existing domain knowledge and physics models within learning schemes and accelerated computing platforms.''
\end{quoteit}
This definition has been published in a 2024 workshop abstract \cite{icerm_brown}, while the following one is part of a 2022 booklet \cite{Schilders}:
\begin{quoteit}
``SciML is a combination of Scientific Computing and ML, combining mathematical models with data--based reasoning, presented as a unified set of abstractions and a high-performance implementation.
In this new area of research, many successes have already been found, with tools like physics--informed neural networks, universal differential equations, deep backward stochastic differential equation solvers for high dimensional partial differential equations, and neural surrogates, showcasing how deep learning can greatly improve scientific modelling practice.''
\end{quoteit}

The next one, reported on the web page of Oden Institute at Austin (Texas) \cite{Sciml-OdenInstitute}, highlights the range of applicability and the power of SciML to tackle the social challenges of the third millennium: 
\begin{quoteit}
Scientific Machine Learning brings together the complementary perspectives of computational science and computer science to craft a new generation of machine learning methods for complex applications across science and engineering. In these applications, dynamics are complex and multiscale, data are sparse and expensive to acquire, decisions have high consequences, and uncertainty quantification is essential. The greatest challenges facing society (clean energy, climate change, sustainable urban infrastructure, access to clean water, personalized medicine and more) by their very nature require predictions that go well beyond the available data. Scientific machine learning achieves this by incorporating the predictive power, interpretability and domain knowledge of physics--based models.
\end{quoteit}

We know that many problems can be tackled and solved by AI alone (image recognition, text generation, recommendation systems, ...) whereas many engineering problems can be effectively tackled by digital models alone (the computation of the response of an electric circuit, or that of the flow rate of a water network, problems from structural mechanics, fluid dynamics, just to mention a few). 

There are, however, very complex problems that until recently have been solved only thanks to digital models, which today can take advantage of the 
mutual interaction between Machine Learning (ML) algorithms and the knowledge of physical processes underpinning digital models to drive the learning process of ML algorithms. 

We believe the following may be a correct, concise, and effective definition:
\begin{quoteit}
``Scientific Machine Learning is an interdisciplinary field empowered by the
synergy of physics-based computational models with machine-learning algorithms
for scientific and engineering applications.''
\end{quoteit}

Generally speaking, SciML emerges when physics exhibits intricate and multiscale behaviours, datasets are sparse and costly to obtain, decisions carry significant consequences,
or we need to quantify uncertainty.

\begin{figure}
    \centering
    \includegraphics[width=0.9\linewidth]{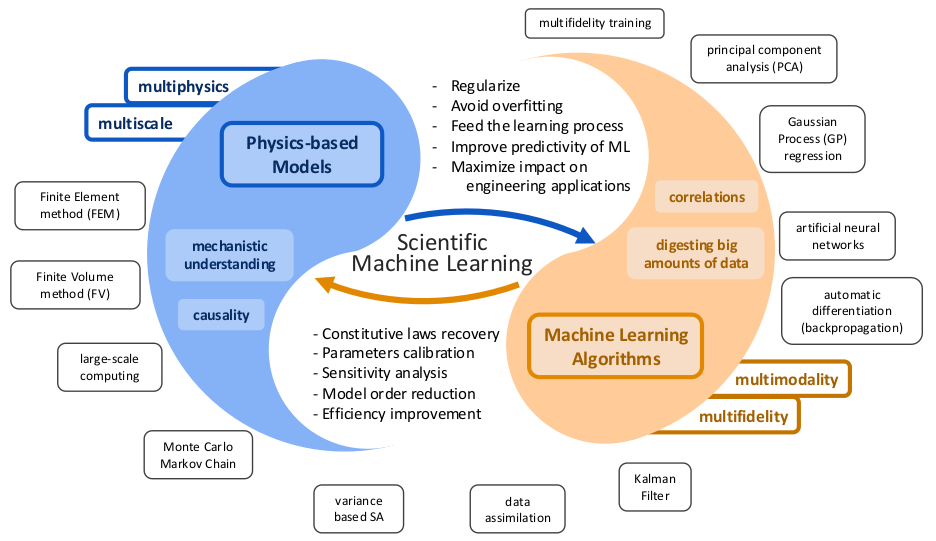}
    \caption{Cooperation between Digital Models and Machine Learning algorithms}
    \label{fig:yin-yang}
\end{figure}

%
%
There are in fact many reasons to regard SciML as a successful cooperative game between Digital Models (DM) and ML algorithms (see Fig. \ref{fig:yin-yang}).
For instance, we can use ML algorithms in favour of DMs:
\begin{itemize}[noitemsep]
    \item to recover (unknown) constitutive laws that are not known yet by using discrete input--output functions, and to enhance descriptive properties of DMs,
    \item to estimate and calibrate parameters by solving multi--query and inverse problems, 
    \item to understand which parameters are more/less sensitive,
    \item to achieve model reduction,
    \item to improve the efficiency of DMs by accelerating solver methods with cheaper (quicker) ML algorithms (curse of dimensionality).

\end{itemize}
On the other side, we can use DMs to empower ML applications with a deeper awareness of the problem structure, for instance, to accomplish the following tasks:
\begin{itemize}[noitemsep]
\item to regularize ML algorithms by penalizing the cost function and avoid overfitting, 
\item to augment data for the (otherwise scarce) training set, by adding solutions obtained from numerical algorithms,
\item to feed the operator learning process,
\item to improve the analytical and predictive power of ML and maximize its impact on science and engineering applications.
\end{itemize}

%
%
\null\textbf{Domain knowledge.} 
Scientific knowledge is essential to describe phenomena of a certain complexity and we believe that data--driven models alone cannot yet completely replace digital models in these cases. For instance, in the case of the digital heart modelling \cite{QDMV2019} that we will briefly review in Sec. \ref{sec:SML-IHM}, multi--physics simulations are needed to unravel the complex interactions among electrophysiology, heart mechanics and blood fluid dynamics and they are complex amalgams of biology, physical laws, mathematics, statistics, and computer science. Scientists (here mathematicians, engineers, and doctors) utilize their deep subject expertise, extensive experience, and finely--tuned intuition to formulate hypotheses and design analytical methodologies aimed at either confirming or disproving them. 
As noticed in Sect. \ref{sec:DDM}, data--driven models may fit training data very well, but at the same time, they might not be able to generalize (i.e., predict outcome values for previously unobserved inputs) with comparable accuracy. Rather, they might provide results on previously unseen data that would be physically unacceptable and meaningless in the absence of outside oversight by scientists or, more in general, physical principles. 

Awareness of \emph{domain knowledge} can enrich the quality of domain--agnostic data, reduce the amount of data needed while accelerating both training and predictive tasks, enhance the robustness and generalization properties in ML, and improve the accuracy, interpretability, and defensibility of SciML models.
By leveraging insights that transcend the limitation of available data, SciML can amalgamate the predictive capacity, interpretative abilities, and domain expertise intrinsic to physics--based models to make predictions that overstep the limitation of existing data.

Domain knowledge manifests in various forms, including physical principles (e.g., ab initio or first--principles physics), constraints (e.g., symmetries, invariances),
and structural configurations (e.g., discrete, graph--like, non-smooth data).

Various theoretical frameworks and computational tools, including solvers and simulations, are available to streamline the incorporation of domain knowledge, enhancing efforts in SciML.
Leveraging domain knowledge aids both supervised and unsupervised ML, as well as synthetic data generation (e.g., through constrained generative adversarial networks (GANs)) and reinforcement learning. Although scientific data may adhere, albeit imperfectly, to underlying physical laws, the direct application of domain knowledge empowers the learning process to address more intricate and computationally demanding phenomena using fewer labelled data points.

Many mathematical and, more in general, scientific areas are involved in SciML. Some, like functional analysis, geometry, and linear algebra, are more theoretical and are invoked to formulate and study the well--posedeness of the mathematical model and to accomplish parameter identification. Others are more applicative: with physics, we formulate constitutive laws, numerical analysis lets us discretize models, probability and statistics allow us to deal with data assimilation, control problems and optimization help us find optimal parameters providing the required output, and computer science is the cornerstone for implementing computer simulations and designing efficient neural networks and training algorithms that deliver the desired results.

%
%
We can identify some of the main \emph{challenges of SciML} as follows:
\emph{(i)} empower ML applications with a deeper understanding of the problem structure,
\emph{(ii)} exploit the analytical and predictive power of ML to enhance descriptive properties and efficiency of physics--based (digital) models and maximize its impact on science and engineering applications,
\emph{(iii)} operate in physical environments marked by large--scale, 3D, multi--modal data streams that are confounded with noise, sparsity, irregularities and other complexities that are common with machines and sensors interacting with the real world.

\null\textbf{Domain knowledge vs data exploitation.}
As remarked in \cite{KLKP}, knowledge of physics and availability of data play a complementary role and can guide our approach towards SciML. More precisely, we can distinguish three regimes: 
\begin{description}[noitemsep]
    \item[(i)] we know enough about the physics
and we need very few data to solve the problem. An example is given when we want to model the dynamics of water in a channel: the model is provided by the incompressible Navier--Stokes equations and we only require to know the geometry of the channel, a set of initial conditions, the boundary data, and the density and viscosity of the water. In this case, digital models are typically used without the need to exploit ML algorithms;
\item[(ii)] we have partial knowledge of physics and/or the corresponding data. An example is given when we want to describe the ionic dynamics of the cardiac cells in order to model cardiac electrophysiology and only noisy or deficient data are available to estimate the parameters involved in the equations. The physical model is given by a set of non--linear differential equations whose parameters are not known and need to be suitably calibrated on a patient--specific basis. At the same time in realistic scenarios the needed data are not easily measurable \cite{Regazzoni-lincei-2021} and the underlying process itself is not fully understood (only phenomenological models are available). This is the most frequent situation in practice, in this case, digital models are combined with ML algorithms to predict the results;
\item[(iii)] we know very little or nothing about physics but we dispose of a huge amount of data. In this case data--driven methods can help us understand the process, or more precisely, the input--output function that maps data into solution.
\end{description}

\begin{figure}
    \centering
    \includegraphics[width=0.9\linewidth]{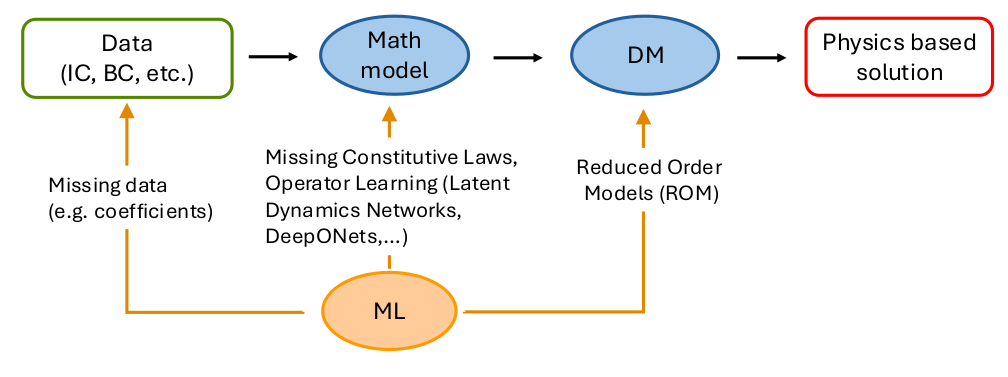}
    \caption{Physics--based approaches, how ML algorithms can improve digital models}
    \label{fig:ML2DM}
\end{figure}

\begin{figure}
    \centering
    \includegraphics[width=0.9\linewidth]{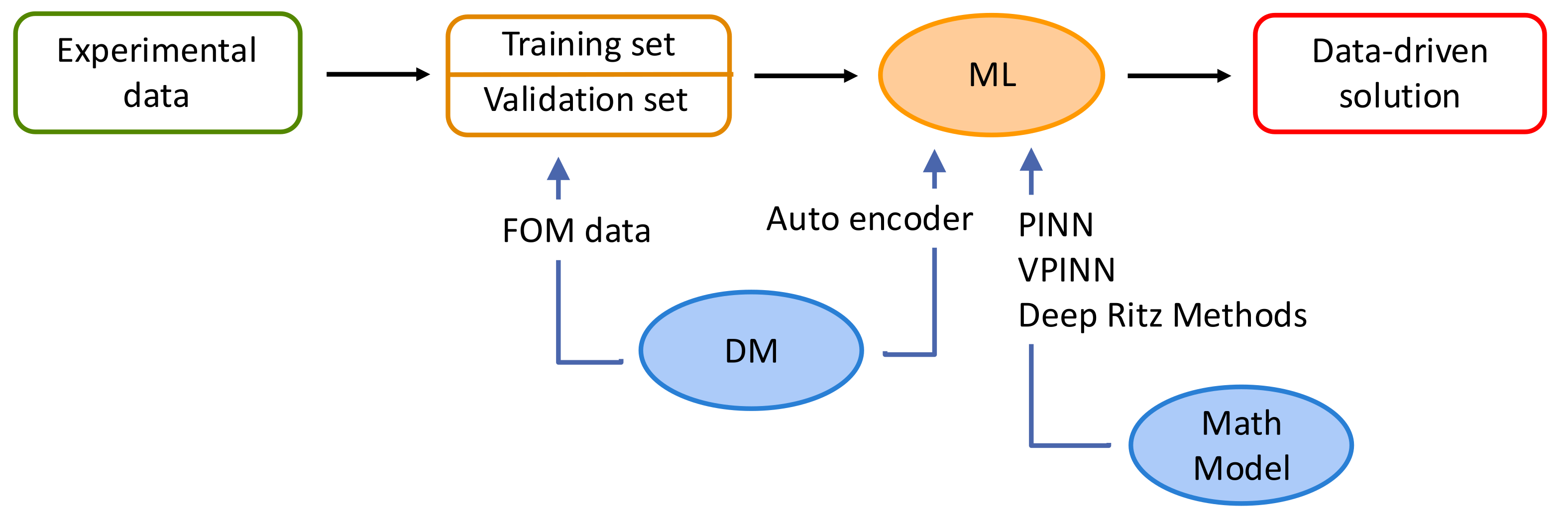}
    \caption{Data--driven based approaches, how digital models can improve ML algorithms}
    \label{fig:DM2ML}
\end{figure}

\null\textbf{How physics--driven digital models and data--driven models interact.}
The learning process of Neural Networks (NNs) used in SciML is typically supervised and the training is almost governed by physics and/or digital models (DM).
Although there are many possibilities for combining ML, physics, and digital models, for the sake of exposition, we identify below two possible paradigms summarized in Fig. \ref{fig:ML2DM} and \ref{fig:DM2ML}.

In \emph{physics--based approaches}, the main process is ruled by the mathematical model, fed by suitable data (like initial conditions (IC), boundary conditions (BC), forcing terms, etc.), which in turn feeds the digital model to obtain the physics--based solution. ML algorithms can provide missing data and constitutive laws, help to learn differential or integral operators, identify the latent dynamics of phenomena, and provide Reduced Order Models to improve the efficiency of digital models. We notice that the interplay between ML and DMs could involve either the whole simulation or only specific components of the process.

In \emph{data--driven approaches}, the main process is ruled by a ML algorithm whose learning process is fed by training and validation datasets. 
DMs can be used to enrich datasets for learning by generating pairs of high--fidelity input--output data, 
so that the NN embeds physical principles.
This approach has been considered since the 90s, see for example \cite{Hoole1995, Ratner1996, Rao1996}, and later adopted by several authors.

DMs can be also used to replace some components of a ML algorithm, once again incorporating physical principles. For example, inside autoencoders, a NN can be replaced by a reduced order model whose snapshots have been obtained by a high--fidelity digital model, see \cite{Tenderini2022}.

The loss function of the ML algorithm can incorporate information provided by the mathematical model of a physical problem in the form of equations residuals. Initially proposed in \cite{Lagaris1998,Lagaris2000}, this approach has been reconsidered and widely developed more recently under the name of Physics Informed Neural Network (PINN) \cite{Raissi2019, KLKP, Shin2024}, see Sec. \ref{sec:PINN}. Suitable modifications are Variational PINN, see Sec. \ref{sec:VPINN}, and Deep Ritz Methods, see Sect. \ref{sec:DeepRitz}.

We warn the reader that this classification is not rigid and should be taken with a grain of salt, the literature is filled with ad--hoc approaches where the interplay between DMs and ML is variably implemented.

\null\textbf{Biases.} 
The concept of bias was introduced by Mitchell in \cite{Mitchell-bias-1980} as follows: \emph{the term bias is used to refer to any basis for choosing one generalization over another, other than strict consistency with the observed training instances}. 
Following \cite{KLKP}, we can distinguish among observational, inductive and learning biases.
\begin{mydescription}
\item[Observational biases] are the input--output pairs provided by a high--fidelity DM and used to enrich the training set of a ML algorithm. They are weak mechanisms for embedding physical principles into ML models during training. The drawback is that a lot of data is typically needed.
\item[Inductive biases] drive a learning algorithm to prioritize one solution over another, independent of the observed data \cite{battaglia2018, Mitchell1997}. They can be introduced through the choice of the hypothesis space, or by designing specialized ML architectures to embed tailored building blocks. The most widespread examples are convolutional NNs \cite{Lecun2015}, Graph NNs \cite{Gori2005, Scarselli2009} and Neural ODEs \cite{chen2018neural}. We find earlier architectures that mimic particular discretization approaches like Finite Elements for PDEs or Finite Differences for ODEs in \cite{Takeuchi1994, Ramuhalli2005}.  
\item[Learning biases] are introduced into the loss function and soft penalty constraints and steer the learning process towards predictions that weakly (approximately) satisfy the physical principles. We find them in the ``physics--informed learning'' paradigm presented in the previous paragraph. 
\end{mydescription}


When we solve Partial Differential Equations (PDE) with numerical methods, we typically ask how reliable an algorithm is for a given type and quantity of data, how robust a solution is to slight variations in data or the addition of noise, and how rigorously the assumptions and underlying theories have been defined and validated. For classical techniques, these questions lead to familiar concepts, including well--posedness, stability, numerical approximation, and uncertainty quantification (see Sect. \ref{sec:DM}). 
However, this is not the case when ML algorithms are involved.

Although SciML interests more and more researchers and the number of published papers on the subject grows exponentially \cite{Herrmann2024}, a clear and representative description of the state of the art is almost impossible.
In the next subsections, we present some of the most widespread paradigms in the field of SciML for PDEs, without pretending to be comprehensive of all the literature devoted to the topic.

\subsection{Surrogate modelling of high--fidelity DM}\label{sec:DM2ML} 

Let us consider the general formulation of a boundary--value problem 
\begin{eqnarray*}
    \left\{\begin{array}{ll}
    \mathcal{P}(z,u)=0 & \mbox{\emph{(physics--based model)}}\\
    y=\mathcal{O}(z) & \mbox{\emph{(observation)}},
    \end{array}\right.
\end{eqnarray*}
where 
$u\in U$ is the \emph{input} which can include constant parameters and/or space dependent functions representing boundary conditions, forcing terms, material coefficients, viscosities, conductivity, and so on, 
$z\in Z$ is the \emph{state}, i.e., the solution, and $y\in Y$ is the output, which depends on the state. For instance, when we solve Navier--Stokes equations around the wing of an aeroplane, we typically compute drag and lift coefficients (output) which depend on the velocity and pressure (state). 

We can represent the \emph{data--to--solution map} as a function $\mathcal{M}: U\to Y$ which maps the input space into the output one such that $\mathcal M:u\mapsto y$. Then, we discretize it by a numerical model 
$\mathcal{M}_h: U\to Y_h$ named \emph{full--order model (FOM)}, such that
$\mathcal M_h: u\mapsto y_h \simeq y$. The approximate solution $y_h$ is typically named \emph{high--fidelity solution} because, potentially, it can be as much accurate as needed.
However, FOM are often associated with large computational costs, which often represent serious obstacles in practical applications. This motivates the construction of fast (albeit less accurate) models, termed surrogate models (SM) or emulators. These surrogate models are learned based on a set of precomputed simulations, obtained thanks to the FOM itself.

Specifically, given $N$ instances of the input $\{u_i\}_{i=1}^N \subset U$ that are chosen by sampling $U$ with a goal--driven strategy, we compute the corresponding solutions
$(y_h)_i =\mathcal{M}_h(u_i)$ by the numerical FOM to form the training set $S=\{(u_i,(y_h)_i),\ i=1,\ldots, N\}$.
Then, we choose the hypothesis space
$\mathcal{H}\subset\{f:U\to Y\}$ of the surrogate model $\mathcal M_{SM}$
and we look for the optimal surrogate model
\begin{equation}\label{eq:surrogate_model}
\mathcal{M}_{SM}^* =
\argmin{\mathcal{M}_{SM}\in \mathcal{H}}\sum_{i=1}^N\|(y_h)_i-\mathcal{M}_{SM}(u_i)\|_Y^2
\end{equation}
that minimizes the discrepancy between the FOM targets $(y_h)_i$
and the outputs $\mathcal{M}_{SM}(u_i)$ of the surrogate model. 

\null\textbf{Computational costs.}
Once the surrogate model has been trained, typically, the time $T_{SM}$ required to compute the SM solution is much smaller than the time $T_{FOM}$ required to compute the FOM solution. However, we must keep in mind that the training phase is usually very expensive and exceeds the cost of a simple FOM resolution.
Thus, should we have to solve a single direct PDE, at the moment, high--fidelity DMs have no competitors in the ML world. If, instead, our interest is in solving inverse problems, having real--time solutions, or estimating parameters, the use of SM can be advantageous. Indeed, in this case, the time spent to solve many FOM problems can be balanced by that required to train the NN. We could find ourselves in the following situations:
\begin{description}
    \item[Many--query scenarios.] If we need to solve $n_{query}\gg 1$ problems that differ in the choice of parameters and apply FOM to compute their solutions, the total time required would be $n_{query}\cdot T_{FOM}$. If instead we train the surrogate model with $N$ full order model solutions and then solve the surrogate model $n_{query}$ times, denoting by $T_{training}$ the time required to train the surrogate model, the SM solution is convenient whenever
    \begin{equation} \label{eqn:many-query}
        N\cdot T_{FOM} + n_{query}\cdot T_{SM} + T_{training} <
        n_{query}\cdot T_{FOM}.
    \end{equation}
    As the leading term of the left-hand side is typically $N\cdot T_{FOM}$, invoking a SM is advantageous whenever $n_{query}>N$.
    
Many--query scenarios can occur in sensitivity analysis, scenario analysis (what happens if we change coefficients or boundary conditions, for instance), parameter estimation (in particular Bayesian parameter estimation), uncertainty quantification (forward and/or backward), multiscale problems, and more. Sometimes, the SM can be used together with the FOM, in a \emph{multifidelity} framework. 
\item[Real--time solutions.] Many real-world applications require real-time solutions. Notable examples are robotic-assisted surgery and continuous patient monitoring in intensive care units, where system responsiveness is crucial for timely interventions. In this case, the expensive training can be performed \emph{offline}, while the cheaper SM solution is achieved \emph{online}, to respond rapidly to real-time demands.
\end{description}

In addition to real--time solutions and many--query scenarios, SM can also be useful in other contexts.

One important use of SM is in creating differentiable or interpretable approximations of black--box solvers. For example, when working with an ``oracle'' or proprietary software that provides solutions without exposing internal calculations or logic, surrogate models that approximate the input--output behaviour of the black--box solver can in fact deliver a differentiable model, which is essential for optimization tasks and sensitivity analysis. This also enables further insights into the solution dynamics that the original solver does not reveal.

Surrogate models are generated similarly to data--driven models (see Sect. \ref{sec:DDM}), they both are fed by data and trained. However, to train classic data--driven models only pure experimental data are considered. Moreover, the physics behind the input--output map is unknown. Instead, surrogate models are trained using input--output pairs provided by FOM which express an a--priori known physics.

\begin{figure}
    \centering
    \includegraphics[width=0.5\linewidth]{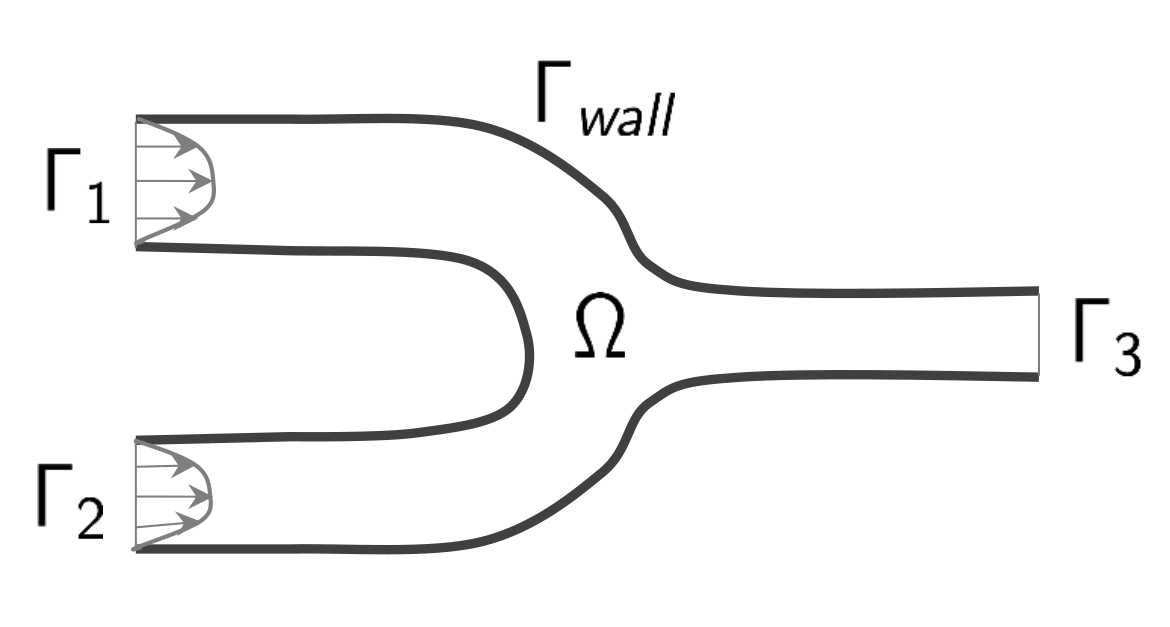}
    \caption{The computational domain and setting for problem \ref{eq:NS-parametricPDE}}
    \label{fig:NS-parametric}
\end{figure}

\null\textbf{Example: a parametric PDE.}
Let us consider the steady Navier--Stokes equations in a converging channel $\Omega\subset \mathbb{R}^{3}$ to predict the pressure $p$ at any point $\mathbf x\in \overline\Omega$ (the output), given the parameters $\mathbf u=(\mu_1,\mu_2)\in\mathbb{R}^2$  (the input)
that modulate the inflow velocities from two inlets, see Fig. \ref{fig:NS-parametric}.  The 
physical problem reads
\begin{eqnarray}\label{eq:NS-parametricPDE}
\left\{\begin{array}{ll}
-\nu\Delta \mathbf{v} +(\mathbf{v}\cdot\nabla)\mathbf{v}+\frac{1}{\rho}\nabla p =\mathbf{0} & \mbox{ in }\Omega\\[1mm]
\nabla\cdot\mathbf{v}=0 & \mbox{ in }\Omega\\[1mm]
\mathbf{v}=\mu_1\mathbf{v}^1_{in} & \mbox{ on }\Gamma_1\\[1mm]
\mathbf{v}=\mu_2\mathbf{v}^2_{in} & \mbox{ on }\Gamma_2\\[1mm]
-\nu\nabla\mathbf{v}\cdot\mathbf{n}+\frac{1}{\rho}p\mathbf{n}=\mathbf{0} & \mbox{ on }\Gamma_3\\[1mm]
\mathbf{v}=\mathbf{0} & \mbox{ on }\Gamma_{wall}.
\end{array}\right.
\end{eqnarray}

The vector $z=(\mathbf v, p)$ represents the state. We choose a Finite Element Method (FEM) solver as FOM, so that we can expand the discretized pressure 
\begin{equation*}    p_h(\mathbf{x})=\sum_{j=1}^{N_h}p_j\varphi_j(\mathbf{x})
\end{equation*}
with respect to the finite--dimensional Lagrange basis function and we denote by
$\mathbf{p}=[p_j]\in\mathbb{R}^{N_h}$ the vector of the nodal values of $p$, i.e., $p_j=p_h(\mathbf{x}_j)$.

Chosen $N$ training input samples $\{\mathbf u_i=(\mu_{1,i},\mu_{2,i})\}_{i=1}^N$, for every $i$ we compute the corresponding FOM solution and we define the discrete input--output map $\mathcal M_h:\mathbf u_i \mapsto \mathbf{p}_i$. In this case, the training set is defined by $S=\{(\mathbf u_i ,\mathbf{p}_i), \ i=1,\ldots, N\}$.

\begin{figure}
    \centering
    \includegraphics[width=0.8\linewidth]{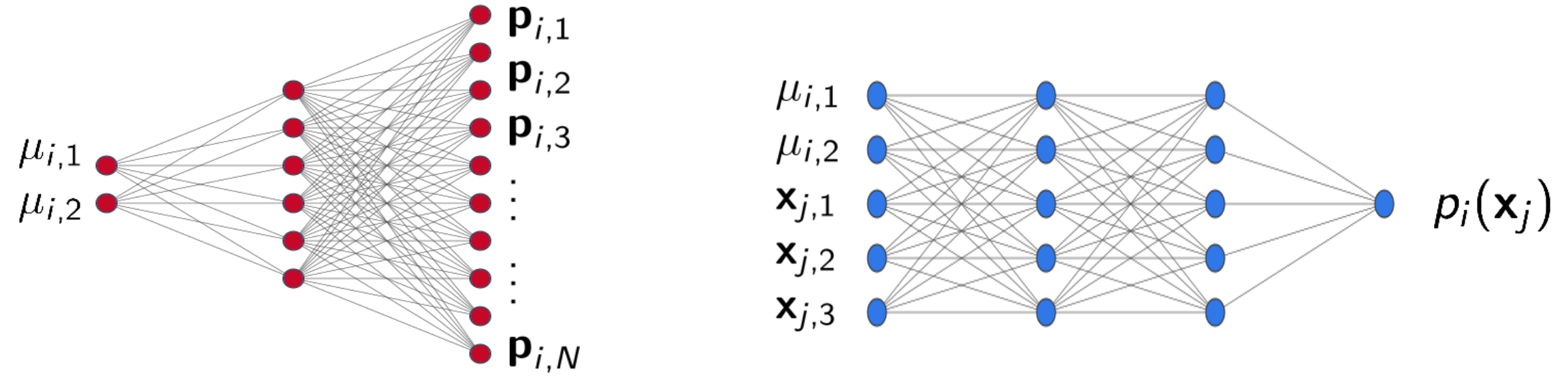}
    \caption{Two possible FFNNs for the parametric Navier--Stokes equation (\ref{eq:NS-parametricPDE})}
    \label{fig:NNs-parametric-NS}
\end{figure}

This is an example of supervised learning, we are using a DM to feed a NN and we are introducing observational biases. We can follow two different strategies to design a NN to compute the SM solution, as drawn in Fig. \ref{fig:NNs-parametric-NS}.

 The first approach consists of taking inputs of $n=2$ values (the physical parameters), $m=N_h$ outputs that are the values of the pressure at all FEM mesh nodes, and $M$ trainable parameters, then we design the neural network
$\NNgeneric_1: \mathbb{R}^{2+ M}\to \mathbb{R}^{N_h}$ such that
\begin{equation}\label{eq:NN1-parametricPDE}
\mathbf{w}^*=\argmin{\mathbf{w}\in\mathbb{R}^M}\sum_{i=1}^N\|\NNgeneric_1(\mathbf u_i;\mathbf{w})-\mathbf{p}_i\|^2_2,   
\end{equation} 
see Fig. \ref{fig:NNs-parametric-NS}, left.

The second approach consists of choosing $N_{po}$ points  $\mathbf{x}_j \in \Omega\subset\mathbb{R}^3$ and evaluating the pressure at those points. More precisely, for any sample $i=1,\ldots,N$  and any point $\mathbf x_j$ with $j=1,\ldots,N_{po}$, we take inputs of $n=5$ values given by the two physical parameters and the $3$ coordinates of the point $\mathbf x_j$. The output of size $m=1$ is the value of the pressure at the chosen node $\mathbf{x}_j$ while the trainable parameters are $M$. Then we design the FFNN
$\NNgeneric_2: \mathbb{R}^{5+ M}\to \mathbb{R}$ such that (see Fig. \ref{fig:NNs-parametric-NS}, right):
\begin{equation}\label{eq:NN2-parametricPDE}
\mathbf{w}^*=\argmin{\mathbf{w}\in\mathbb{R}^M}\sum_{i=1}^N\sum_{j=1}^{N_{po}}|\NNgeneric_2(\mathbf u_i,\mathbf x_j;\mathbf{w})-p_h(\mathbf x_j)|^2.
\end{equation} 
In comparison with \eqref{eq:NN1-parametricPDE}, the approach \eqref{eq:NN2-parametricPDE} has the advantage of using a resolution--invariant and mesh--less representation, and thus is particularly flexible to handle cases where training data are available at different resolutions. In addition, the NN obtained in \eqref{eq:NN2-parametricPDE} is typically very lightweight since input and output are of much lower dimensionality than \eqref{eq:NN1-parametricPDE}. This results in a much smaller number of trainable parameters, which often leads to better generalization for the same amount of data available for training \cite{regazzoni2024ldnets}.

In Sec.~\ref{sec:operator-learning} a more complex scenario will be addressed in which the solution of the differential problem does not depend simply on scalar parameters but also on functions. However, before we will discuss how to embed the knowledge of the physics governing the problem at hand into the learning process.

%
%
\subsection{Physics--Informed learning}\label{sec:PI-ML}

The approaches presented in the previous section leverage DMs with the sole purpose of generating training data for ML algorithms. The physics is incorporated as observational biases. Note that using observational biases does not require knowing the underlying physical problem. Indeed, data could be generated by a black--box solver.

In contrast, Physics--Informed Machine Learning methods aim to incorporate physical principles into the learning process. This is achieved by embedding the governing equations of the physical model into the loss function of the ML algorithm, i.e., we are adding learning biases into the system.

Physics--informed learning can serve a variety of purposes, with the flexibility to adapt the same method to different use cases by making minor adjustments to the list of trainable variables and the terms comprising the loss function. For instance, these methods can be applied to solve PDEs (forward solutions), estimate parameters, address inverse problems, fit experimental measurements, or enhance the resolution in case of noisy data. In this discussion, we will begin by examining the simplest case (solving PDEs) and then extend the approach to other applications. To illustrate this general framework, we consider a simple PDE as an example.

\null\textbf{FEM approximation of a PDE.} 
Let $\Omega\subset{\mathbb R}^d$, with $d=2,3$, be an open bounded domain, with Lipschitz boundary.
Given $f:\Omega\to\mathbb{R}$ and $g:\partial\Omega\to \mathbb{R}$ regular enough and a real coefficient $\gamma>0$, we look for
the solution $u:\Omega\to {\mathbb R}$ of the second--order elliptic equation with Neumann boundary conditions
\begin{eqnarray}\label{eq:strong-elliptic}
\left\{\begin{array}{ll}
-\Delta u+\gamma u=f & \mbox{ in }\Omega\\[2mm]
\nabla u\cdot \mathbf{n}=g &\mbox{ on }\partial\Omega.
\end{array}\right.
\end{eqnarray}
After introducing the Sobolev space $H^1(\Omega)=\{v:\Omega\to{\mathbb R}\ s.t.\ 
\int_\Omega v^2+|\nabla v|^2 < +\infty\}$, its weak form reads
\begin{equation}\label{eq:weak_elliptic0}
\mbox{find}\  u\in V:\ \int_{\Omega}\left(\nabla u\cdot \nabla v + \gamma u v\right)d\Omega=
\int_\Omega fv d\Omega+\int_{\partial\Omega}gv\,d\partial\Omega\qquad \forall v\in V
\end{equation}
or, equivalently,
\begin{equation}\label{eq:weak_elliptic}
\mbox{find}\  u\in V:\ a(u,v)=F(v)\qquad \forall v\in V,
\end{equation}
where $a:V\times V\to \mathbb R$ defined by $a(u,v)=\int_{\Omega}\left(\nabla u\cdot \nabla v + \gamma u v\right)d\Omega$ is a coercive, symmetric, and continuous bilinear form, and $F:V\to \mathbb R$ defined by $F(v)=\int_\Omega fv d\Omega+\int_{\partial\Omega}gv\,d\partial\Omega$ is a linear and continuous functional. If $f\in L^2(\Omega)$and $g\in L^2(\partial\Omega)$, this problem admits a unique solution \cite{Quarteroni-Valli-1994}.

Classical methods to discretize problem (\ref{eq:weak_elliptic}) are
the Galerkin methods \cite{Quarteroni-Valli-1994}, for which, after introducing a discretization parameter $h>0$ and the finite--dimensional space $V_h\subset V$, with $N_h=\dim(V_h)$, 
we look for $u_h\in V_h$ such that
\begin{equation}\label{eq:galerkin}
a(u_h,v_h)=F(v_h)\qquad \forall v_h\in V_h.
\end{equation}
The Galerkin solution $u_h$ is an approximation of $u$. 

The space $V_h$ is the \emph{test space} (i.e. the space to which the test functions $v_h$ belong) and the \emph{trial space} (i.e. the space to which the solution $u_h$ belongs) at the same time.
After choosing a basis of functions 
$\{\varphi_j(\mathbf{x})\}_{j=1}^{N_h}$ in $V_h$, the approximate solution $u_h$ can be written as
$u_h(\mathbf{x})=\sum_{j=1}^{N_h}u_j\varphi_j(\mathbf{x})$ and (\ref{eq:galerkin}) becomes
\begin{equation}\label{eq:galerkin2}
\mbox{find } u_1,\ldots,u_{N_h}:\ 
\sum_{j=1}^{N_h}u_j a(\varphi_j,\varphi_i)= F(\varphi_i)\qquad \forall j=1,\ldots,N_h.
\end{equation}
Thus, by setting $A_{ij}=a(\varphi_j,\varphi_i)$, $b_i=
F(\varphi_i)$, we have $A=[a_{ij}]\in \mathbb{R}^{N_h\times N_h}$, $\mathbf{u}=[u_j]\in\mathbb{R}^{N_h\times 1}$, and
$\mathbf{b}=[b_i]\in\mathbb{R}^{N_h\times 1}$ and the algebraic formulation of (\ref{eq:galerkin2}) reads
$A\mathbf{u}=\mathbf{b}$. The unknowns $u_j$ are named \emph{degrees of freedom}, while the array $\mathbf b $ is the right--hand side.

\begin{figure}
    \centering
    \includegraphics[width=0.3\linewidth]{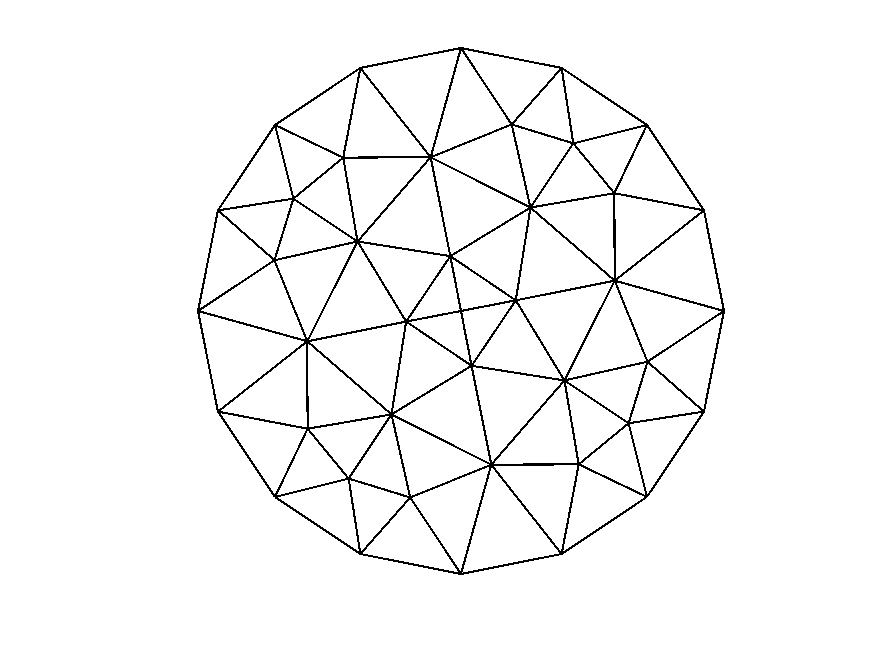}\quad
    \includegraphics[width=0.3\linewidth]{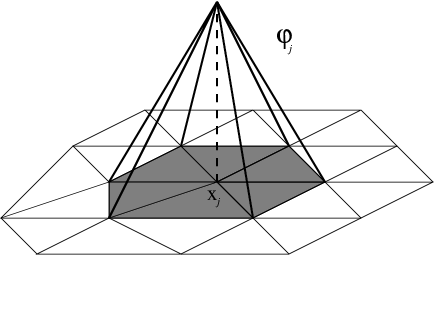}\quad
    \includegraphics[width=0.3\linewidth]{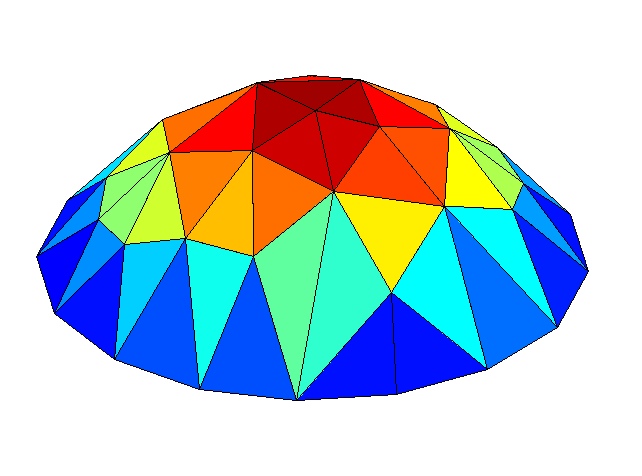}
    \caption{A computational mesh made of triangles (left), a $\mathbb{P}_1-$FEM hat function (centre), a function of $V_h$ (right)}
    \label{fig:FEM-mesh-basis}
\end{figure}

Finite Elements are a special case of Galerkin methods. Given the parameter $h>0$, a mesh\footnote{A mesh is a conforming partition of $\Omega$ in triangles/quadrilaterals when $d=2$, or in tetrahedra/hexahedra when $d=3$, all with diameter less than $h$ (the \emph{mesh--size}), see Fig. \ref{fig:FEM-mesh-basis}, left.} $\mathcal{T}_h$ on $\Omega$ is built. 
When simplexes are chosen and $V_h$ is defined by
\begin{equation}\label{eq:Vh}
    V_h=\left\{v_h\in\mathcal{C}^0(\overline\Omega)\  s.t.\ 
v|_{T_k}\in\mathbb{P}_1\ \forall T_k\in\mathcal{T}_h\right\},
\end{equation}
we are implementing linear $\mathbb{P}_1-$FEM, and the basis functions $\varphi_j$ are hat functions like that drawn in the central picture of Fig. \ref{fig:FEM-mesh-basis}. The functions of $V_h$ are linear combinations of hat functions, the one in the right picture of Fig. \ref{fig:FEM-mesh-basis} is an example.

\null

\null\textbf{Towards NNs.} 
We notice that the trial space $V_h$ of Galerkin methods and the space $V_{NN}$ of feed--forward neural networks (FFNN) with a fixed architecture play a similar role: $V_h$ is the space in which we look for the FEM solution, $V_{NN}$ is the space in which we look for the solution provided by FFNNs. Moreover, we can think that the trainable parameters $\mathbf w$ which characterize a FFNN are analogous to the unknown degrees of freedom $\mathbf u$ that identify the FEM solution (see Tab. \ref{tab:FEM-NN}). 

Thus, we could think of using the space $V_{NN}$ as the trial solution space for PDEs. This could be explained as follows. First of all, FFNNs are very efficient in approximating functions, we recall the universal approximation theorem \ref{thm:Cybenko}. Second, 
one of the most expensive parts of a Finite Elements code is the construction of the mesh, while FFNNs are mesh--free. Nevertheless, although the PDE is linear, FFNNs are not linear models and this brings complications: $V_{NN}$ is not a linear space and it is not a good candidate to be the test space for a Galerkin formulation.

\begin{table}
\begin{center}
\begin{tabular}{cccc}
\toprule
space & solution & input & unknown parameters\\
\midrule
$V_h$ & $\displaystyle u(\mathbf{x})=\sum_{j=1}^{N_h}u_j\,\varphi_j(\mathbf{x})$ & $\mathbf x$ & $\mathbf u=[u_1,\ldots, u_{N_h}]^t\in\mathbb R^{N_h}$\\
$V_{NN}$ & $ u_{NN}(\mathbf x;\mathbf w)$ & $\mathbf x$ & $\mathbf w\in\mathbb R^N$\\
\bottomrule
\end{tabular}
\end{center}
    \caption{Analogy between FEM and FFNNs}
    \label{tab:FEM-NN}
\end{table}

Recently, many approaches have been proposed that follow alternative paradigms to the Galerkin one, we cite the three most widespread ones.
The first one consists of working on the strong form of the PDE and leads to Physics--Informed Neural Networks (PINNs) \cite{Raissi2019,KLKP, Cuomo-rozza-2022,Shin2024}. However, when the coefficients of the PDE are discontinuous and/or pointwise loads are given, solutions in a strong sense of PDE do not always exist.  A possible alternative consists of considering the weak form of the PDE and applying a Petrov--Galerkin approach (with different trial and test spaces) leading to Variational PINNs (VPINNs) \cite{kharazmi2019-vpinn,Kharazmi2021}. A  third alternative makes use of the energy form of the PDE and provides the Deep Ritz method \cite{E2018}.
We refer to \cite{Shin2024} for an overview of further extensions of PINNs.

%
%
\subsubsection{Physics--Informed Neural Networks (PINNs)}
\label{sec:PINN}

PINNs allow us to reconstruct physical fields modelled by PDEs, by integrating the information from the PDE (learning biases) into the loss function.
In principle, many types of PDEs can be approximated by PINNs: integer--order PDEs, integro--differential equations, fractional PDEs, and stochastic PDEs.
PINNs can be used as direct solvers of PDEs, especially to address the case when we miss data (e.g., boundary or initial conditions, forcing terms, operator coefficients) which would be essential to approximate the problem in a classical sense. In addition, they can be even used to solve inverse problems.

\begin{figure}
    \centering
    \includegraphics[width=0.8\linewidth]{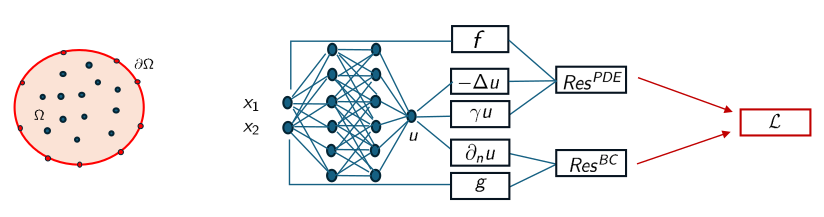}
    \caption{The collocation nodes (left) and the PINN (right). We follow the branch $Res^{PDE}$ if the input $\mathbf x_i$ belongs to the interior of the domain and the branch $Res^{BC}$ if $\mathbf x_i$ is on the boundary}
    \label{fig:pinn1}
\end{figure}

%
%
\null\textbf{PINNs as PDE solvers.} 
We start by examining the case where PINNs are used to solve a PDE.
Let us consider, as an example, the PDE (\ref{eq:strong-elliptic}) in strong form, PINNs are implemented by following these steps \cite{Shin2024}:
\begin{enumerate}
\item define the residuals of the differential equation (\ref{eq:strong-elliptic})$_1$ and of the boundary condition (\ref{eq:strong-elliptic})$_2$:
\begin{eqnarray}\label{eq:pinn-residuals}
\begin{array}{ll}
Res^{PDE}[u(\mathbf{x})]=-\Delta u(\mathbf{x}) +\gamma u(\mathbf{x})-f(\mathbf{x}) & \mathbf x\in\Omega\\[1mm]
Res^{BC}[u(\mathbf{x})]=\nabla u(\mathbf{x})\cdot \mathbf n-g(\mathbf{x}) & \mathbf x\in\partial\Omega,
\end{array}
\end{eqnarray}

\item choose a set of collocation points
$\mathbf{x}_i^{PDE}\in\Omega$ (the blue ones in the left picture of Fig. \ref{fig:pinn1}) and another set of points
$\mathbf{x}_i^{BC}\in\partial\Omega$ (the red ones), set
$\{\mathbf{x}_i\}=\{\mathbf{x}_i^{PDE}\}\cup
\{\mathbf{x}_i^{BC}\}$, and provide the value $u_i$  of the solution at the collocation node $\mathbf x_i$ for any $i=1,\ldots,N$, computed, for instance, by FEM.

\item design the FFNN like, e.g., in the right picture of Fig. \ref{fig:pinn1},

\item train the NN: given the pairs $(\mathbf{x}_i, u_i)$ for $i=1,\ldots, N$,
look for the minimizer
\begin{equation}\label{eq:pinn-minimzer}
\mathbf{w}^*=\argmin{\mathbf{w}\in\mathbb{R}^M}
\mathcal{L}(\mathbf{w}),
\end{equation}
where the loss function is 
\begin{equation*}
\mathcal{L}(\mathbf{w})
=
\mathcal{L}_{PDE}(\mathbf{w})+
\alpha_{BC}\mathcal{L}_{BC}(\mathbf{w})
\end{equation*}
with
\begin{eqnarray}
\label{eq:loss-pinn-respde}
\mathcal{L}_{PDE}(\mathbf{w})&=&\frac{1}{2N_{PDE}}\sum_{i=1}^{N_{PDE}}\left(Res^{PDE}[u_{NN}(\mathbf{x}_i^{PDE};\mathbf{w})]\right)^2\\
\label{eq:loss-pinn-resbc}
\mathcal{L}_{BC}(\mathbf{w})&=&\frac{1}{2N_{BC}}\sum_{i=1}^{N_{BC}}\left(Res^{BC}[u_{NN}(\mathbf{x}_i^{BC};\mathbf{w})]\right)^2,
\end{eqnarray}
and $\alpha_{BC}\in\mathbb R$ is a suitable hyper--parameter.
\item take $u_{NN}(\mathbf x;\mathbf w^*)$ as the PINN solution.
\end{enumerate}

\begin{figure}
    \centering
    \includegraphics[width=0.8\linewidth]{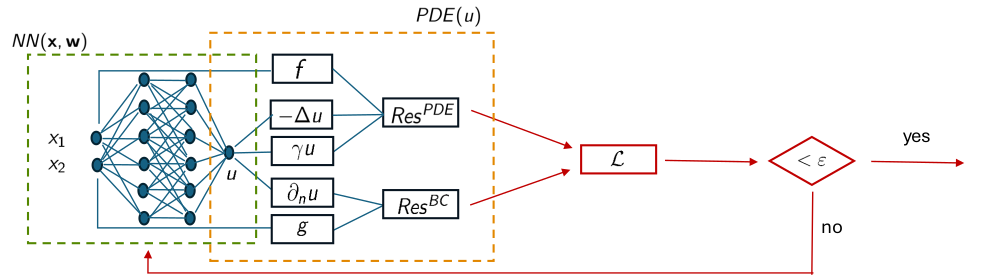}
    \caption{The minimization process in PINNs}
    \label{fig:pinn2}
\end{figure}

In Fig. \ref{fig:pinn2}, the iterative process to minimize the loss function is drawn.

Provided that the activation function is regular (notice that, if $\sigma\in {\mathcal C}^\infty(\mathbb R)$, then $ u_{NN}\in\mathcal{C}^\infty$ on the input space $\mathcal X\subset \mathbb R^n$) and noticing that $u_{NN}$ can be written as a composition of functions like in (\ref{eq:f_NN}), partial derivatives in the residuals are evaluated through automatic differentiation similarly to what done in Sec. \ref{sec:backpropagation}. Then, the loss function can be minimized via gradient--based methods as seen in Sec. \ref{sec:optimization}.

The summations in (\ref{eq:loss-pinn-respde}) and (\ref{eq:loss-pinn-resbc}) can be interpreted as approximations of integrals, i.e.
\begin{equation*}
\underbrace{\frac{1}{N}
\sum_{i=1}^{N}\left(Res[u_{NN}(\mathbf{x}_i;\mathbf{w})]\right)^2}_{\displaystyle I_N} \simeq 
\underbrace{\int_\Omega \left(Res[u_{NN}(\mathbf{x};\mathbf w)]\right)^2 d\Omega}_{\displaystyle I_{ex}}.
\end{equation*}
If the nodes are randomly chosen, this corresponds to using the Monte Carlo
quadrature rule, which provides an error $|I_{ex}-I_N|\approx  N^{-1/2}$.
Alternatively, if we choose the nodes more carefully and replace the arithmetic average with a weighted sum
$I_{ex}\simeq I_{N,\omega}=\sum_{i=1}^N \omega_i \left(Res[u_{NN}(\mathbf{x}_i;\mathbf w)]\right)^2$ \cite[Ch. 4]{Quarteroni-CS-2014}, we might obtain higher accuracy
$|I_{ex}-I_{N,\omega}| \leq C N^{-s}$
with $s\geq 1$ depending on the regularity of both data and solution and possibly on the space dimension.

%
%

One of the main strengths of PINNs is their mesh--free nature. This makes them particularly appealing because they bypass the often lengthy and tedious process of mesh generation. Indeed, PINNs can be effectively applied to irregular and complex domains \cite{Lagaris2000,KLKP}. Additionally, the availability of many open--source machine learning platforms, which are highly optimized for modern hardware, has made the implementation of PINNs very straightforward. Full utilization of the computational power offered by GPUs allows for realising very efficient PINN architectures.
Automatic differentiation makes PINNs easily applicable to non--linear PDEs. 
More precisely, while solving a non--linear PDE with FEM requires modifying heavily the code designed for linear problems (we have to compute the residual of the PDE and apply Newton--like methods), PINNs are straightforward to implement and only those lines of the code performing the calculation of the residual have to be changed.

On the other hand, computing the optimal parameters $\mathbf{w}^*$ involves solving an optimization problem, making the computational cost of PINNs typically higher than that of FEMs, even for linear PDEs. Thus,  traditional numerical methods usually outperform PINNs in solving well--posed forward problems. However, in cases involving PDEs in high--dimensional spaces, where mesh construction becomes challenging and is affected by the curse of dimensionality, PINNs offer a compelling alternative. This is due to the neural networks' ability to efficiently approximate functions in high--dimensional spaces.

A critical aspect of PINNs is that, when the solution of the PDE contains high--frequencies or multiscale features, PINNs using fully connected architectures struggle to converge during the training \cite{Wang2022}.
This is a consequence of the so--called \emph{spectral bias}, a phenomenon where neural networks tend to represent low--frequency components more effectively than high--frequency ones \cite{Wang2022}. This limitation poses challenges in accurately solving problems with solutions characterized by high--frequency components, such as wave--like behaviours or sharp discontinuities. To address this issue, several solutions have been proposed in the literature. These include using periodic activation functions, such as Fourier Features \cite{tancik2020fourier}, which enhance the network's ability to capture high frequencies, and adaptive or weighted loss functions \cite{Wang2021-SISC}, which balance the representation across different spectral scales. These approaches have shown promise in mitigating the effects of spectral bias, thereby improving the capability of PINNs to solve complex problems.

Since PINNs represent a method that has spread to the scientific community relatively recently, its theoretical understanding  is less mature than for other methods, although much progress has been made recently. Some aspects remain to be clarified, and in practical contexts very often empirical tests are needed to find the best configuration. In the following section, we discuss some of the theoretical results that pose a mathematically sound basis for using PINNs.

%
%
\null\textbf{A posteriori error bound.} 
Approximation properties of PINNs have been studied in \cite{Shin2020,molinaro2023,mishra2022-inverse,mishra2023}.
We report here the results proved in \cite{molinaro2023,mishra2023} about the generalization error for feed--forward PINNs.

Let us consider a regular domain $\Omega\subset\mathbb R^d$, a Banach space $V$ (for second--order elliptic PDEs, a typical choice is $H^2(\Omega)$), a boundary operator $\mathcal{B}:V\to L^2(\partial\Omega)$, $f\in L^2(\Omega)$, $g\in L^2(\partial\Omega)$, and the problem of finding $u\in V$ solution of 
\begin{eqnarray}\label{eq:pde_abstract}
\left\{
\begin{array}{ll}
\mathcal{D}(u)=f & \mbox{ in }\Omega\\
\mathcal{B}(u)=g & \mbox{ on }\partial\Omega.
\end{array}\right.
\end{eqnarray}
Let us assume that (\ref{eq:pde_abstract}) admits a unique solution. We define the loss function 
\begin{equation}\label{eq:cost-pinn-res}
{\mathcal{L}}(\mathbf{w})=
\mathcal{L}_{PDE}(\mathbf{w})+\alpha_{BC}\mathcal{L}_{BC}(\mathbf{w}),
\end{equation}
where $\alpha_{BC}\in\mathbb R$ is a suitable weight, while 
$\mathcal{L}_{PDE}$ and $\mathcal{L}_{BC}$ are defined by
\begin{eqnarray*}
\begin{array}{l}
\displaystyle
\mathcal{L}_{PDE}(\mathbf{w})=\frac{1}{2}\sum_{i=1}^{N_{PDE}}
\omega_i^{PDE}|\mathcal{D}(u_{NN}(\mathbf{x}_i^{PDE};\mathbf{w}))-f(\mathbf{x}_i^{PDE})|^2\\[1mm]
\displaystyle
\mathcal{L}_{BC}(\mathbf{w})=\frac{1}{2}\sum_{i=1}^{N_{BC}}
\omega_i^{BC}|\mathcal{B}(u_{NN}(\mathbf{x}_i^{BC};\mathbf{w}))-g(\mathbf{x}_i^{BC})|^2
\end{array}
\end{eqnarray*}
with $\omega_i^{PDE}$ and $\omega_i^{BC}$ suitable quadrature weights. In the case of the Monte Carlo formula, we have $\omega_i^{PDE}=1/N_{PDE}$, $\omega_i^{BC}=1/N_{BC}$.
We assume that, when increasing the number of collocation nodes, we keep constant the ratio between the number of nodes in the interior and on the boundary, so that we can define a single parameter $N$ such that $N_{PDE} \approx N$ and $N_{BC} \approx N$.

Then, we train the network by looking for $\mathbf w^*$ solution of \begin{equation*} 
\mathbf{w}^*=\argmin{\mathbf{w}\in\mathbb{R}^M}
{\mathcal{L}}(\mathbf{w}),
\end{equation*}
 and define $u^*=u_{NN}(\cdot; \mathbf{w}^*)$ as the PINN solution.

The following theorem is adapted from \cite{mishra2023, molinaro2023}:
\begin{theorem}
Provided that the following stability assumption
\begin{eqnarray*}
\|v_1-v_2\|_V\leq &C_{PDE}(\|v_1\|_V,\|v_2\|_V)\cdot\Big[
\|\mathcal{D}(v_1)-\mathcal{D}(v_2)\|^2_{L^2(\Omega)}\\
& +\alpha_{BC}
\|\mathcal{B}(v_1)-\mathcal{B}(v_2)\|^2_{L^2(\partial\Omega)}\Big]^{1/2},
\end{eqnarray*}
holds for any $v_1,\ v_2$ belonging to a closed subspace of $V$,
then there exist two positive constants $c_1$ and $c_2$ such that
\begin{equation}\label{eq:PINN_error}
\|u^*-u\|_V\leq c_1 {\mathcal{L}}(\mathbf{w}^*)^{1/2}+c_2 N^{-s/2}, 
\end{equation}
where $s$ is the accuracy order of the quadrature formula.
\end{theorem}
The term $\|u^*-u\|_V$ represents the generalization error, $\mathcal{L}(\mathbf{w}^*)^{1/2}$ is the training error, while $c_2N^{-s/2}$ bounds the quadrature error.
The constant $c_1$ depends on the exact solution, the NN solution and the PINN architecture (i.e., the number of layers, neurons and training collocation points), while $c_2$ depends on the exact solution and the quadrature formula.

The error estimate (\ref{eq:PINN_error}) establishes that minimizing the residual of the PDE (including the residual of boundary conditions) leads to the control of the generalization error, provided that a sufficient number of collocation nodes is employed. We refer to \cite{mishra2023,molinaro2023} for a more in--depth analysis of this topic.

\null\textbf{PINNs for inverse problems.}
As anticipated, the scope of PINNs is not limited to the solution of PDEs. They can be effectively employed in a variety of scenarios, including the solution of inverse problems. In the following, we illustrate how PINNs can be used to address the identification of parameters of a PDE, starting from some suitable measures on the solution.
Let us consider the problem (\ref{eq:strong-elliptic}), we know the functions $f$ and $g$ and want to estimate the parameter $\gamma\in{\mathbb R}$ starting from $N_{obs}$ observations $u_i^{obs}=u(\mathbf x_i^{obs})$, with $i=1,\ldots,N_{obs}$. This is an instance of \emph{inverse problem} and, in particular, a \emph{parameter identification problem}.

Should we apply a classical approach based on FEM, we should face the \emph{constrained minimization problem}
\begin{eqnarray}\label{eq:inverse-problem}
\begin{array}{l}
    \displaystyle \min_{\gamma>0,\ u_h\in V_h} \sum_{i=1}^{N_{obs}}|u_i^{obs}-u_h(\mathbf x_i^{obs})|^2\\
    \displaystyle s.t.\ \int_\Omega\nabla u_h \cdot \nabla v_h +\int_\Omega \gamma u_h v_h =\int_\Omega f v_h +\int_{\partial\Omega} gv_h \quad \forall v_h\in V_h,
    \end{array}
\end{eqnarray}
for a convenient finite dimensional subspace $V_h$ of the Sobolev space $H^1(\Omega)$.
Then, we could invoke standard algorithms for constrained optimization (like, e.g., those based on the ``trust region'' or the ``interior point'' methods \cite{Nocedal}) or resort to the optimal control theory \cite{Lions-1971-OCS} and solve the related optimality system. In the latter case, if we adopted an iterative solver to solve the optimality system, at each iteration we should solve two FEM problems: one primal (i.e. problem (\ref{eq:inverse-problem})$_2$) and one adjoint of (\ref{eq:inverse-problem})$_2$, see \cite{Lions-1971-OCS, Manzoni-Quarteroni-Salsa-book}.

Instead, the PINN approach consists of defining the loss function
\begin{equation}\label{eq:pinn-loss-obs}
    \mathcal L(\mathbf w,\gamma)= \mathcal L_{obs}(\mathbf w)+
    \alpha_{PDE}\mathcal L_{PDE}(\mathbf w,\gamma) +\alpha_{BC}\mathcal L_{BC}(\mathbf w)
\end{equation}
with 
\begin{eqnarray*}
    \mathcal L_{obs}(\mathbf w)&=&\frac{1}{N_{obs}}\sum_{i=1}^{N_{obs}}|u_i^{obs}-u_{NN}(\mathbf x_i^{obs};\mathbf w)|^2,\\
    \mathcal L_{PDE}(\mathbf w,\gamma)&=&\frac{1}{2}\sum_{i=1}^{N_{PDE}}\omega_i^{PDE}|\mathcal D(u_{NN}(\mathbf x_i^{PDE};\mathbf w);\gamma)-f(\mathbf x_i^{PDE})|^2,
\end{eqnarray*}
 and $\mathcal L_{BC}(\mathbf w)$ as in (\ref{eq:loss-pinn-resbc}),
and computing
\begin{equation*}
    \mathbf w^*, \gamma^* =\argmin{\mathbf w\in \mathbb R^M, \gamma\in \mathbb R^+} \mathcal L(\mathbf w,\gamma).
\end{equation*}
The advantages of PINNs compared with the classical approach are that only a minimal modification with respect to the direct problem is required (we have added a new parameter to the array $\mathbf w$ of the NN parameters and the term $\mathcal L_{obs}$ to the loss function.). Moreover, PINNs also work in case of defective data, for instance, when some of the boundary conditions are missing. The flexibility of PINNs in seamlessly 
addressing different types of problems is a key feature that makes them particularly appealing in the context of inverse problems.

In general terms, the problem of parameter estimation with PINNs can be framed as follows:
\begin{equation} \label{eqn:PINN-inverse}
    \mathbf w^*, \boldsymbol{\gamma}^* = \argmin{\mathbf w\in \mathbb{R}^M, \boldsymbol{\gamma}\in \mathbb{R}^P} 
    \left[\frac{1}{N_{obs}} \sum_{i=1}^{N_{obs}} \big| u_i^{obs} - u_{NN}(\mathbf x_i^{obs}; \mathbf{w}) \big|^2
    + \mathcal{L}_{phys}(u_{NN}(\cdot; \mathbf{w}), \boldsymbol{\gamma})\right],
\end{equation}
where the physics--based regularization term is defined as
\begin{equation}\label{eq:pinn-loss-phys}
    \begin{split}
        \mathcal{L}_{phys}(u, \boldsymbol{\gamma}) =
        &\frac{\alpha_{PDE}}{2N_{PDE}} \sum_{i=1}^{N_{PDE}} \big( Res^{PDE}[u(\mathbf{x}_i^{PDE}); \boldsymbol{\gamma}] \big)^2 \\
        &+ \frac{\alpha_{BC}}{2N_{BC}} \sum_{i=1}^{N_{BC}} \big( Res^{BC}[u(\mathbf{x}_i^{BC}); \boldsymbol{\gamma}] \big)^2.
    \end{split}
\end{equation}

\null\textbf{Multifidelity PINNs.}
Inverse problems are often ill--posed due to data scarcity. Sometimes a surplus of boundary conditions is given on a small subset $\Gamma$ of the boundary $\partial \Omega$, while no information is available on $\partial\Omega\setminus \Gamma$. Other times, boundary data are completely missing, while only measurements inside the domain are available. Consequently, multiple local minima and challenging optimization landscapes occur. 
In these circumstances, despite the \emph{physics--informed regularization}, achieving convergence and stability remains difficult. To address this issue, a multi--fidelity approach can be introduced \cite{meng2020composite}, leveraging a combination of low-- and high--fidelity models. The purpose is to enhance convergence of the optimization process and improve robustness with respect to noise.
Two different datasets are typically used. The first one $(\mathbf x_i^{obs}, u_i^{obs})$ with $i=1,\ldots, N_{obs}$, feeds the high--fidelity model. The second one, $(\mathbf x_i^{low}, u_i^{low})$, for $i = 1, \dots, N_{low}$, feeds the low--fidelity model and is typically characterized by a larger number of data but lower fidelity, such as those collected from numerous low--quality sensors. 

{Typically,} the low--fidelity model $u_L(\mathbf x;\mathbf w_L)$ is a standard FFNN.
The high--fidelity model $ u_H $ can be expressed generally as:
\begin{equation} \label{eqn:PINN_multifidelity}
    u_H(\mathbf{x}; \mathbf{w}) = \mathcal{L}_H(\mathbf{x}, u_L(\mathbf{x}; \mathbf{w}_L); \mathbf{w}_H)
    + \mathcal{NN}_H(\mathbf{x}, u_L(\mathbf{x}; \mathbf{w}_L); \mathbf{w}_H),
\end{equation}
where $ \mathcal{L}_H $ is a linear transformation of $u_L$, and $ \mathcal{NN}_H $ is a neural network correcting the low--order model $u_L$. 
The rationale for this linear/non--linear splitting is that, if the low--fidelity model is well--designed, one would expect a strong linear correlation between $u_H$ and $u_L$; by explicitly incorporating a linear correlation term, the learning process is guided to more effectively capture this type of relationship.
Here, $ \mathbf{w} = [\mathbf{w}_L, \mathbf{w}_H] $ encapsulates the parameters of both low-- and high--fidelity models. 

The optimization problem in the multi--fidelity setting becomes:
\begin{equation*}
    \begin{split}
        \mathbf w^*, \boldsymbol{\gamma}^* = \argmin{\mathbf w\in \mathbb{R}^M, \boldsymbol{\gamma}\in \mathbb{R}^P} 
        &\Big[\frac{1}{N_{obs}} \sum_{i=1}^{N_{obs}} \big| u_i^{obs} - u_H(\mathbf x_i^{obs}; \mathbf{w}) \big|^2 
        + \mathcal{L}_{phys}(u_H(\cdot; \mathbf{w}), \boldsymbol{\gamma}) \\
        & +\frac{\alpha_{low}}{N_{low}} \sum_{i=1}^{N_{low}} \big| u_i^{low} - u_L(\mathbf x_i^{low}; \mathbf{w}_L) \big|^2\Big], 
    \end{split}
\end{equation*}
where $u_H$ is given by \eqref{eqn:PINN_multifidelity}.
Here, the low--fidelity term $ \frac{\alpha_{low}}{N_{low}} \sum |u_i^{low} - u_L|^2 $ incorporates additional information into the optimization, improving the model's ability to generalize even with sparse high--fidelity data.

A more advanced multi--fidelity approach  \cite{regazzoni2021physics} involves leveraging a surrogate model trained on data generated through a numerical solver, by following e.g. the methods described in Sec.~\ref{sec:DM2ML}. 
This allows us to capture parametric dependencies in the solution. The corresponding optimization problem for training $ u_L $ is:
\begin{equation*}
    \mathbf w_L^* = \argmin{\mathbf w_L \in \mathbb{R}^{M_L}} 
    \left[\frac{1}{N_{low}} \sum_{i=1}^{N_{low}} \big| u_i^{low} - u_L(\mathbf{x}_i^{low}, \boldsymbol{\gamma}_i^{low}; \mathbf{w}_L) \big|^2 + R(\mathbf{w}_L)\right],
\end{equation*}
where $ R(\mathbf{w}_L) $ is a regularization term that may include physics--based constraints.
Once trained, the parametrized low--fidelity model $ u_L $ is fixed (i.e., its parameters are kept equal to $\mathbf{w}_L^*$) and used as a prior for the high--fidelity model:
\begin{equation*}
    u_H(\mathbf{x}, \boldsymbol{\gamma}; \mathbf{w}_H) = 
    \mathcal{L}_H(\mathbf{x}, u_L(\mathbf{x}, \boldsymbol{\gamma}; \mathbf{w}_L^*); \mathbf{w}_H)
    + \mathcal{NN}_H(\mathbf{x}, u_L(\mathbf{x}, \boldsymbol{\gamma}; \mathbf{w}_L^*); \mathbf{w}_H).
\end{equation*}
where one often sets $\mathcal{L}_H(\mathbf{x}, u_L; \mathbf{w}_H) \equiv u_L$ to enforce the low-fidelity model as a prior. We remark that, unlike in \eqref{eqn:PINN_multifidelity}, the multi--fidelity expression accounts for the parametric dependence of the solution, that is the NN explicitly depends on the parameters $\boldsymbol \gamma$.

The final optimization for the inverse problem becomes:
\begin{eqnarray} \label{eqn:MPINN-enhanced}
\begin{array}{ll}
    \displaystyle \mathbf w_H^*, \boldsymbol{\gamma}^* = \argmin{\mathbf w_H \in \mathbb{R}^{M_H}, \boldsymbol{\gamma} \in \mathbb{R}^P} 
    &\displaystyle \Big[\frac{1}{N_{obs}} \sum_{i=1}^{N_{obs}} \big| u_i^{obs} - u_H(\mathbf{x}_i^{obs}, \boldsymbol{\gamma}; \mathbf{w}_H) \big|^2\\[2mm]
    & \displaystyle + \mathcal{L}_{phys}(u_H(\cdot, \boldsymbol{\gamma}; \mathbf{w}_H), \boldsymbol{\gamma})\Big].
    \end{array}
\end{eqnarray}

This approach leverages the strengths of both low-- and high--fidelity models, balancing computational efficiency with accuracy. It is particularly advantageous in settings where parametric variations significantly influence the solution (see e.g. the application reported in Sec.~\ref{sec:mpinn-ionic}). 

\null\textbf{PINNs for data fitting.}
Very frequently, observations are affected by noise, so we might be interested in estimating the denoised solution of a PDE, starting from noising data.
This represents another type of problem that can be effectively addressed by PINNs.
Let us consider the example of estimating the blood velocity $\mathbf u$ in a vessel. Let $\Omega\subset \mathbb R^3$ be the computational domain (modelled by a curved cylinder), with $\Gamma_{wall}$, $\Gamma_{in}$ and $\Gamma_{out}$ representing impermeable, inflow, and outflow boundaries, respectively. For any time $t\in(0,T)$ of the simulation, the velocity $\mathbf u=\mathbf u(\mathbf x, t)$ satisfies the Navier--Stokes equation
\begin{eqnarray}\label{eq:NS-denoising}
\left\{\begin{array}{ll}
L_M(\mathbf u,p):=\rho\frac{\partial\mathbf u}{\partial t}-\mu\Delta \mathbf{u} +\rho(\mathbf{u}\cdot\nabla)\mathbf{u}+\nabla p =\mathbf{0} & \mbox{ in }\Omega\\[1mm]
L_D(\mathbf u):=\nabla\cdot\mathbf{u}=0 & \mbox{ in }\Omega\\[1mm]
\mathbf{u}=\mathbf{0} & \mbox{ on }\Gamma_{wall},
\end{array}\right.
\end{eqnarray}
and we dispose of the noisy observation $\mathbf u_i^{obs}\simeq \mathbf u(\mathbf x_i^{obs},t_i^{obs})$ for $i=1,\ldots,N_{obs}$, with $\mathbf x_i^{obs}\in \Omega$ and $t_i^{obs}\in(0,T)$.
We are interested in estimating the velocity of the blood at certain points of the inflow and outflow boundaries and at given times, by using a PINN.

We might design a single NN that, given in input the independent variables $\mathbf x\in\mathbb R^3$ and $t\in\mathbb R$, provides both the velocity $\mathbf u$ and the pressure $p$ as outputs; otherwise, we might use two different NNs, one providing the velocity and the other the pressure. In both cases, we denote by $\mathbf u_{NN}(\mathbf x, t; \mathbf w)$ and $p_{NN}(\mathbf x, t; \mathbf w)$ the solution computed by the PINN, and denote by $\mathcal I_M$ and $\mathcal I_D$ two (not necessarily disjoint) subsets of the indices set $\mathcal I_{obs}=\{1,\ldots, N_{obs}\}$. Then we define the losses
\begin{eqnarray*}
    \mathcal L_{obs}(\mathbf w)&=&\frac{1}{N_{obs}}\sum_{i\in \mathcal I_{obs}} \|\mathbf u_i^{obs}-\mathbf u_{NN}(\mathbf x_i^{obs}, t_i^{obs};\mathbf w)\|^2,\\
    \mathcal L_{M}(\mathbf w)&=&\frac{1}{|\mathcal I_M|}\sum_{i\in \mathcal I_M}\|L_M(\mathbf 
    u_{NN}(\mathbf x_i^{M}, t_i^{M};\mathbf w), 
    p_{NN}(\mathbf x_i^{M}, t_i^{M};\mathbf w))\|^2,\\
    \mathcal L_{D}(\mathbf w)&=&\frac{1}{|\mathcal I_D|}\sum_{i\in \mathcal I_D} (L_D(\mathbf 
    u_{NN}(\mathbf x_i^{D}, t_i^{D};\mathbf w)))^2,
\end{eqnarray*}
and look for
\begin{equation*}
    \mathbf w^*=\argmin{\mathbf w\in \mathbb R^M} \left[\mathcal L_{obs}(\mathbf w)+
    \alpha_{M}\mathcal L_{M}(\mathbf w) +\alpha_{D}\mathcal L_{D}(\mathbf w)\right],
\end{equation*}
with $\alpha_M$ and $\alpha_D$ suitable weights.
The term $\mathcal{L}_{obs}$ of the loss function incorporates data, while $\mathcal{L}_{M}$ and $\mathcal{L}_{D}$ encode the learning biases.
One way to interpret this approach is to view it as a fitting problem, where the terms incorporating physics act as regularization terms, favouring solutions that adhere to known physical principles.

We observe that also in this case the formulation of the problem is very similar to that of the direct one. 

Thanks to their flexibility, PINNs can be applied even when boundary conditions or initial conditions are missing or not completely known. PINNs can be effective for ill--posed problems and in the small data regime (indeed when data are scarce, learning biases supersede missing measurements).
%
%

\null\textbf{Enforcement of Dirichlet boundary conditions.}
We now discuss how to enforce Dirichlet boundary conditions in the PINN framework. The approaches discussed here can be applied to the different types of problems we have discussed so far, including the resolution of PDE, inverse problems and data fitting.
Let us consider the second--order elliptic equation (\ref{eq:strong-elliptic}) but with Dirichlet boundary conditions instead of the Neumann ones:
\begin{eqnarray}\label{eq:strong-elliptic-dir}
\left\{\begin{array}{ll}
-\Delta u+\gamma u=f & \mbox{ in }\Omega\\[2mm]
u=g &\mbox{ on }\partial\Omega.
\end{array}\right.
\end{eqnarray}
We may work following different approaches. 

The first way consists of defining the residual corresponding to the boundary condition as
\begin{equation}\label{eq:res_dirchlet}
Res^{BC}[u(\mathbf x)]=u(\mathbf x)-g(\mathbf x) \qquad \mathbf x\in \partial \Omega
\end{equation}
and using this one instead of \ref{eq:pinn-residuals}$_2$ inside the definition of (\ref{eq:loss-pinn-resbc}) \cite{Lagaris2000}. 

This approach could be interpreted as analogous of 
 the \emph{penalty method} of Babu\v{s}ka \cite{Babuska1973} in the Galerkin framework, for which we look for $u\in V=H^1(\Omega)$ such that
\begin{equation*}
    \int_\Omega \left(\nabla u\cdot \nabla v + \gamma u v -fv\right) +\alpha_{BC} \int_{\partial\Omega} (u -g)v=0\quad \forall v\in V,
\end{equation*}
where $\alpha_{BC}>0$ is a suitable penalization parameter.

\medskip
A second way \cite{Lagaris1998} consists of defining the solution as
\begin{equation}\label{eq:unn_dir}
u_{NN}(\mathbf x; \mathbf w) = G(\mathbf x)+A(\mathbf x)\tilde u_{NN}(\mathbf x; \mathbf w)
\end{equation}
where $G(\mathbf x)$ is a \emph{lifting} of $g(\mathbf x)$ to $\Omega$, i.e. a function $G:\overline\Omega\to \mathbb R$ such that $G|_{\partial\Omega}=g$ and contains no adjustable parameters; $A(\mathbf x)$ is a function named \emph{mask} defined on $\overline\Omega$ such that $A|_{\partial \Omega}=0$ (for instance, if $\Omega=(0,1)^2$, we could take $A(\mathbf x)=x_1x_2(1-x_1)(1-x_2)$); while $\tilde u_{NN}(\mathbf x; \mathbf w)$ is a single--output FFNN with parameters $\mathbf w$ and input $\mathbf x$.

The optimal values of the parameters are computed by minimizing the loss function 
\begin{equation}
    \mathcal L(\mathbf w)=\mathcal L_{PDE}(\mathbf w)=\frac{1}{2N_{PDE}}\sum_{i=1}^{N_{PDE}}(Res^{PDE}(u_{NN}(\mathbf x_i^{PDE};\mathbf w))^2.
\end{equation}
Other strategies to impose Dirichlet conditions rely on formulating the problem as a constrained optimization problem via the augmented Lagrangian method (see \cite{Shin2024} and the references therein).

\subsubsection{Variational PINNs (VPINNs)}
\label{sec:VPINN}
 
Let us consider the elliptic problem (\ref{eq:weak_elliptic}) with Neumann boundary conditions, the NN space $V_{NN}=\big\{u(\mathbf{x})=u_{NN}(\mathbf{x};\mathbf{w}),\ 
\mbox{ with }\mathbf{w}\in\mathbb{R}^{M}\big\}$ and the FEM space $V_h =\big\{u(\mathbf{x})=\sum_{j=1}^{N_h}u_j\,\varphi_j(\mathbf{x})\big\}$.
The Petrov--Galerkin formulation of (\ref{eq:weak_elliptic}) reads:
find 
\begin{equation}\label{eq:Petrov-Galerkin}
u_{NN}\in V_{NN}\quad s.t.\quad a(u_{NN},v_h)=F(v_h)\qquad \forall v_h\in V_h.
\end{equation}
Looking for the solution of (\ref{eq:weak_elliptic}) with VPINN means to minimize the residual of the Petrov--Galerkin equation (\ref{eq:Petrov-Galerkin}), that is we look for \cite{kharazmi2019-vpinn, Kharazmi2021}
\begin{equation}
    \mathbf{w}^*=\argmin{\mathbf{w}\in\mathbb{R}^M}
\frac{1}{2N_h}\sum_{i=1}^{N_h}|F(\varphi_i)-a(u_{NN}(\cdot;\mathbf{w}),\varphi_i)|^2.
\end{equation}
This weak formulation has many advantages with respect to the strong one. First of all, reducing the derivative order by integration by parts leads to less computationally expensive algorithms. Moreover, in some particular cases, the loss function can be expressed analytically, making it possible to obtain more accurate estimates than those obtained for PINNs. Finally, ad--hoc quadrature formulas require a number of quadrature nodes noticeably lower than those required by PINNs to compute the residuals in strong form.

We notice that standard PINNs can be seen as a special kind of VPINNs, with basis functions equal to Dirac deltas centred at collocation points.

\subsubsection{Deep Ritz Method (DRM)}
\label{sec:DeepRitz}

When the bilinear form $a$ associated with the differential operator is symmetric like the one of the problem (\ref{eq:strong-elliptic}), the weak form of the PDE can be interpreted as the Euler--Lagrange equation of an energy--functional. In particular, the weak equation (\ref{eq:weak_elliptic0}) is the Euler--Lagrange equation of the energy--functional
\begin{equation}
    \mathcal E(v)=\frac{1}{2}\int_\Omega |\nabla v|^2 +\frac{\gamma}{2}\int_\Omega |v|^2 -\int_\Omega fv -\int_{\partial\Omega} gv \quad \forall v\in V=H^1(\Omega), 
\end{equation}
so that $u=\argmin{v\in V}\mathcal E(v)$.

\begin{figure}
    \centering
    \includegraphics[width=0.7\linewidth]{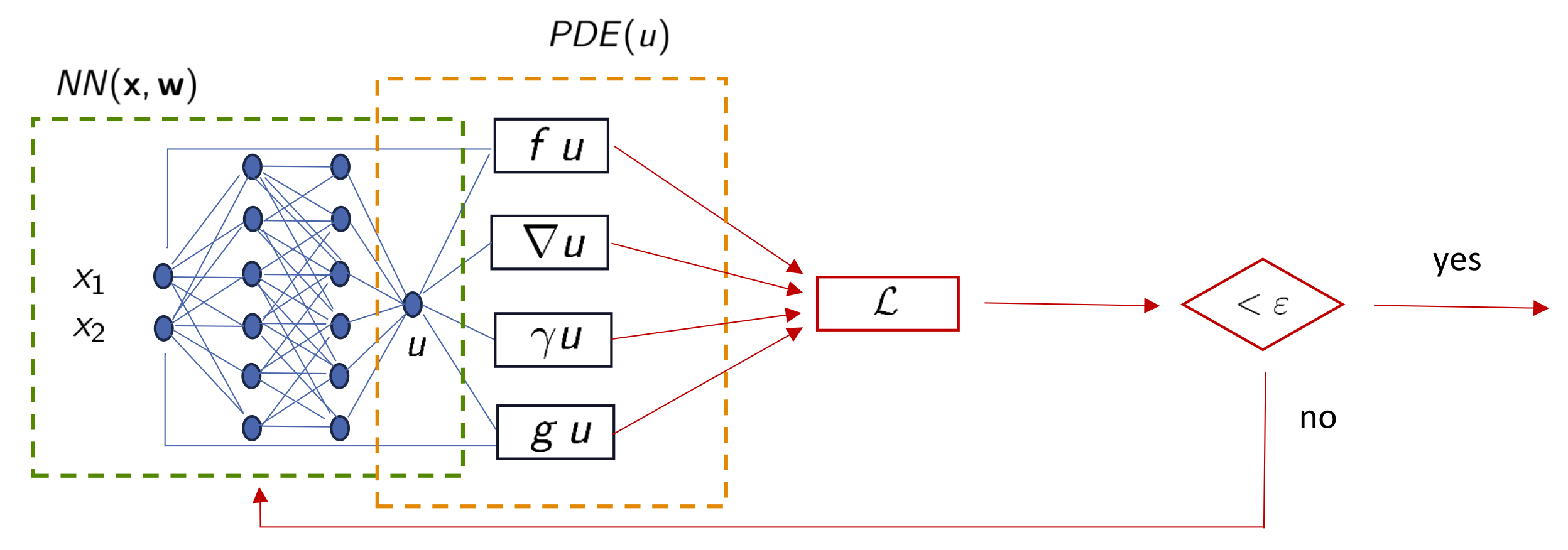}
    \caption{Deep Ritz Models}
    \label{fig:DRM}
\end{figure}

Chosen $N_{PDE}$ points $\mathbf x_i^{PDE}$ in $\Omega$  and $N_{BC}$ points $\mathbf x_i^{BC}$ on $\partial\Omega$ randomly sampled, the Deep Ritz Method \cite{E2018} looks for the parameters
\begin{equation}
    \mathbf w^*=\argmin{\mathbf w\in \mathbb R^M}\mathcal L(\mathbf w)
\end{equation}
where the loss function is defined by 
\begin{eqnarray*}
\begin{array}{lll}
\mathcal L(\mathbf w)=&\displaystyle\frac{1}{N_{PDE}}\sum_{i=1}^{N_{PDE}}
&\displaystyle \Big[\frac{1}{2}|\nabla u_{NN}(\mathbf x_i^{PDE};\mathbf w)|^2\\[2mm]
&&\displaystyle +\frac{\gamma}{2}|u_{NN}(\mathbf x_i^{PDE};\mathbf w)|^2
-f(\mathbf x_i^{PDE})u_{NN}(\mathbf x_i^{PDE};\mathbf w)\Big]\\[2mm]
&\multicolumn{2}{l}{\displaystyle
-\frac{1}{N_{BC}}\sum_{i=1}^{N_{BC}}g(\mathbf x_i^{BC})u_{NN}(\mathbf x_i^{BC};\mathbf w).}
\end{array}
\end{eqnarray*}
The solution $u_{NN}$ is a deep neural network and it is used to evaluate the loss function (see Fig. \ref{fig:DRM}). As for PINNs, partial derivatives are computed by automatic differentiation recalling that $u_{NN}$ can be written as a composition of functions like in (\ref{eq:f_NN}). In \cite{E2018} the authors designed a NN with an arbitrary number of blocks (it will be a hyperparameter), each one composed of two layers with the corresponding activation function and an identity map (like in ResNet \cite{He-resnet-2016}) put in parallel.

As said before, not all PDEs admit an energy formulation, so DRM has limited applicability. Similarly to PINNs, also DRM leads to non--linear problems (in $\mathbf w$) even when the PDE is linear (in $u$), so in such cases, classical Galerkin methods typically outperform DRM. Moreover, even when the energy functional $\mathcal E$ is convex, the loss function is non--convex and can feature local minima.

Comparing DRM with PINNs, we notice that:
\emph{(i)} DRM does not include hyperparameters (weights of the residuals for PINNs), so it is not subject to the tuning issue; 
\emph{(ii)} the treatment of Dirichlet boundary conditions is less easy than for PINNs, however, the formulation (\ref{eq:unn_dir}) with lifting and mask does work;
\emph{(iii)} the extension to defective and inverse problems is not easy.

\subsubsection{Optimizing a DIscrete Loss (ODIL)}

The Optimizing DIscrete Loss (ODIL) method introduced in \cite{Karnakov2024} computes the solution of a PDE as the minimizer of a loss function which is the sum of discrete residuals of the PDE at a set of points. While PINNs evaluate residuals only by using automatic differentiation and do not discretize the PDE, ODIL first discretizes the PDE by a grid method (a Full Order Method), e.g., finite differences or finite volumes, so that the residuals of the discrete system are algebraic expressions. To minimize the loss function, ODIL exploits algorithms typically used in ML, like the Adam method (see Sect. \ref{sec:optimization}), or Gauss--Newton \cite{Zhu1997, Liu1989} and quasi--Newton methods \cite{Nocedal, Quarteroni-CS-2014}. In all cases, gradients are evaulated by automatic differentiation. To model unknown coefficient functions of the PDE, feed forward neural networks are employed.

ODIL shares common features with the least--squares finite element method and the discretize--then--differentiate approach for PDE--constrained optimization problems \cite{Bochev1998, Gunzburger2002}. 

By minimizing the discrete loss function, ODIL maintains the accuracy and conservation properties of the full order discretization method employed. 
It has been used to solve both forward and inverse problems, like field reconstruction and body shape inferring, both starting from noisy data \cite{Karnakov2024}. 

ODIL is formulated as follows. After discretizing the PDE with a grid method, let $\mathbf u\in\mathbb R^{N_u}$ be the array of the unknown degrees of freedom and $\mathbf w\in \mathbb R^{N_w}$ the array of unknown parameters. $\mathbf w$ could include coefficients of PDEs, like e.g., the constant $\gamma$ of (\ref{eq:strong-elliptic}), or parameters of a neural network modelling coefficient functions of PDEs, like conductivity function in the heat equation (see (\ref{eq:heat})).

Let $F_i(\mathbf u, \mathbf w)=0$ for $i=1,\ldots, N_F$ be the residuals associated with the discretized PDE, boundary conditions, initial conditions evaluated at some specified (either internal or boundary) grid--points, or other constraints. We define the loss function as
\begin{equation}\label{eq:odil-loss}
\mathcal L(\mathbf u;\mathbf w)=\frac{1}{N_F}\sum_{i=1}^{N_F}[F_i(\mathbf u;\mathbf w)]^2
\end{equation}
and we look for
\begin{equation}\label{eq:odil-min}
\mathbf u^*,\,\mathbf w^*=\argmin{\mathbf u\in\mathbb R^{N_u},\ \mathbf w\in \mathbb R^{N_w}} \mathcal L(\mathbf u; \mathbf w).
\end{equation}

As an example, let us consider the heat equation 
\begin{eqnarray}\label{eq:heat}
\left\{\begin{array}{ll}
\displaystyle \frac{\partial u}{\partial t}-\frac{\partial }{\partial x}\left(k(u)
\frac{\partial u}{\partial x}\right)=0 & \mbox{ in }(0,L)\times(0,T)\\[1mm]
u(0,t)=u(L,t)=0 & \mbox{ in }(0,T)\\
u(x,0)=u_0(x) & \mbox{ in }(0,L)
\end{array}\right.
\end{eqnarray}
modelling the temperature $u(x,t)$ in a bar of length $L$ on the time interval $(0,T)$, starting from the initial temperature $u_0(x)$.
We want to infer the conductivity $k(u)$ of the medium, provided that we know the temperature $U_k$ at a finite set of points $(\tilde x_k,\tilde t_k)$ for $k=1,\ldots, N_{obs}$, i.e., $u(\tilde x_k, \tilde t_k)=u_{0,k}.$
We split the space interval $(0,L)$ into $N_x$ cells of length $h$ and the time interval $(0,T)$ in $N_t$ time--intervals of length $\Delta t$. Then, we discretize (\ref{eq:heat})$_1$ with centred finite deferences in space and the Crank--Nicolson method in time. Moreover, we model the unknown conductivity function $k$ with a feed--forward neural network with one input (temperature $u$ at a point), one output (value of conductivity at the same point) and parameters $\mathbf w$, i.e., $k=k(u;\mathbf w)$. 

We denote by $\mathbf u=[u_j^n]_{j=0,\ldots,N_x,\ n=0, \ldots,N_t}$ the global  array of degrees of freedom and by $\mathbf u^n =[u_j^n]_{j=0,\ldots,N_x}$ the array of degrees of freedom at time $t_n$.

For any $j=1,\ldots, N_x-1$  and $n=1,\ldots, N_t$, we define the discrete residuals
\begin{equation}\label{eq:odil-residuals}
F_j^n(\mathbf u;\mathbf w)=\frac{u_j^{n+1}-u_j^n}{\Delta t}-
\frac{1}{2}[G_j(\mathbf u^n; \mathbf w)+G_j(\mathbf u^{n+1}; \mathbf w)],
\end{equation}
where
\begin{equation*}
    G_j(\mathbf u^n;\mathbf w)=\frac{1}{h^2}\left[
    k\left(\frac{u_{j+1}^n+u_j^n}{2};\mathbf w\right)(u_{j+1}^n-u_j^n) -
    k\left(\frac{u_{j}^n+u_{j-1}^n}{2};\mathbf w\right)(u_{j}^n-u_{j-1}^n)\right].
\end{equation*}
The boundary and initial conditions $u_0^n=0$, $u_{N_x}^n=0$ for $n=1,\ldots,N_t$ and $u_j^0=u_0(x_j,0)$ for $j=1,\ldots,N_x-1$ are exploited in (\ref{eq:odil-residuals}). Finally, the loss function is defined by
\begin{equation}\label{eq:odil-heat-loss}
\mathcal L(\mathbf u; \mathbf w)=\frac{1}{N_x N_t}\sum_{j=1}^{N_x}\sum_{n=1}^{N_t}[F_j^n(\mathbf u;\mathbf w)]^2+\frac{w_{obs}}{N_{obs}}\sum_{k=1}^{N_{obs}}(u_{j(k)}^{n(k)}-u_{0,k})^2,
\end{equation}
where $j(k)$ and $n(k)$ are the space and time indices of the cell containing the measurement point $(\tilde x_k,\tilde t_k)$, while $w_{obs}\in \mathbb R$ is a suitable weight.

In \cite{Karnakov2024} the authors compare ODIL and PINNs on three forward problems (1D wave equation,2D Poisson equation, and 2D lid--driven cavity problem) and one inverse problem (inferring conductivity from temperature, 1D problem) showing that, on these tests, ODIL outperforms PINNs in terms of accuracy and computational costs.

%
%
\subsection{Operator learning}\label{sec:operator-learning}

ML algorithms have primarily been developed for learning mappings between finite--dimensional spaces, specifically the 
input and output spaces $\mathcal{X} \subseteq \mathbb{R}^n$ and $\mathcal{Y} \subseteq \mathbb{R}^m$ within the supervised learning framework of Sec.~\ref{sec:supervised_learning}.
However, in problems arising in the realm of Scientific Computing, one has often to deal with maps between infinite dimensional functional spaces. Pivotal examples are represented by the data-to-solution maps underlying PDEs, where -- for example -- one is interested in the map from a forcing term to the solution of an elliptic PDE, or from a time-dependent boundary condition to the solution of a hyperbolic PDE.
Traditional machine learning algorithms (e.g., FFNNs) do not generalize straightforwardly, as they are designed for finite--dimensional spaces.
The emerging field of \textit{Operator learning} focuses on ML algorithms capable of learning \textit{operators}, understood as mappings between function spaces.

Operator learning represents an approach used in the framework described in Sec.~\ref{sec:DM2ML} to develop surrogate models of high-fidelity DMs, by leveraging training data that are obtained as solutions of the high--fidelity DM itself. Once trained, these surrogate models can replace high--fidelity DMs to deliver rapid approximations of physical problems, particularly in many-query scenarios where computational cost is a concern.

We are therefore interested in approximating either a linear or a non--linear operator that can be written as a mapping from one functional space to another. To outline the key--points of this topic, let us start with a simple example.
Let us suppose we want to learn the differential operator that governs a linear second--order elliptic PDE given on the domain $\Omega\subset \mathbb R^d$, with $d=1,2,3$,  with homogeneous Dirichlet boundary conditions. Given the functions $\mu_1:\Omega\to\mathbb R^+$ and $\mu_2:\Omega\to\mathbb R$, the strong form of the PDE reads: find the solution $y:\Omega\to \mathbb R$ such that
\begin{eqnarray}\label{eq:Elliptic}
\left\{\begin{array}{ll}
-\nabla\cdot(\mu_1(x)\nabla y(x))=\mu_2(x) & \forall x\in\Omega\\[1mm]
y(x)=0 & \forall x\in\partial\Omega.
\end{array}\right.
\end{eqnarray}
Let the parameter functions be sufficiently regular, i.e.,
\begin{equation}
u=(\mu_1,\mu_2)\in V_u:=\{\mu_1\in L^\infty(\Omega):\ \mu_1\geq \overline{\mu_1}>0 \ a.e. \ \Omega\}\times L^2(\Omega),
\end{equation}
then we can write the weak form of (\ref{eq:Elliptic}):
find $y\in V_y:=H^1_0(\Omega)$ such that
\begin{equation}
 \int_\Omega \mu_1\nabla y\ \nabla \hat y\, d\Omega=\int_\Omega \mu_2 \hat y\, d\Omega\quad \forall \hat y\in V_y.
\end{equation}
This solution exists and is unique (\cite{Quarteroni-Valli-1994}).

The goal now is to approximate the operator 
\begin{equation}
 \mathcal{G}:V_u\to V_y:\quad u\mapsto y=\mathcal{G}(u)
\end{equation}
that maps the infinite--dimensional \emph{input space} $V_u$ to the infinite--dimensional \emph{solution space} $V_y$, and not only the values of the solution $y$ at some given points. The operator learning theory aims at generalizing feed--forward neural networks, which typically operate on finite--dimensional spaces, to infinite--dimensional spaces. 

These new NNs are particularly useful in contexts where one wants to solve a PDE many times for different types of boundary conditions or different discretizations, avoiding solving the PDE numerically from scratch or training a new FFNN each time.
Indeed, standard NN architectures strongly depend on the discretization and the type of boundary conditions used to produce the training data.

%
%
\subsubsection{Deep Operator Networks (DeepONet)} \label{sec:deeponet}

An early paper on operator learning is that of Chen and Chen in 1995 \cite{chen1995}. Therein, the authors proved that, under regularity assumptions on both parameters and solution, for any continuous operator $\mathcal G$ (either linear or non--linear) and any required precision $\varepsilon$ we can find a neural network $\mathcal G_{\mathbf w}$ approximating $\mathcal G$ and depending on the parameters $\mathbf w$. The authors also provide the NN architecture of $\mathcal G_{\mathbf w}$. 

\begin{theorem}[Universal approximation theorem for non--linear operators]\cite{chen1995}\label{thm:chen1995}
Let $X$ be a Banach space, $K_1\subseteq X$ compact, $K_2\subseteq \mathbb{R}^d$ compact, $V_u\subset\mathcal{C}^0(K_1)$ compact, 
 and $\mathcal{G}:V_u\to \mathcal{C}^0(K_2)$ a continuous operator (either linear or non--linear). Let $\sigma$ be a sigmoid function (i.e. a bounded function $\sigma:\mathbb R \to \mathbb R$ such that $\lim_{x\to -\infty}\sigma(x)=0$ and $\lim_{x\to+\infty}\sigma(x)=1$). Then for any $\varepsilon>0$, there are positive integers $M,\ N,\ m$, real numbers $c_i^k,\ \zeta_k,\ \xi_{ij}^k$, and points $\boldsymbol\omega_k\in\mathbb R^d$, $\mathbf x_j\in K_1$, $i=1,\ldots,M$, $k=1,\ldots, N$, $j=1,\ldots, m$, such that 
\begin{equation}\label{eq:chen_thesis}
\left|\mathcal G(u)(\mathbf z)-\underbrace{\sum_{k=1}^N\left\{ \left[\sum_{i=1}^M c_i^k \sigma 
\left(\sum_{j=1}^m \xi_{ij}^k u(\mathbf x_j)+\theta_i^k\right)\right]
\sigma(\boldsymbol \omega_k\cdot \mathbf z+\zeta_k)\right\}}_{\mathcal G_{\mathbf w}(u)(\mathbf z)}\right|< \varepsilon
\end{equation}
holds for all $u\in V_u$ and $\mathbf z\in K_2$.
\end{theorem}

\begin{figure}
    \hskip 2.cm   
\scalebox{0.6}{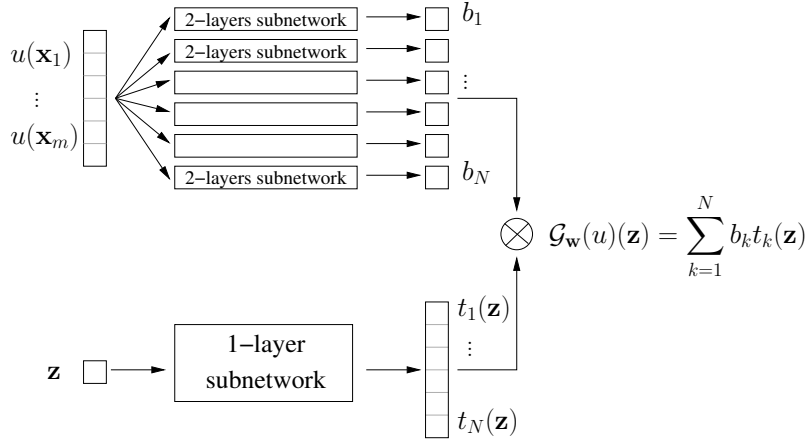}
    \caption{The NN by Chen and Chen approximating any (non)linear continuous operator \cite{chen1995}}
    \label{fig:chen_no}
\end{figure}

Although the input of the operator consists of a function, to feed the neural network we have to evaluate it at a finite set of points that are named \emph{sensors}, namely the $m$ points $\mathbf x_j\in K_1$. 
The point $\mathbf z\in \mathbb R^d$ instead, used to evaluate the solution $y=\mathcal G(u)\in \mathcal C^0(K_2)$, is not fixed, and can be varied after the training stage.
The array parameter $\mathbf w$ contains the real numbers $c_i^k,\ \zeta_k,\ \xi_{ij}^k$.

The neural network is composed of two subnetworks, exploiting the same activation function $\sigma$. The first subnetwork (the term between square brackets in (\ref{eq:chen_thesis})) is a stack of $N$ two--layers networks acting on the values of the parameter functions at the sensors, the second one ($\sigma(\boldsymbol\omega_k\cdot \mathbf z+\zeta_k)$) is a one--layer network acting on the variable $\mathbf z$ at which we want to evaluate the solution function. The results of the two subnetworks are combined to provide the final output which is the value of $\mathcal G_{\mathbf w}(u)(\mathbf z)\simeq y(\mathbf z)$.

The universal approximation theorem \ref{thm:chen1995} guarantees a small approximation error for a sufficiently rich network requiring a large number of sensors. However, as we have seen in Sect. \ref{sec:supervised_learning}, the approximation error alone does not guarantee the efficacy of a NN, since also optimization and generalization errors contribute to making a NN effective.

In \cite{Lu2021}, the authors propose \emph{Deep Operator Networks (DeepONets)} which improve the NN proposed by Chen and Chen in \cite{chen1995}. More precisely, the two shallow subnetworks of \cite{chen1995} are replaced in DeepONet by deep neural subnetworks, so that a limited number of sensors is required to ensure low generalization and optimization errors and low computational costs. 

\begin{figure}
    \begin{center}    
\scalebox{0.5}{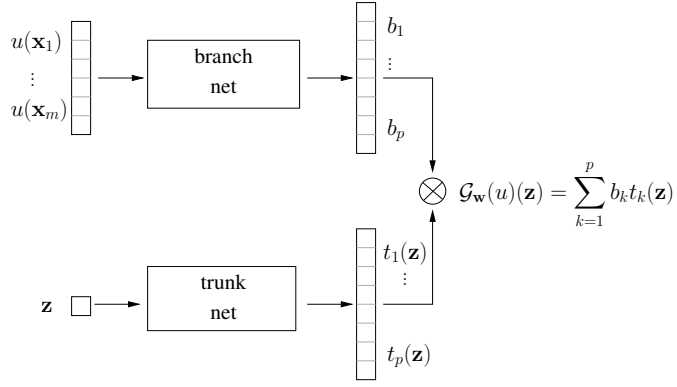}
    \end{center}
    \caption{A Deep Operator Network (DeepONet)}
    \label{fig:deeponet}
\end{figure}

The subnetwork acting on $\mathbf z$ is named \emph{trunk net} and predicts the values of $p$ shape functions $\{t_1, t_2,\ldots, t_p\}$ at the point $\mathbf z$, while the subnetwork acting on the values $[u(\mathbf x_1),\ldots, u(\mathbf x_m)]$ is named \emph{branch net} (it is a unique unstacked subnetwork) and predicts the values of the problem--specific
coefficients $b_k$ such that $y(\mathbf z)\simeq\mathcal G_{\mathbf w}(u)(\mathbf z)=\sum_{k=1}^p b_k t_k(\mathbf z)$, see Fig. \ref{fig:deeponet}. 

The operator $\mathcal G_{\mathbf w}$ does not provide the analytic expression of an approximation of the solution $y$ but approximates the value of $y$ at given points $\mathbf z$, leveraging the shape functions $t_k(\mathbf z)$ that replace high--fidelity (e.g. FEM) shape functions.  

The following \emph{generalized universal approximation theorem for operators}, that is an immediate consequence of Theorem \ref{thm:chen1995}, has been stated in \cite{Lu2021}.
\begin{theorem}
Suppose that $X$ is a Banach space, $K_1\subset X$, $K_2\subset \mathbb R^d$ are two compact sets in $X$ and $\mathbb R^d$, respectively, $V_u$ is a compact set in $\mathcal C^0(K_1)$. Assume that $\mathcal G:V_u\to \mathcal C^0(K_2)$ is a non--linear continuous operator. Then, for any $\varepsilon>0$, there exist positive integers $m,\ p$, continuous vector functions $\mathbf g:\mathbb R^m\to \mathbb R^p$, $\mathbf f:\mathbb R^d \to \mathbb R^p$, and $\mathbf x_1,\, \mathbf x_2,\ldots, \mathbf x_m \in K_1$ such that
 \begin{equation*}
     \left| \mathcal G(u)(\mathbf z)-\langle\underbrace{\mathbf g(u(\mathbf x_1),u(\mathbf x_2),\ldots,u(\mathbf x_m))}_{branch}, 
     \underbrace{\mathbf f(\mathbf z)}_{trunk}\rangle\right|<\varepsilon
 \end{equation*}
holds for all $u\in V_u$ and $\mathbf z\in K_2$ , where $\langle\cdot,\cdot\rangle$ denotes the dot product in $\mathbb R^p$. Furthermore, the functions $\mathbf g$ and $\mathbf f$ can be chosen as diverse classes of neural networks, which satisfy the classical universal approximation theorem of functions, for example, (stacked/unstacked) fully connected neural networks, residual neural networks and convolutional neural networks.
\end{theorem}

As stated in the previous theorem, DeepONet is a high--level network architecture, and the architectures of its inner trunk and branch networks are not defined a--priori. 
Finally, we notice that the branch net constrains the input to a fixed location (otherwise all the parameters will have to be recalculated), however, it is possible to make DeepONet discretization invariant \cite{deHoop2022}.

%
%
\subsubsection{Neural Operators}\label{sec:Neural-Operators}

A FFNN transforms a vector $\mathbf{x} \in \mathbb{R}^n$ into a vector $\mathbf{y} \in \mathbb{R}^m$. This process operates sequentially, by transforming at each step a vector into another vector (namely the vector $\mathbf{a}^{[\ell-1]}$ into $\mathbf{a}^{[\ell]}$, see Algorithm~\ref{alg:FFNN}). Neural Operators generalize this vector--to--vector transformation process to a function--to--function transformation \cite{Kovachki2022}.
The key observation underlying this generalization is that a vector $\mathbf{a}^{[\ell]} \in \mathbb{R}^{N_\ell}$ can be interpreted as the values of a function from the discrete domain $\Omega_\ell = \{1,2,\dots,N_\ell\}$ into $\mathbb{R}$. 
This suggests defining a Neural Operator as a generalized FFNN, where we set instead $\Omega_\ell = \mathbb{R}^d$ in each layer (i.e. for any $\ell$). Thus, the Neural Operator will map an input function $u \colon \mathbb{R}^d \to \mathbb{R}^{n_u}$ into an output function $v \colon \mathbb{R}^d \to \mathbb{R}^{n_v}$, by sequentially transforming it through a series of intermediate functions $a^{[\ell]}\colon \mathbb{R}^{d} \to \mathbb{R}^{M_\ell}$. This enables neural operators to naturally handle problems with inputs and outputs defined on continuous domains, such as those encountered in physics-informed modelling or solutions to partial differential equations. The sequential mapping of functions at each layer mirrors the structure of FFNNs but extends their applicability to infinite-dimensional spaces.

To achieve this generalization, the standard FFNN update equation,
$\mathbf{a}^{[\ell]}=\sigma(W^{[\ell]}\mathbf{a}^{[\ell-1]}+\mathbf{b}^{[\ell]})$, is replaced with its infinite-dimensional counterpart. Specifically, the matrix-vector multiplication $W^{[\ell]}\mathbf{a}^{[\ell-1]}$ is replaced by a kernel convolution $\int_\Omega\kappa(x,\hat x) a^{[\ell-1]}(\hat x) d\hat x$. 
Notably, in the discrete case ($\Omega_\ell = \{1,2,\dots,N_\ell\}$), the kernel integral reduces to a matrix multiplication, ensuring consistency between discrete and continuous formulations.
In practice, using the same idea from the ResNet architecture \cite{He-resnet-2016}, an additional term is summed before applying the non--linear activation function, yielding the expression:
\begin{equation}\label{eq:layer-neural-operator}
a^{[\ell]}(x)=\sigma_\ell\left(\int_\Omega\kappa(x,\hat x) a^{[\ell-1]}(\hat x) d\hat x + b^{[\ell]}(x) + W^{[\ell]} a^{[\ell-1]}(x)\right),  
\end{equation}
where $W^{[\ell]} \in \mathbb{R}^{M_{\ell} \times M_{\ell - 1}}$ represents a matrix.
Finally, a lifting operator ($\mathcal P\colon \mathbb{R}^{n_u} \to \mathbb{R}^{M_0}$) is applied to the input to map it into a higher--dimensional space suitable for the model. Similarly, a projection operator ($\mathcal Q\colon \mathbb{R}^{M_L} \to \mathbb{R}^{n_v}$) is applied to the output of the last layer to map it back to the original space. Both operators act pointwise.

The complete Neural Operator $\mathcal G_{\mathbf w}$ is therefore expressed as:
\begin{equation}\label{eq:Gw-neuraloperator}
\mathcal G_{\mathbf w}=\mathcal Q \circ \sigma_L \circ T^{[L]} \circ \ldots \circ  \sigma_\ell \circ T^{[\ell]} \circ \ldots
\circ\sigma_1\circ T^{[1]}\circ \mathcal P
\end{equation}
where, for any $\ell=1,\ldots, L$, $b^{[\ell]}$ is the bias function and $\sigma_\ell$ the activation function of the layer $\ell$. The linear transformation $T^{[\ell]}(a^{[\ell-1]})$ is defined as $T^{[\ell]}(a^{[\ell-1]})=W^{[\ell]}a^{[\ell-1]}+\mathcal K^{[\ell}(a^{[\ell-1]})+b^{[\ell]}$, where $\mathcal K^{[\ell]}$ is a kernel operator associated with $\kappa(x,\hat x)$.

The most notable properties of these networks are that they are discretization--invariant and guarantee universal approximation. 

Discretization invariance means that neural operators keep the same model parameters even if the discretization of the underlying functional spaces is varied. More precisely, a \emph{discretization--invariant model} with a fixed number of parameters satisfies the following rules \cite{Kovachki2022}: \emph{(i)} acts on any discretization of the input function, i.e. accepts any set of points in the input domain;
\emph{(ii)} can be evaluated at any point of the output domain; \emph{(iii)} converges to an infinite--dimensional operator as the discretization is refined.
The last property guarantees consistency in the limit with the continuous PDE as the discretization is refined.  

Universal approximation guarantees that any continuous operator can be approximated by a neural operator up to a given tolerance. This property is ensured even though the internal layers only perform linear (actually, affine) operations, thanks to the presence of non--linear activation functions.

The parameter vector $\mathbf w$ is the concatenation of the parameters of $\mathcal P$, $\mathcal Q$, $W^{[\ell]}$, $\mathcal K^{[\ell]}$, and $b^{[\ell]}$ for any $\ell$ and it is computed by standard gradient--based minimization algorithms. We notice that both inputs and all the functions entering the loss function need to be
discretized, but, in view of the discretization invariance, the learned parameters $\mathbf w$ may be used with other discretization.

The kernel $\kappa(x,\hat x)$ takes care of the non--local properties of the solution that the differential equation, by its definition, would not be able to describe. We refer to \cite{Kovachki2022} for a detailed description of neural operators.

The integrals in (\ref{eq:layer-neural-operator}) need to be approximated, and different strategies can be adopted, from Monte Carlo method, to averaging aggregation, to more sophisticated quadrature formulas.
In \cite{Kovachki2022}, four classes of efficient
parametrizations are described: graph neural operators, multi--pole graph neural operators, low--rank neural
operators, and Fourier neural operators. 

\begin{figure}
    \begin{center}    
\scalebox{0.65}{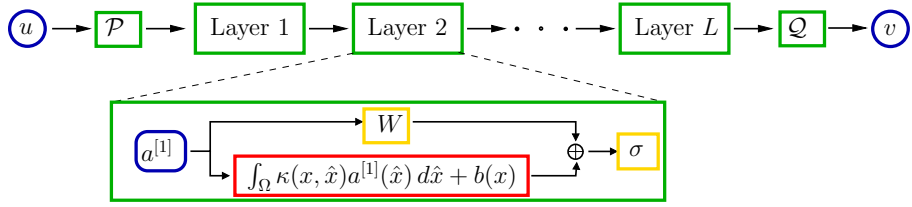}
    \end{center}
    \caption{A Neural Operator}
    \label{fig:neural-operator}
\end{figure}

%
%
\null\textbf{Fourier Neural Operators (FNO).}
FNO were proposed in \cite{Li2021}. Instead of working in the physical domain, here the kernel is parametrized in the frequency space, and a Fourier transform and its inverse are applied on the right and the left of the kernel, respectively.

Each layer, named \emph{Fourier layer} takes the form
\begin{eqnarray}\label{eq:Fourier-layer}
a^{[\ell]}(x)=\sigma_\ell\left(W^{[\ell]} a^{[\ell-1]}(x)+(\mathcal F^{-1}(K^{[\ell]}\mathcal F(a^{[\ell-1]}))(x) + b^{[\ell]}(x)\right), 
\end{eqnarray}
where $K^{[\ell]}$ is a tensor (representing the kernel in the frequency space) which is learnable like the linear operators $W^{[\ell]}$ and the bias function $b^{[\ell]}$, while $\mathcal F$ denotes the Fourier transform along the spatial dimension $\mathbf x$ and $\mathcal F^{-1}$ its inverse. 

Integrals in the Fourier transform are computationally expensive, thus a Fourier Fast Transform FFT can be employed, provided that the kernels are invariant to translations.
A drawback of FNO is that the number of parameters grows very fast unless the Fourier coefficients of the output decay very fast, but fortunately, this is typically the case in PDEs.

FNO enjoys the universal approximation property \cite{Kovachki2021-FNO} and, are discretization invariant.  
FNO is significantly more capable than convolutional neural networks of learning low--frequencies and demonstrates larger errors in higher frequencies than in lower frequencies \cite{qin2024}. To improve the capability of FNO in capturing even high frequencies, in \cite{qin2024} a second FNO is trained to predict the residual of the first FNO’s predictions.
\begin{figure}
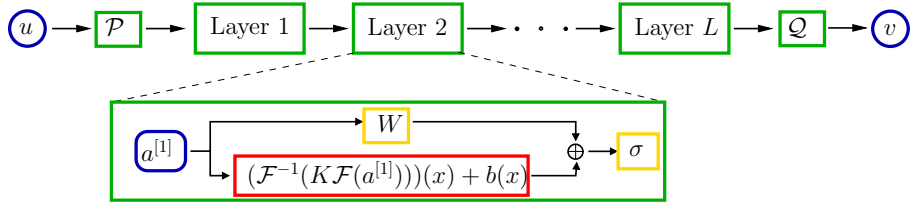

    \begin{center}    
\scalebox{0.65}{}
    \end{center}
    \caption{A Fourier Neural Operator}
    \label{fig:FNO}
\end{figure}

\subsubsection{Operator learning for time--dependent problems} \label{sec:operator-learning-time-dependent}

Operator learning methods apply to operators between functional spaces. Often, the independent variable is the spatial coordinate, as in stationary PDEs. However, in many cases, the independent variable is time, as in evolutionary PDEs (e.g. parabolic or hyperbolic). 
A straightforward approach to {handling} time-dependent problems with operator learning is treating the time variable as a generic independent variable, akin to the $(d+1)$--th spatial coordinate. In this case, DeepONets or Neural Operators can be used, although they do not account for the sequential nature of time. Other methods, often based on recurrent (RNN), and autoregressive architectures (see Sect. \ref{sec:architectures}), have been developed ad--hoc to incorporate temporal sequentiality. These methods are the main focus of this section.

Rather than listing every single method as done thus far, in Secs. \ref{sec:operator-learning-time-dependent}--\ref{sec:operator-learning_hidden-dynamics-discovery}--\ref{sec:space-time learning} we will instead focus on addressing different approaches to specific tasks. For each approach, we will discuss advantages and disadvantages and refer to the most significant contributions in the literature. For the various tasks, we will highlight some of the many possible methods that leverage diverse combinations of fundamental components such as the time advancement scheme, {the} loss function, and {the} neural network architecture. 

\null\textbf{Neural ODEs to learn dynamical systems.}
Let us consider a generic time--dependent model, in the form:
\begin{equation} \label{eqn:dynODEbase}
    \left\{
    \begin{aligned}
        & \frac{d}{dt} \dynOut(t) = \dynRHS(\dynOut(t), \dynInp(t)) & & t \in (0, T) \\
        & \dynOut(0) = \dynOut_0, & &
    \end{aligned}
    \right.
\end{equation}
where $\dynOut \in \mathbb{R}^{\dynNumOut}$ is the vector of state variables, which in this case coincides with the output of interest, while $\dynInp \in \mathbb{R}^{\dynNumInp}$ contains the input variables. {The latter are either} time--dependent or constant (in which case, they are typically called \textit{parameters}). The {(generally non--linear)} term $\dynRHS$ governs the temporal evolution of the system. Suppose we {know} a set of temporal trajectories $\widehat\dynOut^j(t)$, which we call \textit{samples}, {each one corresponding} to a given input $\widehat\dynInp^j(t)$ and an initial condition $\widehat\dynOut_0^j$, where $j = 1, \dots, \dynNumSamples$ represents the sample index. 
Operator learning for {this kind of} time--dependent problems consists of learning, from the training samples, an operator that maps a given initial condition and an input trajectory to an output trajectory.

In practice, we suppose that both the input and the output are sampled at time intervals ($0 = t_0, t_1, \dots, t_{\dynNumTimes} = T$). Using the notation $\widehat\dynInp^j_i = \widehat\dynInp^j(t_i)$ and $\widehat\dynOut^j_i = \widehat\dynOut^j(t_i)$, to denote {input and output} training trajectories, our training dataset consists of $\{(\widehat\dynInp^j_i , \widehat\dynOut^j_i)\}_{i = 0, \dots, \dynNumTimes}^{j = 0, \dots, \dynNumSamples}$.

A successful approach to tackle this operator learning problem is to use neural networks of recurrent nature (see Sec.~\ref{sec:architectures}). Here we consider the case of \textit{Neural ODEs} \cite{chen2018neural}, which {approximate} the dynamics through the following system of ODEs:
\begin{equation} \label{eqn:dynODE_neuralODE}
    \left\{
    \begin{aligned}
        & \frac{d}{dt} \dynOut(t) = \dynNNdyn(\dynOut(t), \dynInp(t); \dynWdyn) & & t \in (0, T) \\
        & \dynOut(0) = \dynOut_0 & &
    \end{aligned}
    \right.
\end{equation}
where $\dynNNdyn$ is a neural network (typically a FFNN), which receives as input the vector of state variables $\dynOut(t)$ concatenated with the vector of input variables $\dynInp(t)$, and returns an approximation of the vector of temporal derivatives $\frac{d}{dt} \dynOut(t)$. The trainable parameters of the neural network are collected in $\dynWdyn$.

The system of ODEs \eqref{eqn:dynODE_neuralODE} is approximated by means of suitable numerical methods. A common choice is the Forward Euler method, due to its computational convenience, which results in the following scheme:
\begin{equation} \label{eqn:dynODE_neuralODE_discrete}
    \left\{
    \begin{aligned}
        & \dynOut^{(k+1)} = \dynOut^{(k)} + \Delta t \dynNNdyn(\dynOut^{(k)}, \dynInp(k \Delta t); \dynWdyn) & & k = 0, \dots, T/\Delta t - 1 \\
        & \dynOut^{(0)} = \dynOut_0. & &
    \end{aligned}
    \right.
\end{equation}
Note that the integration time step $\Delta t$ does not need to match the sampling time step  (i.e. $\Delta t^{\text{sampling}}_j = t_{j+1} - t_{j}$). Should the input $\dynInp$ be unavailable at specific time points $k \Delta t$, interpolation of the input values can be performed during both training and testing phases.
For applications requiring high numerical accuracy, more advanced methods such as Runge--Kutta schemes are often preferred \cite{chen2018neural}.

\null\textbf{Derivative--based vs trajectory--based training.}
To obtain the parameters $\dynWdyn$, a supervised learning strategy is used, given that we have a labelled dataset (i.e., the data $\widehat\dynInp^j_i$ and $\widehat\dynOut^j_i$) There are mainly two possible approaches for training.
Note that the two approaches differ solely on the way the training is carried out.

The first approach called \textit{derivative--based training} minimizes the difference between the temporal derivative {computed} by the neural network and the temporal derivative {approximated} directly from {training} data. {The latter} can be obtained through finite differences, e.g. defining 
\begin{equation*}
    \delta{\dynOut}^j_i := \frac{\widehat\dynOut^j_{i+1} - \widehat\dynOut^j_i}{t_{i+1} - t_i},
\end{equation*}
or by using higher--order formulas.
The minimization problem then becomes:
\begin{equation}
    \dynWdyn^*=\argmin{\dynWdyn}
    \sum_{j = 1}^{\dynNumSamples} \sum_{i = 0}^{\dynNumTimes - 1} \left\| \delta{\dynOut}^j_i - \dynNNdyn(\widehat\dynOut^j_i, \widehat\dynInp^j_i; \dynWdyn) \right\|^2.
\end{equation}
The second approach called \textit{trajectory--based training} {instead minimizes} the difference between the temporal trajectory {${\dynOut}^j(t)$} approximated by the Neural ODE through numerical integration (like, e.g., in (\ref{eqn:dynODE_neuralODE_discrete}) and the temporal trajectory available from training data. The minimization problem then becomes:
\begin{equation}\label{eq:loss2-neuralODE}
    \dynWdyn^*=\argmin{\dynWdyn}
    \sum_{j = 1}^{\dynNumSamples} \sum_{i = 0}^{\dynNumTimes} \left\| \widehat\dynOut^j_i - {\dynOut}^j(t_i) \right\|^2.
\end{equation}
When the sampling time step $\Delta t_j^{\text sampling}$ is not a multiple of the integration time step $\Delta t$, the values ${\dynOut}^j(t_i)$ are computed by interpolating those provided by (\ref{eqn:dynODE_neuralODE_discrete}).
Note that the loss function {in (\ref{eq:loss2-neuralODE}) implicitly depends} on $\dynWdyn$, as ${\dynOut}^j_i$ {itself} depends on $\dynWdyn$ through the numerical solution of \eqref{eqn:dynODE_neuralODE}. This approach is also known as \textit{end--to--end training}, as the numerical approximation of the ODE system is incorporated into the training process.

Trajectory--based training is computationally more challenging than derivative--based training, as it requires the numerical solution of an ODE system for each sample. 
Additionally, this process is incorporated into the calculation of the loss function, and therefore the automatic differentiation engine must be able to compute the gradient of the loss function with respect to the parameters $\dynWdyn$ across the numerical solution of an ODE system. This process requires particular attention, especially if implicit advancement schemes, which require the solution of non--linear equation systems, are used. For this reason, explicit solvers, such as Forward Euler or Runge-Kutta, are often preferred in this context. Note, however, that modern automatic differentiation tools, contained in the most popular ML libraries, can handle these complexities transparently, although this potentially leads to deep computational graph{s}, which can require a lot of memory in the case of many data. Alternatively to black--box automatic differentiation, the adjoint method proposed in \cite{chen2018neural} significantly reduces the memory footprint of the gradient calculation, as it avoids the need to store the entire computational graph, but reconstructs the gradient through the solution of an adjoint ODE system.

On the other hand, derivative--based training is characterized by a very shallow computational graph, as each term that {contributes to} the loss function does not depend on the other ones. This makes the gradient calculation much more efficient, {although less accurate,} especially if the data are noisy. In such cases, the finite difference operation can amplify the noise present in the data, and therefore the estimate of the temporal derivative can be very unstable. Furthermore, at testing time, when the model is used to make predictions by numerically integrating the Neural ODE, such errors can accumulate unpredictably, leading to solutions that diverge significantly from the expected trajectories. The trajectory--based approach performs the integration of the ODE at training time, hence it allows controlling error accumulation and is therefore more robust to the presence of noise in the data. To mitigate the effect of noise within the derivative--based approach, more robust temporal derivative reconstruction techniques, such as total--variation regularization \cite{chartrand2011numerical}, can be used.
However, in the presence of noisy data, and especially for long--term predictions, trajectory--based training is typically the more robust choice.

\null\textbf{Comparison with methods that do not consider time--dependence.}
Now that we have introduced operator learning for time--dependent problems using autoregressive models such as Neural ODEs, it is useful to compare them with techniques that do not account for the temporal nature of the problem. Indeed, the operator that maps the time--dependent function $\dynInp(t)$ to the time--dependent function $\dynOut(t)$ can also be learned using other operator learning methods, such as DeepONets and Neural Operators. However, as anticipated above, these methods do not account for the sequential nature of time and, therefore, cannot capture the temporal dependence of the variables. In particular, these methods do not guarantee consistency with the \textit{arrow of time}, meaning that the temporal evolution of a system is unidirectional.
This implies that, although these methods can approximate the dynamics of a system, they cannot guarantee that the prediction of $\dynOut(t)$ at a given time $t \in (0,T)$ solely depends on the input $\dynInp(t')$ for $t' \leq t$. This property is {instead} guaranteed by Neural ODEs, which can capture the temporal dependence of the variables and, therefore, ensure consistency with the arrow of time.

Another advantage of using autoregressive models, compared to general--purpose operator learning methods, is {their} ability to handle time series of variable length. In particular, they can perform time extrapolation, i.e., predict the temporal evolution of the system even outside the time range in which the data were sampled. This is possible due to the autoregressive nature of the model, which allows integrating the temporal dynamics of the system from an initial condition for arbitrary times. Clearly, this capability is limited by the quality of the training data and the type of dynamics under consideration. Training data {ought be} representative of the system's dynamics over long periods. Indeed, due to the data-driven nature of the approach used, the predictions are reliable as long as they do not deviate too far from the regions in the state space covered by the training data. For example, let us suppose we are learning in a data--driven manner the dynamics of a robotic arm, as a function of external stimuli and of the power applied to the joints. The learned data-driven model might be able to predict the temporal evolution of the arm even for longer times than those used to train the model. However, if the long--term dynamics cause the arm to move to positions far from those covered by the training data, the model likely provide{s} inaccurate predictions. This is an intrinsic limitation of all data--driven methods, and therefore it is important that the training data are representative of the system's dynamics for which predictions are to be made.

It should be noted that the above observations are not limited to Neural ODEs but also apply to other neural network architectures such as RNNs, LSTM, and GRU (see Sec.~\ref{sec:architectures}) that account for the temporal nature of the problem {and} can be used as alternatives to Neural ODEs. One advantage of Neural ODEs, however, is their flexibility when time sampling is non-uniform, whereas other autoregressive methods may require pre--processing of data to uniform temporal sampling.

\null\textbf{Sparse Identification of non--linear Dynamics (SINDy).}
Sparse Identification of non--linear Dynamics (SINDy) is a method aimed at data--driven identification of a dynamical system, with the goal of obtaining a parsimonious model that captures the fundamental phenomena of the system \cite{brunton2016discovering}. Like Neural ODEs, SINDy can capture the temporal dependence of variables, thus ensuring consistency with the arrow of time. Unlike Neural ODEs, however, the right--hand side of the ODE system is not approximated with a neural network but through a linear combination of a set of candidate functions, which can be polynomials, trigonometric functions, exponentials, etc. The learned model is therefore interpretable, as the identified terms retain physical meaning (e.g., dissipation, reaction, forcing terms). The method is designed so that the resulting combination is \emph{parsimonious}, meaning that only a few significant terms are present.

SINDy requires a library of candidate functions, representing the possible terms that may appear in the right--hand side of ODE system. A vector $\boldsymbol{\theta}(\dynOut,\dynInp)$ is defined, containing $N_\theta$ candidate functions (polynomials, trigonometric functions, exponentials, etc.) in terms of $\dynOut$ and $\dynInp$. For example, if $\dynOut=(y_1,y_2)$ and $\dynInp=(u_1)$, a library of functions can be:
\begin{equation*}
        \boldsymbol{\theta}(\dynOut, \dynInp) = \left[\, 1,\; y_1,\; y_2,\; u_1,\; y_1^2,\; y_1y_2,\; \sin(y_2),\; e^{-y_1}u_1,\; \dots\right].
\end{equation*}
Multiplicity and type of candidates determine both expressive capacity and risk of overfitting.

Next, a model is sought in the form of
\begin{equation} \label{eqn:dynODE_SINDy}
    \left\{
    \begin{aligned}
        & \frac{d}{dt} \dynOut(t) = \mathbf{W} \boldsymbol{\theta}(\dynOut(t) , \dynInp(t) ) & & t \in (0, T) \\
        & \dynOut(0) = \dynOut_0 & &
    \end{aligned}
    \right.
\end{equation}
where $\mathbf{W} \in \mathbb{R}^{\dynNumOut \times N_\theta}$ is a matrix of unknown coefficients that defines the active terms of $\boldsymbol{\theta}$ in each row of the ODE system. The problem of identifying the dynamical system thus becomes a generalized linear regression problem. Specifically, sparse regression techniques are used to promote solutions where many elements of $\mathbf{W}$ are zero. This {yields} a parsimonious and interpretable model, where only a few terms are significant, drawing inspiration from the field of compressed sensing \cite{baraniuk2007compressive,baraniuk2010model}. An approach that exploits the derivative--based loss consists of solving the following optimization problem:
\begin{equation}
    \mathbf{W}^* = \argmin{\mathbf{W}} 
    \sum_{j = 1}^{\dynNumSamples} \sum_{i = 0}^{\dynNumTimes - 1} \left\| \delta{\dynOut}^j_i - \mathbf{W} \boldsymbol{\theta}(\widehat\dynOut^j_i, \widehat\dynInp^j_i) \right\|_2^2 + \lambda \left\| \mathbf{W} \right\|_0,    
\end{equation}
for a suitable hyperparameter $\lambda$, where $\| \cdot \|_0$ is the norm that counts the number of non-zero elements in a matrix. Since the problem would require an intractable combinatorial search, $\| \cdot \|_1$ norm regularization or a sequential thresholding procedure is used in place of the zero norm. \cite{brunton2016discovering,Brunton-Kutz-2022}.

Once the coefficients $\mathbf{W}^*$ are identified, a parsimonious model is obtained that uses only the columns of $\boldsymbol{\theta}$ corresponding to the non-zero coefficients. To evaluate the quality of this model, the predicted dynamics are compared with those measured from the test dataset. It is essential to verify the model's robustness to noise and small variations in the data, as well as its extrapolation capability on trajectories not included in the training set.

The SINDy method is particularly effective when the true dynamics are sparse, meaning that only a few terms in the governing equations are significant. By leveraging sparsity, SINDy can identify interpretable models that provide insights into the underlying physical processes.
Although lacking any NN component, this method is however driven by data.
This is why we deem it appropriate to present it under the SciML heading.

\subsubsection{Intrinsic or hidden dynamics discovery}
\label{sec:operator-learning_hidden-dynamics-discovery}

The methods described in Sec.~\ref{sec:operator-learning-time-dependent} are aimed at learning the dynamics of a system from a dataset of trajectories of the system state. However, very often this state has a huge dimensionality, {so that} instead of learning the dynamics in the state space, it is preferred to learn an intrinsic or essential dynamics that takes place in a low--dimensional space. In other cases, it is needed to discover hidden dynamics that are not directly observable in the data. We will start by considering the first case and then move on to the second.

When the system state has a very high dimensionality (typically when $\dynNumOut$ is of the order of $10^2$ or higher), instead of learning the dynamics in the state space, one {could} reduce the dimensionality of the variables at hand and find a low--dimensional representation that {nonetheless} captures the intrinsic dynamics of the system. This can {unveil} relationships between state variables and identify the most relevant variables, finding an essential and therefore more interpretable representation of the system, ultimately improving {our understanding of the system dynamics. An extra benefit is the reduction of the computational cost associated with the numerical approximation of a high--dimensional system, especially in contexts where real--time predictions are necessary. Notable instances are system control and creation of digital twins \cite{hesthaven2022reduced,peirlinck2021precision}. A paradigmatic example is when the dynamical system under consideration arises from semi-discretization in space of a mathematical model based on PDEs, as described in Sec.~\ref{sec:numerical_models} (see e.g. \eqref{eq:IBVP-semidiscrete}), where the state variable is a high--dimensional vector representing the discretized solution of the PDE.

\null\textbf{Projection based methods.} 
Projection--based ROMs (Reduced Order Models) consist of projecting the system{'s} dynamics into a low--dimensional subspace \cite{benner2015survey,hesthaven2022reduced}. The latter is defined as a linear subspace of dimension $\dynNumLat \ll \dynNumOut$ of the state space of the full--order model (FOM).
Let $W \in \mathbb{R}^{\dynNumOut \times \dynNumLat}$ be a matrix whose columns contain an orthonormal basis of the reduced space. The state variable of the ROM, denoted by $\dynLat\in \mathbb{R}^{\dynNumLat}$, is related to the state variable of the FOM $\dynOut\in \mathbb{R}^{\dynNumOut}$ through the relations $\dynLat = W^T \dynOut$ and $\dynOut \simeq W \dynLat$. By projecting the FOM \eqref{eqn:dynODEbase} into the reduced subspace, we obtain the {corresponding} ROM:
\begin{equation} \label{eqn:dynODE_ROM}
    \left\{
    \begin{aligned}
        & \frac{d}{dt} \dynLat(t) = W^T \dynRHS(W \dynLat(t), \dynInp(t)) & & t \in (0, T) \\
        & \dynLat(0) = W^T \dynOut_0. & &
    \end{aligned}
    \right.
\end{equation}
An example in the case where FOM arises from the semi-discretization in space of a differential model is provided in \eqref{eq:IBVP-ROM1-semidiscrete}.

As anticipated in Sec.~\ref{sec:numerical_models}, to define the matrix $W$, the most common approaches are greedy algorithms and POD (Proper Orthogonal Decomposition) \cite{Hesthaven-2016,Quarteroni-Manzoni-Negri}. The latter consists of collecting the snapshots of the training dataset into a matrix
\begin{equation*}
    Y = \left[
        \widehat\dynOut_1^1 |
        \widehat\dynOut_2^1 |
        \dots |
        \widehat\dynOut_{\dynNumTimes}^1 |
        \widehat\dynOut_1^2 |
        \widehat\dynOut_2^2 |
        \dots |
        \widehat\dynOut_{\dynNumTimes}^2 |
        \dots |
        \widehat\dynOut_1^{\dynNumSamples} |
        \widehat\dynOut_2^{\dynNumSamples} |
        \dots |
        \widehat\dynOut_{\dynNumTimes}^{\dynNumSamples} \right],
\end{equation*}
where $N_S$ represents the number of training samples, while $N_T$ is the number of time--steps,
and then defining $W$ as the matrix whose columns are the first $\dynNumLat$ eigenvectors associated with the $\dynNumLat$ largest eigenvalues of the covariance matrix $YY^T$. From a computational point of view, {these eigenvectors} can be obtained in an efficient and numerically stable manner through the Singular Value Decomposition (SVD) of the matrix $Y$ {(see, e.g., \cite[Ch. 4]{Quarteroni-CS-2014})}.

Projection--based approaches perform well {when} the system{'s} dynamics {is} dominated by a few modes that can be captured in a low--dimensional space. However, many systems are not suitable to be compactly described through a linear projection.
A useful concept to quantify the {suitability} of a linear projection to approximate the solution manifold is the \textit{Kolmogorov n-width} {\cite{kolmogorov}}.
Let $M \subset \mathbb{R}^{\dynNumOut}$ be the solution manifold associated with the dynamical system.
In this case, the Kolmogorov $n-$width is defined as}
\begin{equation}
    d_n(M) = 
    \inf_{\substack{W \in \mathbb{R}^{\dynNumOut \times n}}} 
    \sup_{\dynOut \in M} \inf_{\dynLat\in \mathbb{R}^{n}} \|\dynOut - W \dynLat\|.
\end{equation}
In words, the Kolmogorov $n-$width is the minimum distance of the solution manifold from linear subspaces of dimension $n$. Intuitively, the Kolmogorov $n-$width measures the best approximation error that can be obtained by approximating the solution manifold with a linear subspace of dimension $n$.

In many cases, the Kolmogorov $n-$width decreases slowly with $n$, indicating that a very high--dimensional subspace is needed to obtain a good approximation, making projection--based ROMs ineffective. { This happens for} systems characterized by complex and non--linear dynamics, turbulent flows, and systems with multiscale dynamics, where linear dimensionality reduction may not be sufficient to capture all the essential features of the system. Moreover, even linear equations can present a solution manifold with a slowly decreasing Kolmogorov $n$--width. An example is given by travelling wave solutions like wave equation.

Furthermore, when the function $\dynRHS$ is non--linear, numerically approximating the equation \eqref{eqn:dynODE_ROM} requires projecting the reduced state back to the full space each time, thus losing any computational advantage. To overcome such a drawback, hyper--reduction techniques are necessary, aiming to reduce computational cost while maintaining acceptable accuracy. Hyper--reduction techniques include methods such as the \textit{Discrete Empirical Interpolation Method (DEIM)} \cite{chaturantabut2010nonlinear}, which seeks to approximate non--linear terms using a reduced number of sampling points. However, finding a good compromise between accuracy and computational cost is challenging, as excessive reduction in the number of sampling points can lead to significant loss of accuracy, while too many sampling points can {annihilate} the benefits of dimensionality reduction.

\null\textbf{Autoencoder--based methods.}
{To overcome} limitations of projection--based dimensionality reduction techniques, system {reduction} {can be achieved via} non--linear transformations \cite{kashima2016nonlinear,hartman2017deep,lee2020model}.
These approaches {rely on} replacing the transformations $\dynLat = W^T \dynOut$ and $\dynOut \simeq W \dynLat$ {by successive application of two non--linear maps, an} encoder $\dynLat = \dynNNenc(\dynOut; \dynWenc)$ and a decoder $\dynOut \simeq \dynNNdec(\dynLat; \dynWdec)$, where $\dynNNenc$ and $\dynNNdec$ are neural networks with trainable parameters $\dynWenc$ and $\dynWdec$, respectively. {Decoders and encoders are typically made of} FFNNs, or, when the variable $\dynOut$ consists of values of a spatial field sampled on a Cartesian grid convolutional neural networks (CNNs).

{The sequential application} of networks $\dynNNenc$ and $\dynNNdec$ constitutes an autoencoder (see Sec.~\ref{sec:architectures}): the dimension of the reduced state variable $\dynLat$ is lower than that of $\dynOut$, so that autoencoders learn a compact representation of data. Thanks to {their} ability to represent non--linear relationships, autoencoders can capture compact representations much more efficiently than linear dimension reduction methods. As a matter of fact, linear dimensionality reduction underlying projection--based ROMs can be seen as a special case of autoencoders, where {both} the encoding function ($\dynLat = W^T \dynOut$) and decoding function ($\dynOut \simeq W \dynLat$) are linear.

Autoencoders are trained to minimize the reconstruction error, i.e., the discrepancy between original data and data reconstructed starting from encoded reduced state variable. In other words, an autoencoder is trained by simultaneously training the two NNs, through the following optimization problem:
\begin{equation}\label{eqn:autoencoder-manifold}
    \dynWenc^*, \dynWdec^* = \argmin{\dynWenc, \dynWdec} 
    \sum_{j = 1}^{\dynNumSamples} \sum_{i = 0}^{\dynNumTimes} \left\| \widehat\dynOut^j_i - \dynNNdec(\dynNNenc(\widehat\dynOut^j_i; \dynWenc); \dynWdec) \right\|^2.
\end{equation}
Once autoencoder is trained and thus able to extract a compact representation $\dynLat \in \mathbb{R}^{\dynNumLat}$, often called \textit{latent variable}, of the state $\dynOut \in \mathbb{R}^{\dynNumOut}$, it is necessary to learn dynamics in the reduced space $\mathbb{R}^{\dynNumLat}$ (also called \textit{latent space}). This can be done, similarly to projection--based ROMs, by projecting dynamics onto the tangent space of the non--linear manifold, although this operation can frustrate benefits of reduction due to high computational costs of this operation \cite{lee2020model}. A generally more effective alternative is to learn dynamics directly in the reduced space, by training a time--dependent model (e.g., based on Neural ODEs, RNNs, or LSTM \cite{liu2022hierarchical,maulik2021reduced, vlachas2022multiscale}, {see Sect. \ref{sec:operator-learning-time-dependent}}) {as follows}:
\begin{equation} \label{eqn:dynODE_ROM_autoencoder}
    \left\{
    \begin{aligned}
        & \frac{d}{dt} \dynLat(t) = \dynNNdyn(\dynLat(t), \dynInp(t); \dynWdyn) & & t \in (0, T) \\
        & \dynLat(0) = \dynNNenc(\dynOut_0; \dynWenc). & &
    \end{aligned}
    \right.
\end{equation}
This model is trained to reproduce trajectories obtained by encoding training data, i.e., the trajectories defined by $\dynLat^j_i = \dynNNenc(\dynOut^j_i; \dynWenc)$. Similarly to the case considered in Sec.~\ref{sec:operator-learning-time-dependent}, the model {\eqref{eqn:dynODE_ROM_autoencoder}} can be trained using either a derivative--based or trajectory--based approach, and the trade--offs described above apply. 
Additionally, the SINDy method can also be used in combination with autoencoders to obtain parsimonious and interpretable models of the system's intrinsic dynamics \cite{champion2019data}.
Whatever model is used to describe dynamics in the reduced space, the predicted output of interest is then obtained by decoding the reduced state variable, i.e., $\dynOut(t) = \dynNNdec(\dynLat(t); \dynWdec)$.

Typically, the autoencoder model is trained prior to the time--dependent model, and the latter {uses} optimal parameters of encoder and decoder. An alternative is to perform training jointly, the so--called end--to--end approach, where autoencoder and time--dependent model are trained simultaneously to minimize reconstruction error and dynamics prediction error, respectively. The end--to--end approach can lead to better generalization of the model, as the autoencoder can learn a representation more suited to the system's dynamics.
The price to pay is {greater} computational complexity, as end--to--end training is characterized by a deeper computational graph and a larger parameter space.

\null\textbf{Model learning with latent variables.}
The methods described in the previous paragraph explicitly construct a compact representation of the state variable $\dynOut(t)$ using an encoder and a decoder. An alternative approach is to learn the system's dynamics directly in a latent space, without explicitly constructing the encoding function. This approach, introduced in \cite{RDQ-JCP-2019} under the name of \textit{model learning}, consists of learning a model like
\begin{equation} \label{eqn:dynODE_ROM_latent}
    \left\{
    \begin{aligned}
        & \frac{d}{dt} \dynLat(t) = \dynNNdyn(\dynLat(t), \dynInp(t); \dynWdyn) & & t \in (0, T) \\
        & \dynOut(t) = \dynNNdec(\dynLat(t); \dynWdec) & & t \in (0, T) \\
        & \dynLat(0) = \dynLat_0 & &
    \end{aligned}
    \right.
\end{equation}
where $\dynNNdyn$ and $\dynNNdec$ are neural networks with trainable parameters $\dynWdyn$ and $\dynWdec$, respectively. 
The latent variable $\dynLat$ provides a compact and meaningful representation of the data, but it is not explicitly constructed through an encoding function, as in the case of autoencoder--based methods. As we will see later, this approach features several advantages over autoencoder--based methods.

The model is trained in an end--to--end manner, minimizing the prediction error of the observed variable $\dynOut$, through the following optimization problem:
\begin{equation}
    \dynWdyn^*, \dynWdec^* = \argmin{\dynWdyn, \dynWdec} 
    \sum_{j = 1}^{\dynNumSamples} \sum_{i = 0}^{\dynNumTimes} \left\| \widehat\dynOut^j_i - \dynNNdec(\dynLat^j(t_i); \dynWdec) \right\|^2,
\end{equation}
where $\dynLat^j(t_i)$ represents the latent variable obtained by solving the ODE \eqref{eqn:dynODE_ROM_latent} for sample $j$ at time $t_i$.
The loss function {implicitly depends} on $\dynWdec$, as $\dynLat^j(t_i)$ depends on $\dynWdyn$ through the numerical solution of \eqref{eqn:dynODE_ROM_latent}. 
We note that the approach followed here uses a trajectory-based loss function, combined with the decoder. We remark that, clearly, the derivative--based loss function is not compatible with the presence of latent variables, since their trajectories are not known a priori.

Special attention must be given to the treatment of initial conditions, due to the lack of an encoder that maps the initial condition $\dynOut_0$ to the corresponding latent variable $\dynLat_0$. In many cases, the samples start from a common state: in this case, it can be assumed without loss of generality that $\dynLat_0 = \mathbf{0}$. Otherwise, the initial condition associated with each training sample can be treated as an additional parameter to be optimized, along with the model parameters. During inference, the initial condition is then estimated through data assimilation from the observation of some time instances of the system outputs (see \cite{regazzoni2021combining,ziarelli2024model} for more details).

Model learning with latent variables is particularly effective because the latent space is not fixed a priori, as it is when using a pre--trained autoencoder, but is discovered simultaneously with the dynamics. This leads to discovering a latent space that can not only accurately reconstruct the original system variable but also predict the intrinsic dynamics. A two--step training, on the other hand, tends to favour an encoding that is functional to the state reconstruction but may overlook features that are essential for capturing the dynamics. Compared to end--to--end training with autoencoders, model learning with latent variables is generally computationally convenient, as it does not require an encoder and thus has fewer trainable parameters \cite{regazzoni2024ldnets}. Moreover, the model learning approach with latent variables shows great flexibility, for example, when the variable to be reconstructed has a non--uniform spatial sampling (a situation where it is not possible to apply an encoder) as we will show in Sec.~\ref{sec:space-time learning}. Additionally, the model learning approach with latent variables can naturally capture hidden dynamics and handle non--Markovian systems (i.e., systems in which the state at time $t_{i+1}$ may depend on other previous states other than that at time $t_i$), as we will see in the next paragraph.

\null\textbf{Hidden dynamics discovery.}
Another context in which the methods described in Sec.~\ref{sec:operator-learning-time-dependent} may not be sufficiently descriptive is when the observable variables available for training do not contain enough information to fully capture the system's dynamics. 
In particular this is the case of non-Markovian systems, where the system's evolution depends not only on the observed variable $\dynOut(t)$ and the input $\dynInp(t)$, but also on additional unobserved internal variables which encode the system's history.
In these contexts, methods capable of discovering hidden dynamics, not directly observable in the data, are necessary.

To exemplify, let us consider a system whose dynamics are governed by the following system of ODE
\begin{equation} \label{eqn:dynODEstate}
    \left\{
    \begin{aligned}
        & \frac{d}{dt} \dynState(t) = \dynRHS(\dynState(t), \dynInp(t)) & & t \in (0, T) \\
        & \dynOut(t) = \dynOBS(\dynState(t)) & & t \in (0, T) \\
        & \dynState(0) = \dynState_0, & &
    \end{aligned}
    \right.
\end{equation}
where $\dynState(t) \in \mathbb{R}^{\dynNumState}$ is an internal state variable not observable, $\dynOut(t) \in \mathbb{R}^{\dynNumOut}$ is the observable variable, and $\dynInp(t) \in \mathbb{R}^{\dynNumInp}$ is the input. The evolution of the state variable $\dynState(t)$ is governed by the function $\dynRHS$, while the observable variable $\dynOut(t)$ is obtained through the observation function $\dynOBS$. In this context, we talk about \textit{hidden dynamics discovery} because the goal is to reconstruct the system's dynamics from observations of pairs $(\dynOut(t), \dynInp(t))$, without access neither to FOM equations \eqref{eqn:dynODEstate} nor observations of the state variable $\dynState(t)$.

Interestingly, the model learning method with latent variables can discover hidden dynamics of the system without requiring any modifications. In fact, the model \eqref{eqn:dynODE_ROM_latent} can discover compact representations of the data, which can include hidden information, and learn the system's dynamics directly in this latent space. In this way, the model can capture the hidden dynamics of the system, even if not directly observable, thus capturing non--Markovian effects.

Autoencoder--based methods, on the other hand, are not suitable, as the reduced state $\dynLat$ is defined by encoding the observed variable $\dynOut$, and therefore may not contain all the information necessary to capture the hidden dynamics. In this context, one approach is to augment the state variable obtained through the encoder with latent variables, aimed at capturing non--Markovian effects \cite{vlachas2022multiscale}.
Another approach is to use a time--delay embedding of the observed variable $\dynOut(t)$, defined as $\widetilde{\dynOut}(t):= (\dynOut(t), \dynOut(t + \Delta t^{\text{sampling}}), \dots, \dynOut(t + (q-1) \Delta t^{\text{sampling}})]$ and applying the autoencoder to this delay--augmented space of dimensions $q \cdot \dynNumOut$ \cite{bakarji2023discovering}. This approach is rooted in Takens' theorem \cite{takens1981dynamical}, which states that under certain assumptions, the delay--augmented state has an attractor that is diffeomorphic (i.e., deformable in a continuously differentiable manner) to the attractor of the underlying, although unobserved, state. Clearly, this requires that the observations are sampled at constant intervals $\Delta t^{\text{sampling}}$.

\subsubsection{Space-time operator learning} \label{sec:space-time learning}
The methods described in Secs.~\ref{sec:operator-learning-time-dependent} and \ref{sec:operator-learning_hidden-dynamics-discovery} aim to learn the temporal dynamics of a system from a dataset of trajectories.
In many cases, systems of practical interest exhibit not only temporal but also spatial dependence.
Examples include models of fluids, structures, and populations, which are typically described by PDEs whose solutions are functions of both space and time. In these cases, it is necessary to learn the relationship between the spatial and temporal fields of the system and problem data. This task is known as \textit{space-time operator learning}.

Consider the following abstract spatio--temporal problem, defined in the spatial domain $\Omega \subset \mathbb{R}^d$ and in the time interval $(0, T)$:
\begin{equation} \label{eqn:dynPDEstate}
    \left\{
    \begin{aligned}
        & \frac{d\dynPDEState}{dt} (\mathbf{x}, t) = \dynRHSop(\mathbf x,\dynPDEState, \dynInp)
        & & \mathbf{x} \in \Omega, \, t \in (0, T) \\
        & \dynPDEOut(\mathbf{x}, t) = \dynOBSop(\mathbf x,\dynPDEState)
        & & \mathbf{x} \in \Omega, \, t \in (0, T) \\
        & \dynPDEState(\mathbf{x},0) = \dynPDEState_0(\mathbf{x}) & &
    \end{aligned}
    \right.
\end{equation}
where $\dynPDEState(\mathbf{x}, t) \in \mathbb{R}^{\dynNumState}$ is the state variable of the problem, which we generally assume to be unobservable, while $\dynPDEOut(\mathbf{x}, t) \in \mathbb{R}^{\dynNumOut}$ is the observable output.
Clearly, the case when the internal state is fully observable is encompassed in this framework, as the particular case when the observation operator $\dynOBSop$ is the identity operator.
Typically, we assume that, at each time, the solution of the problem belongs to a certain functional space that we denote by $\dynPDEOutSpace$, {typically} $L^2(\Omega)$ or $H^1(\Omega)$. Thus, we have $\dynPDEOut(\cdot, t) \in \dynPDEOutSpace$ for every $t \in (0, T)$. Additionally, $\dynInp(t) \in \mathbb{R}^{\dynNumInp}$ is a vector that collects the inputs (which for simplicity we assume to depend on time but not on space), which can {represent} forcing terms, coefficients of the PDE, possibly including the boundary conditions. Note that this also includes the case of inputs that are constant in time, i.e., scalar parameters that influence the system's dynamics. We denote by $\dynRHSop$ a differential operator that may include spatial derivatives with respect to the space variable $\mathbf x$ and incorporate in its definition suitable boundary conditions. The goal of space-time operator learning is to learn the operator that maps $\dynPDEState(\mathbf{x}, t)$ and $\dynInp(t)$, for $ \mathbf{x} \in \Omega$ and $t \in (0, T)$, to the spatio--temporal dynamics of $\dynPDEOut(\mathbf{x}, t)$.

It is worth noting that space--time operator learning can be addressed by methods such as DeepONets and Neural Operators, treating the time variable as if it were an additional spatial variable. However, as already noted in Sec.~\ref{sec:operator-learning-time-dependent}, the time variable has a different nature compared to spatial variables, and it can be advantageous to treat it differently, in order to learn models that are inherently consistent with the arrow of time and invariant with respect to temporal translations. Additionally, using methods that account for the special nature of the time variable can allow for time extrapolation and, in general, {enjoy} better generalization properties.

\null\textbf{Discretize in space and learn in time.}
A widely used approach for performing space--time operator learning involves first semi--discretizing the solution in space, thus encoding the output $\dynPDEOut(\cdot, t)$ at a fixed time into a vector (typically of very high dimensionality) $\dynPDEOutVEC(t) \in \mathbb{R}^{N_h}$. In this way, the space--time problem \eqref{eqn:dynPDEstate} is reduced to a time--dependent problem, which can be addressed using the methods described in Sec.~\ref{sec:operator-learning-time-dependent} and Sec.~\ref{sec:operator-learning_hidden-dynamics-discovery}.

Space discretization is performed through a discretization operator $\dynPDEOutDiscretize \colon \dynPDEOutSpace \to \mathbb{R}^{N_h}$, where $N_h \gg 1$ is the dimension of the vector $\dynPDEOutVEC(t)$. The discretization operator can be constructed in various ways, depending on the nature of the problem. Typically, it is constructed through the pointwise evaluation of the solution on a spatial grid of nodes. Alternatively, it can be defined through the coefficients associated with the expansion with respect to a set of basis functions, such as the Finite Element expansion reported in \eqref{eq:uN}. An alternative is to consider the coefficients associated with the expansion with respect to a Fourier basis, which involves performing the discrete Fourier transform of the solution and considering the coefficients associated with the transform. In any case, the result is a high-dimensional vector of dimension $N_h$ that represents the spatial solution at a fixed time. The subscript $h$ refers to the characteristic size associated with the discretization (e.g., the mesh element size for a grid--based discretization; the sampling period in the case of the discrete Fourier transform).

Once the discretization operator $\dynPDEOutDiscretize$ is defined, the space--time problem \eqref{eqn:dynPDEstate} is reduced to learning the dynamics of the vector $\dynPDEOutVEC(t) = \dynPDEOutDiscretize(\dynPDEOut(\cdot, t))$, using e.g. the methods described in Sec.~\ref{sec:operator-learning-time-dependent} and Sec.~\ref{sec:operator-learning_hidden-dynamics-discovery}.
Specifically, the discretized vector $\dynPDEOutVEC(t)$ is typically first encoded into a reduced variable $\dynLat(t)$ using an autoencoder \cite{lee2020model,vlachas2022multiscale,oommen2022learning} or POD \cite{wang2018model}, and then the dynamics in the reduced space are learned using a time--dependent model. The output of the model is then decoded back into the high--dimensional space, obtaining the predicted solution at a fixed time. Note that, in the case of autoencoder--based methods, the autoencoder can be chosen in a way reminiscent of the spatial structure of the problem, for example using a convolutional network \cite{oommen2022learning,vlachas2022multiscale,lee2020model}. However, this approach is confined to Cartesian grids on simple--shaped domains. Generalizations to arbitrarily shaped domains can instead be based on graph convolutional neural networks \cite{Wu20214,pichi2024graph}. The dynamics in the reduced space can be learned by using either RNNs \cite{liu2022hierarchical}, LSTMs \cite{maulik2021reduced,vlachas2022multiscale,wang2018model}, neural ODEs \cite{linot2022data}, or the SINDy method, should one be interested in learning an interpretable dynamics \cite{champion2019data}. As an alternative, methods that do not account for the arrow of time can be used to learn the dynamics in the reduced space, such as simple FFNNs \cite{fresca2020deep} or DeepONets \cite{oommen2022learning}.

The \textit{discretize in space and learn in time} approach thus allows reducing the space--time problem to a temporal problem, significantly simplifying the problem. However, this simplification {could ignore} important information, particularly those related to the spatial structure of the problem. Additionally, the semi--discretization in space can lead to high--dimensionality problems, which can be challenging to handle {by} standard methods.

\null\textbf{Latent Dynamics Networks (LDNets).}
The Latent Dynamics Networks (LDNets) {method}, proposed in \cite{regazzoni2024ldnets}, is based on the idea of learning the spatio--temporal dynamics directly in a latent space, without performing a spatial semi--discretization. 
Instead of discretizing the output field with respect to a fixed grid and storing it as a high--dimensional vector, LDNets represent the field as a continuous function. This function consists of a neural network (typically, a FFNN), which takes spatial coordinates (e.g. $x$, $y$ and $z$) as input while the output is the corresponding value of the field at those coordinates. The field is evaluated at any point in space and is not limited to a predefined grid resolution. This eliminates the trade-off between resolution and memory usage which affects grid--based methods. Additionally, the neural network takes as additional input the latent state: the field is in this way \textit{conditioned} on (that is subject to) the latent state, meaning that for a fixed latent state, the output field is represented as a continuous function of the spatial coordinates, but when the latent state varies, the output field changes accordingly. This allows the latent state to capture the essential dynamics of the system in a compact form, reducing the dimensionality of the problem while still enabling accurate representation of complex spatio--temporal patterns.
The LDNet model is defined as follows:
\begin{equation} \label{eqn:LDNet}
    \left\{
    \begin{aligned}
        & \frac{d}{dt} \dynLat(t) = \dynNNdyn(\dynLat(t), \dynInp(t); \dynWdyn) 
        & & t \in (0, T) \\
        & \dynOut(\mathbf{x}, t) = \dynNNdec(\mathbf{x}, \dynLat(t); \dynWdec)
        & & \mathbf{x} \in \Omega, \, t \in (0, T) \\
        & \dynLat(0) = \dynLat_0 & &
    \end{aligned}
    \right.
\end{equation}
where $\dynNNdyn$ and $\dynNNdec$ are neural networks with parameters $\dynWdyn$ and $\dynWdec$, respectively. The latent variable $\dynLat \in \mathbb{R}^{\dynNumLat}$ represents a compact and meaningful representation of the spatial field. 
The LDNet method can be seen as a generalization to the case of spatio--temporal problems of model learning with latent variables (see Sec.~\ref{sec:operator-learning_hidden-dynamics-discovery}): also in this case, the latent variable is not fixed a priori, as it is when using a pre--trained autoencoder, but is discovered simultaneously with the dynamics, so it can also discover hidden information that is not directly observable in the data. In this way, the LDNet model can capture non--Markovian dynamics, as the latent variable can represent hidden information that influences the system's dynamics.

Training is performed through the following optimization problem, which simultaneously trains the two neural networks:
\begin{equation} \label{eqn:LDNet_loss}
    \dynWdyn^*, \dynWdec^* = \argmin{\dynWdyn, \dynWdec} 
    \sum_{j = 1}^{\dynNumSamples} \sum_{i = 0}^{\dynNumObservation} \left\| \widehat\dynOut^j_i - \dynNNdec(\mathbf{x}_i, \dynLat^j(t_i); \dynWdec) \right\|^2,    
\end{equation}
where $\dynLat^j(t_i)$ represents the latent variable obtained by solving the ODE \eqref{eqn:LDNet}{$_1$} for sample $j$ at time $t_i$, and $\widehat\dynOut^j_i$ is the observed output corresponding to the pair $(\mathbf{x}_i, t_i)$, for $i = 0, \dots, \dynNumObservation$. 

The continuous spatial representation allows capturing the spatial structure of the problem without the need for semi-discretization, thus obtaining a mesh-less and resolution-invariant method of great flexibility. Indeed, as can be seen from the loss function \eqref{eqn:LDNet_loss}, the training does not require the training data to be sampled on a fixed grid, but data with variable sampling between time instances and even between samples can be used. Additionally, the ability to obtain the output at any point, on a query basis, allows the use of stochastic training methods, thus lightening the computational burden associated with training, which would not be the case if the decoder returned the entire batch of evaluations. It is also important to note that the decoder--only nature of LDNets, combined with the continuous spatial representation, allows LDNets to perform training working only in low--dimensional spaces: both $\dynNNdyn$ and $\dynNNdec$ have typically low-dimensional inputs and outputs, and in any case independent of the spatial resolution of the problem (associated with $N_h$, in the case where the data are generated through, e.g., Finite Element simulations). This makes LDNets computationally efficient, easily scalable, and especially -- thanks to their parsimonious structure -- very efficient at generalizing. Indeed, in the test cases considered in \cite{regazzoni2024ldnets}, LDNets have shown greater accuracy compared to autoencoder--based methods, with orders of magnitude fewer trainable parameters.

\subsubsection{{Foundation models for operator learning}}

\emph{Foundation models} are generalist models that are pretrained on vast amounts of data extracted from heterogeneous distributions. Subsequently, they can optionally be finetuned on a few task--specific samples. Pretraining allows them to exploit available non--specific data, and finetuning makes them adapt to several downstream tasks.
Foundation models can find fertile ground in the field of PDEs because \emph{(i)} different PDEs share features or kernels that can be learned in the pretraining phase, \emph{(ii)} often, the high--quality samples needed to train a specific PDE are few and these could be only used in the finetuning phase, \emph{(iii)} foundation models show a high capacity to generalize to unseen and (apparently) unrelated PDEs.

%
%
\null\textbf{Poseidon.}
Poseidon \cite{herde2024poseidon} is a  \emph{foundation model}, designed to learn PDE solution operators. Let us consider the Initial Boundary Value Problem (\ref{eq:IBVP}), let $\mathcal X\subset L^p(\Omega;\mathbb R^m)$ for some $1\leq p<\infty$ be a functional space and $u\in \mathcal{C}([0,T],\mathcal X)$ be the solution of (\ref{eq:IBVP}). We are interested in looking for the solution operator $\mathcal S:[0,T]\times \mathcal X \to \mathcal X$ such that $u(t)=\mathcal S(t, u_0)$ is the solution of (\ref{eq:IBVP}) at any time $t\in[0,T]$. Given a data distribution $\mu\in \text{Prob}(\mathcal X)$, the Operator Learning Task reads:
\begin{quoteit}
    Given any initial datum $u_0\sim \mu$, find an approximation $\mathcal S^*\approx \mathcal S$ to the solution operator $\mathcal S$ of (\ref{eq:IBVP}), in order to generate the entire solution trajectory $\{\mathcal S^*(t, u_0)\}$ forall $t\in[0,T]$.
\end{quoteit}

\begin{figure}
    \centering
    \includegraphics[width=0.9\linewidth]{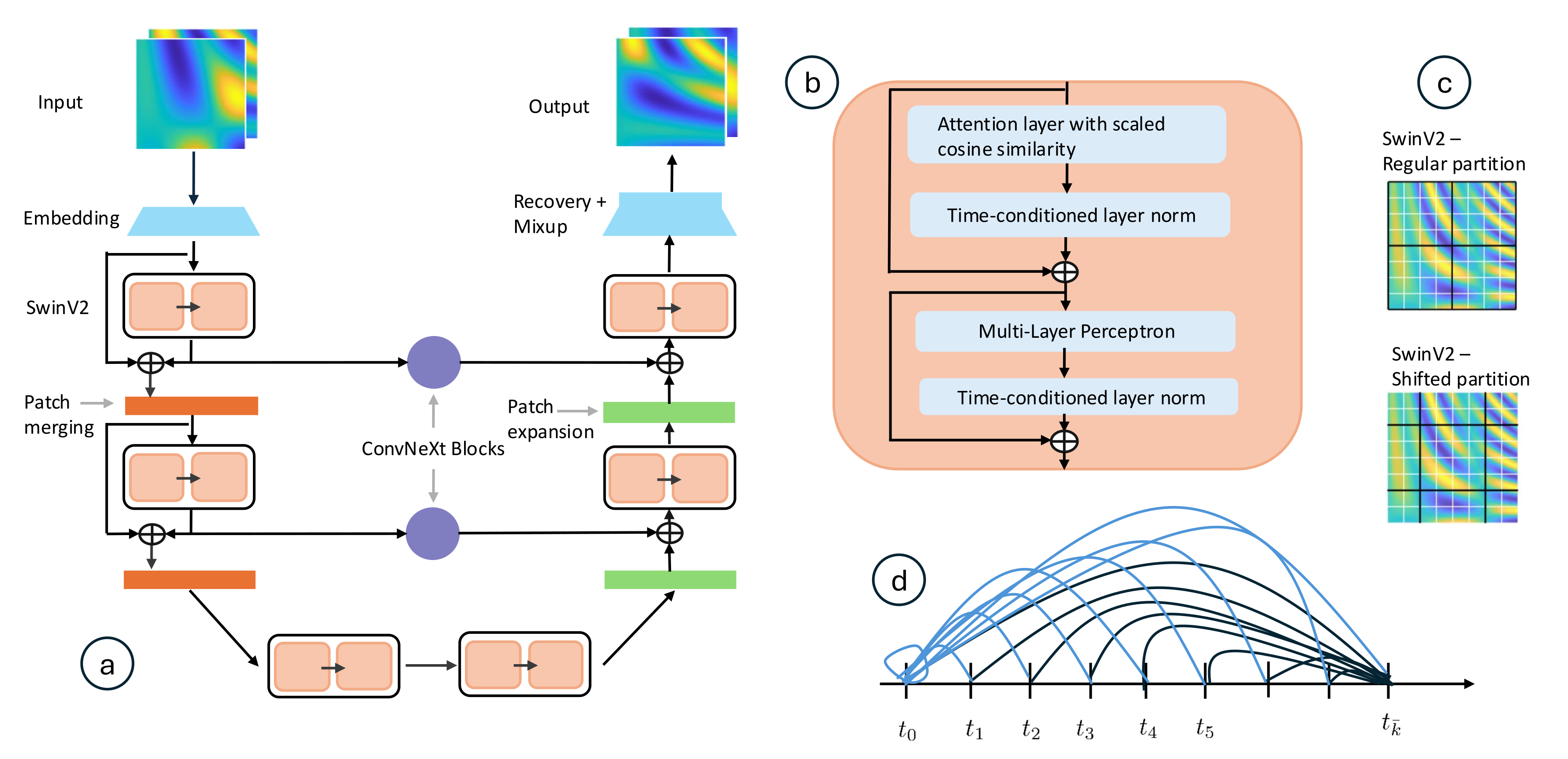}
    \caption{(a) The scalable Operator Transformer scOT that is the backbone of Poseidon, (b) one of the two blocks of the SwinV2 transformer, (c) regular and shift window partitions, (d) the all2all training strategy: light--blue lines refer to the pairs with $k=0$ and $\bar k\geq k$, black lines refer to the pairs with $\bar k$ given and $k\leq \bar k$.}  
    \label{fig:poseidon}
\end{figure}

The backbone of Poseidon, named \emph{scalable Operator Transformer (scOT)} is a \emph{hierarchical multiscale vision transformer} based on SwinV2 \cite{liu2022swinV2} (see Sect. \ref{sec:ML-models}) and with lead--time conditioning, see Fig. \ref{fig:poseidon}.  
ScOT is a U--Net style \cite{Ronneberger2015-unet, cao2021unet} encoder--decoder architecture whose basic blocks are Shifted Windows Version2 (SwinV2) transformer blocks (see Sect. \ref{sec:ML-models}). 
U--Nets are characterized by a downscaling phase (in the encoder) and an upscaling phase (in the decoder). The downscaling is implemented by patch mergings, while upscaling is by patch expansions. Patches are the basic units into which an image is split up and play the role of tokens in transformer architecture. We recall that SwinV2 considers two different partitions of the image in windows (a window is a cluster of patches), one regular and the other shifted to the regular, and then it applies self--attention models to each window, locally. Alternating between regular and shifted partitions ensures that adjacent patches belonging to different windows can relate to each other during the attention mechanism. Indeed, such an interaction would be absent if only the regular partition were used.
Moreover, ConvNeXT blocks (i.e., deep CNNs with GELU\footnote{$GELU(x)=x\int_{-\infty}^x f(t)\, dt$, where $f(t)$ is the normal probability density function. It is approximated by the formula  $GELU(x)=0.5\,x\,(1+\tanh[\sqrt{2\pi}\ (x+0.044715\,x^3)])$} activation function) connect the encoder and decoder at intermediate points.

Following \cite{Perez2017}, to enable continuous--in--time evaluation, scOT replaces all standard normalization layers inside the SwinV2 blocks with the so--called \emph{lead--time conditioned} layer norm
\begin{equation}
    \text{LN}_{\alpha(t),\beta(t)}(\mathbf v)(x)=\alpha(t)\odot\frac{\mathbf v(x)-\mu_{\mathbf v}(x)}{\sigma_{\mathbf v}(x)}+\beta(t),
\end{equation}
where $\mu_{\mathbf v}$ and $\sigma^2_{\mathbf v}$ are the mean and variance of $\mathbf v$, while $\alpha(t)=\alpha_1 t+\alpha_2$, and $\beta(t)=\beta_1 t+\beta_2$ are affine functions in $t$ with $\alpha_i, \beta_i$ learnable parameters.

Poseidon is based on a novel training strategy, named \emph{all2all}, in which the ground truth training data are given in the form of trajectories and the loss function is defined by exploiting the semi--group property satisfied by the solution operator $\mathcal S$. More precisely, 
let $\{\mathcal S(t_k, u_{0,i})\}$ (with $k=0,\ldots, K$ and $a_i\sim\mu$, $i=1,\ldots, M$) be the ground truth training data. Denoting with $\theta$ the array of the learnable parameters, the loss function is defined by
\begin{equation}\label{eq:poseidon_loss}
    \widehat{\mathcal L}(\theta)=
    \frac{1}{M\widehat K} \sum_{i=1}^M\sum_{\bar k=0}^{K}\sum_{k=0}^{\bar k}\|\mathcal S(t_{\bar k}-t_k,u_i(t_k))-\mathcal S_\theta^*(t_{\bar k}-t_k,u_i(t_k))\|^p_{L^p(\Omega)},
\end{equation}
where $u_i(t_k)=\mathcal S(t_k,u_{0,i})$ (approximately) and $\widehat K=(K+1)(K+2)/2$ is the total number of pairs $(k,\bar k)$ with $k\leq \bar k$. In fact, (\ref{eq:poseidon_loss}) exploits the fact that the solution operator $\mathcal S$ of (\ref{eq:IBVP}) has the semi--group property
\begin{equation}
    u(\bar t)=\mathcal S(\bar t, u_0)=\mathcal S(\bar t - t, u(t))=\mathcal S (\bar t -t, \mathcal S(t, u_0)), 
\end{equation}
for any $t,\ \bar t\in [0,T]$ and for any initial condition $u_0$.  
We notice that, at each time step $t_{\bar k}$, all the previous values $u_i(t_k)$ with $k\leq \bar k$ are used to evaluate the loss.

In \cite{herde2024poseidon}, the data for the pretraining have been produced by solving both incompressible Navier--Stokes and compressible Euler  equations with periodic boundary conditions on $\Omega=[0,1]^2$, along the time interval $(0,1)$. Two and four different types of initial conditions have been considered for the Navier--Stokes and Euler equations, respectively, to represent a large variety of dynamics. Navier--Stokes equations have been solved by Fourier Spectral methods on a grid of $128\times 128$ frequencies and 21 time--steps. Euler equations have been solved 
 by a high--order finite volume scheme on a grid of $512\times 512$ cells, (later downsampled to $128\times 128$) and 21 time--steps.

The pretraining dataset contains 77840 trajectories, 720 of which are for the test set and 1440 for the validation set. Each trajectory is uniformly sampled at eleven time--steps (with every other time--step selected).

Fifteen downstream tasks have been solved, going from Navier--Stokes and compressible Euler equations with periodic boundary conditions (like in the pertaining step), to partially or totally new physics. Among the new physics, they considered: Navier--Stokes equations enriched with an advection--diffusion equation modelling the transport of a passive tracer; Navier--Stokes equation with a forcing term; compressible Euler equation with gravitation; the wave equation with homogeneous Dirichlet boundary conditions; the Allen-Cahan equation for phase transition; the steady state compressible Euler equation around an airfoil by using a non--cartesian mesh; the Poisson equation with homogeneous Dirichlet conditions and Gaussian bumps forcing term; the Helmholtz equation in the frequency domain with Dirichlet boundary conditions.

The downstream tasks have been solved with finite difference schemes, finite volumes, finite elements, or Fourier spectral methods. From 1260 to 20000 trajectories have been produced to finetuning each task.

To finetune the pretrained foundation model for any downstream task, the vector of learnable parameters is split into two parts. The parameters associated with the initial embedding and final recovery are randomly initialized and trained, while those associated with the U--Net (SWin2, patch merging, patch expansion, and ConvNeXt layers) are initialized by transferring the corresponding parameters from the pretrained model.

Numerical results show that Poseidon can learn effective representations from a small set of PDEs during the pretraining and generalize well to unseen and unrelated physics downstream \cite{herde2024poseidon}.

\null\textbf{In--Context Operator Networks (ICON).}
\emph{In--context learning} refers to the ability of a generative language model to learn or perform a specific task without re--training or fine--tuning, but solely by specifying the \emph{context}, which is the description of the task jointly with some examples related to that task \cite{dong-etal-2024-survey}.

In \cite{Yang2023}, the idea of in--context learning has been extended to learn operators that underlie differential equations. 
An \emph{In--Context Operator Newtowrks} ICON model is a transformer encoder--decoder architecture \cite{Vaswani2017} (see Sect. \ref{sec:architectures} and Fig. \ref{fig:transformer}) designed to learn operators acting on mathematical functions. 
Functions are expressed discretely as sets of key--value pairs, where, \emph{keys} are the input (selected values of the independent variables) and \emph{values} are the outputs (corresponding values of the dependent variables) of the function.

Inputs and outputs of the operators to be learned are named \emph{conditions} and \emph{quantities of interest (QoIs)}, respectively.
For instance, in a forward differential problem, conditions are the coefficient functions, or the initial/boundary conditions, while the QoI is the problem solution. In an inverse problem, the QoI could be a parameter function, while the solution of the differential equation, or even an observed variable, is the condition.  

A pair of one condition and the corresponding QoI is called 
\emph{example}, while the condition that we want ICON to return the corresponding QoI is named \emph{question condition}. A collection of examples plus the question condition compose a \emph{prompt} of an ICON model.
Since the QoI is a function, we must provide the model with a set of \emph{queries}, which are the keys with which we want to evaluate the QoI returned by the model. 

\begin{figure}
    \centering
    \includegraphics[width=0.8\linewidth]{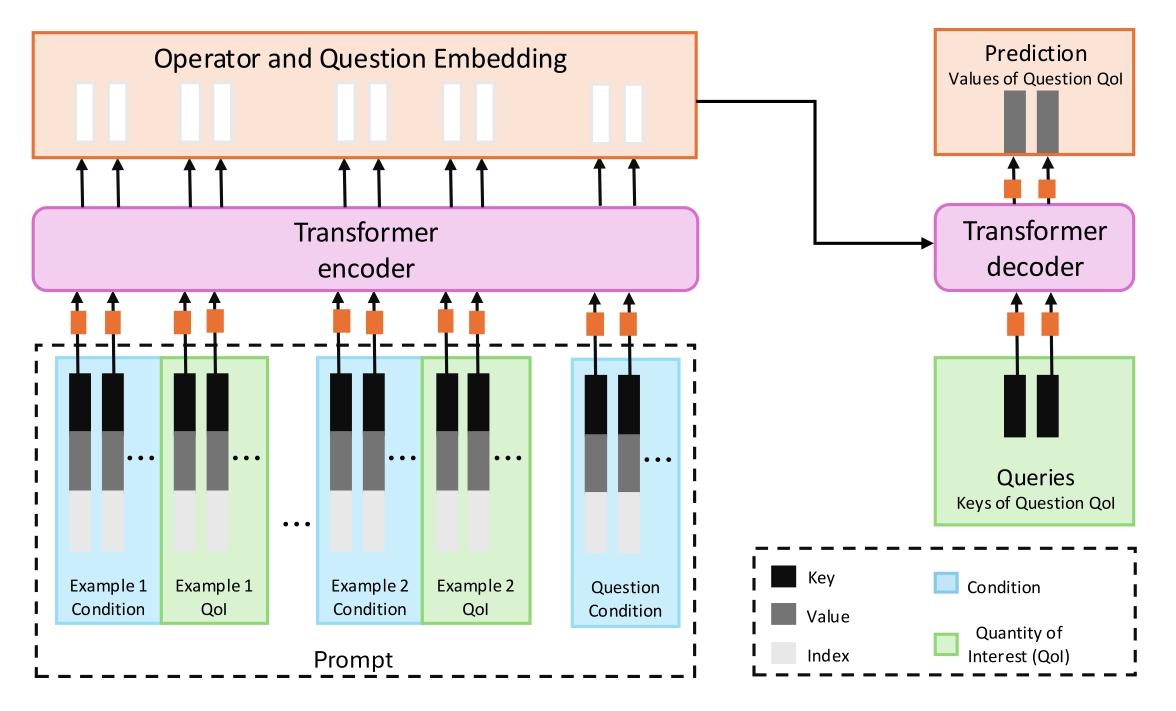}
    \caption{The neural network architecture of ICON}
    \label{fig:ICON}
\end{figure}

The transformer encoder takes the prompt as an input and returns the operator and question embedding 
(i.e. the internal representation of both the learned operator and question), which, jointly with the queries, are the inputs for the transformer decoder. Finally, the decoder provides the values of the question QoI (see Fig. \ref{fig:ICON}).

The ICON architecture can process input sequences of any length and a variable number of examples. Moreover, it accepts key--value pairs of different lengths for each condition/QoI function.

Training data can contain numerical solutions to different kinds of differential problems. At each iteration of the training process, a random batch of prompts, queries, and labels is built, and the batch can contain data relative to different operators. The loss function is the Mean Square Error between the NN's answers and the labels.
The training set presented in \cite{Yang2023} includes linear first--order ODEs, a damped oscillator, the Poisson equation, and linear and non--linear reaction--diffusion equations, all as forward and inverse problems in one variable, as long as a mean--field control problem in one and two dimensions. 

After the training, the NN is used to predict the question QoI, based on the question condition and a few examples describing the operator. In a forward pass, ICON learns the operator from the examples and applies it to return the question QoI without updating its weights or fine--tuning. 

Numerical results reported in \cite{Yang2023} show that
ICON can generalize to either finer or smaller meshes than those used with respect to the size of the key--value pairs used for the training, as well as to operators outside the distribution of operators observed during the training. Moreover, it can learn operators of new forms that were never seen in training data (for instance, obtained by adding a linear term to known operators).

\color{black}

\subsection{Topics related to Scientific Machine Learning not covered in this paper}
\label{sec:other_SciML_topics}

The field of SciML is rapidly expanding, with new research directions continuously emerging. We have decided to make a selection of the topics covered in this document, without any intention of {labelling these} topics as more important than others, but rather following a criterion of thematic coherence. As done in Section \ref{sec:other_ML_topics}, here we also provide a quick glance at topics not covered in this document, so that interested readers can explore them independently.
Indeed, the strategies that enable the combination of data--driven and physics--based methods are virtually limitless, both from a methodological perspective and in relation to the specific thematic or application.

Among the emerging topics in the field of SciML, we mention the data--driven discovery of constitutive laws, which aims to identify complex material relationships from data \cite{kirchdoerfer2017data,oishi2017computational,flaschel2021unsupervised,tac2022data,linka2023new}. 
A field of great interest for many application area is learning the solution operators of PDEs accounting for the variability of the computational domain \cite{kashefi2021point,regazzoni2022usmnets,yin2024dimon}.
Another research direction concerns the acceleration of traditional scientific computing algorithms, where machine learning techniques are used to improve computational efficiency, e.g. by learning artificial viscosity models or optimal stabilization parameters \cite{ling2016reynolds,xiao2021using,caldana2025discovering}. Additionally, SciML for linear algebra and domain decomposition, \cite{heinlein2019machine,klawonn2024machine,antonietti2023accelerating,caldana2024deep}, represents an area of {emerging} interest. Finally, we mention the learning of correction terms with respect to existing physical models aimed at improving the accuracy of models through the integration of empirical data \cite{duraisamy2019turbulence}, and the interplay between machine learning and control theory \cite{ruiz2023neural,ruiz2024control}.

\section{SciML for the iHeart simulator}\label{sec:SML-IHM}

In this section, we apply some of the SciML techniques presented in this work to computational cardiology, a field of great socio-economic relevance which poses numerous challenges for both the mathematical and computational sides.
In Sec.~\ref{sec:IHM} we first present the integrated heart model that describes in detail the functioning of the human heart, and show how it can be used to simulate the heartbeat and its interactions with the rest of the body. In Sec.~\ref{sec:SML-IHM} we then show how SciML techniques can be used to accelerate its numerical resolution and support sensitivity analysis and patient--specific calibration procedures, ultimately addressing questions of research and clinical relevance.
The numerical simulations presented in this section were obtained thanks to \lifex{}, a high--performance software library for solving multiphysics and multiscale problems, particularly specialized in computational cardiology \cite{africa2022lifex,africa2024lifex,africa2023lifex,bucelli2024lifex}.

\subsection{The integrated heart model}\label{sec:IHM}


The heart is a muscular organ consisting of four chambers (two atria on top and two ventricles below) and functions as the body's pump, making blood circulate throughout the body via the circulatory system. 
Essential for sustaining life, it works tirelessly, contracting and relaxing in a rhythmic cycle to deliver oxygen--rich blood to tissues and organs and to expel waste products such as carbon dioxide and metabolic waste.

\begin{figure}
\begin{center}
\includegraphics[width=0.7\textwidth]{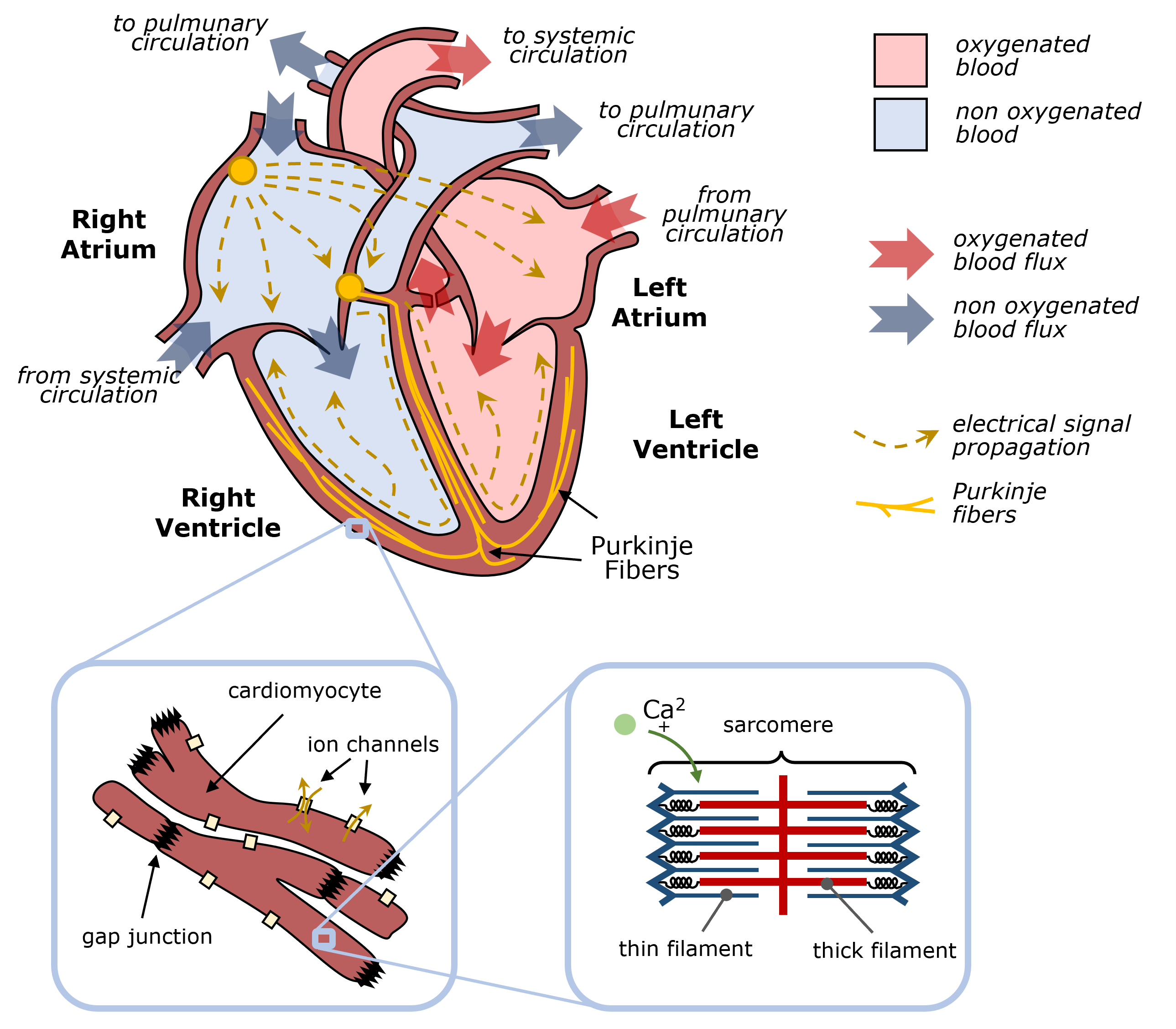}
\end{center}
\caption{The heart, the cardiac muscle's cells (cardiomyocytes), and a sarcomere, the fundamental unit of a muscle's striated tissue, responsible for its contraction}
\label{fig:fig1-QDR}
\end{figure}

For every second of our life, a spontaneous electrical stimulus starts at the sinoatrial node, a small cluster of specialized cells located in the right atrium of the heart that serves as a natural pacemaker (see Fig. \ref{fig:fig1-QDR}). This stimulus generates an electric current that propagates in every cell of the heart, the cardiomyocytes, first spreading in the upper part, the atria, and then, after reaching the atrioventricular node, located between atria and ventricles, in the lower part, the ventricles. 
The electric current is due to the movement of several ions (mainly sodium, calcium, and potassium) across cells' membrane and the whole phenomenon, known as action potential, is studied by the branch of \emph{cardiac electrophysiology}.

Among the various ions that come into play, calcium  induces a force (named \emph{active force}) in every single cardiomyocyte,  making the latter relax and contrive synchronously. Cardiomyocytes can be thought of as one--dimensional microscopic structures which deform longitudinally and create a train of waves that propagate from one cardiomyocyte to the other. 
Thanks to a transmission process across different space--time scales, the force generated in the cardiomyocytes induces a deformation of the entire myocardial tissue providing the contraction of the ventricles and atria. This third important process is named \emph{tissue mechanics}.

Due to the presence of blood in the four chambers, tissue mechanics induces \emph{blood dynamics}: from the right atrium the blood moves into the right ventricle crossing the tricuspid valve and, thanks to the pulmonary circulation, it reaches the lungs to be purified; then it re--enters the left atrium, moves into the left ventricle crossing the mitral valve, and finally, it is injected into the aorta to reach every cell of the body through the systemic circulation. 
The \emph{valve dynamics} is regulated by the blood dynamics, in particular by the different values the blood pressure assumes inside the chambers and the large arteries, and by the tissue mechanics.

Besides reaching all the cells of the body, the oxygenated blood must also reach the cells of the cardiac muscle to make it function properly and maintain its continuous pumping action.
This process, known as \emph{heart's perfusion}, occurs thanks to the presence of the coronary arteries which are supplied by ascending aorta.

Finally, cardiac electrophysiology is combined with the propagation of the electric signal through the torso to compute high--quality electrocardiograms, with the ultimate goal of validating the calibrated models.

\begin{figure}[h!]
\begin{center}
\includegraphics[trim=0 2cm 0 1cm, width=0.8\textwidth]{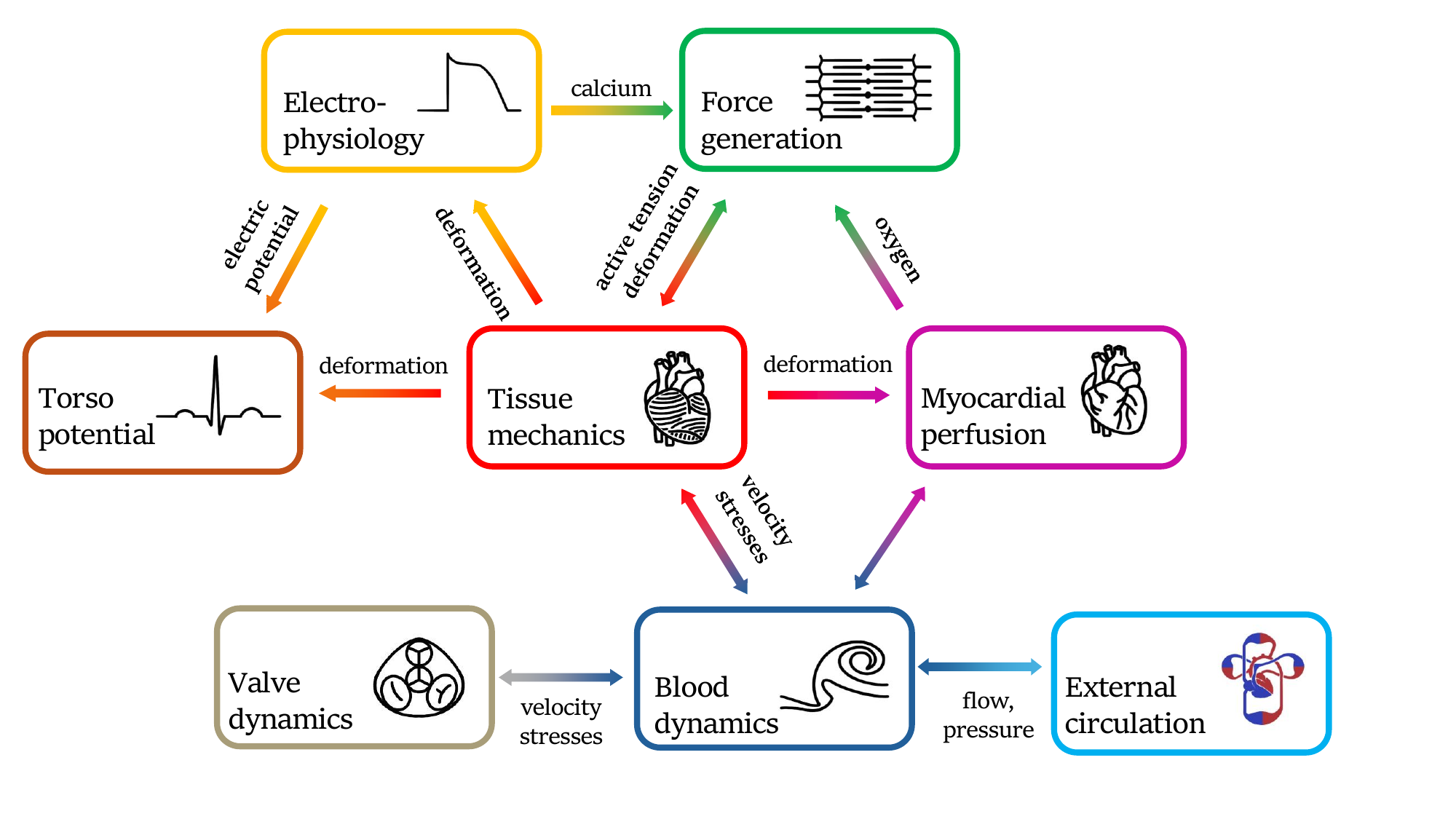}
\end{center}
\caption{Interactions among fundamental processes occurring in the cardiac function}
\label{fig:heart-physiology}
\end{figure}

The complex interactions among all the processes listed above are schematically represented in  Fig. \ref{fig:heart-physiology}. 

\null\textbf{Modelling cardiac electrophysiology.}
The driver of the cardiac function is electrophysiology, i.e., the result of chemical and electrical processes taking place at diﬀerent spatial scales, from the subcellular to the whole organ scale. 

The propagation of action potentials and the excitation of cardiac cells throughout the whole cardiac tissue is modelled by either the monodomain or bidomain partial differential equations. These equations are complemented by appropriate ionic models that provide information on the ionic currents and sub--cellular mechanisms responsible for excitation. 

The \emph{bidomain model} complemented with suitable \emph{ionic models} \cite{ColliFranzone-2014,QDMV2019} is the richest mathematical model describing cardiac electrophysiology and reads as follows, where $\widehat\Omega\subset{\mathbb R}^3$ is the space occupied by the heart muscle at rest.
We look for the extracellular transmembrane potential $v_\text{e}$, the transmembrane potential $v$ (that is the difference between the intracellular potential $v_\text{i}$ and the extracellular one) and the  ionic variables $\mathbf{z}_{\text{ion}}$ (describing ion channels, concentrations, or phenomenological variables) such that
\begin{eqnarray}\label{eq:EP}
\left\{\begin{array}{ll}
    \displaystyle \chi C_\text{m} \frac{\partial v}{\partial t} - \nabla\cdot\left(J\mathbf F^{-1}\mathbf{D}_\text{i}\mathbf F^{-T} \nabla (v+v_\text{e})\right) +\chi I_\text{ion}(v, \mathbf z_{\text{ion}}) = I_\mathrm{app}(\mathbf x, t)  
    & \text{in } \widehat\Omega\times(0,T) \\[2mm] 
    -\nabla\cdot\left(J\mathbf F^{-1}\mathbf{D}_\text{i}\mathbf F^{-T} \nabla v\right) -\nabla\cdot\left(J\mathbf F^{-1}(\mathbf{D}_\text{i}+\mathbf{D}_\text{e})\mathbf F^{-T} \nabla v_\text{e}\right)=0 
    & \text{in } \widehat\Omega\times(0,T) \\[2mm] 
    \displaystyle \frac{\partial\mathbf z_{\text{ion}}}{\partial t} = \mathbf {\boldsymbol\Phi}_\mathrm{ion}(v, \mathbf z_{\text {ion}}) & \text{in } \widehat\Omega\times(0,T) \\[0.5em]
    v(\mathbf x, 0) =v_0(\mathbf{x}),\ 
    \mathbf z_{\text{ion}}(\mathbf x, 0) =\mathbf{z}_{\text{ion},0}(\mathbf{x}),\ & \text{in }\widehat\Omega.
\end{array}\right.
\end{eqnarray}
The parabolic and elliptic equations in (\ref{eq:EP}) are supplemented with homogeneous Neumann boundary conditions to guarantee uniqueness of solution.
The coefficient $C_\text{m}$ is the capacitance per unit area, while $\chi$ is the surface--to--volume ratio of the membrane.  $\mathbf F=\mathbf I + \nabla\mathbf d$ is the deformation gradient tensor depending on the displacement $\mathbf d$ of the myocardium (which is determined by solving the mechanics model) with $J=\text{det}(\mathbf F)$, while $\mathbf{D}_\text{i}$ and $\mathbf{D}_\text{e}$ are the intracellular and extracellular conductivity tensors, respectively, that 
characterize the electrical properties of the tissue and are responsible for the correct propagation of the action potential. They are anisotropic to take into account that the electrical signal propagates with a greater velocity along the direction of muscle fibers compared to the other orthogonal directions \cite{Piersanti_2021, Africa_2022_fibers} and are defined by
\begin{equation}\label{eq:diffusivity}
\mathbf{D_\star}=\sigma_f^\star \frac{\mathbf{Ff}_0\otimes \mathbf{Ff}_0}{|\mathbf{Ff}_0|^2}
+\sigma_s^\star \frac{\mathbf{Fs}_0\otimes \mathbf{Fs}_0}{|\mathbf{Fs}_0|^2}
+\sigma_n^\star \frac{\mathbf{Fn}_0\otimes \mathbf{Fn}_0}{|\mathbf{Fn}_0|^2}, \qquad \star\in\{\text{i},\text{e}\}.
\end{equation}
$\mathbf{f}_0$, $\mathbf{s}_0$, and $\mathbf{n}_0$ are the fibers, sheets, and cross--fibers directions, while $\sigma_f^\star,\ \sigma_s^\star,$ and $\sigma_n^\star$ are the corresponding conductivity coefficients along the three directions.

The mathematical expressions of the function ${\boldsymbol\Phi}_\mathrm{ion}(v, \mathbf z_{\text{ion}})$ (which describes the dynamics of gating variables and ionic concentrations) and the ionic current  $I_\text{ion}(v, \mathbf z_{\text{ion}})$ strictly depend on the choice of the ionic model \cite{TenTusscherPanfilov-2006, courtemanche1998}. Finally, $I_{\text{app}}$ is the external current triggering the action potential.

One of the most important ionic species in $\mathbf{z}_{\text{ion}}$ is the intracellular calcium concentration $\Cai$, a pivotal factor in the excitation--contraction coupling that triggers contractile force generation by muscle cells.

\null\textbf{Modelling active force generation.}
 Cardiomyocytes are organized in sarcomeres, cylinder--shaped elements composed of two families of thick and thin filaments that are made of proteins (see Fig. \ref{fig:fig1-QDR}, (c)).  A particular role is played by the proteins myosin and actin: they belong to the two different types of filaments and, as a result of the action of calcium, they interact making the filaments mutually slide and causing contractions and elongations of sarcomeres \cite{JenkinsKemnitzTortora2007, Huxley1954,Huxley1954a}.

At any time $t$ during the heartbeat, force generation models can be written in the general form
\begin{equation}\label{eq:4.1-QDR}
    \frac{d\actState}{dt} = \actRHS\left(\actState, \Cai,SL, \frac{d\,SL}{dt}
    \right), \qquad
    \actTension= \actOBS(\actState),
\end{equation}
where $\actState$ denotes a vector collecting the state variables associated with the dynamics of the contractile proteins, $\Cai$ is the local calcium concentration, $SL$ is the sarcomere length, while $\actTension$ denotes the generated active force per unit area. In \cite{Regazzoni-phd, RDQ-active-contraction-2018, RDQ-bio-2020, RDQ-active-force-2021} a highly detailed model that explicitly describes sub--cellular mechanisms at the microscale level driven by both calcium concentration and sarcomere length has been proposed and studied extensively.

\null\textbf{Modelling cardiac mechanics.}
The cardiomyocytes’ contraction and the pressure exerted by the blood onto the endocardium induce during each heartbeat large deformations of the heart muscle, which can be up to a few centimeters. 
Denoting by $\widehat\Omega\subset{\mathbb R}^3$ the space occupied by the heart muscle at rest, the myocardium deformation is modelled by solving the momentum balance equation whose unknown is the displacement $\mathbf{d}$ of the tissue with respect to the rest position: 
\begin{equation}\label{eq:mechanics}
\rho\frac{\partial^2 \mathbf{d}}{\partial t^2}-\nabla\cdot\mathbf{P}(\mathbf{d},\actState)=\mathbf{0}\qquad \textrm{in }\widehat\Omega\times(0,T),
\end{equation}
where $\rho$ is the mass density of the
body.
The \emph{ﬁrst Piola–Kirchhoﬀ stress tensor} $\mathbf{P}$ models the stresses acting on the body and embeds both active and passive mechanical properties
\begin{equation}\label{eq:PK-tensor}
\mathbf{P}(\mathbf{d},\actState)=
\mathbf{P}^{\textrm{act}}(\mathbf{d},\actState)
+\mathbf{P}^{\textrm{pass}}(\mathbf{d}).
\end{equation}
The active part $\mathbf{P}^{\textrm{act}}(\mathbf{d},\actState)$ directly depends on the active force $\actTension$ and must be constitutively defined by suitably upscaling the microscopically generated stress at the level of the sarcomere.
The passive part $\mathbf{P}^{\textrm{pass}}(\mathbf{d})$ is related to the passive mechanical response of the heart, it is obtained as the derivative of a strain energy function that characterizes the mechanical properties of the material and incorporates orthotropic hyperelastic constitutive laws \cite{Fedele2023-comprehensive, RSAFDQ2022, Usyk2002, Guccione1991}.

The displacement $\mathbf d$ of the myocardium affects cardiac electrophysiology through the so--called \emph{electro--mechanical feedback} by entering the definition of the conductivity tensors $\mathbf D_\text{i}$ and $\mathbf D_\text{e}$ (indeed, the deformation gradient tensor $\mathbf F$ involved in the definition (\ref{eq:diffusivity}) refers to the displacement $\mathbf d$ solution of (\ref{eq:mechanics})), as well as the blood dynamics through the fluid--structure interaction between the muscle and the blood.

To guarantee the well--posedness of the momentum balance equation (\ref{eq:mechanics}), suitable initial conditions must be defined in  $\widehat\Omega$ and boundary conditions must be supplemented on the external epicardium and internal endocardium surfaces of the heart muscle. On the epicardium, a generalized visco--elastic Robin boundary condition can be imposed to take into account a frictionless contact with the pericardium and surrounding tissues \cite{Bucelli2023-mathematical-model}, while the boundary condition on the endocardium arises from the interaction with the blood dynamics inside the four chambers of the heart.

\null\textbf{Modelling blood and valves dynamics.}
The heartbeat is the result of the coordinated interaction of the four chambers, involving two primary phases: the ventricular systole and the ventricular diastole. The ventricular systole includes an isochoric phase and a subsequent ventricular contraction phase, which leads to the opening of aortic and pulmonary valves when ventricular pressure surpasses aortic and pulmonary pressures and facilitates blood ejection into the systemic circulation. The ventricular diastole, or relaxation phase, includes a second isochoric phase due to the closure of all heart valves, followed by the opening of mitral and tricuspid valves, ventricular expansion, and gradual filling. This phase also encompasses atrial systole, contributing to ventricular filling. 
These processes are described using the incompressible Navier--Stokes equations in a moving domain by following the Arbitrary Lagrangian--Eulerian (ALE) framework \cite{Bucelli2023-mathematical-model}, while the Resistive Immersed Implicit Surface (RIIS) method is used to incorporate cardiac valves \cite{Zingaro2022-geometric-multiscale, Zingaro2023-modeling-isovolumetric, Fedele2017-patient-specific, Fumagalli2020-image-based, Fumagalli2022-image-based}. 

More precisely, let $\widehat\Omega_{\text{f}}$ denote the fluid dynamics domain in its reference configuration (including the four chambers of the heart as well as the first tract of the main vessels connected to the heart) and $\widehat\Sigma$ be the interface between the reference fluid domain (chambers) and the structure domain (muscle). 
 First, we define the displacement $\mathbf{d}_{\text{ALE}}=\mathbf{d}_{\text{ALE}}(t)$ of the fluid domain at time $t\in(0,T)$ as the solution of the non--linear equation
\begin{eqnarray*} 
\begin{array}{ll}
        -\nabla\cdot \mathbf{P}_\text{ALE}(\mathbf d_\text{ALE}) = \mathbf 0 &\text{in } \widehat\Omega,
\end{array}
\end{eqnarray*}
in which an elastic constitutive law is used to define the tensor $\mathbf{P}_\text{ALE}$ \cite{Hoffman2011} and which is
supplemented with the adherence condition between blood and muscle  $\mathbf d_\text{ALE} = \mathbf d$  on $\widehat\Sigma$ ($\mathbf d$ is the solution of (\ref{eq:mechanics})).  
Then, denoting by
$\displaystyle
    \mathbf u_\text{ALE} = \partial \mathbf d_\text{ALE}/\partial t$
the rate of deformation of the fluid domain and by $\Omega_{\text{f}}^t$ the fluid dynamics domain at time $t$, obtained as the image of the reference domain $\widehat\Omega_\text{f}$ through the map $\mathcal{L}_\text{f}:(\hat{\mathbf x},t)\to \mathbf{x}=\hat{\mathbf x}+\mathbf{d}_\text{ALE}$, the Navier--Stokes equations
\begin{eqnarray}\label{eq:NS-ALE}
\left\{\begin{array}{ll}
\displaystyle \rho_\text{f}\left[\frac{\partial\mathbf u}{\partial t} + \left(\left(\mathbf u - \mathbf u_\text{ALE}\right) \cdot\nabla\right)\mathbf u\right] - \nabla\cdot\sigma_\text{f}(\mathbf u, p) + \boldsymbol{\mathcal{R}}(\mathbf u, \mathbf u_\text{ALE}) = \mathbf 0 & \text{ in } \Omega_\text{f}^t\times(0,T) \\
\nabla\cdot \mathbf u = 0 & \text{ in } \Omega_\text{f}^t\times(0,T)
\end{array}\right.
\end{eqnarray}
provide the velocity
$\mathbf u$ and the pressure $p$ of the fluid, where $\rho_\text{f}$ is the density of the fluid, and $\sigma_\text{f}$ is the Cauchy stress tensor. 
The term
\begin{eqnarray*}
\boldsymbol{\mathcal{R}(\mathbf u, \mathbf u_\text{ALE})} = \sum_{\mathrm{k} \in \mathcal{V}} \frac{R_\mathrm{k}}{\varepsilon_\mathrm{k}} \delta_{\varepsilon_\mathrm{k}}\left(\varphi_\mathrm{k}^t(\mathbf x)\right)\left(\mathbf u - \mathbf u_\text{ALE} - \mathbf u_{\Gamma_\mathrm{k}}\right)
\end{eqnarray*}
is a resistive term that aims at modelling the effects of valves' dynamics on the blood at the macroscopic level. More precisely, $\mathcal V$ is the set of indices of heart valves, $R_k$ are resistance coefficients,   $\varepsilon_\mathrm{k}$ valve half--thicknesses, $\delta_{\varepsilon_\mathrm{k}}$ smoothed Dirac functions with $\varphi_\mathrm{k}^t(\mathbf x)$ the distance with a sign from the surface $\Gamma_\mathrm{k}$ of the valve, and $\mathbf u_{\Gamma_\mathrm{k}}$ the velocity of the valve with respect to the moving domain $\Omega_\text{f}^t$.

The coupling between the blood and the muscle is completed by imposing the kinematic condition 
$$\mathbf u\circ {\mathcal L}_\text{f}= \frac{\partial \mathbf{d}}{\partial t} \qquad \mbox{ on }\widehat\Sigma,$$
and the dynamic condition
$$\sigma_\text{s}(\mathbf{d},\mathbf{z}_{\text{act}})=\sigma_\text{f}(\mathbf{u},p) \qquad \mbox{ on }\widehat\Sigma,$$
where $\sigma_\text{s}$ is the Cauchy stress tensor related to the first Piola--Kirchhoff tensor $\mathbf{P}$. This last condition models the stress exerted on the myocardium by the blood contained in the chambers and it is used as a boundary condition to close the system (\ref{eq:mechanics}).

Inlet and outlet conditions for the Navier--Stokes equations (\ref{eq:NS-ALE})  can be provided by solving the blood dynamics in the circulatory system external to the heart.

\null\textbf{Modelling blood circulation.}
Several blood circulation models, with diﬀerent degrees of accuracy, have been proposed in the literature. These range from three--dimensional ﬂuid--structure
interaction models \cite{Peskin1977, QDMV2019, Tagliabue-complex-2017, Tagliabue-fluid-2017, Vigmond2008} to zero-dimensional models, also known as lumped parameters models, whose
variables only depend on time, but not on spatial coordinates \cite{Caruel2014, Gerbi-2019-monolithic, BlancoFeijoo2010, Hirschvogel2017, QDMV2019, RSAFDQ2022}.
In general terms, these models, which are derived by applying the principles of conservation of mass and momentum, are written as systems of Ordinary Differential Equations, in the form
\begin{equation}\label{eq:6.1-QDR}
    \frac{d\mathbf{z}_{\textrm{circ}}}{dt}={\boldsymbol\Phi}_{\textrm{circ}}(\mathbf{z}_{\textrm{circ}},t),
\end{equation}
where the vector $\mathbf{z}_{\textrm{circ}}(t)\in{\mathbb R}^{N_{\textrm{circ}}}$ collects the pressure and volumes in different compartments of the circulatory system.

The coupling between the zero--dimensional model of blood circulation and the
three--dimensional model of the heart is achieved by imposing that the volume enclosed by the ventricular cavity, when the domain $\widehat\Omega$ is moved by the displacement $\mathbf{d}$, equals the one predicted by the zero--dimensional circulation model \cite{RSAFDQ2022,Piersanti2022-3d-0d}.

\null
\textbf{Modelling the perfusion.}
The cardiac perfusion process represents blood and oxygen supply via blood flow in the coronary tree, from epicardial coronary arteries to capillaries. The main coronary arteries (those that are visible at the clinical level by imaging) are comparable to the organ scale and are modelled by 3D incompressible Navier-Stokes equations, while the intramyocardial coronaries and microcirculation, which are invisible at the imaging level, are comparable to the cellular scale and are modelled as a porous media by a multi--compartment Darcy model \cite{Hyde2014, DiGregorio2021, Zingaro2023-perfusion}. 

Let $\Omega_\text{C}\subset\mathbb{R}^3$ denote the coronaries domain, which we suppose non--deformable and composed of $J$ main coronaries, $\Omega_\text{M}\subset{\mathbb R}^3$ the myocardial muscle, which we decompose in $J$ non--overlapping subdomains $\Omega_\text{M}^j$ such that each main coronary feeds exactly one of them, and $\Gamma^j$ the outflow boundary of the coronary $j\in\{1,\ldots, J\}$. The multi--compartment Darcy model \cite{Cookson2012} consists of three coexisting Darcy equations in the whole myocardium $\Omega_\text{M}$, each of one characterized by a different permeability tensor $\mathbf{K}_\text{i}$ and corresponding to a different length scale of capillaries. The first compartment is  the one upstream and is characterized by a volumetric source $g_1$ that, by enforcing mass conservation, should be provided by the outgoing coronary flow rate:
$$g_1=\sum_{j=1}^J \frac{\chi_{\Omega_{\text M}^j}}{|\Omega_{\text M}^j|} 
\int_{\Gamma_j}\mathbf{u}\cdot\mathbf{n},$$
where $\chi_{\Omega_{\text M}^j}$ denotes the characteristic function of the subregion $\Omega_{\text M}^j$. 
Always by enforcing mass conservation, the intermediate compartments feature null sources and sinks, while the farthest away compartment (modelling the microvascolature) shows a sink volumetric term that models the coronary venous return. By imposing a prescribed flow rate at the coronaries inflow and defective outflow conditions on $\Gamma^j$ for the Navier-Stokes equations \cite{Vergara2011-defective}, the coupled system reads:
\begin{eqnarray*}
\left\{\begin{array}{ll}
\displaystyle \rho_{\text f}\left(\frac{\partial \mathbf u}{\partial t}+(\mathbf u \cdot \nabla)\mathbf u\right)+\sigma_{\text f}(\mathbf{u},p)=\mathbf{0} & \mbox{ in }\Omega_{\text C}\\[1mm]
\displaystyle (\sigma_{\text f}(\mathbf{u},p)\mathbf{n})\cdot \mathbf n +
\frac{1}{\alpha^j}\int_{\Gamma^j}\mathbf u \cdot \mathbf n=
\frac{1}{|\Omega_\text{M}^j|}\int_{\Omega_\text{M}^j}p_\text{M,1} & \mbox{ on }\Gamma^j\\[1mm]
\mathbf{u}_{\text {M,i}}+\mathbf K_\text{i}\nabla p_{\text{M,i}}=\mathbf{0},\qquad i=1,2,3 &  \mbox{ in }\Omega_{\text M}\\
\nabla\cdot \mathbf{u}_\text{M,1}=g_1 -\beta_{1,2}(p_\text{M,1}-p_\text{M,2}) & \mbox{ in }\Omega_{\text M}\\[1mm]
\displaystyle \nabla\cdot \mathbf{u}_\text{M,2} =
-\beta_{2,1}(p_\text{M,2}-p_\text{M,1})
-\beta_{2,3}(p_\text{M,2}-p_\text{M,3}) & \mbox{ in }\Omega_\text{M}\\[1mm]
\displaystyle \nabla\cdot \mathbf{u}_\text{M,3} =
-\gamma(p_\text{M,3}-p_\text{veins})
-\beta_{2,3}(p_\text{M,3}-p_\text{M,2}) & \mbox{ in }\Omega_\text{M},\\[1mm]
\mathbf u_\text{M,i}\cdot\mathbf n =0,\qquad i=1,2,3 &\mbox{ on }\partial\Omega_{\text M},
\end{array}\right.
\end{eqnarray*}
where $\mathbf u$ and $p$ are the velocity and pressure of the blood in the coronaries, $\mathbf{u}_{\text {M,i}}$ and $p_{\text{M,i}}$ are Darcy velocities and pressures in the compartments $\text{i}\in\{1,2,3\}$. The coefficients $\alpha^j$ are the conductances (supposed to be dependent on the perfusion region), $\beta_{ik}\geq 0$ represent the inter--compartment pressure--coupling coeﬃcients, and in the last equation we have accounted for the coronary venous return through the venous pressure $p_\text{veins}$, with $\gamma$ a suitable coeﬃcient \cite{DiGregorio2021}.

\null\textbf{Modelling torso potential.}
Physics--based computational models of the cardiac function must accurately simulate the patient's electrophysiology to offer predictive support for clinical use. Validation of the calibrated model on a specific patient also includes comparing the calculated electrocardiogram (ECG) and the measured ECG.
For this purpose, it is mandatory to reproduce the ECG waveform to a high--quality standard and such a goal can be achieved by considering myocardial displacement in ECG generation, i.e., by combining a cardiac electromechanical model with a simulation of the action potential in the torso, and even, with the simulation of torso deformation induced by the myocardial displacement \cite{Boulakia2010, Gillette2023, Zappon2024}.

Modelling the generation of electrophysiological clinical output typically involves an electrophysiology solver (like, e.g., the bidomain model (\ref{eq:EP})) and a model, like an elliptic equation, to calculate the ECG starting from the electrical signal propagation through the torso. Let us denote by $\Omega_{\text{H}}\subset{\mathbb R}^3$ the computational domain corresponding to the heart muscle, by $\Omega_\text{T}$ the surrounding torso domain, and by $\Gamma=\partial\Omega_\text{H}\cap\partial\Omega_\text{T}$ their common interface. 
Denoting by $\mathbf D_\text{T}$ the isotropic diffusion tensor in the torso, at each time $t\in(0,T)$ the electric potential $v_\text{T}$ in the torso is computed by solving the coupled system
\begin{eqnarray}\label{eq:ep-torso}
\left\{\begin{array}{ll}
\text{electrophysiology  model (\ref{eq:EP})} &\text{in }\Omega_H\\
-\nabla \cdot(\mathbf D_\text{T}\nabla v_\text{T}) = 0 &\text{in }\Omega_{\text{T}}\\
v_\text{T} = v_\text{e} &\text{on }\Gamma\\
(\mathbf D_\text{T}\nabla v_\text{T}) \cdot \mathbf{n} = (\mathbf D_\text{e}\nabla v_\text{e})\cdot\mathbf{n} &\text{on }\Gamma,
\end{array}\right.
\end{eqnarray}
where $v_\text{e}$ is the extracellular potential in the heart muscle and $\mathbf D_\text{e}$ the extracellular diffusivity tensor introduced in (\ref{eq:EP}). Whether the displacement of the torso should be considered, (\ref{eq:ep-torso}) should be enriched with an elastic model that provides the displacement $\mathbf{d}_\text{T}$ of the torso, which, in turns, would allow computing the corresponding deformation gradient $\mathbf F_\text{T}$ (with $J_\text{T}=\text{det}\mathbf F_\text{T}$) to replace $\mathbf D_\text{T}$ in (\ref{eq:ep-torso}) with 
$J_\text{T}\mathbf F_\text{T}^{-1}\mathbf D_\text{T}\mathbf F_\text{T}^{-T}$ \cite{Zappon2024}.


\subsection{Multifidelity PINNs for the estimation of ionic parameters}
\label{sec:mpinn-ionic}

The theoretical framework of multifidelity PINNs, introduced in Section~\ref{sec:PINN}, provides a robust methodology for parameter estimation in systems governed by Ordinary Differential Equations (ODEs). This approach is particularly relevant in cardiac electrophysiology, where estimating ionic parameters is crucial for accurately modelling the electrical behaviour of cardiac cells. In this subsection, we present an application of multifidelity PINNs to the Bueno-Orovio (BO) ionic model, which captures the dynamics of cardiac cellular electrophysiology \cite{regazzoni2021physics}.

\null\textbf{The Problem.}
\newcommand{\taufi}{\tau_{\text{fi}}}
The BO model describes the temporal evolution of the transmembrane potential $u(t)$ and the gating variables $v(t)$, $w(t)$, and $s(t)$ through a system of non--linear ODEs:
\begin{equation} \label{eqn:BuenoOrovio}
    \begin{aligned}
        &\dfrac{du}{dt} = - ( {I}_{\mathrm{fi}} + {I}_{\mathrm{so}} + {I}_{\mathrm{si}}) + {I}_{\mathrm{app}} \\
        &\dfrac{dv}{dt} = \dfrac{(1 - H(u - \theta_v)) (v_\infty - v)}{\tau_v^-} - \dfrac{H(u - \theta_v) v}{\tau_v^+}  \\
        &\dfrac{dw}{dt} = \dfrac{(1 - H(u - \theta_w)) (w_\infty - w)}{\tau_w^-} - \dfrac{H(u - \theta_w) w}{\tau_w^+}  \\
        &\dfrac{ds}{dt} = \dfrac{1}{\tau_s} \left( \dfrac{1 + \tanh(k_s (u - u_s))}{2}-s \right). \\
        \end{aligned}
\end{equation}
where the ionic currents are defined as:
\begin{equation}
    \begin{aligned}
    {I}_{\mathrm{fi}} &= -\frac{v H(u - \theta_v) (u - \theta_v) (u_u - u)}{\taufi},\\
    {I}_{\mathrm{so}} &= \frac{(u - u_o) (1 - H(u - \theta_w))}{\tau_o} + \frac{H(u - \theta_w)}{\tau_{\mathrm{so}}},\\
    {I}_{\mathrm{si}} &= -\frac{H(u - \theta_w) w s}{\tau_{\mathrm{si}}}.
\end{aligned}
\end{equation}
Here, $H(\cdot)$ is the Heaviside function, and the parameters $\taufi$, $\tau_v^\pm$, $\tau_w^\pm$, $\tau_s$, and others govern the dynamics of the ionic currents and gating variables. Accurately estimating $\taufi$ is critical, as it directly influences the action potential shape and duration.
Note that, in this context, the dynamics depend solely on time, as spatial variability is not considered. 

\null\textbf{Data generation.}
To generate the training dataset (whose input is the independent variable $t$ and the outputs are the values of both the transmembrane potential and gating variables at time $t$), we employ a numerical solver based on the Finite Difference method to approximate the solution of the BO model. The solver uses a time step of $\Delta t = 0.1$ ms and simulates up to a final time of $T = 0.8$ s. The resulting transmembrane potential $u(t)$ is then subsampled at intervals of 25 ms, and Gaussian noise with zero mean and standard deviation $\sigma$ is artificially added to the data. In addition to measurements of the transmembrane potential, we include data for the gating variables $v(t)$, $w(t)$, and $s(t)$. However, since acquiring these measurements is typically more challenging in practice, they are sampled much less frequently, with one data point collected every 0.3 s.
Further information on the dataset generation are available in \cite{Regazzoni-lincei-2021}.

\null\textbf{Methodology.}
We address the problem of identifying the parameter $\taufi$ starting from measurements of the transmembrane potential $u(t)$ and the gating variables $v(t)$, $w(t)$, and $s(t)$. According to the notation of Sec.~\ref{sec:PINN}, we have $P = 1$, that is a single unknown physical parameter, namely $\gamma = \taufi$.
We solve this parameter estimation problem using both the standard PINN approach (Eq.~\eqref{eqn:PINN-inverse}) and the MPINN approach (Eq.~\eqref{eqn:MPINN-enhanced}). 
As the differential problem at hand has four unknowns ($u,v,w,s$), we employ a fully connected neural network with four output neurons, each corresponding to an unknown. The input consists instead of the time variable $t$, namely the independent variable of the differential problem \eqref{eqn:BuenoOrovio}.
Hence, we will denote by 
$u_{NN}(t;\mathbf w)$,
$v_{NN}(t;\mathbf w)$,
$w_{NN}(t;\mathbf w)$ and
$s_{NN}(t;\mathbf w)$
the four output neurons of the fully connected NN, where $\mathbf w$ represents the network trainable parameters.
Having four unknowns, the data-fitting term in the loss function is also given by the sum of four terms:
\begin{eqnarray*}
    \mathcal L_{obs}(\mathbf w)
    &=&\displaystyle
    \frac{1}{N_{obs, u}}\sum_{i=1}^{N_{obs, u}}|u_i^{obs}-u_{NN}(t_i^{obs, u};\mathbf w)|^2
    \\&&+\displaystyle
    \frac{1}{N_{obs, gat}}\sum_{i=1}^{N_{obs, gat}}|v_i^{obs}-v_{NN}(t_i^{obs, gat};\mathbf w)|^2 
    \\&&+\displaystyle
    \frac{1}{N_{obs, gat}}\sum_{i=1}^{N_{obs, gat}}|w_i^{obs}-w_{NN}(t_i^{obs, gat};\mathbf w)|^2 
    \\&&+\displaystyle
    \frac{1}{N_{obs, gat}}\sum_{i=1}^{N_{obs, gat}}|s_i^{obs}-s_{NN}(t_i^{obs, gat};\mathbf w)|^2  
\end{eqnarray*}
where
$t_i^{obs, u}$ for $i = 1, \dots, N_{obs, u}$
are the time instants where observations of the transmembrane potential are available, while 
$t_i^{obs, gat}$ for $i = 1, \dots, N_{obs, gat}$
are the times where observations of the gating variables are available.
Similarly, the physics--informed term $\mathcal L_{phys}(\mathbf w,\gamma)$  incorporates the residuals of all the equations of the system \eqref{eqn:BuenoOrovio}, while the BC contribution is missing.
For additional details, we refer the interested reader to \cite{regazzoni2021physics}.

For the PINN setup, we use a feed--forward neural network (FFNN) with three hidden layers, containing 32, 24, and 16 neurons, respectively. In the MPINN framework, the high-fidelity correction network $\mathcal{NN}_H$ adopts the same architecture, while the low--fidelity parametric emulator uses three hidden layers with 32, 16, and 8 neurons. Note that the configuration of decreasing layer sizes, commonly used in practice as an alternative to constant--sized layers or layers that first increase and then decrease, is not critical in this context, as the results would not vary significantly. The low--fidelity dataset consists of 75 numerical simulations, subsampled at 10 ms intervals.

Given the non--linear nature of the BO model, a one-shot training approach starting from a random initialization of the NN parameters often leads to solutions far from the global minimum. If the initial guess for the solution is significantly incorrect, the physics--based term in the loss function can inadvertently drive the parameter estimation away from its true value. 

To fix this shortcoming, both the PINN and MPINN methods employ a staged training strategy, adjusting the relative contributions of the loss components over different phases.
Initially, we perform 500 epochs with $\alpha_{PDE} = 10^{-6}$,
effectively treating the problem as a data--fitting task while using a small physics weight to prevent the residuals from diverging. Once the data is reasonably well-fitted, we increase the contribution of the physics-informed term, setting $\alpha_{PDE} = 10^{-3}$, and perform 10 000 additional training epochs. For the MPINN method, an intermediate step is introduced between these phases: 500 iterations are dedicated solely to optimizing the unknown parameter $\tau_\mathrm{fi}$, leveraging the low-fidelity parametric map. This intermediate optimization step is computationally inexpensive, as it involves adjusting a single variable, and provides a more informed initial guess for the parameter.

\null
\textbf{Results and Discussion.}
To evaluate the performance of the PINN and MPINN methods, we conduct 10 independent training runs for each approach, using different random initializations of the NN parameters. The experiments are repeated for two noise levels: $\sigma = 0.05$ and $\sigma = 0.025$. Table~\ref{tab:MPINN:comparison} reports the relative errors in estimating $\tau_\mathrm{fi}$ for the two methods under both noise conditions. 

\begin{table}[h!]
\centering
\begin{tabular}{l|rr}
\textbf{Noise level ($\sigma$)} & \textbf{0.05} & \textbf{0.025} \\
\hline
\textbf{Error (PINN)} & 0.216 & 0.108 \\
\textbf{Error (MPINN)} & 0.013 & 0.005 \\
\end{tabular}
\caption{Comparison of average relative errors in the estimation of $\tau_\mathrm{fi}$ between PINN and MPINN methods for different noise levels.}
\label{tab:MPINN:comparison}
\end{table}

The results clearly demonstrate the superior performance of the MPINN approach, which achieves significantly lower estimation errors {for} both noise levels. By leveraging the low-fidelity model as a prior, the MPINN method not only improves accuracy but also accelerates convergence compared to standard PINNs.

This study highlights the potential of multifidelity PINNs in parameter estimation problems, particularly in scenarios with noisy and sparse data. The ability to integrate information from multiple fidelities, and to leverage the prior given by a surrogate model, offers a robust approach to tackling highly non--linear parameter estimation problems, such as those encountered in cardiac modelling.

\subsection{Physics--aware NNs for the inverse problem of electrocardiography}

Electrocardiographic imaging (ECGI) is a novel technique to measure cardiac electric activity by exploiting both body surface signals and thoracic CT--scans. ECGI is more accurate than solely plain body surface signals like electrocardiograms (ECG) or body surface potential maps. In particular, the imaging component allows us to determine the patient--specific organs' geometries.

\begin{figure}
     \centering
     \includegraphics[width=0.25\textwidth] {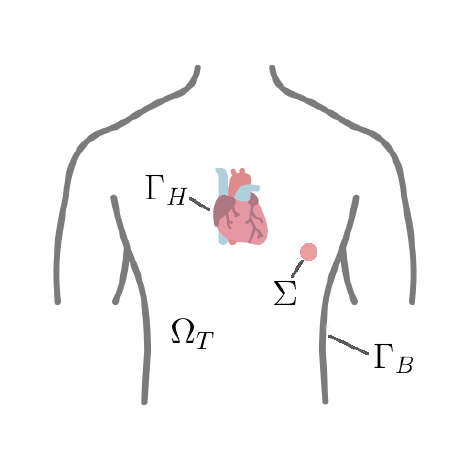}\quad
     \includegraphics[width=0.7\textwidth] {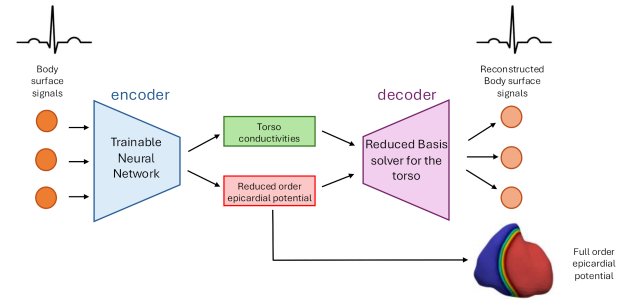}\quad
     \caption{On the left, the computational domain of the inverse problem of electrocardiography (Image elaborated starting from icons by https://www.flaticon.com/authors/bsd and https://www.flaticon.com/authors/smashicons). On the right, scheme of the architecture of ST--RB--DNN}
     \label{fig:ecgi-torso}
 \end{figure}

\null\textbf{The problem.}
The mathematical model of ECGI is the so--called inverse problem of electrocardiography, its final goal is to look for the epicardial potential field $v$ which generates body surface signals $y(v)$ as close as possible to some target ones $z$ in the least--squares sense. More precisely, let $\Omega_T\subset \mathbb R^3$ be the computational domain corresponding to the torso, $\Gamma_B$ the body surface, $\Sigma\subset \Gamma_B$ the portion of the body surface where measurements $z$ are available, and $\Gamma_H$ the epicardium surface (see Fig. \ref{fig:ecgi-torso}). Given the epicardial potential $v\in H^{1/2}(\Gamma_H)$ at a certain time $t\in(0,T)$, let $y(v)\in H^1(\Omega_T)$ be the torso potential, i.e., the solution of the forward elliptic problem
 \begin{eqnarray}\label{eq:ecgi-forward}
 \left\{\begin{array}{ll}
 -\nabla\cdot( \mathbf D_T\nabla y(v))=0 & \mbox{ in }\Omega_T\\
 y(v)=v & \mbox{ on } \Gamma_H\\
 \nabla y(v)\cdot \mathbf n_B=0 & \mbox{ on }\Gamma_B,
 \end{array}\right.
 \end{eqnarray}
 where $\mathbf D_T$ is the electrical conductivity tensor in $\Omega_T$, and $\mathbf n_B$ is the outward normal unit vector to $\Gamma_B$.
 The inverse problem of electrocardiography reads: given a function $z\in L^2(\Sigma)$ representing the measured electric potential on the body surface, look for the epicardial potential $u_H\in H^{1/2}(\Gamma_H)$ such that
 \begin{equation}\label{eq:ecgi-inverse}
 u_H=\argmin{v\in H^{1/2}(\Gamma_H)}\left[\frac{1}{2}\int_{\Sigma}|y(v)-z|^2d\sigma+R(v)\right]
 \end{equation}
and $y(v)$ is the solution of (\ref{eq:ecgi-forward}).
$R(v)$ is a suitable regularization term which ensures the problem (\ref{eq:ecgi-inverse}) is well--posed and includes a penalization coefficient $\alpha$ which is hard to tune. $\alpha$ should be large enough to guarantee that $R(v)$ is effective, but at the same time, too large values provide inaccurate solutions.
Moreover, solving (\ref{eq:ecgi-forward}) requires evaluating the electric conductivity tensor $\mathbf D_T$ which can be inferred only through thoracic CT--scans, implying radiations on patients.

\null\textbf{Methodology.}
In \cite{Tenderini2022}, the authors designed a physically informed deep learning model of autoencoder type, named \emph{Space Time--Reduced Basis--Deep Neural Network (ST--RB--DNN)}, as an alternative to solving (\ref{eq:ecgi-forward})--(\ref{eq:ecgi-inverse}). The \emph{input} is constructed starting from body surface signals $z$ on $\Sigma\times(0,T)$, e.g., the ECG. The \emph{encoder} is a FFNN which takes body surface signals in input and produces two kinds of outputs. The first output is a reduced order spatio--temporal representation of the epicardial potential $v$ on $\Gamma_H\times(0,T)$. It is named latent potential and is responsible for the signal $z$. The second output is the set of the electric conductivities $\mathbf D_T$ of the different organs considered in the torso.
The \emph{decoder} is a deterministic (not trainable) tensorial reduced order model which takes electric conductivities, as well as the latent spatio--temporal representation of the extracellular potential $v_e$ as Dirichlet datum on $\Gamma_H$, expands the space--time reduced epicardial potentials along the time coordinate, and solves the problem (\ref{eq:ecgi-forward}) for any $t\in(0,T)$ to compute the reduced--in--space and full--order--in--time torso potential. Then it reconstructs the body surface signals $y(v)|_{\Sigma\times(0,T)}$ .

We note that the output of the whole ST--RB--DNN model is not limited to body surface signals, but it also contains the reconstructed full order epicardial potential field.

The first way of training ST--RB--DNN treats it as a pure autoencoder (AE) architecture. Denoting by $\mathbf w$ the array containing the trainable parameters of the encoder, the loss function is defined as
\begin{eqnarray}\label{eq:ecgi-loss-pureAE}
\mathcal L^{AE}(\mathbf w)=\mathcal L_{sig}(\mathbf w)+
\mathcal L_{reg}(\mathbf w)
\end{eqnarray}
where
$\mathcal L_{sig}$ is the Mean Square Error (MSE) (see Sect. \ref{sec:performance}) between the measured body surface signals and those computed by the autoencoder, and $\mathcal L_{reg}$ is a regularization term which helps in preventing overfitting. The training set is the collection of body surface signals.

Alternatively, ST--RB--DNN can be trained by minimizing the loss function
\begin{equation}\label{eq:ecgi-loss}
\mathcal L(\mathbf w)=w_{BC}\mathcal L_{BC}(\mathbf w)+\mathcal L_{sig}(\mathbf w)+
\mathcal L_{reg}(\mathbf w),
\end{equation}
where $\mathcal L_{BC}$ is the Mean Absolute Error (MAE) (see Sect. \ref{sec:performance}) between the target space--time reduced epicardial potential field and the approximated one provided by the NN, while the parameter $w_{BC}\in\mathbb R^+$ represents a weight. Now, the training set includes both body surface signals and the corresponding epicardial potential fields. Because acquiring epicardial maps from clinical data requires invasive procedures, training and test datasets can be built \emph{in silico}, i.e., by computing them numerically as the solution of a mathematical model.

In \cite{Tenderini2022}, after introducing scalar parameters characterizing the problem setting, epicardial potentials have been computed by solving the parametric form of the bidomain equation (\ref{eq:EP}). Parameters have been used to set initial activation patterns and cardiac conductivities, while the geometry of the heart was fixed and the heart has been assumed to be isolated from the torso. Considering the parametric nature of the problem, the authors employed the spatio--temporal Reduced Basis solver proposed in \cite{Choi2019} which allows encoding space--time--dependent fields into a very low number of coefficients, almost independent of the grid refinement along both the spatial and temporal coordinates.

The isolated heart assumption expressed by the homogenous Neumann condition on the epicardium offers the advantage of getting a one--way coupling between the bidomain equations and the torso problem (\ref{eq:ecgi-forward}). This means that currents are not required to be continuous at the epicardial surface and torso potential has no effect on the heart electrical activity. Although not physically supported, this assumption has negligible impact on the activation pattern at the epicardium.

Two different ways have been considered to organize input body surface signals: in the form of time series and as low--frequency Discrete Fourier Transform (DFT) coefficients. 
Consequently, the NN architecture of the encoder has been designed differently in the two cases. 
Notice that time--series--based ST--RB--DNN model can only work at a fixed acquisition frequency, while DFT--based model does not have this limitation.

\null
\textbf{Results and discussion.}
Numerical results presented in \cite{Tenderini2022} refer to a benchmark test case where both the train and test datasets have been generated numerically. The bidomain equation, coupled with the Aliev Panvilov ionic model \cite{AlievPanfilov1996}, has been solved on a reference biventricular geometry, while the torso has been modelled as a homogeneous and isotropic medium. Snapshots for Reduced Basis spaces have been computed by the finite element method in space and finite differences in time.

400 data points were employed to form the dataset. 
Hyperparameters of the encoder and the weights in the loss function were selected using a grid search process (see Sect. \ref{sec:hyperparameters}).
Numerical results show that, after selecting the best values of hyperparameters, the first training procedures (\ref{eq:ecgi-loss-pureAE}) provides slightly lower errors on the 12--lead ECG signals compared to the second training strategy (\ref{eq:ecgi-loss}).  On the contrary, the latter is more accurate in predicting epicardial activation maps.  
When the model is trained with 12--lead ECG signals, the signal preprocessing via DFT improves model performances. 
More details on numerical results are given in \cite{Tenderini2022}.

\subsection{Learning the microscopic dynamics in the framework of a coupled multiscale problem} \label{sec:mor-sarcomeres}
We now illustrate how operator learning methods for time--dependent problems, described in Sec.~\ref{sec:operator-learning-time-dependent}, can be used in a multiscale context where a PDE model at the macroscale is coupled to an ODE model at the microscale. The microscopic model must be solved a large number of times (typically, once for each node of the mesh describing the macroscopic scale), besides being characterized by a large number of unknowns and very fine time scales. This results in a high computational cost, often unsustainable.

This represents a typical situation where operator learning can be used to learn efficient surrogates of the microscopic model which will be then coupled with the original, high-fidelity model. This hybrid model faithfully captures the multiscale dynamics while being computationally efficient.

\null\textbf{The problem.}
The specific example that we consider {is} the multiscale problem of cardiac electromechanics described in Sec.~\ref{sec:IHM}. In this case, the macroscopic dynamics are described by the electrophysiology model \eqref{eq:EP} coupled with the mechanics model \eqref{eq:mechanics}. The microscopic dynamics are instead described by an active force generation model, which accounts for subcellular biochemical mechanisms, a system of ODEs that involves a large number of unknowns and very fine time scales. For the reader's convenience, we report here the general expression of the active force generation model, which can be written in the form
\begin{equation}\label{eq:force_generation_FOM}
    \left\{
        \begin{aligned}
            & \frac{d\actState}{dt}(t)= \actRHS\left(\actState(t), \Cai(t), SL(t), \frac{d\,SL}{dt}(t)\right)
            & & \qquad t \in (0,T] \\
            & \actTension(t) = \actOBS(\actState(t))
            & & \qquad t \in (0,T] \\
            & \actState(0) = \mathbf{z}_{\textrm{act}, 0}
        \end{aligned}
    \right.
\end{equation}
where $\actState$ is the vector that collects the state variables associated with the dynamics of the contractile proteins, $\Cai$ is the local calcium concentration, $SL$ is the sarcomere length, while $\actTension$ denotes the active tension generated, and represents the output of interest of the microscopic model, which allows it to be linked to the macroscopic dynamics.
Biophysically detailed force generation models can have many state variables and require a very fine time step to capture the rapid dynamics associated with biochemical mechanisms.
For example, the RDQ18 model \cite{RDQ-active-contraction-2018} is characterized by 2176 state variables and a characteristic time step of $\SI{2.5e-5}{\second}$. This means that if the macroscopic mesh is composed of, say, $10^6$ nodes, the number of variables used to describe the microscopic dynamics would be on the order of $10^9$.

\null\textbf{Methodology.}
To reduce the computational cost associated with approximating a multiscale cardiac model where the fine scale is described by models of the type \eqref{eq:force_generation_FOM}, in \cite{RDQ-ML-2020} it was proposed to use an operator learning method that allows for the construction of an efficient surrogate of the microscale model, trained from a dataset of precomputed simulations using numerical approximation of the FOM \eqref{eq:force_generation_FOM}. Note that, since we only want to surrogate the microscale model here, the generation of the training data can be done by considering the model \eqref{eq:force_generation_FOM} uncoupled from the macroscale model. In this way, the training dataset can be generated with a much lower computational effort compared to the cost of a full multiscale simulation.

The training dataset consists of a set of $\dynNumSamples$ trajectories $(\widehat\dynInp^j(t), \widehat\actTension^j(t))$, $j=1,\dots,\dynNumSamples$, where the input $\dynInp(t) \in \mathbb{R}^2$ is composed of $\Cai(t)$ and $SL(t)$, since the RDQ18 model considered here does not explicitly depend on $\frac{d\,SL}{dt}$ \cite{RDQ-active-contraction-2018}. 
All samples of the training dataset start from equilibrium initial conditions, associated with an input value $\dynInp_0$ corresponding to the conditions present at the beginning of the cardiac cycle (tele-diastolic condition). 
To generate the training dataset, we consider: 50 double step inputs, from $\dynInp_0$ to a random value, and then back to $\dynInp_0$; 45 oscillatory inputs, where range and frequency are randomly sampled; 60 random walk inputs, for a total of $\dynNumSamples = 155$ training samples (see \cite{RDQ-ML-2020} for more details). These inputs are chosen to cover different response ranges of the system to different types of inputs. 
The trajectories are sampled at a time step of $\Delta t = \SI{1e-2}{\second}$, defining the time instants $t_i = i\Delta t$, $i=0,\dots,\dynNumTimes^j$, where $\dynNumTimes^j = T_j / \Delta t$ is the number of time instants for the $j$-th sample and $T_j$ the corresponding duration.

Due to the strong nonlinearity of the FOM dynamics, projection-based reduced order modelling methods are ineffective in this case, as shown in \cite{Regazzoni-phd}. Therefore, to achieve a drastic dimensionality reduction of the model, we use an operator learning approach based on model learning with latent variables, as described in Sec.~\ref{sec:operator-learning_hidden-dynamics-discovery}. In particular, the surrogate model is constructed as:
\begin{equation}\label{eq:force_generation_ROM}
    \left\{
        \begin{aligned}
            & \frac{d\dynLat}{dt}(t) = \dynNNdyn\left(\dynLat(t), \dynInp(t); \dynWdyn \right)
            & & \qquad t \in (0,T] \\
            & \actTension(t) = \dynLat(t) \cdot \mathbf{e}_1
            & & \qquad t \in (0,T] \\
            & \dynLat(0) = \dynLat_0
        \end{aligned}
    \right.
\end{equation}
where $\dynLat(t) \in \mathbb{R}^{\dynNumLat}$ represents the latent state, whose dimension $\dynNumLat$ is tuned together with the other hyperparameters, and $\dynNNdyn$ represents a neural network with trainable parameters $\dynWdyn$. The term $\mathbf{e}_1$ denotes the first vector of the canonical basis of $\mathbb{R}^{\dynNumLat}$, ensuring that the output $\actTension$ of the surrogate model coincides with the first component of the latent state. This approach, referred to as \textit{output-inside-the-state} in \cite{RDQ-JCP-2019}, is particularly useful when the number of state variables is low (as in this case, where there is only one output variable), and it allows for training a single neural network instead of two, as required in the standard formulation \eqref{eqn:dynODE_ROM_latent}. 

The initial condition of the latent variable is set to $\dynLat_0 = (\bar\actTension, 0, \dots, 0)^T$, where $\bar\actTension = \actOBS(\mathbf{z}_{\textrm{act}, 0})$ represents the active tension associated with the tele--diastolic condition from which the samples originate.

For training the surrogate model, we use a semi-physical approach, leveraging the physical knowledge of the model to guide the training of the neural network. In particular, the loss function used for training the surrogate model consists of two terms. First, we have a data-fitting term that penalizes the discrepancy between the output of the surrogate model and the output of the FOM:
\begin{equation*}
    \mathcal L_{data}(\dynWdyn) = \frac{1}{\dynNumSamples}\sum_{j=1}^{\dynNumSamples}\frac{1}{\dynNumTimes^j}\sum_{i = 0}^{\dynNumTimes^j} \left| \widehat\actTension^j(t_i) - \dynLat(t_i; \dynWdyn) \cdot \mathbf{e}_1 \right|^2.
\end{equation*}
To this term, we add a term that forces the latent state variable to return to the initial condition $\dynLat_0$ at the end of the trajectories for which the FOM also returns to the tele--diastolic condition. Denoting by $J_{\text{cycle}}$ the set of samples for which the trajectory returns to the tele--diastolic condition (associated with the 50 double-step inputs mentioned above), the additional loss term is defined as: 
\begin{equation*}
    \mathcal L_{cycle}(\dynWdyn) = \frac{1}{| J_{\text{cycle}}|}\sum_{j \in J_{\text{cycle}}}
    \frac{ \| \dynLat(T_j) - \dynLat_0 \|^2}
    {\frac{1}{\dynNumTimes^j} \sum_{i = 0}^{\dynNumTimes^j} \| \dynLat(t_i) - \dynLat_0 \|^2}
\end{equation*}
This term enhances the stability of the surrogate model over a large number of cycles, as is of interest in the context of cardiac mechanics, and prevents the neural network-based model from deviating too far from the FOM dynamics. 

Finally, we introduce another bias into the learning process, namely that the system's initial condition constitutes an equilibrium point. Formally, this translates for the system \eqref{eq:force_generation_ROM} into the constraint $\dynNNdyn(\dynLat_0, \dynInp_0; \dynWdyn) = \mathbf{0}$ that can be enforced through an additional loss term, of the type 
\begin{equation*}
    \mathcal L_{eq}(\dynWdyn) = 
    \|\dynNNdyn(\dynLat_0, \dynInp_0; \dynWdyn)\|^2.
\end{equation*}
Alternatively to this weak enforcement of the constraint, we can modify the model architecture to strongly impose the constraint directly in the neural network, replacing the right-hand side of \eqref{eq:force_generation_ROM} by
\begin{equation}\label{equilibrium_strong_imposition}
    \dynNNdyn(\dynLat, \dynInp; \dynWdyn) =  \widetilde\NNgeneric_{\text{dyn}}(\dynLat, \dynInp; \dynWdyn) - \widetilde\NNgeneric_{\text{dyn}}(\dynLat_0, \dynInp_0; \dynWdyn),
\end{equation}
where $\widetilde\NNgeneric_{\text{dyn}}$ represents a neural network with trainable parameters $\dynWdyn$. This approach forces the output to coincide with $\mathbf{0}$ when the input is $\dynLat_0$ and $\dynInp_0$. 

In conclusion, the training consists of the following minimization problem:
\begin{equation*}
    \dynWdyn^* = \argmin{\dynWdyn} \left[
    \mathcal L_{data}(\dynWdyn) 
    + \Wratio \mathcal L_{cycle}(\dynWdyn)
    + \Weq \mathcal L_{eq}(\dynWdyn)
    \right],
\end{equation*}
where $\Wratio$ and $\Weq$ are hyperparameters that balance the contributions of the different loss terms, whose tuning is done along with the other hyperparameters. The equilibrium condition can be imposed weakly by setting $\Weq > 0$, or strongly as in \eqref{equilibrium_strong_imposition} and setting $\Weq = 0$.

\null\textbf{Results and discussion.}
The numerical tests presented in \cite{RDQ-ML-2020} show that the operator learning approach can accurately and reliably learn the dynamics of the microscale model with only $\dynNumLat = 2$ variables in the reduced model and 2 hidden layers with 6 neurons each. Once the training is completed, the surrogate model is validated by comparing its predictions with those of the FOM model on a set of tests, including different types of tests from those used for training, such as steady-state tests, isometric and shortening twitches with physiological calcium transients, and long-term simulations to test the stability of the surrogate model over a large number of cycles. The results show that the surrogate model can faithfully capture the dynamics of the FOM model and generalize well on test data not used during training, with an average error of less than $2\%$ compared to the FOM model. Remarkably, the surrogate model allows for numerical resolution over 1000 times faster than the corresponding FOM model.

An ablation test conducted in \cite{RDQ-ML-2020} highlights the effect of the different loss terms and the choice between the weak and strong imposition of the equilibrium condition on the generalization ability of the surrogate model (see Tab.~\ref{tab:MOR-sarcomeres_errors}). The results show that adding both terms enforcing prior knowledge (cycle condition and equilibrium condition) in the loss function yields a more stable and accurate surrogate model than the one trained only with the data-fitting term, and the best results are obtained with the strong imposition of the equilibrium condition.

\newcommand{\e}[1]{\times 10^{#1}} 
\begin{table}
    \begin{tabularx}{\textwidth}{llrr}
        \toprule
        Equilibrium condition & Cycle condition & Training error & Test error \\
        \midrule
        -                       & -                   & $1.62\e{-2}$ & $2.66\e{-2}$ \\
        weak ($\Weq = 10^{-1}$) & -                   & $1.52\e{-2}$ & $2.10\e{-2}$ \\
        strong                  & -                   & $1.70\e{-2}$ & $3.10\e{-2}$ \\
        weak ($\Weq = 10^{-1}$) & $\Wratio = 10^{-1}$ & $1.48\e{-2}$ & $2.35\e{-2}$ \\
        strong                  & $\Wratio = 10^{-1}$ & $1.44\e{-2}$ & $1.97\e{-2}$ \\
        \bottomrule  
    \end{tabularx}  
    \caption{Training and test relative errors obtained by training the surrogate model \eqref{eq:force_generation_ROM} with or without imposition of the equilibrium condition (in either weak or strong form) and with or without imposition of the cycle condition.}
    \label{tab:MOR-sarcomeres_errors}
\end{table}

\begin{figure}
    \centering
    \includegraphics[width=0.9\textwidth]{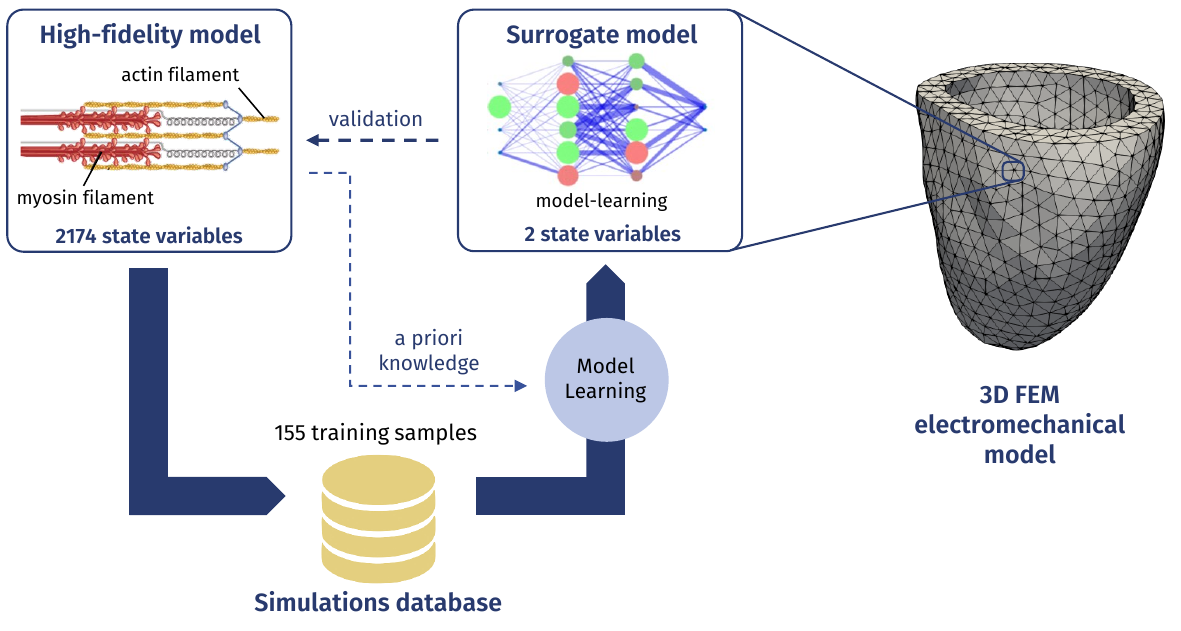}
    \caption{Using model learning with latent variables to reduce the computational cost associated with multiscale cardiac electromechanics models. From the high-fidelity model, a dataset of simulations is constructed, which is then used -- along with prior knowledge of the problem's physics -- to learn a surrogate model. Once validated against the FOM, the model is deployed at each node of the computational mesh to obtain a hybrid multiscale model based on classical FEM modelling techniques and neural network-based models.}
    \label{fig:MOR-sarcomeres}
\end{figure}
The surrogate model can be used to calculate the active tension generated at each node of the macroscale mesh, significantly reducing the computational cost associated with solving the multiscale model (see Fig.~\ref{fig:MOR-sarcomeres}).
An interesting aspect is that the stabilizing nature of the interaction between the scales in cardiac mechanics causes force fluctuations to be dampened thanks to feedback mechanisms (see \cite{RDQ-ML-2020} for a detailed mathematical treatment). This phenomenon has the beneficial effect of also damping the errors due to the approximation of the surrogate model compared to the FOM. 
Indeed, when the surrogate model is coupled with the macroscale model, the overall error of the multiscale model in calculating the main biomarkers is even lower than that of the standalone force generation model simulations, as shown in Tab.~\ref{tab:MOR-sarcomeres_cardiac_indicators} for a left ventricle simulation, where errors are on the order of $10^{-3}$.  The obtained speed--up makes the cost associated with solving the activation model negligible compared to that of other physics, and the dimensionality reduction of the state space significantly reduces the number of variables in the multiscale model: from 2198 variables for each node of the macroscale mesh (18 ionic variables, the transmembrane potential, 2176 activation variables, 3 components of the displacement) to only 24 variables (18 + 1 + 2 + 3).

\begin{table}
	\begin{tabularx}{\textwidth}{lrrr}
		\toprule
		Biomarker & FOM \eqref{eq:force_generation_FOM} & Surrogate model \eqref{eq:force_generation_ROM} & Relative error \\
		\midrule
		Stroke volume (\si{\milli\liter}) & 56.64 & 56.39 & $4.33\e{-3}$ \\
		Ejection fraction (\si{\percent}) & 44.48 & 44.29 & $4.33\e{-3}$ \\
		Maximum pressure (\si{\mmHg})     & 108.94 & 109.10 & $1.52\e{-3}$ \\
		Work (\si{\milli\joule})          & 662 & 659 & $4.85\e{-3}$ \\
		\bottomrule  
	\end{tabularx}  
	\caption{Cardiac biomarkers obtained through multiscale simulations where active force generation is modelled either by the FOM \eqref{eq:force_generation_FOM} or by the ROM \eqref{eq:force_generation_ROM}, and associated relative errors.}
	\label{tab:MOR-sarcomeres_cardiac_indicators}
\end{table}

\subsection{Time-dependent operator learning for a multiphysics coupled problem} \label{sec:EMROM}

In the previous section, we saw how operator learning methods can be used to surrogate the fine-scale dynamics in multiscale problems. Now we will see the application of operator learning in a different context, namely for the coupling of multiphysics systems, through a use case proposed in \cite{regazzoni2022machine}. Specifically, we consider the problem of coupling the cardiac electromechanics model, described by equations \eqref{eq:EP}\--\eqref{eq:4.1-QDR}--\eqref{eq:mechanics}, with the blood circulation model, described by \eqref{eq:6.1-QDR}. This problem is characterized by a strong interaction between the two scales, as the contraction of the cardiac muscle influences the flows and pressures in the circulatory network, which in turn influences the deformation of the cardiac muscle. The ultimate goal is to accelerate the evaluation of the input-output map associated with the multiphysics model, in order to quickly explore the dependence between the multiple parameters of the models and the outputs of clinical interest.

\null\textbf{The problem.}
We consider the 3D cardiac electromechanics model described by equations \eqref{eq:EP}\--\eqref{eq:4.1-QDR}--\eqref{eq:mechanics}, coupled with the 0D blood circulation model described by \eqref{eq:6.1-QDR}. For simplicity, we consider the case of a single cardiac chamber, specifically the left ventricle, but the extension to the case of multiple chambers is possible, as discussed in \cite{salvador2024whole}. The resulting coupled model can be written in the following compact form:
\begin{equation} \label{eqn:EMROM_FOM-coupled}
    \left\{
        \begin{aligned}
            & \frac{d \mathbf{z}_{\textrm{em}}}{dt} (\mathbf{x}, t) = \dynRHSop(\mathbf{z}_{\textrm{em}}, p_{LV}, \mathbf{x}, t, \mathbf{p}_{\textrm{em}}) 
            & & \qquad \mathbf{x} \in \Omega, \, t \in (0, T) \\
            & \frac{d\mathbf{z}_{\textrm{circ}}}{dt}(t)={\boldsymbol\Phi}_{\textrm{circ}}(\mathbf{z}_{\textrm{circ}}, p_{LV},t, \mathbf{p}_{\textrm{circ}}),
            & & \qquad t \in (0,T] \\
            & V_{LV}^{\textrm{em}}(\mathbf{z}_{\textrm{em}}) = V_{LV}^{\textrm{circ}}(\mathbf{z}_{\textrm{circ}})
            & & \qquad t \in (0,T] \\
            & \mathbf{z}_{\textrm{em}}(\mathbf{x},0) = \mathbf{z}_{\textrm{em},0}(\mathbf{x}) 
            & & \qquad \mathbf{x} \in \Omega \\
            & \mathbf{z}_{\textrm{circ}}(0) = \mathbf{z}_{\textrm{circ},0} 
            & & \\
        \end{aligned}
    \right.
\end{equation}
where $\mathbf{z}_{\textrm{em}}$ represents the vector that collects the state variables of the electromechanics model, i.e. ionic variables, transmembrane potentials, activation variables, tissue displacement, and velocity. The vector $\mathbf{z}_{\textrm{circ}}$ represents instead the state vector of the blood circulation model.
The two models are coupled through the kinematic compatibility condition $V_{LV}^{\textrm{em}} = V_{LV}^{\textrm{circ}}$, where the two terms represent the left ventricle volume calculated by the electromechanics model and by the blood circulation model, respectively. This constraint is enforced by means of a Lagrange multiplier, namely the pressure of the blood contained in the left ventricle $p_{LV}(t)$, which therefore appears as an input for both models. Both models depend on a set of physical parameters, represented respectively by $\mathbf{p}_{\textrm{em}}$ and $\mathbf{p}_{\textrm{circ}}$. For example, $\mathbf{p}_{\textrm{em}}$ includes the electrical conductivity of the tissue, the angle formed by the fibers, the contractility of the cardiomyocytes, while $\mathbf{p}_{\textrm{circ}}$ includes parameters such as the resistance and compliance of the blood vessels.

In clinical applications, the pressure $p_{LV}(t)$ and the volume $V_{LV}(t)$ of the left ventricle represent the quantities of fundamental importance for the diagnosis and therapy of cardiac pathologies, such as heart failure, and for the evaluation of the effectiveness of pharmacological treatments or medical devices. Developing computational tools capable of predicting both of them accurately and quickly in response to variations in the physical parameters $\mathbf{p}_{\textrm{em}}$ and $\mathbf{p}_{\textrm{circ}}$ is essential. 

\null\textbf{Methodology.}
The approach proposed in \cite{regazzoni2022machine} is based on the idea of surrogating the computationally expensive part of the multiphysics system \eqref{eqn:EMROM_FOM-coupled}, i.e., the 3D electromechanical model, while keeping the 0D circulation model in its full-order form, as it is computationally inexpensive.

To train the surrogate model, we first collect a dataset of transients of volumes and pressures obtained through the numerical approximation of the FOM model \eqref{eqn:EMROM_FOM-coupled}, by randomly sampling the physical parameters $\mathbf{p}_{\textrm{em}}$ and $\mathbf{p}_{\textrm{circ}}$. The dataset consists of $\dynNumSamples$ transients, each associated with the corresponding physical parameters: $(\widehat V_{LV}^j(t), \widehat p_{LV}^j(t), \widehat{\mathbf{p}}_{\textrm{em}}^j, \widehat{\mathbf{p}}_{\textrm{circ}}^j)$, for $j=1,\dots,\dynNumSamples$ and for $t \in (0, T)$.

The surrogate model is then defined by relying on the model learning approach with latent variables (see Sec.~\ref{sec:operator-learning_hidden-dynamics-discovery}), as the following system of ODEs:
\begin{equation}  \label{eqn:EMROM_ROM}
    \left\{
    \begin{aligned}
        & \frac{d}{dt} \dynLat(t) = \dynNNdyn\left(
            \dynLat(t), 
            p_{LV}(t), 
            \cos\left(\frac{2\pi t}{T_{HB}}\right),
            \sin\left(\frac{2\pi t}{T_{HB}}\right),
            \mathbf{p}_{\textrm{em}}; 
            \dynWdyn\right) 
            & & t \in (0, T) \\
        & V_{LV}(t) = V_{LV}^{NN}(\dynLat(t)) := \dynLat(t) \cdot \mathbf{e}_1 & & t \in (0, T) \\
        & \dynLat(0) = \dynLat_0 & &
    \end{aligned}
    \right.
\end{equation}
where $\dynLat(t) \in \mathbb{R}^{\dynNumLat}$ represents the latent variables, $\dynNNdyn$ represents a neural network with trainable parameters $\dynWdyn$, and $\mathbf{e}_1$ denotes the first vector of the canonical basis of $\mathbb{R}^{\dynNumLat}$. We are thus using the \textit{output-inside-the-state} approach to reduce the number of neural networks to be trained, as in Sec.~\ref{sec:mor-sarcomeres}. Consistently, the initial condition of the latent variables is set to $\dynLat_0 = (V_{LV}^{em}(\mathbf{z}_{\textrm{em},0}), 0, \dots, 0)^T$. Additionally, to leverage our physical knowledge of the system, since we know that the non-autonomous component in \eqref{eqn:EMROM_FOM-coupled} corresponds to the electrical stimulus that periodically activates the cardiac tissue, we introduce two additional inputs to the neural network $\dynNNdyn$, consisting of the coordinates of a point rotating on the unit circle with the same period $T_{HB}$ as the heartbeat, thus informing the surrogate model of the stimulus' periodicity. Finally, the model also receives the values of the physical parameters $\mathbf{p}_{\textrm{em}}$ as input, so that during training it can learn how the system dynamics depend on these parameters.

Training is then performed with a trajectory-based loss function (see Sect. \ref{sec:operator-learning-time-dependent}) through the following minimization problem:
\begin{equation*}
    \dynWdyn^* = \argmin{\dynWdyn} \left[
        \frac{1}{\dynNumSamples\,\dynNumTimes}\sum_{j=1}^{\dynNumSamples}\sum_{i = 0}^{\dynNumTimes} \left| \widehat V^j(t_i) - V^{NN}_{LV}(\dynLat^j(t_i)) \right|^2
        + \lambda \|\dynWdyn\|^2,
    \right],
\end{equation*}
where $\dynLat^j$ is obtained by numerically integrating the system \eqref{eqn:EMROM_ROM} from the initial condition $\dynLat_0$, by considering the physical parameters $\widehat{\mathbf{p}}_{\textrm{em}}^j$ and by imposing the input pressure equal to $\widehat p_{LV}^j(t)$. Note that during training, the model \eqref{eqn:EMROM_ROM} is solved independently of the circulation model to which it will later be coupled: in particular, the pressure transient is provided here as input and is not determined by the interaction with the circulation model. This strategy allows constructing a surrogate model of the 3D components of the model (i.e., the dynamics of the variables $\mathbf{z}_{\textrm{em}}$) without having to simultaneously learn the dynamics of the variables $\mathbf{z}_{\textrm{circ}}$.

To mitigate overfitting, we use Tikhonov regularization (see Sec.~\ref{sec:regularization}), where $\lambda$ is a hyperparameter that controls the regularization of the neural network weights.

Once the NN--based model has been trained, i.e., the optimal values $\dynWdyn^*$ have been determined, it can be coupled with the blood circulation model to obtain a hybrid model that can be used to evaluate the input-output map associated with the multiphysics model in a computationally efficient manner, by numerically approximating the following system of equations:
\begin{equation} \label{eqn:EMROM_ROM-coupled}
    \left\{
        \begin{aligned}
            & \frac{d}{dt} \dynLat(t) = \dynNNdyn\left(
            \dynLat(t), 
            p_{LV}(t), 
            \cos\left(\frac{2\pi t}{T_{HB}}\right),
            \sin\left(\frac{2\pi t}{T_{HB}}\right),
            \mathbf{p}_{\textrm{em}}; 
            \dynWdyn^*\right), 
            & & \ t \in (0, T) \\
            & \frac{d\mathbf{z}_{\textrm{circ}}}{dt}={\boldsymbol\Phi}_{\textrm{circ}}(\mathbf{z}_{\textrm{circ}}, p_{LV}(t),t, \mathbf{p}_{\textrm{circ}}),
            & & \ t \in (0,T] \\
            & V_{LV}^{NN}(\dynLat(t)) = V_{LV}^{\textrm{circ}}(\mathbf{z}_{\textrm{circ}}(t))
            & & \ t \in (0,T] \\
            & \dynLat(0) = \dynLat_0 
            & & \\
            & \mathbf{z}_{\textrm{circ}}(0) = \mathbf{z}_{\textrm{circ},0} 
            & & \\
        \end{aligned}
    \right.
\end{equation}

\null\textbf{Results and discussion.}
In \cite{regazzoni2022machine}, two test cases were considered. In the first one, the variability with respect to a single parameter of the electromechanical model, namely the contractility called $a_{XB}$, is considered, i.e., setting $\mathbf{p}_{\textrm{em}} = (a_{XB})$. In the second one, four parameters relevant to different aspects of the electromechanical model dynamics are considered instead, namely the electrical conductivity of the tissue in the fiber direction $\sigma_f$, the transmural fiber rotation angle $\alpha$, the tissue stiffness $C$, and again the contractility $a_{XB}$, i.e., setting $\mathbf{p}_{\textrm{em}} = (\sigma_f, \alpha, C, a_{XB})^T$. The number of samples used for training is $\dynNumSamples = 30$ for the 1-parameter case, $\dynNumSamples = 40$ for the 4-parameter case. In all cases, simulations containing 5 heartbeats each are considered.

In Tab.~\ref{tab:EMROM:errors_LV}, we present the errors quantifying the accuracy of the surrogate model coupled with the circulation model \eqref{eqn:EMROM_ROM-coupled} in reproducing the results of the full-order model (FOM) \eqref{eqn:EMROM_FOM-coupled}, for both the 1-parameter and 4-parameter cases. The table reports errors calculated from simulations lasting 5 heartbeats (the same duration as in the training set) and from simulations lasting 10 heartbeats, to assess the surrogate model's ability to generalize to longer simulations. The results show no significant difference between the two cases. This ability to extrapolate in time is particularly notable, as the simulations in the training dataset did not reach a limit cycle, meaning the model did not have the opportunity to observe full transients for most parameter values. Nonetheless, the surrogate model coupled with the circulation model is able to provide accurate results even in this context. This capability is important in clinical applications, where it is essential to evaluate the model’s behaviour under limit-cycle conditions (the only ones of clinical relevance). However, achieving such a limit cycle can often require many heartbeats, presenting significant computational challenges \cite{augustin2021computationally}.

\begin{table}
    \begin{center}
           \begin{tabular}{ rlrrrrrr }
           \toprule
           & & \multicolumn{6}{c}{\textbf{5 heartbeats}} \\
           & & $\PLV(t)$ & $\VLV(t)$ & $\pminLV$ & $\pmaxLV$ & $\VminLV$ & $\VmaxLV$
           \\
           \multirow{ 2}{*}{1-parameter} & \multicolumn{1}{l|}{relative error } &
           0.0336 & 0.0090 & 0.0097 & 0.0046 & 0.0139 & 0.0035
           \\
           & \multicolumn{1}{l|}{\Rtwo} &
                  &        & 0.9969 & 0.9986 & 0.9990 & 0.9995
           \\
           \multirow{ 2}{*}{4-parameter} & \multicolumn{1}{l|}{relative error } &
           0.0620 & 0.0285 & 0.0517 & 0.0272 & 0.0471 & 0.0127
           \\
           & \multicolumn{1}{l|}{\Rtwo} &
                  &        & 0.9437 & 0.9530 & 0.9594 & 0.9706
           \\
           \midrule
           & & \multicolumn{6}{c}{\textbf{10 heartbeats}} \\
           & & $\PLV(t)$ & $\VLV(t)$ & $\pminLV$ & $\pmaxLV$ & $\VminLV$ & $\VmaxLV$
           \\
           \multirow{ 2}{*}{1-parameter} & \multicolumn{1}{l|}{relative error} &
           0.0293 & 0.0071 & 0.0113 & 0.0037 & 0.0096 & 0.0031
           \\
           & \multicolumn{1}{l|}{\Rtwo} &
                  &        & 0.9992 & 0.9998 & 0.9985 & 0.9994
           \\
           \multirow{ 2}{*}{4-parameter} & \multicolumn{1}{l|}{relative error} &
           0.0631 & 0.0265 & 0.0442 & 0.0147 & 0.0382 & 0.0122
           \\
           & \multicolumn{1}{l|}{\Rtwo} &
                  &        & 0.9223 & 0.9996 & 0.9923 & 0.9906
           \\
           \bottomrule
           \end{tabular}
           \caption{Accuracy of the surrogate model coupled with the circulation model \eqref{eqn:EMROM_ROM-coupled} in reproducing the results of the FOM \eqref{eqn:EMROM_FOM-coupled}, in 5 heartbeats long (top) and 10 heartbeats long (bottom) simulations. For $PLV(t)$ and $\VLV(t)$ transients we report the relative $L^2$ error in time, while for scalar biomarkers, such as minimum and maximum pressure (i.e. $\pminLV$ and $\pmaxLV$) and volume (i.e. $\VminLV$ and $\VmaxLV$), we report the relative error and the $R^2$ coefficient of determination. 
           We consider both the 1-parameter and the 4-parameter cases, as indicated in the first column.}
           \label{tab:EMROM:errors_LV}
    \end{center}
\end{table}

A common approach in computational cardiology for making rapid predictions of the input-output map that links the parameters of the electromechanical model to variables of interest (such as $\pminLV$, $\pmaxLV$, $\VminLV$, and $\VmaxLV$) is to train a \textit{static} emulator that directly learns the mapping between the parameters $\mathbf{p}_{\textrm{em}}$ and $\mathbf{p}_{\textrm{circ}}$ and the variables of interest, typically using a FFNN \cite{Longobardi2020,Cai2021surrogate}. However, since such emulators must simultaneously capture the dependence of the output on both the electromechanical and circulation model parameters, they often require a significantly large number of training samples (on the order of 1000 or more). Moreover, they fail to capture the temporal dynamics of the system. In contrast, the approach proposed in \cite{regazzoni2022machine}, which leverages time-dependent operator learning, allows the model to surrogate only the computationally expensive part, drastically reducing the number of training simulations and providing accurate predictions even for long simulations.

Once trained, the reduced-order model (ROM) in \eqref{eqn:EMROM_ROM-coupled} offers a substantial computational advantage compared to the full-order model (FOM) in \eqref{eqn:EMROM_FOM-coupled}: the FOM takes approximately 4 hours of computation on 32 cores for each heartbeat, while the ROM requires less than a second on a single core of a standard laptop. This speed--up -- over a factor of $10^5$ -- is made possible by the surrogate model's ability to learn the essential dynamics of the FOM from the training data and generalize to new data, enabling accurate results in dramatically reduced computational times. However, when evaluating the computational advantages, the cost of generating the training dataset must be considered. Despite the approach's ability to achieve accurate predictions with relatively few training samples, the computational cost of generating the training data is significant. Additionally, the training time should be accounted for.

Therefore, the surrogate model proves particularly advantageous for applications requiring the exploration of the model's behaviour across a wide range of parameters, such as sensitivity studies and parameter estimation, in order to recompense the computational effort required by the construction of the surrogate model, as outlined in \eqref{eqn:many-query}. The next two sections delve into two such applications.

\subsection{NN--based surrogate models for global sensitivity analysis}

\textbf{The problem.}
In this section, we focus on the application of surrogate models for global sensitivity analysis (GSA) of cardiac electromechanics models.
GSA is a fundamental tool for understanding the behaviour of complex systems, such as cardiac models, by quantifying the influence of the model's parameters on {some} output quantities of interest (QoIs). 
However, performing GSA on cardiac models is computationally expensive, as it requires evaluating the model for a large number of parameter combinations. Surrogate models can significantly reduce the computational cost by providing fast and accurate predictions of the model's output for a wide range of parameters.

\null\textbf{Methodology.}
GSA is commonly performed by sampling the parameter space and calculating appropriate indicators, such as Borgonovo indices \cite{plischke2013global}, Morris elementary effects \cite{morris1991factorial}, Sobol indices \cite{sobol1990sensitivity,homma1996importance}, and Kucherenko indices \cite{kucherenko2012estimation}. Here, we focus on variance-based sensitivity analysis, which utilizes a probabilistic framework, however, this approach about surrogate models can also be applied in combination with other sensitivity analysis methods.

Sobol indices measure the sensitivity of a quantity of interest (QoI), denoted by $\qoi_j$, to a specific parameter, denoted by $\param_i$. The \textit{first-order Sobol index} (denoted $\SobolFirst{i}{j}$) quantifies the effect when the parameter varies independently over the changes in the other parameters:
\begin{equation*}
\SobolFirst{i}{j} = \frac{\variance_{\param_i}\left[ \expected_{\param_{\sim i}}\left[ \qoi_j | \param_i\right] \right]}{\variance\left[ \qoi_j \right]},
\end{equation*}
where $\param_{\sim i}$ represents the set of all parameters except for the $i$-th one. 
In this expression, the expected value $\expected_{\param_{\sim i}}\left[ \qoi_j | \param_i\right]$ represents the average value of $\qoi_j$ when $\param_i$ is fixed, while all other parameters $\param_{\sim i}$ vary according to their distributions. Intuitively, this captures how the behaviour of the QoI depends on the fixed value of $\param_i$ while accounting for the variability due to all other parameters.
The variance $\variance_{\param_i}\left[ \expected_{\param_{\sim i}}\left[ \qoi_j | \param_i\right] \right]$ then measures how much this expected value changes as $\param_i$ varies. This quantifies the contribution of $\param_i$ to the overall variability of the QoI. Finally, normalizing by the total variance $\variance\left[ \qoi_j \right]$ provides a relative measure of sensitivity, enabling comparison between different parameters.

To account for the combined effects of parameters, including their interactions, the \textit{total-effect Sobol index} $\SobolTotal{i}{j}$ is used:
\begin{equation*}
    \SobolTotal{i}{j} = \frac{\expected_{\param_{\sim i}}\left[ \variance_{\param_i}\left[ \qoi_j | \param_{\sim i}\right] \right]}{\variance\left[ \qoi_j \right]}
                  = 1 - \frac{\variance_{\param_{\sim i}}\left[ \expected_{\param_i}\left[ \qoi_j | \param_{\sim i}\right] \right]}{\variance\left[ \qoi_j \right]}.
\end{equation*}
This index captures the influence of a parameter when it varies on its own as well as in conjunction with other parameters \cite{sobol1990sensitivity}.

The Saltelli method \cite{homma1996importance,saltelli2002making} provides a way to estimate the Sobol indices by employing Sobol quasi-random sequences to approximate the necessary integrals. In practice, this method involves evaluating the model across a large number of parameter combinations and subsequently processing the resulting QoIs to estimate the Sobol indices.

In \cite{regazzoni2022machine}, variance-based GSA of the electromechanics-circulation model is performed on the surrogate model \eqref{eqn:EMROM_ROM-coupled} to approximate the evaluation of the QoIs, by simultaneously considering the variability with respect to the circulation model parameters $\mathbf{p}_{\textrm{circ}}$ and the electromechanical model parameters $\mathbf{p}_{\textrm{em}}$ (that is, we set $\param = (\mathbf{p}_{\textrm{em}}^T, \mathbf{p}_{\textrm{circ}}^T)^T$).
A set of 20 QoIs is considered, including the maximum and minimum pressures and volumes associated with the cardiac chambers and the arterial systemic circulation along the heartbeat.
For each parameter choice, the surrogate model is used to simulate the electromechanical model for a certain number of heartbeats, until a limit cycle is reached, and the QoIs are calculated with respect to the last heartbeat. The Saltelli method is then applied to estimate the Sobol indices, by evaluating the surrogate model at a set of parameter combinations generated by the Sobol quasi-random sequence. 

\null\textbf{Results and discussion.}
For the sake of space, we do not report the estimated Sobol indices (for more details, the interested reader is referred to \cite{regazzoni2022machine}). Instead, we focus on the computational benefits of using surrogate models for GSA. For this, we consider a realistic scenario to performing GSA on a cardiac electromechanics model using a 160-core cluster.
Computational times were measured on Intel Xeon E5-2640 v4 2.4 GHz CPUs.
In this case, conducting the GSA would require simulating 74 000 parameter sets to achieve accurate results. On average, the system reaches a limit cycle after 10 heartbeats, leading to a total of 740 000 heartbeats that need to be simulated. Running this large number of simulations with the full-order model (FOM) would be practically impossible, as it would take around 68 years of continuous computation on the 160-core cluster.
In contrast, using the surrogate model allows for the estimation of the Sobol indices in just 7.5 days. This includes 6.7 days to generate the training dataset, 18 hours to train the model, and 1 hour and 17 minutes to compute the QoIs required by the Saltelli method. As a result, employing the surrogate model provides a remarkable speed--up of 3 300 times.

\subsection{NN--based surrogate models for Bayesian parameter estimation} \label{sec:EMROM_MCMC}

\null\textbf{The problem.}
Personalizing a cardiac electromechanical model for a specific patient involves more than just using geometry obtained from imaging data. It also requires estimating the model's key parameters based on available clinical measurements. However, these measurements often consist of only a few scalar values, and solving the inverse problem (i.e., determining the physical parameters $\param$ from the observed QoIs $\qoi$) must account for the noise inherent in these measurements, which introduces uncertainty into the parameter estimates.

\null\textbf{Methodology.}
Bayesian methods, such as Markov Chain Monte Carlo (MCMC) \cite{brooks1998markov} and Variational Inference \cite{hoffman2013stochastic}, address these challenges within a robust statistical framework. These methods compute the \textit{likelihood}, representing the probability distribution of the desired parameter values given the observed QoIs ($\qoiObs$). Unlike optimization techniques that yield point estimates, Bayesian approaches provide a probability distribution over the parameter space, reflecting the \textit{credibility} of different parameter combinations.

This credibility assessment integrates measurement uncertainty -- modelled as noise with covariance matrix $\NoiseCov$ -- and prior knowledge about the parameters, encoded in a \textit{prior distribution} $\piPrior(\param)$. Using the parameters-to-QoIs map $\forward \colon \param \mapsto \qoi$, the observed QoIs are expressed as $\qoiObs = \forward(\param) + \error$, where $\error \sim \mathcal{N}( \cdot | \mathbf{0}, \NoiseCov)$ represents Gaussian measurement noise. Bayes' theorem provides the \textit{posterior distribution}, which quantifies the belief in parameter values after observing $\qoiObs$:
\begin{equation*}
\piPost(\param) = \frac{1}{Z}\, \mathcal{N}( \qoiObs | \forward(\param), \NoiseCov) \, \piPrior(\param),
\end{equation*}
where the normalization constant $Z$ is given by:
\begin{equation*}
Z = \int_{\paramSpace} \mathcal{N}( \qoiObs | \forward(\widehat{\param}), \NoiseCov) \, d\piPrior(\widehat{\param}).
\end{equation*}
Computing $\piPost$ is computationally intractable because it involves approximating the integral defining $Z$. MCMC provides an efficient way to approximate $\piPost$ with moderate computational effort. 

Similar to the Saltelli method used in sensitivity analysis, MCMC requires numerous model evaluations for various parameter values. Importantly, this method is non-intrusive, i.e. it does not require knowledge of the underlying mathematical model, but relies solely on evaluations of the map $\forward \colon \param \mapsto \qoi$.
To reduce computational costs, we can therefore replace the FOM with a surrogate model that significantly accelerates the process, enabling efficient computation of the posterior distribution while maintaining accuracy.

The Bayesian framework also enables a rigorous treatment of the approximation error introduced when the high-fidelity model \eqref{eqn:EMROM_FOM-coupled} is replaced by its surrogate \eqref{eqn:EMROM_ROM-coupled}. Let $\forwardRed$ represent the approximated parameters-to-QoIs map defined by the surrogate model \eqref{eqn:EMROM_ROM-coupled}. In this context, the high-fidelity map can be expressed as $\forward(\param) = \forwardRed(\param) + \errorROM$, where $\errorROM$ denotes the surrogate model's approximation error.

Consequently, the observed QoIs can be expressed by $\qoiObs = \forwardRed(\param) + \errorROM + \errorEXP$, where $\errorEXP$ accounts for the experimental measurement error. Assuming that the two error sources are independent, the covariance of the total error $\error = \errorROM + \errorEXP$ satisfies $\NoiseCov = \NoiseCovROM + \NoiseCovEXP$. Here, $\NoiseCovROM$ represents the covariance of the surrogate model's approximation error, which can be estimated using a validation set, while $\NoiseCovEXP$ corresponds to the covariance of the experimental error, determined by the specific measurement protocol. This formulation ensures that the surrogate model's approximation error is appropriately accounted for during parameter estimation.

\null\textbf{Results and discussion.}
To evaluate the surrogate model's ability to accelerate parameter estimation for multiscale cardiac electromechanical models, the following test was performed in \cite{regazzoni2022machine}. A high-fidelity simulation using \eqref{eqn:EMROM_FOM-coupled} was first conducted to generate a pair of QoIs $\qoiObs$ commonly measured in the clinical practice, namely maximum and minimum arterial pressures. Subsequently, the surrogate model \eqref{eqn:EMROM_ROM-coupled} was employed to approximate the high-fidelity model, and a Bayesian parameter estimate was performed to reconstruct the value of two relevant parameters, namely the active contractility and the systemic arterial resistance, assuming the values of the remaining parameters to be known. We note that the two parameters belong to different models: the active contractility is a parameter of the electromechanical model ($\mathbf{p}_{\textrm{em}}$), while the systemic arterial resistance is a parameter of the circulation model ($\mathbf{p}_{\textrm{circ}}$). 
This test case demonstrates this approach ability to estimate parameters of both the surrogated model (the electromechanical one) and the original, high fidelity circulation model. Finally, the estimated parameters were validated against the true values used to produce $\qoiObs$.

To simulate the presence of measurement errors, synthetic noise of varying magnitudes was artificially added to the exact QoIs. Specifically, the noise was sampled from independent Gaussian distributions with zero mean and variance $\NoiseMagnEXP^2$. Three cases were considered: $\NoiseMagnEXP^2 = 0$ (i.e., no noise), $\NoiseMagnEXP^2 = \SI{0.1}{\mmHg\squared}$, and $\NoiseMagnEXP^2 = \SI{1}{\mmHg\squared}$. For both unknown parameters, a non-informative prior was adopted, represented by a uniform distribution over the ranges used to train the ROM. 
The covariance of the total error was thus defined as $\NoiseCov = \NoiseCovROM + \NoiseCovEXP$, where $\NoiseCovEXP = \NoiseMagnEXP^2 \, \mathbb{I}_2$ represents the experimental measurement error covariance ($\mathbb{I}_2$ is the 2-by-2 identity matrix), and $\NoiseCovROM$ is the ROM approximation error covariance, estimated from its empirical statistical distribution on the validation set.

\begin{figure}
    \centering
    \includegraphics[width=.95\textwidth]{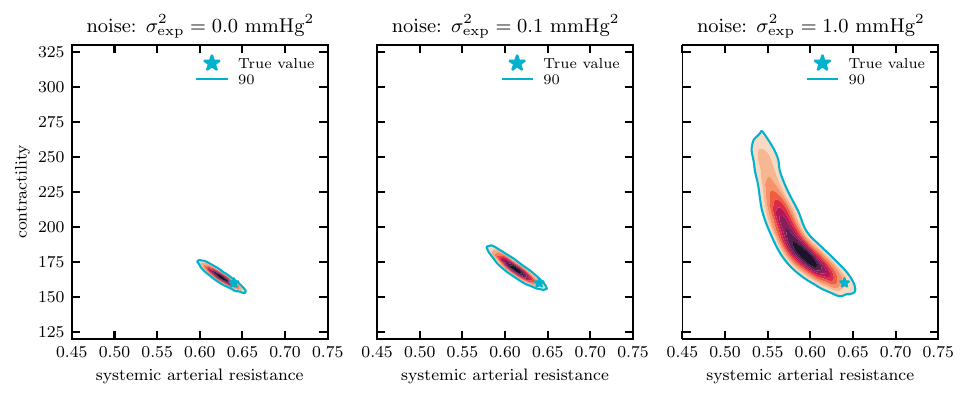}
    \caption{
    {Output of the Bayesian estimation presented in Sec.~\ref{sec:EMROM_MCMC}.}
    We depict the posterior distribution $\piPost$, estimated by means of the MCMC method, for $\NoiseMagnEXP^2 = 0$ (left), $\NoiseMagnEXP^2 = \SI{0.1}{\mmHg\squared}$ (middle) and $\NoiseMagnEXP^2 = \SI{1}{\mmHg\squared}$ (right).
    The blue lines show the 90\% credibility regions, while the red stars represent the exact value of the unknown parameters.    }
    \label{fig:MCMC}
\end{figure}

Figure~\ref{fig:MCMC} displays the posterior distribution $\piPost$ for the parameter pair under the three noise levels considered. The blue line indicates the 90\% credibility region, which is the region of the parameter space with the highest posterior probability that encompasses 90\% of $\piPost$. In each noise scenario, the credibility region includes the true parameter values, denoted by a blue star. As expected, higher noise levels result in larger credibility regions, reflecting greater uncertainty in the parameter estimates.

One notable advantage of Bayesian parameter estimation methods, compared to deterministic approaches, lies in their ability to quantify the uncertainty associated with parameter estimates. Additionally, Bayesian methods capture correlations between parameters, as evidenced by the oblique shape of the credibility regions. This correlation arises because changes in contractility and arterial resistance can produce similar effects on the measured QoIs (minimum and maximum arterial pressures), leading to highly correlated posterior distributions.

Similarly to the application to GSA, the use of NN--based surrogate models for Bayesian parameter estimation results in remarkable computational savings. Specifically, the estimation process requires approximately 960 000 heartbeats to be simulated. Leveraging on the surrogate model, the entire procedure -- including the generation of the training dataset -- can be completed in just 6 days and 8 hours. In contrast, performing the same task with the FOM would take over 87 years on a 160-core cluster, yielding a remarkable speed--up of 5 000x. 
{As in the previous section, computational times were measured on Intel Xeon E5-2640 v4 2.4 GHz CPUs.}
Furthermore, if Bayesian calibration needs to be repeated for different datasets, the NN--based surrogate model does not need retraining. In such cases, only the MCMC algorithm must be re-executed, which, thanks to the surrogate model, takes only 13 hours and 20 minutes -- further emphasizing the efficiency and practicality of the proposed approach.

\subsection{Latent Dynamics Networks to accelerate electrophysiology simulations}\label{sec:LDNet-EP}

Simulating cardiac electrophysiology is critical for understanding and predicting the electrical behaviour of cardiac tissue. Applications range from diagnosing arrhythmias to optimizing defibrillation and ablation strategies. Clinically, these simulations can help personalize treatment plans, improve device design, and support decision--making in emergency scenarios.

From a numerical perspective, electrophysiological models are characterized by steep gradients in transmembrane potentials, nonlinearity, and multi--scale coupling. These features make them computationally expensive, especially when high--resolution meshes or long simulation times are required. This computational burden limits the feasibility of traditional numerical methods in scenarios demanding quick (real time, or nearly real time) responses.

Surrogate models offer approximations of the results of full-scale simulations with significantly lower computational costs, thus enabling real-time predictions.

\null
\textbf{The Problem.}
We consider the Monodomain equation coupled with the Aliev--Panfilov model to describe electrical signal propagation in cardiac tissue. The model equations are:
\begin{equation}\label{eq:AlievPanfilov}
\left\{
\begin{aligned}
    &\frac{\partial v}{\partial t} - \mathbf{D} \Delta v = K v (1 - v) (v - \alpha) - v w + I_\mathrm{app}(\mathbf x, t)  
    && \text{in } \Omega\times(0,T) \\
    &\frac{\partial w}{\partial t} = \left(\gamma + \frac{\mu_1 w}{\mu_2 + v}\right) (- w - K v (v - b - 1))
    && \text{in } \Omega\times(0,T) \\
    & \left(\mathbf{D} \nabla v \right) \mathbf{n} = \mathbf{0}
    && \text{in } \partial\Omega\times(0,T) \\
    &v(\mathbf x, 0) = 0, \, w(\mathbf x, 0) = 0 
    && \text{in } \Omega \\
\end{aligned}
\right.
\end{equation}
where $v(\mathbf x, t)$ is the transmembrane potential, $w(\mathbf x, t)$ is the recovery variable, and $I_\mathrm{app}(\mathbf x, t)$ represents external electrical stimuli applied at specific locations. The problem features steep depolarization fronts and wave collisions, making it a challenging test for space--time operator learning.

We consider two test cases. In the former, we consider the one-dimensional domain $\Omega = (0,L)$, and two stimulation points located at $x = 1/4 \, L$ and $x = 3/4 \, L$, where we apply square impulses to mimic the action of a pacemaker. 
In the second test case, we consider a 2D square domain with a first wavefront propagating to the right, followed by the application of a second circular stimulus, at the centre of the square domain.
Depending on the stimulus radius and the timing of its application, three distinct outcomes can occur:
\begin{itemize}
    \item tissue refractoriness: the circular stimulus fails to trigger a second activation because it is applied while the tissue remains in a refractory state;
    \item focal activation: a single focal activation arises when the circular stimulus is delivered after the vulnerable window;
    \item re--entrant drivers: two self--sustained re--entrant drivers emerge, continuously reactivating the tissue.
\end{itemize}
The second test case features for more intricate spatial patterns, including bifurcating phenomena, as detailed above.

\null\textbf{Methodology.}
For both test cases, in \cite{regazzoni2024ldnets} different space-time operator learning methods were compared.
The methods considered were (see Sec.~\ref{sec:space-time learning}): a projection-based method, consisting of {Proper Orthogonal Decomposition} (POD) combined with {Discrete Empirical Interpolation Method} (DEIM), denoted as POD/DEIM; autoencoder-based methods, wherein the dynamics in the latent space is learned either through a Neural ODE or an {Long Short--Term Memory} (LSTM, {see Sect. \ref{sec:architectures}}), denoted by AE/ODE and AE/LSTM, respectively; and Latent Dynamics Networks (LDNets). Furthermore, we consider the case where after training of the models AE/ODE and AE/LSTM, their parameters are further optimized in an end--to--end manner (see Sec.~\ref{sec:operator-learning_hidden-dynamics-discovery}). These methods are denoted by AE/ODE--e2e and AE/LSTM--e2e, respectively.

For all the methods considered, the same training dataset was generated by solving the governing equations \eqref{eq:AlievPanfilov} using a numerical solver. It includes time series of the applied stimuli at specific locations, and the models were trained to predict the transmembrane potential $v(\mathbf x, t)$. 
For the 1D test case, the governing equations were approximated by means of the finite difference method both in space and time, on a regular grid with $\Delta t = \SI{10}{\micro\second}$ and 800 points in space.
Then, the space--time grid was subsampled by retaining 500 time instants and 100 points in space.
The 2D problem was numerically approximated using the $P1$ finite element method on a structured grid with an element size of $h = \SI{0.5}{\milli\meter}$. A semi-implicit time--stepping scheme was employed, dividing the time domain into steps of $\Delta t = \SI{0.25}{\milli\second}$. To construct the datasets, the space--time grid was subsampled, retaining 2,694 spatial points and 180 time instants.
In the 1D test case, 100 training/validation samples and 100 test samples were considered, while in the 2D test case, 200 training/validation samples and 75 test samples were used.

To ensure a fair comparison, an automatic tuning algorithm was used to select the optimal hyperparameter values for the different methods, with an upper bound of $\dynNumLat \leq 12$ set on the latent space dimension. 
Specifically, the Tree-structured Parzen Estimator (TPE) Bayesian algorithm \cite{bergstra2011algorithms,optuna2019} is used in conjunction with a $k$-fold cross-validation procedure (see Algorithm~\ref{alg:kcross}).
The results reported are obtained using the optimal hyperparameter configuration chosen by the tuning algorithm, independently for each method. 

\null\textbf{Results and Discussion.}
The results obtained in both test cases and {using} the different methods are reported in Tab.~\ref{tab:EP_1D_errors} (1D test case) and Tab.~\ref{tab:EP_2D_errors} (2D test case). For simplicity, in the 2D test case only the methods that performed the best in the 1D case are considered.

Due to the presence of travelling wavefronts, this problem exhibits a slow decline in the Kolmogorov $n-$width (see Sec.~\ref{sec:operator-learning_hidden-dynamics-discovery}), which results in reduced accuracy when reconstructing the electrical potential using the POD-DEIM method with 12 modes. Achieving satisfactory results requires more than 24 modes, but this leads to an increased computational cost during the prediction phase. Moreover, the POD-DEIM method offers minimal speed-up compared to the other methods under consideration. This limitation is due to the need for the POD-DEIM model to be solved using the same temporal discretization as the high-fidelity model for numerical stability reasons. This requirement imposes a significant constraint relative to the other methods considered in the comparison.

Better accuracy is attained by both autoencoder--based methods and LDNets, owing to their capacity to model non--linear relationships between the latent states and the solution. Among these, LDNet provides superior performance, achieving a testing normalized RMSE of approximately $\num{7e-3}$ in both the 1D and 2D test cases. The testing normalized RMSE of the autoencoder--based methods is approximately five times larger than that of LDNets or more. In the 2D case, which adds an extra spatial dimension compared to the 1D case, the advantage of LDNets over the other methods is more pronounced. As shown in \cite{regazzoni2024ldnets}, autoencoder--based methods display various artifacts in their solutions, particularly failing to accurately model scenarios where tissue refractoriness prevents signal propagation. In contrast, LDNet's predictions are nearly indistinguishable from those of the high-fidelity model, effectively capturing the three distinct behaviours exhibited by the system in this test case. The LDNet method, in the test cases considered in \cite{regazzoni2024ldnets}, achieves superior accuracy with a significantly smaller number of trainable parameters: autoencoder--based methods require over ten times more parameters in the 1D case and more than 400 times more in the 2D case.

\begin{table}
    \centering
    \begin{footnotesize}
    \begin{tabular}{l|r|r|r|r|r|r|r|r}
        \toprule
                                        & \multicolumn{2}{c|}{normalized RMSE}           &  Trainable parameters      &  \multicolumn{2}{c}{{Wall time (s)}}            \\
                                        & \multicolumn{1}{c|}{training}         &  \multicolumn{1}{c|}{testing}                    &                                      &  {offline      }&  {online   }       \\ \midrule
        {FOM}                           &                   &                         &                                         &                        &  {37.321   }       \\ 
        POD-DEIM ($\dynNumLat = 12$)     &  \num{4.05e-01}  &  \num{3.92e-01}         &                                         &  \num{797     }&  {5.839    }       \\ 
        POD-DEIM ($\dynNumLat = 24$)     &  \num{3.59e-01}  &  \num{3.47e-01}         &                                         &  \num{799     }&  {7.720    }       \\
        POD-DEIM ($\dynNumLat = 36$)     &  \num{1.71e-01}  &  \num{1.62e-01}         &                                         &  \num{861     }&  {7.442    }       \\
        POD-DEIM ($\dynNumLat = 48$)     &  \num{7.48e-02}  &  \num{7.57e-02}         &                                         &  \num{1124    }&  {7.976    }       \\
        POD-DEIM ($\dynNumLat = 60$)     &  \num{2.97e-02}  &  \num{2.90e-02}         &                                         &  \num{1242    }&  {8.408    }       \\
                AE/LSTM                   &  \num{1.90e-01}  &  \num{1.98e-01}        & \num{17933}                             &  \num{11009}   &  {0.005    }       \\
            AE/LSTM-e2e                   &  \num{2.05e-02}  &  \num{5.87e-02}        & \num{17933}                             &  \num{33851}   &  {0.005    }       \\
                AE/ODE                   &  \num{2.09e-02}  &  \num{4.58e-02}         & \num{22697}                             &  \num{23982}   &  {0.017    }       \\
            AE/ODE-e2e                   &  \num{1.78e-02}  &  \num{3.37e-02}         & \num{22697}                             &  \num{97821}   &  {0.017    }       \\
                LDNet                   &  \num{7.09e-03}   &  \num{7.37e-03}         & \num{ 1708}                             &  \num{22887}   &  {0.014    }       \\
        \bottomrule
    \end{tabular}    
    \end{footnotesize}
    \caption{Training and test errors, the number of trainable parameters, and the computational times for the offline and online phases are reported for the 1D test case described in Sec.~\ref{sec:LDNet-EP}. The computational times were measured on an Intel Xeon Processor E5-2640 2.4GHz. The offline phase corresponds to the model construction process. For the POD/DEIM method, this includes generating the basis for the solution manifold and DEIM approximation. For the other methods, the offline phase involves neural network training. Conversely, the online phase refers to predicting the system's evolution for a new sample after the model has been constructed.}
    \label{tab:EP_1D_errors}
\end{table}

\begin{table}
        \centering
        \begin{footnotesize}
        \begin{tabular}{l|r|r|r|r|r|r|r|r}
            \toprule
                                           & \multicolumn{2}{c|}{normalized RMSE}           &  Trainable parameters     &  \multicolumn{2}{c}{{Wall time (s)}}            \\
                                           & \multicolumn{1}{c|}{training}         &  \multicolumn{1}{c|}{testing}                   &                                     &  {offline}  &  {online   }       \\ \midrule
            {FOM}                          &                  &                             &                                     &                    &  {807.210}       \\ 
                  AE/ODE                   &  \num{6.96e-02}  &  \num{7.83e-02}             & \num{1193732}                       &         \num{75315}&  {0.191}       \\
              AE/ODE-e2e                   &  \num{3.97e-02}  &  \num{4.23e-02}             & \num{1193732}                       &         \num{95479}&  {0.188}       \\
                   LDNet                   &  \num{7.31e-03}  &  \num{7.57e-03}             &    \num{2789}                       &         \num{90349}&  {0.139}       \\
            \bottomrule
        \end{tabular}    
        \end{footnotesize}
        \caption{Training and test errors, the number of trainable parameters, and the computational times for the offline and online phases are reported for the 2D test case described in Sec.~\ref{sec:LDNet-EP}. See caption of Fig.~\ref{tab:EP_1D_errors} for more details.}
        \label{tab:EP_2D_errors}
\end{table}

\section{Some final thoughts and concluding remarks} \label{sec:conclusions}

We would like to conclude this paper with a few considerations regarding both Machine Learning (ML) and Scientific Machine Learning (SciML). 

Our goal is to shed light on the truly distinctive features of ML and SciML models. From this perspective, it will be helpful to use the well--known least--squares algorithm as a reference point, given its popularity in data fitting problems. Some of the conclusions will be quite surprising and go against the common perceptions of artificial intelligence.

An essential component of the learning process is the choice of the model, that is, the function $f$ mapping inputs (training data) to outputs (the answer to our question). The model $f$ is defined by a set of parameters and possibly hyperparameters. The parameters are determined through a process of minimization of an appropriate function $\mathcal L$ (called the loss function), while hyperparameters are specified by the user.

It is worth noting that least--squares methods, widely used in the mathematical community centuries before the advent of Machine Learning, can also be viewed as learning processes.
As a matter of fact, least--squares regression looks for the best parameters of a model $f$ -- typically a polynomial -- by minimizing the Mean Square Error (MSE) between available target data and the values predicted by the model itself. Like in Machine Learning, least--squares minimize a loss function (the MSE) and work on a training dataset (the target data used to evaluate the loss). The choice of the model $f$ can be driven by physical knowledge of the process at hand. This ansatz can be accurate in some circumstances, and inappropriate in others.

For example, suppose we aim to determine the constitutive stress--strain relationship of a material from experimental data. After collecting a set of stress--strain data points we can fit a model -- e.g. a linear one -- to the data. Afterwards, we can use the model to predict the stress associated with a strain value not present in the training set. This is a typical example of a learning process, where the model is trained on a set of data and then used to make predictions on new data. However, should the true stress--strain curve be more complex than linear, e.g., of exponential type, the linear model would fail to generalize effectively, i.e., it would not predict accurate stresses in correspondence of new strains.

The origins of the least-squares method trace back to Carl Friedrich Gauss, who applied it in 1795 to predict the orbit of the asteroid Ceres \cite{stigler1981gauss} after it was no longer observable because hidden by the sun. Assuming an elliptical orbit, Gauss fitted six orbital parameters to observational data \cite{gauss1877theoria}. Remarkably, his predictions differed significantly from those of other astronomers yet proved almost exact when Ceres was rediscovered after passing behind the sun, earning Gauss international acclaim \cite{gustafsson2018scientific}. This represents one of the earliest examples of successful generalization in history.

It is worth noticing that assuming an elliptic orbit was a strong hypothesis that Gauss used to constrain the model. In the same way, fitting a stress--strain relationship with a linear or rather exponential model relies on strong prior assumptions that the user imposes on the model. 
In typical Machine Learning applications, instead, one typically opts for a more complex model, such as a neural network (NN), that can capture a wide range of input--output relationships without making any a--priori assumptions on the law that has generated the training data.

For this reason, Machine Learning is often described as \textit{automatically discovering} the best model for the data. Strictly speaking, the model $f$ is still chosen a priori by the user, but the hypothesis space is so rich that, essentially, the algorithm autonomously discovers a model.

As for the more purely algorithmic aspects, we identify some differences between classical least--squares approaches and Machine Learning algorithms, such as NNs, which are worth emphasizing.

In the case of linear least--squares methods, i.e., when the model $f$ depends linearly on the parameters, the minimization algorithm is deterministic: the gradient of the loss function is computed analytically, and its nullification ($\nabla \mathcal L=0$) leads to a linear system (normal equations), which will then be solved using algebraic techniques. 

In the case of non--linear least--squares methods, $f$ depends non--linearly on the parameters, thus the resulting normal equations are nonlinear and, as such, they should be solved by an iterative method which, typically, performs well on small--size systems.

If the learning model is based on NNs, due to the complex compositional structure of the model function $f$, the minimization algorithm typically leverages backpropagation for gradient computation (as now the analytical evaluation of $\nabla \mathcal L$ is cumbersome to do with paper--and--pen, thus automatic differentiation is preferable) and usually employs an iterative minimization algorithm based on stochastic gradient descent.
Hence, the deterministic nature of the process is lost. Nevertheless, it is worth noting that, once trained via the minimization process, a NN is deterministic. This means that the associated algorithm can be described unambiguously and, given the same input data, will always produce the same output. Therefore, a NN retains the property of \emph{reproducibility}.
Considerations apart deserve the cases of reinforcement learning, adversarial neural networks, and generative AI algorithms, where the determinism is lost.

Another consideration concerns the problem size. Typically, a least--squares approach (whether linear or non--linear) relies on a very limited number of parameters, unlike NN approaches which instead may involve extraordinarily large numbers of parameters, alongside significantly larger training datasets. The many parameters, together with the highly non--linear and compositional structure of NNs, enhance their \emph{capability} to represent highly complex and high--dimensional input--output processes (see, for instance, Fig. \ref{fig:parameters_size}). This is, e.g., the case of convolutional neural networks for image recognition, such as ResNet--50 (with 25 million parameters), and AlexNet (60 million parameters). 

Further, a frequently debated aspect is the presumed lack of \emph{interpretability} of ML algorithms based on NN models. Strictly speaking, just as in the case of least--squares, NN models are in a certain sense interpretable once training is complete: indeed, once parameters and hyperparameters are determined, the input--output transfer function can be represented in finite and unambiguous terms. Of course, it must be acknowledged that the ``readability'' of such a function can be highly problematic due to its non--linear compositional form. In this respect, commonly used least--squares model functions are far more readable and interpretable. To put it differently, identifying the role of each parameter in shaping the NN response is very challenging.

Going back to the stress--strain fitting example, a linear relationship is easily interpretable: the slope of the line is a measure of the stiffness of the material, while the intercept represents the residual stress at zero strain. With a polynomial model, the interpretation becomes more complex but still manageable. However, interpreting the parameters of a NN model is much more challenging, as the relationship between the input and the output is not straightforward and depends on the entire network structure. This is one of the reasons why NNs are often referred to as black--box models.
However, it is important to note that the distinction between interpretable and black--box models is not absolute.
There exists a spectrum of models going from simpler and more interpretable ones to others more complex yet richer and harder to interpret.
In this transition, there is a trade--off between interpretability and model complexity. While simpler models offer clear, understandable relationships between inputs and outputs, they cannot often capture intricate patterns within the data. On the other hand, more complex models, such as neural networks, can handle highly non--linear and high--dimensional relationships, but at the cost of losing interpretability. This balance is a key consideration when choosing the appropriate model to solve a given problem.

Finally, turning to Scientific Machine Learning, a few observations are in order. The first one is that SciML is a relatively new area, much is still a work in progress, and the most significant achievements often result from trial and error rather than being situated within a rigorous and established theoretical framework. Key characteristics of this area include the co--presence, in various forms, of physics--based methods alongside ML algorithms. With some liberty, we could say that this falls within the domain of ``grey boxes'', which result from the fusion of white--box models (physics--based models that are interpretable and explainable) and black--box models (those driven by data and that entirely ignore the context that generated that data). Due to the almost complete arbitrariness with which these combinations can be made, the classifications we have sketched in this paper should be considered a preliminary attempt to outline some identifying elements that may begin to shed light on a still quite obscure (albeit vibrant) field. We do not doubt that this attempt will soon require revision, given the extraordinary interest we are experiencing in these topics. 
\section*{Acknowledgments}
This article draws inspiration from the content of short courses recently conducted by its authors, specifically, by F.R. at Politecnico di Torino in Fall 2023, A.Q. at the University of Maryland in College Park in Spring 2024, and A.Q. at Sorbonne Université in Fall 2024. A.Q. and F.R. acknowledge the grant Dipartimento di Eccellenza 2023–2027 (Dipartimento di Matematica, Politecnico di Milano).
A.Q. acknowledges the financial support of the Brin Mathematics Research Center of the University of Maryland at College Park and the Laboratoire Jacques--Louis Lions of Sorbonne Université at Paris. For the results on the cardiac model, the E.R.C. AdG n. 740132 is gratefully acknowledged. P.G. has received support from the project PRIN, MUR, Italy, CUP 20227K44ME. F.R. has received support from the project PRIN2022, MUR, Italy, 2023-2025, P2022N5ZNP “SIDDMs: shape-informed data--driven models for parametrized PDEs, with application to computational cardiology”, funded by the European Union (Next Generation EU, Mission 4 Component 2). P.G. and F.R. are members of the Gruppo Nazionale Calcolo Scientifico – Istituto Nazionale di Alta Matematica (GNCS INdAM). 

\bibliographystyle{plain}
\bibliography{bibl}

\end{document}